\documentclass{amsart}

\usepackage{amsmath} 
\usepackage{amssymb}
\usepackage{mathrsfs}

\newtheorem{theorem}{Theorem}[section] 
\newtheorem{claim}[theorem]{Claim}
\newtheorem{cc}[theorem]{The Crucial Claim}

\newtheorem{conclusion}[theorem]{Conclusion}
\newtheorem{observation}[theorem]{Observation}

\theoremstyle{definition}
\newtheorem{definition}[theorem]{Definition}

\newtheorem{example}[theorem]{Example}

\newtheorem{thesis}[theorem]{Thesis}
\newtheorem{discussion}[theorem]{Discussion}
\newtheorem{convention}[theorem]{Convention}
\newtheorem{conjecture}[theorem]{Conjecture}

\newtheorem{hypothesis}[theorem]{Hypothesis}

\theoremstyle{remark}
\newtheorem{remark}[theorem]{Remark}
\newtheorem{question}[theorem]{Question}
\newtheorem{notation}[theorem]{Notation}

\newcommand{\rd}{{\rm rd}}
\newcommand{\du}{{\rm du}}
\newcommand{\Av}{{\rm Av}}

\newcommand{\vwt}{{\rm vwt}}
\newcommand{\ntr}{{\rm ntr}}

\newcommand{\mxK}{{\rm mxK}}
\newcommand{\Kind}{{\rm Kind}}
\newcommand{\supp}{{\rm supp}}
\newcommand{\ord}{{\rm ord}}
\newcommand{\otp}{{\rm otp}}

\newcommand{\loc}{{\rm loc}}
\newcommand{\Mod}{{\rm Mod}}
\newcommand{\Levy}{{\rm Levy}}
\newcommand{\comp}{{\rm comp}}
\newcommand{\Reg}{{\rm Reg}}
\newcommand{\Def}{{\rm Def}}

\newcommand{\Th}{{\rm Th}}
\newcommand{\tp}{{\rm tp}}

\newcommand{\EM}{{\rm EM}}
\newcommand{\ai}{{\rm ai}}

\newcommand{\vk}{{\rm vk}}

\newcommand{\rK}{{\rm rK}}
\newcommand{\fs}{{\rm fs}}
\newcommand{\qf}{{\rm qf}}
\newcommand{\tr}{{\rm tr}}
\newcommand{\vK}{{\rm vK}}
\newcommand{\iif}{{\rm if}}
\newcommand{\aK}{{\rm aK}}
\newcommand{\eK}{{\rm eK}}
\newcommand{\tK}{{\rm tK}}

\newcommand{\sK}{{\rm sK}}

\newcommand{\uK}{{\rm uK}}

\newcommand{\qK}{{\rm qK}}

\newcommand{\vwa}{{\rm vwa}}
\newcommand{\nsp}{{\rm nsp}}
\newcommand{\wat}{{\rm wat}}
\newcommand{\uf}{{\rm uf}}
\newcommand{\reg}{{\rm reg}}

\newcommand{\ve}{{\rm v}}

\newcommand{\lin}{{\rm lin}}
\newcommand{\ind}{{\rm ind}}

\newcommand{\pcf}{{\rm pcf}}

\newcommand{\leeft}{{\rm left}}
\newcommand{\riight}{{\rm right}}

\newcommand{\lc}{{\rm lc}}
\newcommand{\pK}{{\rm pK}}
\newcommand{\IND}{{\rm IND}}
\newcommand{\Deef}{{\rm Def}}
\newcommand{\EC}{{\rm EC}}
\newcommand{\AP}{{\rm AP}}
\newcommand{\aut}{{\rm aut}}

\newcommand{\eq}{{\rm eq}}

\newcommand{\ged}{{\rm gd}}

\newcommand{\id}{{\rm id}}

\newcommand{\deef}{{\rm def}}

\newcommand{\Dom}{{\rm Dom}}
\newcommand{\Rang}{{\rm Rang}}

\newcommand{\rest}{{\restriction}}

\newcommand{\set}{{\rm set}}
\newcommand{\wilog}{{\rm without loss of generality}}
\newcommand{\Wilog}{{\rm Without loss of generality}}

\newcommand{\then}{{\underline{then}}}
\newcommand{\when}{{\underline{when}}}
\newcommand{\Then}{{\underline{Then}}}

\newcommand{\Iff}{{\underline{iff}}}
\newcommand{\mn}{{\medskip\noindent}}
\newcommand{\sn}{{\smallskip\noindent}}
\newcommand{\bn}{{\bigskip\noindent}}

\newcommand{\bbB}{{\mathbb B}}
\newcommand{\bbI}{{\mathbb I}}

\newcommand{\ga}{{\mathfrak a}}
\newcommand{\gk}{{\mathfrak k}}
\newcommand{\gC}{{\mathfrak C}}

\newcommand{\cC}{{\mathscr C}}
\newcommand{\bbC}{{\mathbb C}}

\newcommand{\bbD}{{\mathbb D}}
\newcommand{\cD}{{\mathscr D}}
\newcommand{\cW}{{\mathscr W}}
 
\newcommand{\gS}{{\mathfrak S\/}} 

\newcommand{\cH}{{\mathscr H}}

\newcommand{\cF}{{\mathscr F}}
\newcommand{\cE}{{\mathscr E}}

\newcommand{\bbL}{{\mathbb L}}
\newcommand{\bbN}{{\mathbb N}}

\newcommand{\gn}{{\mathfrak n}}
\newcommand{\cN}{{\mathscr N}}

\newcommand{\cP}{{\mathscr P}}

\newcommand{\bbQ}{{\mathbb Q}}
\newcommand{\bbR}{{\mathbb R}}
\newcommand{\bbZ}{{\mathbb Z}}

\newcommand{\cS}{{\mathscr S}}
\newcommand{\cT}{{\mathscr T}}
 
\newcommand{\gx}{{\mathfrak x}} 
\newcommand{\cU}{{\mathscr U}}
\newcommand{\cX}{{\mathscr X}}

\newcommand{\cf}{{\rm cf}}

\newcount\skewfactor
\def\mathunderaccent#1#2 {\let\theaccent#1\skewfactor#2
\mathpalette\putaccentunder}
\def\putaccentunder#1#2{\oalign{$#1#2$\crcr\hidewidth
\vbox to.2ex{\hbox{$#1\skew\skewfactor\theaccent{}$}\vss}\hidewidth}}

\newenvironment{PROOF}[2][\proofname.]
   {\begin{proof}[#1]\renewcommand{\qedsymbol}{\textsquare$_{\rm #2}$}}
   {\end{proof}}

\begin{document}

\title {Dependent dreams: recounting types}
\author {Saharon Shelah}
\address{Einstein Institute of Mathematics\\
Edmond J. Safra Campus, Givat Ram\\
The Hebrew University of Jerusalem\\
Jerusalem, 91904, Israel\\
 and \\
 Department of Mathematics\\
 Hill Center - Busch Campus \\ 
 Rutgers, The State University of New Jersey \\
 110 Frelinghuysen Road \\
 Piscataway, NJ 08854-8019 USA}
\email{shelah@math.huji.ac.il}
\urladdr{http://shelah.logic.at}
\thanks{Comes from F951, first typed 08/Sept/8
\newline
The author thanks Alice Leonhardt for the beautiful typing.  Partially
supported by the BSF and the NSF. Publication 950.}



\date{February 24, 2012}

\begin{abstract}
We investigate the class of models of a general dependent
theory.  We continue \cite{Sh:900} in particular investigating the so
called ``decomposition of types"; our thesis is that
what holds for stable theory and for Th$(\bbQ,<)$ hold for dependent
theories.  Another way to say this is: we have to look at small enough
neighborhood and use reasonably definable types to analyze general
types; also we presently concentrate on complete types over
saturated models (and sometimes just quite saturated models).
We now mention the main results understandable without reading the
paper.  First, a parallel to the ``stability spectrum", we consider the
``(problem of) recounting of types", that is assume 
$\lambda = \lambda^{< \lambda}$
is large enough, $M$ a saturated model of $T$ of cardinality
$\lambda$, let $\gS_{\aut}(M)$ be the set of complete
types over $M$ up to being conjugate, i.e. we identify $p,q$ when
some automorphism of $M$ maps $p$ to $q$.  Whereas for independent $T$
usually the number is 
$2^\lambda$, for dependent $T$ the number is $\le \lambda$
moreover it is $\le |\alpha|^{|T|}$ when $\lambda = \aleph_\alpha$ and
$\lambda$ is not too small, see \S(5B).  Second, for stable $T$,
recall that a model is
$\kappa$-saturated iff it is $\aleph_\varepsilon$-saturated and every
infinite indiscernible set (of elements) of cardinality $< \kappa$ is
not $\subseteq$-maximal.  We prove here an analog in \S(7B).  Third, 
if $M$ is saturated and $p \in \bold S(M)$ 
then $p$ is the average of an indiscernible sequence of length $\|M\|$ inside
the model, see \S(6A).  Fourth, we prove a (weak) relative 
of the existence of indiscernibles, see \S(4A).
Lastly, the so-called generic pair conjecture was proved
in \cite{Sh:900} for $\kappa$ measurable, here it is essentially
proved, i.e. for $\kappa = \kappa^{< \kappa} > |T| + \beth_\omega$,
see \S(7A).
\end{abstract}

\maketitle
\numberwithin{equation}{section}
\setcounter{section}{-1}

\section*{Anotated Content} 
\bn
\S0 \quad Introduction, pg.5
\bigskip

\S(0A) \quad What is done here, pg.5
\bigskip

\S(0B) \quad From higher perspective: the test, pg.7
\bigskip

\S(0C) \quad Basic definition, pg.9
\bigskip

\noindent
\S1 \quad Presenting Questions, Definitions and facts,  pg.11
\bigskip

\S(1A) \quad Recounting type, pg.11
\mn
\begin{enumerate}
\item[${{}}$]   [We count complete types on saturated models up to
conjugacy (= being automorphic), and ask what is 
the spectrum under this definition.]
\end{enumerate}
\bigskip

\S(1B) \quad On the outside definable sets and uf$(p)$, pg.18
\mn
\begin{enumerate}
\item[${{}}$]   [In particular we define bounded/medium/large
directionality; see more in \S(1C), \S(6B).  This is continued in
Kaplan-Shelah \cite{KpSh:946}.]
\end{enumerate}
\bigskip

\S(1C) \quad Indiscernibility, pg.22
\bigskip

\S(1D) \quad Limit models and the generic pair conjecture, pg.26
\mn
\begin{enumerate}
\item[${{}}$]   [For stable $T$ there is a neat characterization of
$\kappa$-saturated models.  What can we say for dependent $T$?  We
also present the generic pair conjecture and comment on
$(\kappa,2)$-limit models.]
\end{enumerate}
\bn
\S2 \quad Decomposition of types, pg.29
\bigskip

\S(2A) \quad Decompositions - the basics, p.29
\mn
\begin{enumerate}
\item[${{}}$]   [We suggest to try to analyze types, 
i.e.  $\tp(\bar d,M)$ where $M$ is $\kappa$-saturated and
$\bar d \in {}^{\theta^+ >}{\gC}$ via decompositions $\bold x \in 
\pK_{\kappa,\mu,\theta},\qK_{\kappa,\mu,\theta}$ and
relevant so-called solutions.]  
\end{enumerate}
\bigskip

\S(2B) \quad Smoothness and $(\bar\mu,\theta)$-sets, pg.37
\mn
\begin{enumerate}
\item[${{}}$]   [For $T$ dependent some sets $A \subseteq \gC_T$
behave like sets $A \subseteq \gC_{T'},T'$ stable, those are the
(smooth) $(\bar \mu,\kappa)$-sets.  Now for enough $\bold x \in \text{
pK}_{\mu,\bar\kappa,\theta}$, i.e. the smoother ones, 
the set $B^+_{\bold x}$ are such sets; so
it is natural that they interest us here.]
\end{enumerate}
\bigskip

\S(2C) \quad Measuring non-solvability and reducts, p.40
\mn
\begin{enumerate}
\item[${{}}$]   [We define $\ntr(\bold x)$ and relatives, a cardinal
measuring how badly solvability fails.  This will help in proving
density of $\tK$ later.]
\end{enumerate}
\bigskip

\noindent
\S3 \quad Strong Analysis: pg.43
\bigskip

\S(3A) \quad Introducing $\rK,\tK,\vK$, pg.43
\mn
\begin{enumerate}
\item[${{}}$]   [We define $\tK_{\kappa,\bar\mu,\theta}$, the most
desriable decomposition and define $\rK_{\kappa,\bar\mu,\theta}$, an
approximation to it and see what is nice about having 
few smooth decompositions.  
We further define $\uK,\vK$ which are 
``poor" relatives of $\qK,\tK$ respectively.  They are needed
as we succeed to prove density for $\vK$ essentially when $\kappa >
\beth_\omega$ and for $\tK$ only in some cases.
We give relevant definitions and basic facts, in particular about 
$\tK_{\kappa,\bar\mu,\theta}$ and review sufficient conditions
for indiscernibility.]
\end{enumerate}
\bigskip

\S(3B) \quad Sequence homogeneity and indiscernibles, pg.48
\mn
\begin{enumerate}
\item[${{}}$]   [Why is $\tK_{\kappa,\bar\mu,\theta}$ so desirable?
First, we prove that for a decomposition $\bold x \in 
\tK_{\kappa,\bar\mu,\theta}$, the model $M_{[\bold x]}$,
i.e. $M_{\bold x}$ expanded by individual constants for every $b \in
B^+_{\bold x}$ and coding the type of $\bar d_{\bold x} \char 94 \bar
c_{\bold x}$, is a $\kappa$-sequence homogeneous model.  We further show it
even for $\bold x \in \vK_{\kappa,\bar\mu,\theta}$; now
$\kappa$-sequence homogeneity implies uniqueness, so this points
the way to their uses in showing that for $M \in \EC_{\kappa,\kappa}(T)$, the
number of $p \in \bold S^\theta(M)$ up to conjugacy is small.

Second, we give some sufficient conditions for indiscernibility
related to $\tK$ and $\vK$.  We shall use some later]
\end{enumerate}
\bigskip

\S(3C) \quad Toward Density of $\tK$, pg.59
\mn
\begin{enumerate}
\item[${{}}$]   [To help prove the density we give sufficient
conditions when the union of an increasing sequence of decompositions
from $\rK_{\kappa,\bar\mu,\theta}$ belongs to
$\tK_{\kappa,\bar\mu,\theta}$ or to $\vK_{\kappa,\bar\mu,\theta}$.  For
the case $\kappa = \mu_2 = \mu_1 = \mu_0$ is a weakly compact cardinal
we prove the $\le_1$-density of $\tK_{\kappa,\bar\mu,\theta}$ in
$\pK_{\kappa,\bar\mu,\theta}$.] 
\end{enumerate}
\bigskip

\noindent
\S4 \quad Density, pg.64
\bigskip

\S(4A) \quad Partition theorems for dependent $T$, pg.64
\mn
\begin{enumerate}
\item[${{}}$]  [We prove two polarized partition theorems, 
showing dependency of $T$
has meaningful implications in this direction; 
they can be looked at as a substitute of the
existence of theorems of indiscernible sets for stable $T$.]
\end{enumerate}
\bigskip

\S(4B) \quad Density of $\tK$ in ZFC occurs, pg.70
\mn
\begin{enumerate}
\item[${{}}$]   [Our goal here is to prove the density of 
$\tK_{\kappa,\kappa,\theta}$ when $\kappa = \mu^+,\mu$ is singular strong limit
of cofinality $> \theta$ and $T$ is countable; also when $\kappa$ is
strongly inaccessible.  For this we prove that
some $\bar e$ universally solves $\bold x \in 
\qK_{\kappa,\kappa,\theta}$.  A crucial point is that instead
of using ``$\kappa$ weakly compact" (as in \S(3C)) we use a
partition theorem for dependent $T$ from \S(4A).]
\end{enumerate}
\bigskip

\noindent
\S5 \quad Stronger Density, pg.80
\bigskip

\S(5A) \quad More density of $\tK$, pg.80
\mn
\begin{enumerate}
\item[${{}}$]    [We prove the density of 
$\tK_{\kappa,\mu,\theta}$ when $\mu$ is as in \S(4A) and $\kappa <
\mu^{+\omega}$, see \ref{p2} so under G.C.H. we can prove a weak
version of the recounting types: for $M \in \EC_{\kappa,\kappa}(T)$
there are $\le \kappa$.  For this we use the partition theorem from \S(4A).
Note that under full G.C.H., this covers all large enough (regular) cardinals
$\kappa$ but only for $\mu$ close enough to $\kappa$.  So the
conclusion concerning the recounting of types are weak; still this 
proves a strong distinction between dependent $T$ and independent $T$.]
\end{enumerate}
\bigskip

\S(5B) \quad Density of $\vK$: Exact recounting of types and $\vK$, pg.88
\mn
\begin{enumerate}
\item[${{}}$]  [Here we use $\vK$, which was only a burden so far.  In the
relevant cases we prove its density (in $\pK$) and conclude the right 
number of types up to conjugacy (for $\kappa = 
\kappa^{< \kappa}$ large enough, $M \in \EC_{\kappa,\kappa}(T)$).]
\end{enumerate}
\bigskip

\S(5C) \quad Exact recounting of types and $\vk$, pg.95
\bigskip

\noindent
\S6 \quad Indiscernibles, pg.98
\bigskip

\S(6A) \quad Indiscernibles materializing $\bold m$, pg.98
\bigskip

\S(6B) \quad Indiscernibles existence from bounded directionality, pg.102
\bigskip

\noindent
\S7 \quad Applications, pg.107
\bigskip

\S(7A) \quad The generic pair conjecture, uniqueness of
$(\kappa,\sigma)$-limit model, pg.107
\mn
\begin{enumerate}
\item[${{}}$]  [Note that the case $(\kappa,2)$ is the generic pair.
We prove it for $\lambda > \beth_\omega + |T|$.]
\end{enumerate}
\bigskip

\S(7B) \quad Criterion for saturativity, pg.109
\bigskip

\noindent
\S8 \quad Concluding Remarks, pg.112
\newpage

\section{Introduction}

\subsection {What is done here} \
\bigskip

This is a step in trying to understand a dependent elementary class
$\Mod_T$.  Our approach is:

\begin{thesis}
\label{y3}  
1) It is fruitful to prove that questions on (first order complete) $T$ and a
cardinal does not depend too much on the cardinal, by finding
syntactical equivalent condition; this suggests it is an interesting
dividing line.

\noindent
2) We should first analyze saturated models (then quite saturated
models and only then general models).

\noindent
3) In particular we should first try to understand complete types 
over saturated models, etc.
\end{thesis}

More specifically:
\begin{thesis}
\label{y5}  
For $M \in \EC_{\kappa,\kappa}(T)$ we shall try
to analyze $p \in \bold S^\varepsilon(M)$ by types of two simple
kinds:
\medskip

\noindent
\underline{Kind A}:  $\Av(D,M),D$ an ultrafilter on ${}^\varepsilon B$ for
some $B \subseteq M$ of cardinality $< \mu$ ($\mu$ a fix cardinal
$\ll \|M\|$).
\medskip

\noindent
\underline{Kind B}:  $\Av(\bold I,M)$ where 
$\bold I = \langle \bar a_\alpha:\alpha <
\lambda\rangle$ an indiscernible sequence (of $\varepsilon$-tuples) inside $M$.
\end{thesis}

\begin{remark}  For stable $T$, if $M$ is $|T|^+$-saturated then
every $p \in \bold S(M)$ is $\Av(\bold I,M)$ for some indiscernible
sequence (so set) $\bold I$ of cardinality $\aleph_0$, so it 
falls under both kinds.

Consider a fixed complete first order theory $T$ which is dependent.
The problem we try to address here is
analyzing a complete type over a saturated model, say $p \in \bold
S^{< \theta^+}(M)$ where $\theta \ge |T|$.  The reader may wonder 
why not $p \in \bold S^{< \omega}(M)$?   The reason is that 
anyhow we are driven to consider infinitely many variables.  

Trying to analyze $p \in \bold S^\theta(M),M \in
\EC_{\kappa,\kappa}(T)$, clearly whatever occurs for some stable
theories may appear, so in the analysis we allow types definable over
small sets (though presently not stable types, just definable in a weak
sense) where any fix bound will be O.K. but as it happens
``small sets" mean a set of cardinality say $< \beth_\omega + |T|^+$.

Also in dense linear order there are cuts defined say by a
sequence of elements of length any regular $\sigma < \kappa$
(e.g. $p(x) \in \bold S(M)$ say that $x$ induces a cut of $M$ whose lower
half has cofinality $\sigma$), we cannot avoid
this so we allow types gotten as averages of indiscernible sequences
of length $\sigma$.  Note that types related to large cofinalities are
not covered by $\Kind \, A$, just as in \cite[\S1]{Sh:877}, where
the cuts with both cofinalities maximal are fine - there
expanding by them preserve saturation.

An approximation to analyzing $p$ is $\bold x \in 
\pK_{\kappa,\mu,\theta}$; a characteristic case is $\kappa =
\kappa^{< \kappa}$ large enough, $\theta = |T| = \aleph_0,\mu =
\beth_\omega$ (actually we use $\bar \mu$ but ignore it in the introduction).  
Now, see Definition \ref{b05}, such $\bold x$
consist of the
model $M= M_{\bold x}$, which is $\kappa$-saturated (and in general
may have larger cardinality), the sequence 
$\bar d = \bar d_{\bold x}$ realizing a complete type $p$ over $M$
which we are trying to analyze, $\bar c = 
(\ldots \char 94 \bar c_i
\char 94 \ldots)_{i \in v(\bold x)}$ an initial segment of the
analysis where $v(\bold x)$ is an ordinal $< \theta^+$ or just a linear
order of cardinality $\le \theta$.  
This means that for each $i \in v(\bold x)$ one of the
following two cases occurs, letting $r_{\bold x,i} = 
\tp(\bar c_i,M_{\bold x} + \Sigma\{\bar c_j:j < i\})$.

In the first case, formally $i \in v_{\bold x} \backslash u_{\bold
x}$, the type $r_{\bold x,i}$ does not split over some
$B_{\bold x,i} \in [M_{\bold x}]^{< \mu}$ (or even is finitely
satisfiable in it).  So
this type is in a suitable sense definable over some small set as in the
stable case, so is the ``stable part" called ``$\Kind \, A$" above.

In the second case, formally $i \in u_{\bold x}$ the type 
$r_{\bold x,i}$ is the
average of an indiscernible sequence ${\bold I}_{\bold x,i} =
\langle \bar a_{\bold x,i,\alpha}:\alpha < \kappa_i\rangle$ where
$\kappa_i = \cf(\kappa_i) \in [\mu,\kappa)$.  

In \cite{Sh:900}
some relatives were used but there $\mu = \kappa$ hence $B^+_{\bold x} =
\cup\{\bold I_{\bold x,i}:i \in u_{\bold x}\} \cup 
\bigcup\{B_{\bold x,i}:i \in v_{\bold x} \backslash u_{\bold x}\}$ here
corresponds to $B_{\bold x}$ there, so there the analysis is by
information of size just smaller than $\kappa$, whereas here it is by
$\le \theta$ indiscernible sequences of length a regular cardinal 
+ information of bounded size, i.e. $< \mu$ a major difference.

How does such $\bold x$ help?  For each $i \in v_{\bold x}$ we define when
$\bold x$ is active in $i$; it is the parallel of forking,
i.e. of ``$\tp(\bar d_{\bold x},M_{\bold x} + \Sigma\{\bar c_j:j \le
i\})$ forks over $M_{\bold x} + \Sigma\{\bar c_j:j < i\}$", this cannot
occur $\theta^+$ times so there is $\bold y$ above $\bold x$ maximal
in this sense; i.e. we cannot increase $v_{\bold x}$ having a ``new" activity
but not changing $M_{\bold x},\bar d_{\bold x},\bar c_i(i \in v_{\bold
x})$ but possibly increasing $v_{\bold x}$.  Moreover,
see \ref{b20}(2) we have further versions, local and/or less demanding, but
we skip this in the introduction.  The class
of maximal such $\bold y$'s is called $\qK'_{\kappa,\mu,\theta}$, see
Definition \ref{b16}(1); for them we can prove:
\mn
\begin{enumerate}
\item[$(*)$]   if $A \subseteq M_{\bold y},|A| < \mu$ \then \, some
$\bar e \in {}^\theta(M_{\bold x})$ solve $(\bold x,A)$ which means
that $\tp(\bar d_{\bold x},\bar c_{\bold x} + \bar e) \vdash 
\tp(\bar d_{\bold x},\bar c_{\bold x} + A)$ and even uniformly, which is
expressed by ``according to $\bar\psi$".
\end{enumerate}
\mn
This is the parallel of: if $M$ is a dense linear form, $p \in \bold
S(M),{\cC} = (C_1,C_2)$ the cut of the linear order $M$ 
which $p(x)$ induces and it has both
cofinalities $\ge \mu$ and $A \subseteq M,|A| < \mu$ then we can
choose $a \in C_1,b \in C_2$ such that $(a,b)_M \cap A = \emptyset$
hence $(a < x < b) \in p$ and $(a < x < b) \vdash p(x) \restriction A$.
\end{remark}

\noindent
All this seems to support:
\begin{thesis}
\label{y11}  
1) The theory of dependent elementary classes is the
combination of what occurs in stable classes and in the theory of dense
linear orders.

\noindent
2) We analyze general types by decompositions to three kinds: one are
   finitely satisfiable in a small set (or just does not split over a
   small set), second are averages of
   indiscernible sequences, third, are like branches of trees (include
   cuts of a linear order) any ``bounded" subset are 
implied by a very small subset.

But we really gain understanding by the density of
$\tK_{\kappa,\mu,\theta} \subseteq \pK_{\kappa,\mu,\theta}$
for some pair $(\kappa,\mu)$, (to cover all relevant cases better use 
$\vK^\otimes$, see \S3).
That is for $\bar d \in {}^\theta{\gC}$, we can find $\bold x \in
\tK_{\kappa,\mu,\theta}$ such that $\bar d \triangleleft \bar
d_{\bold x},M = M_{\bold x}$ and for every $A \subseteq M$ of
cardinality $< \kappa$ we can find $(\bar c',\bar d')$ in $M$ realizing
the same type as $(\bar c_{\bold x},\bar d_{\bold x})$ over $M$ and
$\tp(\bar d_{\bold x},\bar c_{\bold x} + \bar c' + d') \vdash 
\tp(\bar d_{\bold x},\bar c_{\bold x} + A + \bar c'_{\bold x} + \bar
d')$, even uniformly and fixing the type of $\bar c_{\bold x} \char 94
\bar d_{\bold x} \char 94 \bar c' \char 94 \bar d'$.  In a stronger sense
the type of $\bar c_{\bold x} \char 94 \bar d_{\bold x}$ over $M$
really combine parts definable over a small set and one like a
(partial) order.
\end{thesis}

\noindent
Another thesis is (see \cite[\S1]{Sh:783})
\begin{thesis}
\label{y13}  
In dependent (elementary) classes the family of outside definable
sets (Def$_{< \alpha}(M)$, see Definition \ref{a31}) replace the family of
inside definable sets for stable classes.
\end{thesis}
\bigskip

\noindent
\centerline {$ * \qquad * \qquad *$}
\bigskip

This work may be continued \cite{Sh:F973} and as said above it
continues \cite{Sh:900} though does not depend on it.   
More specifically, how are \cite{Sh:900} and the present work related?

In both cases decomposition ($\pK_{\kappa,\mu,\theta}$ here,
$K^1$ there) are central and $\qK',\qK$ here\footnote{In the context
  of \cite{Sh:900}, i.e. $\mu_0 = \kappa$ essentially we get $\qK' =
  \qK$, see \ref{b22}(3),(5)} are parallel to $\mxK$ there and
also $\le_1,\le_2$ are similar here and there.  In both cases the
model $M_{\bold x}$ is $\kappa$-saturated and $\bar d_{\bold x},\bar
c_{\bold x}$ are of cardinality $\le \theta$ (normally $\kappa >
\theta \ge |T|$).  But here we use $|B_{\bold x,i}| <
\mu$ and allow $\mu \ll \kappa$ rather than $|B_{\bold x,i}| < \kappa$, and
instead use indiscernible sequences $\bold I_{\bold x,i}$ for some
$i$'s.  Hence $B^+_{\bold x}$ 
here stands for $B_{\bold x}$ there, both have cardinality
$< \kappa$, but there $B_{\bold x}$ is any set, here \wilog \, $\bold
x$ is smooth so $B^+_{\bold x}$ is a so called 
$(\bar\mu,\theta)$-set, essentially $\le \theta$ sets each of 
cardinality $< \mu$
plus $\le \theta$ mutually indiscernible sets of $(\le
\theta)$-tuples.  Such sets have some affinity to stable ${\gC}$,
e.g. $|\bold S(B^+_{\bold x})| \le 2^{< \mu} + |B^+_{\bold x}|^{|T|}$.

Also $\tK_{\kappa,\mu,\theta}$ here is related to strict
decompositions in \cite{Sh:900}.  But in \cite{Sh:900} we get
existence assuming only $\kappa$ is a measurable cardinal so quite a 
large cardinal, so cannot prove in ZFC that it exists; 
whereas here this is proved for every large
enough regular cardinal provably in ZFC, and the bound is small (at
least for my taste), $\beth_\omega$, well $+ |T|$, of course.

All this is a good point in favor of large cardinals by the
criterion (first suggested by G\"odel): we can first prove things assuming
them, this helps us to find the way to really sort out things.
\bigskip

\subsection {From Higher Perspective: The Test} \
\bigskip

\noindent
What questions do we address here?
\begin{question}
\label{y31}
\underline{The serious/dull question}  
1) Is the equation dependent/stable = groups/Abelian groups true?

That is, is dependence a better dividing line than stable (among say
elementary classes), but we have been (and are) just too dim to see it?

\noindent
2) The use of cardinals $(> \aleph_0)$ in model theory: 
has it passed its time OR is it the key to dependent classes and will
continue to be central.
\end{question}

Alas, most (relevant) people already know the answers, unfortunately not all
of them know the same answer.

In more serious mode, we suggest here to put dependent theories to
``end of first level examination".  Trying to be objective we ask: do we have
a good analog to what is in the first paper on stable $T$, \cite{Sh:1} 
(and \cite{Sh:12}),
essentially equivalently at the time of stability being three years old.

So here is the test composed of four questions (as presented in a
lecture in MAMLS, Fall 2008 Meeting in honor of Gregory Cherlin) and a
fifth question (as urged by the audience):

\begin{question}
\label{y50}
\underline{Question/Test}  
Find parallels of (1)-(4) and answer (5) for dependent $T$.

\noindent
1)The stablility spectrum Theorem (for 
stable theory $T$  on a model of cardinality $\lambda$ there
are $ \le \lambda$ completer 1-types).

\noindent
2) Strong partition theorems, i.e. existence of 
indiscernibles: for stable $T$ , if $a_\alpha \in \gC$
for $\alpha < \lambda^+$ are given, $\lambda = \lambda^{|T|}$ then for some 
unbounded, even stationary subset $S$ of $\lambda^+$ the sequence
$\langle a_\alpha:\alpha \in S \rangle$ is indiscernible.

\noindent
3) ``Understanding" complete types over models and 
indiscernible sequences (for stable
$T$, the finite equivalence relation theorem which was somewhat later).

\noindent
4) Characterize saturated models by indiscernible sequences, 
(for stable  $T,M$ is $\kappa$ saturated iff it is 
$\aleph_\epsilon$-saturated and every infinite indiscernible
set of cardinality $< \kappa$).

\noindent
5) The generic pair conjecture,
a major question from \cite{Sh:900} and more generally the existence
of $(\lambda,\kappa)$-limit models ($\kappa = 2$ is the generic pair case).

We did not mention two problems having been answered earlier:
majority on indiscernibles (see \S(1C)) and definability of types 
(as we may consider the following theorem as an answer: 
expansion by outside definable sets preserved the theory of the model
being dependent, by \cite[\S1]{Sh:783}).

We will present the questions in \S1 and present solutions to (1),(4)
and the first part of (5) in \S5,\S7.  
Unfortunately we do not solve the original interpretation of questions
(2),(3) as we hoped, but, not surprisingly, we think we 
have excellent excuses.  Now the
answer to the parallel of (3) we considered, i.e. ``no case of high
directionality" that is bounding the number
of ultrafilters $D$ on $M$ such that $\Av(D,M)=p$, has already been 
known to be false for many years, proved by Delon.

As for the existence of indiscernibles, i.e. \ref{y50}(2) and actually
also (3), subsequently Kaplan-Shelah \cite{KpSh:946}, proved that the
premature assertion in the Rutgers lecture is false, this is 
the negative half of the excuse, i.e. this version cannot be proved
being false.

However on the positive side, we believe we have reasonable
substitutes, i.e. reasonable parallels of parts (2),(3) of \ref{y50}
for dependent $T$.

For part (3):
\mn
\begin{enumerate}
\item[$\boxplus_1$]  if $M \in \EM_{\kappa,\kappa}(T)$ and $p
\in \bold S(M)$ \then \, $p$ is the average of an indiscernible sequence
in $M$ of length $\kappa$, see \ref{e1}, (more in \S(6A) and the 
results of \S(6B)).
\end{enumerate}
\mn
About the existence of indiscernibles, i.e. part (2) of \ref{y50}, 
by \S6 we have
\mn
\begin{enumerate}
\item[$\boxplus_2$]  existence for $T$ with low or medium
directionality (introduced in \S(1B)).
\end{enumerate}
\end{question}

Probably this is not convincing:  but a true answer for \ref{y50}(2) 
is another relative (or you may say
a weak version) of the existence of indiscernibles
\medskip

\noindent
\begin{enumerate}
\item[$\boxplus_3$]  if $\kappa = \cf(\kappa) > \aleph_0$ and
$\Delta \subseteq \bbL(\tau_T)$ is finite and
$a_{\alpha,n} \in \gC$ for $\alpha < \kappa,n < n(*) < \omega$
\then \, we can find stationary $\mathscr{S}_n \subseteq \kappa$ for
$n < n(*)$ such that: for $\bar\alpha  \in \prod\limits_{\ell < n}
\cS_\ell$, the $\Delta$-type of $\langle
a_{\alpha_0},\dotsc,a_{\alpha_{n(*)}-1}\rangle$ depends just
on the truth values of $\alpha_{\ell(1)} < \alpha_{\ell(2)}$ for
$\ell(1),\ell(2) < n(*)$.
\end{enumerate}
\medskip

\noindent
This holds by \ref{d21}, (note that we can apply 
it for any permutation of
$\{\langle 0,\dotsc,n(*)-1\rangle\}$ and the formulation here is simpler
because we use the club filter on $\chi$, i.e. use
diagonal intersection of clubs).  Note that for $T$ any completion of
Peano arithmetic (or any 2-independent $T$) this holds only for (some)
large cardinal.

There has been work on dependent theories in the previous century, 
see e.g. in the introductions of \cite[\S1]{Sh:715}, \cite[\S0]{Sh:783}, \cite[\S0]{Sh:900}; there
was much activity in the first decade of the present century, 
but in different directions; on indiscernibility see \S(1C) here.
\bigskip

\subsection {Basic Definitions} \
\bigskip

We assume basic knowledge in model theory.

\begin{convention}
\label{z0} 
1) ${\gC} = {\gC}_T$ is a monster
model of the complete first order $T$.

\noindent
2) The vocabulary of $T$ is $\tau_T$.

\noindent
3) $\bbL(\tau)$ is the set of first order formulas in the vocabulary $\tau$.
\end{convention}

\begin{definition}
\label{z3}  
1) Let $\EC_{\lambda,\kappa}(T)$ be
the class of $\kappa$-saturated models of $T$ of cardinality
$\lambda$; if $\kappa =1$ this means that we omit the 
$\kappa$-saturation; we may omit $\kappa$ when $\kappa = \lambda$.

\noindent
2) Let $\bar x_{\bar a} = \langle x_{a_t}:t \in \Dom(\bar a)
\rangle$ where $\bar a \in {}^I{\gC}$ for some index set $I =
\Dom(\bar a)$, usually $I$ an ordinal.  Let 
$\bar x_{[\alpha]} = \langle x_\beta:\beta
 <\alpha\rangle$, similarly $\bar x_{[u]}$ for $u$ a set or linear
order.  Generally we allow infinite sequence of variables
   but the formulas are finitary so only finitely many variables are
mentioned. 

\noindent
2A) Let $\bar x'_{\bar a} = \langle x'_{a_t}:t \in \Dom(\bar a)
\rangle$, etc.; note $\bar x_{\bar a \rest u} = \bar x_{\bar a} \rest u$.

\noindent
2B) If $\eta \in {}^I \Dom(\bar a)$ 
then: $\bar x_{\bar a,\eta} = \langle x_{a_{\eta(s)}}:
s \in I\rangle$ and $\bar a_\eta = \langle a_{\eta(s)}:s \in 
\Dom(\eta)\rangle$; see \ref{q3}.

\noindent
2C) Let $\ell g(\bar a) = \Dom(\bar a)$.  
Note $\ell g(\bar x_{\bar a}) = \Dom(\bar x_{\bar a})$ and
$\ell g(\bar x_{[u]}) = u$.

\noindent
3) Let $\varphi(\bar x)$ be the pair $(\varphi,\bar x)$, where
\mn
\begin{enumerate}
\item[$\bullet$]  $\varphi$ is a first order formula (in
$\bbL(\tau_T),T$ the first order theory understood from the content
\sn
\item[$\bullet$]  $\bar x$ is a sequence without repetition of
variables, including all the variables occuring in $\varphi$ freely.
\end{enumerate}
\mn
We normally use $\varphi(\bar x,\bar y)$ as a different object
than $\varphi(\bar x \rest u,\bar y \rest v)$ and $\varphi$ may
   stand for such object, e.g. $\langle \psi_\varphi(\bar y,\bar
   z):\varphi = \varphi(\bar x,\bar y) \in \bbL(\tau_T)\rangle$. This is
   ambiguous in principle but clear in practice.  See more in
   Definition \ref{a6}(4).

\noindent
4) We may use $A+B$ instead of $A \cup B$ and $\sum\limits_{t \in I}
A_t$ for $\cup\{A_t:t \in I\}$.
\end{definition}

\begin{observation}
\label{z4}  
The number of formulas $\varphi(\bar x_{\bar c},
x_{\bar d}) \in \bbL(\tau_T)$ is $|T| + |\ell g(\bar c)| + |\ell
g(\bar d)|$ so $\ge \aleph_0$ and maybe $> |T|$.
\end{observation}

\begin{definition}
\label{z5}  
1) For $M \prec {\gC}$ and $B \subseteq {\gC}$ let $M_{[B]}$ be $M$ expanded
by relations definable in ${\gC}$ with parameters from $B$, as in
\cite[\S1]{Sh:783}. 

\noindent
2) Similarly $M_{[p(\bar x)]}$ for $p(\bar x) \in \bold S^\varepsilon(M)$.
\end{definition}

\begin{convention}
\label{z8}  E.g. saying ``$\bar c \char 94 \bar d$
realizes $\tp(\bar c_{\bold x} \char 94 \bar d_{\bold x},A)"$ we may
forget to say $\ell g(\bar c) = \ell g(\bar c_{\bold x}),\ell g(\bar d) =
\ell g(\bar d_{\bold x})$.
\end{convention}

\begin{notation}
\label{z23}  
1) $\tp_\varphi(\bar d,\bar c \dotplus A)$ for
$\varphi = \varphi(\bar x_{\bar d},\bar x_{\bar c},\bar y)$ is
$\{\varphi(\bar x_{\bar d},\bar c,\bar a):\bar a \in {}^{\ell g(\bar
y)}A$ and ${\gC} \models \varphi[\bar d,\bar c,\bar a]\}$.

\noindent
2) Similarly $\tp_\Delta(\bar d,\bar c \dotplus A)$ where $\Delta \subseteq
   \{\varphi:\varphi = \varphi(\bar x_{\bar d},\bar x_{\bar c},\bar y)
   \in \bbL(\tau_T)\}$.

\noindent
3) $\tp_{\pm \varphi}$ means $\tp_{\{\varphi,\neg \varphi\}}$.

\noindent
4) Let $\Gamma_{[\zeta]} = \{\varphi:\varphi = \varphi(\bar
   x_{[\zeta]},\bar y) \in \bbL(\tau_T)\}$; similarly
$\Gamma_{[\zeta],n} = \{\varphi:\varphi = \varphi(\bar
   x^0_{[\zeta]},\dotsc,\bar x^{n-1}_{[\zeta]},\bar y)\}$.

\noindent
5) Let $(\forall^\kappa t \in I)\vartheta(t)$ means: for all but $<
   \kappa$ members $t \in I$ we have $\vartheta(t)$ (but may use
   $(\forall^\infty n)$ instead $(\forall^{\aleph_0} n \in
   \bbN)\vartheta(n))$.  Similarly $(\exists^\kappa t\in I)$ means: there
   are $\ge \kappa$ members $t$ of $I$ such that $\vartheta(t)$.
\end{notation}

\begin{definition}
\label{z25}
1) We say that a model $M$ is a $\kappa$-sequence homogeneous \when
   \,: if $f$ is a partial one-to-one function from $M$ to $M$ of
   cardinality $< \kappa$, i.e. $|\Dom(f)| < \kappa$ and $f$
is elementary in $M$ t hen: for every $a \in M$ for some $b \in M$ the
   function $f' = f \cup \{\langle a,b\rangle\}$ is elementary in $M$,
   where

\noindent
1A) We say the function $f$ is elementary in $M$ when: $\Def(f)
\subseteq M,\Rang(f)$ and if $M \models \varphi[a_0,\ldots]$ and
   $a_0,\ldots \in \Dom(f)$ then $M \models \varphi[f(a_0),\ldots]$.

\noindent
2) We say that a model $M$ is strongly $\kappa$-sequence homogeneous
   \when \,: if $f$ is as in part (1) then $f$ can be extended to an
   automorphism of $M$.

\noindent
3) We say that a model $M$ is strongly $\kappa$-saturated \when \, $M$
   is $\kappa$-saturated and strongly $\kappa$-sequence homogeneous.
\end{definition}

\begin{convention}
\label{z27}  
1) Generally (i.e. from \S2 on if not said otherwise) in this
work, $I$ vary on $K_{\lin}$, the class of linear orders which are endless.
\end{convention}
\newpage

\section {Presenting questions, definitions and facts}

We here recall and make some definitions and questions related to the
family of dependent theories and say some easy things to clarify,
mostly those questions are dealt with later in this work.
\bigskip

\subsection {Recounting types} \
\bigskip

We define the new version of the number of types, i.e. up to
automorphisms, considering saturated model and generalizations.  We
then have a ``first look at them".  First, about the function
$f^{\aut}_T$, counting the types up to automorphisms, see Definition \ref{a3}:
\mn
\begin{enumerate}
\item[$\boxplus$]  $(a) \quad$ if $T$ is stable, the function
$f^{\aut}_T(\lambda)$ is constant, $\le 2^{|T|}$, if $T$ is countable

\hskip25pt the constant value belongs to $\{2,3,\ldots\} \cup
\{\aleph_0,2^{\aleph_0}\}$, see \ref{a7}(1),(2)
\sn
\item[${{}}$]  $(b) \quad$ in (a), for countable $T$ every one of the 
values occurs even for 

\hskip25pt superstable $T$, see \ref{a6}
\sn
\item[${{}}$]  $(c) \quad$ in (a), if $T$ is $\aleph_0$-stable then except
$2^{\aleph_0}$ every one of the values occurs
\sn 
\item[${{}}$]  $(d)(\alpha) \quad$ if $T$ is independent then
$f^{\aut}_T(\lambda) = 2^\lambda$ when $(\exists \mu)(\lambda =
\lambda^{< \lambda} =$

\hskip25pt $2^\mu > |T|)$, see \ref{a9}
\sn
\item[${{}}$]  $\qquad (\beta) \quad$ if $T$ is independent, $\lambda =
\lambda^{<\lambda} > |T|$ but not as in $(\alpha)$ then

\hskip25pt  still $f^{\aut}_T(\lambda) \ge \lambda$
\sn
\item[${{}}$]  $(e) \quad$ if $T$ is dependent and unstable then
$f^{\aut}_T(\aleph_\zeta) \ge |\zeta +1|$, see \ref{a7}(4),(5).
\end{enumerate}
\mn
This explains that the problem is about dependent (unstable) $T$.
Note that the case of independent $T$ and strongly inaccessible
$\lambda > |T|$ is not resolved here, see on it \cite{Sh:F1124}.

The rest of this subsection is devoted to looking at relatives of
$f^{\aut}_T$ motivated by a desire not to use instances of G.C.H.

\begin{definition}
\label{a3} 
1) Let $\bold C := \{\lambda:\lambda = 
\lambda^{< \lambda}\}$ and
$\bold C _{> \mu} = \bold C(> \mu)$ be $\bold C \backslash \mu^+$.

\noindent
2) For $T$ a complete first order theory and $\theta \ge 1$ we define
the function $f^{\aut}_{T,\theta}:\bold C \rightarrow$ Card by
$f^{\aut}_{T,\theta}(\lambda) = 
|\gS^{\theta}_{\aut}(M_\lambda)|$ for 
$M_\lambda \in \EC_{\lambda,\lambda}(T)$, i.e. 
a saturated model of $T$ of cardinality $\lambda$, where

\noindent
3) $\gS^\theta_{\aut}(M) = (\bold
S^\theta(M)/\equiv_{\aut})$ where $\equiv_{\aut}$ or
more\footnote{We can define also when $p_\ell \in 
\bold S^\theta(M_\ell)$ are equivalent = conjugate for 
$\ell=1,2$ as in \cite{Sh:E46} which deal in a non-first order but for a
stable context.}
fully $\equiv^{\aut}_M$ is the 
following equivalence relation: $p,q \in \bold S^\theta(M)$ are
$\equiv_{\aut}$- equivalent iff they are conjugate, i.e. 
there is an automorphism of $M$ mapping $p$ to $q$.

\noindent
4) If we omit $\theta$ we mean $\theta =1$, if we write 
``$< \aleph_0$" we mean ``for any finite $n > 0$".
\end{definition}

\begin{example}
\label{a6} 
1) Assume $T = \Th(\bbQ,<)$, the theory
of dense linear orders with neither first nor last element.  \Then \,
$f^{\aut}_T(\aleph_0)$ is equal to 6, yes, six.

\noindent
2) If $T = \Th(\bbC)$, or $T$ is the theory of some algbraically
closed field of characteristic $p,p$ prime or zero, \then \, 
$f^{\aut}_T(\lambda) = \aleph_0$, for $\lambda \ge \aleph_0$.

\noindent
3) In part (1), in general, 
$f^{\aut}_T(\aleph_\alpha) = 6 + 2|\alpha|$ for $\aleph_\alpha
\in \bold C$.  

\noindent
4) Let $\tau = \{P_i:i < \alpha\},P_i$ a unary predicate and $T$ says
 that each $P_i$ is infinite and $T$ says that each $P_i$ is infinite,
 they are pairwise disjoint, and if $\alpha$ is finite then 
$\{x:\bigwedge\limits_{i < \alpha} \neg P_i(x)\}$ is infinite.
Then $T$ is stable (even totally transcendental so $\aleph_0$-stable
 if $\alpha$ is countable) and $f^{\aut}_T(\lambda) = 2(|\alpha +1|)$ 
for $\lambda \ge \aleph_0 + |\alpha|$.  If $\alpha$ is finite $>0$ and 
$\beta \le \alpha$ and above we demand $P_\ell$ is a singleton when
 $\ell < \beta$, infinite when $\ell \ge \beta$ then we get 
$f^{\aut}_\tau(\lambda) = 2|\alpha - \beta| + |\beta| +2$.

\noindent
5) Let $T = \Th(M)$ where $M = ({}^\omega 2,P^M_n)_{n <
   \omega}$ and $P^M_n = \{\eta \in {}^\omega 2:\eta(n)= 1\}$ for $n
   \in \bbN$.  \Then \, $T$ is countable
   superstable and $f^{\aut}_\tau(\lambda) = 2^{\aleph_0}$
   for $\lambda \ge 2^{\aleph_0}$.

\noindent
6) Let $T = \Th({}^\omega \omega,E_n)_{n < \omega}$ where $E_n =
 \{(\eta,\nu):\eta,\nu \in {}^\omega \omega$ and $\eta \rest n = \nu
 \rest n\}$.  So $T$ is countable, stable not superstable and 
$f^{\aut}_T(\lambda) = 2,f^{\aut}_{T,2}(\lambda) = 
\aleph_0$ for every $\lambda$.
\end{example}

\begin{observation}
\label{a7}  
1) If $T$ is stable, 
\then \, $f^{\aut}_T(\lambda)$ is constant and is
$\le 2^{|T|}$ for every $\lambda \in
\bold C_{> |T|}$ (or just $T$ has a saturated model of cardinality
$\lambda$, e.g. $\lambda = \lambda^{|T|}$).  
Similarly $f^{\aut}_{T,\theta}(\lambda) \le 2^{|T|+\theta}$ and is constant. 

\noindent
2) If $T$ is countable and stable and
e.g. $\lambda = \lambda^{\aleph_0}$ \then \,
$f^{\aut}_T(\lambda)$, is countable so 
$\le \aleph_0$ or is constantly $2^{\aleph_0}$.

\noindent
3) If $T$ is $\aleph_0$-stable \then \, $f^{\aut}_T(\lambda) \le
   \aleph_0$.

\noindent
4) If $T$ is unstable and is dependent, 
\then \, $f^{\aut}_T(\aleph_\zeta) \ge |\zeta + 1|$ for
$\aleph_\zeta \in \bold C$ which is $> |T|$.

\noindent
5) If $T$ is independent, $\lambda > |T|$ is inaccessible then
   $f^{\aut}_T(\lambda) \ge \lambda$.
\end{observation}

\begin{PROOF}{\ref{a7}}
  1) Assume $M$ is saturated of cardinality $>|T|$ or just
a strongly $|T|^+$-sequence homogeneous (see Definition \ref{z25}).  
Every $p \in \bold S^m(M_\lambda)$ is
definable, in fact there is a sequence $\langle \psi_\varphi(\bar y,\bar
z_\varphi):\varphi = \varphi(\bar x,\bar y) \in \bbL(\tau_T)\rangle$
so $\bar x = \bar x_{[m]}$ such that for every $p \in \bold S^m(M)$ there is a
sequence $\bar c^p := \langle \bar c^p_\varphi:\varphi = 
\varphi(\bar x,\bar y) \in \bbL(\tau_T)\rangle$, of sequences from
$M$ such that $\ell g(\bar c^p_\varphi) = \ell g(\bar z_\varphi)$ and
$\varphi(\bar x,\bar b) \in p$ iff $\bar b \in
{}^{\ell g(\bar y)}M$ and $M \models 
\psi_\varphi[\bar b,\bar c^p_\varphi]$, see \cite[Ch.II]{Sh:a}.  
Now the number of complete types of sequences of the form
$\langle \bar c_\varphi:\varphi = \varphi(\bar x,\bar y) 
\in \bbL(\tau_T)\rangle$ in $M$ with $\bar c^p_\varphi \in
{}^{\ell g(\bar z_\varphi)}M$ is $\le 2^{|T|}$.  But $M$ is
strongly $|T|^+$-sequence homogeneous, see Definition \ref{z25}(3), 
so this piece of information suffices, that
is, if $p,q \in \bold S^m(M)$ and tp$(\bar c^p,\emptyset,M) = 
\tp(\bar c^q,\emptyset,M)$ then there is an automorphism $f$ of $M$
which maps $\bar c^p$ to $\bar c^q$ hence $f$ maps $p$ to $q$.  
Of course, this works for $\gS^\zeta_{\aut}(M)$ too, only
the bound is $2^{|\zeta|+|T|}$, so for $\zeta \ge |T|$ we get even
equality.

\noindent
2) As in part (1), but this constant value
is the number of equivalence class of
some Borel relation hence by a theorem of Silver is $\le \aleph_0$
or is $2^{\aleph_0}$, see e.g. \cite{HrSh:152}, \cite{Sh:202}. 

\noindent
3) By the proof of part (1) and the definition of being
   $\aleph_0$-stable.

\noindent
4)  Recall $T$ has the strict order property (by \cite[Ch.II]{Sh:c}) 
hence some formula $\varphi(x,\bar y_n)$
has the strict order property.  We fix such $\varphi$; and any
$M \in \EC_{\aleph_\zeta,\aleph_\zeta}(T)$ for any 
regular $\kappa \le \aleph_\zeta$ we can find an
indiscernible sequence $\bold I_\kappa = \left< \langle
b_{\kappa,\alpha}\rangle \char 94 \bar a_{\kappa,\alpha}:\alpha <
 \kappa\right>$ in $M$ such that:
\mn
\begin{enumerate}
\item[$(*)$]   $(a) \quad {\gC} \models \varphi[b_{\kappa,\beta},
\bar a_{\kappa,\alpha}]$ iff $\alpha < \beta$
\sn
\item[${{}}$]   $(b) \quad \varphi(x,\bar a_{\kappa,\alpha}) \vdash \varphi(x,
\bar a_{\kappa,\beta})$ if $\alpha < \beta$.
\end{enumerate}
\mn
Let $p_\kappa = \Av(\langle b_{\kappa,\alpha}:\alpha
<\kappa\rangle,M)$, so it is enough to prove that for regular
$\kappa_1 \ne \kappa_2$, the types $p_{\kappa_1},p_{\kappa_2}$ are not
conjugate;p (\wilog \, $\kappa_1 < \kappa$).  
For this it is enough to prove $p_{\kappa_1} \ne
p_{\kappa_2}$ (as the assumptions in the choice of $\bold
I_\kappa,p_\kappa$ are preserved by automorphisms of $M$).  Toward
contradiction assume $p_{\kappa_1} = p = p_{\kappa_2}$ and \wilog \,
$\kappa_1 < \kappa_2$.  For $\ell=1,2$ we have
$\alpha < \kappa_\ell \Rightarrow (\forall^{\kappa_\ell} \beta <
\kappa_\ell)\varphi(b_{\kappa_\ell,\beta},\bar a_{\kappa_\ell,\alpha})
\Rightarrow \varphi(x,\bar a_{\kappa_\ell,\alpha}) \in p \Rightarrow
(\forall^{\kappa_{3 - \ell}} \beta < \kappa_{3 - \ell})
\varphi(b_{\kappa_{3-\ell},\beta},
\bar a_{\kappa_\ell,\alpha})$ so applying this to $3 - \ell$ we
have $\alpha < \kappa_{3 -\ell} \Rightarrow 
(\forall^{\kappa_\ell}\beta)[\varphi(b_{\kappa_\ell,\beta},
\bar a_{\kappa_{\ell - \alpha,\alpha}})]$.  So necessarily 
there are an unbounded $u \subseteq \kappa_2$ and
$\alpha_* < \kappa_1$ such that $\alpha_* \le \alpha < \kappa \wedge
\beta \in u \Rightarrow M \models \varphi[b_{\kappa_2,\beta},\bar
a_{\kappa_1,\alpha}]$.  Renaming, \wilog \, $\alpha_* = 0,u = \kappa_2$.

First, assume $\kappa_2 < \aleph_\zeta$.  Let $q(\bar y) = 
\{\neg \varphi(b_{\kappa_1,i},\bar y):i < \kappa_1\} \cup
\{\varphi(b_{\kappa_2,j},\bar y):j < \kappa_2\}$.  If $\bar a \in
{}^{\ell g(\bar y)}M$ realizes
$q(\bar y)$ we get $\neg \varphi(x,\bar a) \in p_{\kappa_1} \wedge
\varphi(x,\bar a) \in p_{\kappa_2}$, contradiction.

But if $r \subseteq q(\bar y)$ is finite and $i_* = \sup\{i:b_{\kappa_1,i}$
appear in $r\}$ then $\bar a_{\kappa_1,i_* +1}$ realizes $r$ so
$q(\bar y)$ is a type in $M$ but we are
assuming $|q|= \kappa_2 < \aleph_\zeta$ and $M$ is saturated 
so $q$ is realized in $M$, contradiction.  

Second, assume $\kappa_2 = \aleph_\zeta$; we could have chosen $p_{\kappa_2}$
using a linear order $I = I'_2 + I''_2$, isomorphic to
$(\kappa_2 + \kappa^*_2)$ such that $I'_2 =\{s_\alpha:\alpha <
\kappa\},I''_\alpha = \{t_\alpha:\alpha < \kappa\}$ and $\alpha <
\beta < \kappa_2 \Rightarrow s_\alpha <_I s_\beta <_I t_\beta <_I t_\alpha$.

We choose $\langle b_s,\bar a_s:s \in I\rangle$ in $M$ such that $M
\models \varphi[b_t,\bar a_s]$ if $s <_I t$.  Also \wilog \, for every
$A \subseteq M$ of cardinality $< \aleph_\zeta$ for some $\alpha <
\aleph_\zeta$ the set $\langle \bar b_s \char 94 \bar a_s:s \in
\{s_\beta,t_\beta\}:\beta \in (\alpha,\aleph_\zeta)\rangle$ is
indiscernible over $A$.  

Lastly, \wilog \, $\beta < \kappa_2 \Rightarrow 
b_{s_\beta} = b_{\aleph_\zeta,\beta}$ so 
$p_{\kappa_2}(\lambda) = \{\psi(x,\bar c):\bar c \subseteq M$
and $M \models \psi[b_{s_\alpha},\bar c]$ for every $\alpha < \kappa_2$
large enough$\}$.  Now for any $\alpha < \kappa_2$ we have
$[\varphi(x,\bar a_{s_\alpha}) \wedge \neg
\varphi(x,\bar a_{t_\alpha})) \in p_{\kappa_2}$ hence
for some $\gamma(\alpha) = \gamma_\alpha < \kappa_1$ we have ${\gC}
\models \varphi[b_{\kappa_1,\gamma(\alpha)},\bar a_{s_\alpha}]
\wedge \neg \varphi[b_{\kappa_1,\gamma(\alpha)},\bar a_{t_\alpha}]$ so
for some $\gamma < \kappa_1$ the set $u = \{\beta <
\kappa_2:\gamma(\beta) = \gamma\}$ is unbounded in $\kappa_2 =
\aleph_\zeta$.  So choose above $A = \{b_{\kappa_1,\gamma}\}$ and get
a contradiction.

\noindent
5) See more in \cite{Sh:F1124}, still we state \ref{a9} below.  
\end{PROOF}

\begin{observation}
\label{a9}  Assume $T$ is independent, \then \,:

$f^{\aut}_T(\lambda) = 2^\lambda$ for $\lambda = 2^\mu \in \bold C_{> |T|}$
\end{observation}

\begin{PROOF}{\ref{a9}}
Because there are $M_0 \in \EC_{\lambda,1}(T)$ 
such that $A \subseteq M_0,|A|=\mu$ such that $\cP
= \{p \in \bold S(M):p$ finitely satisfiable in $A\}$ has cardinality
$2^\lambda$, but $\cP_q = \{p \in \cP,p$ conjugate to $q\}$ has
cardinality $\le \lambda^\mu = \lambda$ for each $q \in \cP$.
\end{PROOF}

Dealing with saturated models, for unstable $T$, force us to have the
suitable cardinality with $(\kappa = \kappa^{< \kappa})$! so our
restriction to such cardinals is natural, that is recall
\begin{claim}
\label{a10}  
If $M \in \EC_{\kappa,\kappa}(T)$ but $T$ is
unstable and $\kappa > \aleph_0$ \then \, $\kappa = \kappa^{< \kappa}$.
\end{claim}
\bigskip

\begin{PROOF}{\ref{a10}}  
By \cite[Ch.III]{Sh:c}.  
\end{PROOF}

\begin{conjecture}
\label{a11} 
1) If $T$ is dependent, \then \,
$f^{\aut}_T(\aleph_\alpha) \le |\alpha|^{2^{|T|}}$ for
$\aleph_\alpha \in \bold C$.

\noindent
2) If $T$ is dependent unstable, \then \, for some $\kappa^+(T) \le
|T|^+$ we have $f_T(\aleph_\alpha) = |\alpha|^{<\kappa^+(|T|)}$ when
   $\aleph_\alpha \in \bold C$ is large enough (see \cite[Ch.III]{Sh:c}
on number of independent orders).
\end{conjecture}

\begin{discussion}
\label{a13} 
1) During a try to improve \cite{Sh:900}, raising this Conjecture
changes my outlook and leads to this work.

\noindent
2) We may like to eliminate the use of G.C.H. or weak relatives,
though \ref{a10} show this is not straight.  We may consider the 
following relatives,
$f^{\wat}_{T,\theta}(-)$ and $f^{\vwa}_{T,\theta}(-)$,
those are not further dealt with in this work, i.e. after \S(1A).  
\end{discussion}

\begin{definition}
\label{a16} 
1) For $\lambda \ge |T|$ let $f^{\wat}_{T,\theta}(\lambda) =  
\min\{\mu$: for every $M \prec {\gC}$ of cardinality 
$\lambda$ there is $N \prec {\gC}$ of  
cardinality $\lambda$ extending $M$ such that 
$|\gS_{\aut}(N)| \le \mu\}$.

\noindent
2) Let $f^{\vwa}_{T,\theta}(\lambda) = \min\{\mu$: for
every $M \prec {\gC}$ of cardinality $\lambda$ there is $N \prec 
{\gC}$ of cardinality $\lambda$ extending $M$ and function $g:
\bold S(M) \rightarrow \bold S(N)$ such that $p \in \bold S(M)
\Rightarrow g(p) \rest M = p$ and $|\{g(p)/\equiv^{\aut}_N:p
\in \bold S(M)\}|\le \mu\}$; so $f^{\vwa}_\lambda(T) \le f^{\wat}_T(\lambda)$. 

\noindent
3) Omitting $\theta$ means $\theta =1$, writing ``$< \theta$" means we use
$\bold S^{< \theta}(-)$.
\end{definition}

\begin{discussion}
\label{a19}  
Let us consider $T = T_{\ord} := \Th(\bbQ,<)$,
we concentrate on $f^{\vwt}_T(\lambda)$, the case
$f^{\wat}_\lambda(T)$ can be analyzed similarly.  For any $\lambda$
letting $\Theta^{\tr}_\lambda = \{\kappa:\kappa = \cf(\kappa) \le
\lambda$ and $\lambda^{<\kappa>_{\tr}} > \lambda\}$, see
Definition \ref{a20d} below, so for some
$M \in \EC_{\lambda,1}(T)$ for each $\kappa \in 
\Theta^{\tr}_\lambda$ it has a set $\kappa$ of $> \lambda$ 
cuts of cofinality $(\kappa,\kappa)$.  Now if we consider $N,M \prec N \in 
\EC_{\lambda,1}(T)$, some of these 
will not be filled, hence $f_T(\lambda) \ge |\Theta^{\tr}_\lambda|$.  

Concerning the size of $\Theta^{\tr}_\lambda$ note that by
Easton forcing (using a not necessarily increasing function $f$ from
RCard to Car), if $\mu = \min\{\mu:2^\mu \ge \lambda\}$ then
$\Theta^{\tr}_\lambda \cap [\mu,\lambda)$ is quite arbitrary.
However, by pcf theorems $\Theta^{\tr}_\lambda \cap \mu$ is quite
small, see \cite{Sh:460}, \cite{Sh:829} and maybe even is provably
always finite.

Given $M \in \EC_{\lambda,1}(T)$ there is $N \in 
\EC_{\lambda,1}(T)$ extending it which is strongly $\aleph_0$-saturated
(equivalently, $2$-transitive), filling as many cuts as we can.  Now all the
cuts of $N$ of cofinality $(\aleph_0,\aleph_0)$ are conjugate; also the types
corresponding to cuts ${\cC}$ with cofinality $(\kappa^1_{\cC},
\kappa^2_{\cC})$ such that $\kappa^1_{\cC} \ne \kappa^2_{\cC} \vee 
\kappa^1_{\cC} = \kappa^2_{\cC} \notin \Theta^{\tr}_\lambda \backslash
\{\aleph_0\}$ are easy to handle; because their number is $\le \lambda$,
\cite[Ch.VIII,\S0]{Sh:a} and we fill the cut ${\cC}$ by $J_{\cC}$
such that $J_{\cC}$ has both cofinalities $\aleph_0$ as well as
treating increasing sequences leading to the cuts from both sides; 
in fact we can choose $N$ such that this occurs to 
any cut of $M$ filled by some member of $N \backslash M$.  

But when $\kappa^1_{\cC} = \kappa^2_{\cC} \notin \Theta^{\tr}_\lambda$ call it
$\kappa_{\cC}$ and it $\in \Theta^{\tr}_\lambda \backslash \{\aleph_0\}$
it is not immediately clear whether all such cuts can be treated to
ensure uniqueness up to conjugacy.

Let $\langle(a_{{\cC},i},b_{{\cC},i}):i < \kappa_{\cC}\rangle$ be a 
decreasing sequence of intervals converging to the cut 
${\cC}$; now the isomorphism type of ${\cC}$ can be handled \when \,:
\mn
\begin{enumerate}
\item[$\boxplus_M$]   the following set contains a club
of $\kappa_{\cC},\{i < \kappa_{\cC}$: the cut of $M$ with
lower half $\{a:\bigvee\limits_{j <i} a <_M a_{{\cC},j}\}$ is
filled in $N$ and the cut of $M$  with upper half 
$\{b:\bigvee\limits_{j<i} b_{{\cC},j} < b\}$ is filled in $N\}$.
\end{enumerate}
\mn
Now as classically known we can find a tree ${\cT}$ of cardinality
$\lambda$ with $\le \lambda$ levels and $\le \lambda$
nodes, with nodes intervals of
$I$ and cuts correspond to branches.  So clearly we can ensure
$\boxplus_M$ and this is clearly enough.  So we can understand
$f^{\vwt}_\lambda(T)$ for $T = \Th(\bbQ,<)$.  We may
formalize \ref{a19} as a claim.  (Note that computing 
$f^{\wat}_{T,\theta}(\lambda),f^{\vwt}_{T,\theta}(\lambda)$ for
$\theta > 1$ is easy from the case $\theta =1$.  We use $\alpha(*) \ge
\omega$ below to simplify.
\end{discussion}

\begin{claim}
\label{a20}
Let $T = T_{\ord} := \Th(\bbQ,<)$.  For any cardinal
$\lambda = \aleph_{\alpha(*)} \ge \aleph_\omega$ we have
$f^{\wat}_T(\lambda) = |\alpha(*)|,f^{\vwt}_T(\lambda) =
|\Theta| + 1 = |\Theta|$ where 
$\Theta = \Theta^{\tr}_\lambda := 
\{\theta:\theta = \cf(\theta) \le
\lambda$ and $\lambda^{<\theta>_{\tr}} > \lambda\}$, see below.
\end{claim}

\begin{definition}
\label{a20d} 
$\lambda^{<\theta>_{\tr}} =
\sup\{|\lim_\theta(\cT)|:\cT \subseteq {}^{\theta >}\lambda$ is closed
under initial segments and has cardinality $\le \lambda\}$ where
lim$_\theta(\cT) = \{\eta \in {}^\theta \lambda:\eta \rest i \in \cT$
for every $i < \theta\}$.
\end{definition}

\begin{PROOF}{\ref{a20}}
Let $M \in \EC_{\lambda,1}(T)$ be given, \wilog \, $M$ is  such that:
\mn
\begin{enumerate}
\item[$(*)$]  for every $\theta \in \Reg \cap \lambda^+$, in $M$
\begin{enumerate}
\item[$(a)$]  there is an increasing sequence $\langle
a^1_{\theta,\alpha}:\alpha \le \theta\rangle$
\sn
\item[$(b)$]  there is an decreasing sequence $\langle
a^2_{\theta,\alpha}:\alpha \le \theta\rangle$
\sn
\item[$(c)$]  if $\theta \in \Theta$ there is a tree $\cT_\theta
\subseteq {}^{\theta >} \lambda$
exemplifying $\lambda^{<\theta>_{\tr}} > \lambda$ and members
$\langle b_{\theta,\eta},c_{\theta,\eta}:\eta \in \cT_\theta\rangle$ such
that $\nu \triangleleft \eta \in \cT_\theta \Rightarrow
b_{\theta,\nu} <_M b_{\theta,\eta} <_M c_{\theta,\eta} <_M
c_{\theta,\nu}$ and $[\eta \char 94 \langle \alpha \rangle,\eta \char
94 \langle \beta\rangle \in \cT_\theta,\alpha < \beta \Rightarrow
c_{\theta,\eta \char 94 \langle \alpha \rangle} <_M b_{\theta,\eta \char
94 \langle \beta \rangle}]$.
\end{enumerate}
\end{enumerate}
\mn
Assume $M \prec N \in \EC_{\lambda,1}(T)$ and let $N^+$ be such
that $N \prec N^+$ and $N^+$ is $\lambda^+$-saturated.  For $\ell \le 4$ choose
$d_\ell \in N^+$ be such that: $d_0 \in M,(\forall a \in N)(d_1 < a <
d_2),(\forall a \in N)(a < d_0 \rightarrow a < d_3 < d_0),(\forall a
\in N)(d_0 < a \rightarrow d_0 < d_4 < a)$.  For $\theta \in \Theta$,
let $\eta = \eta_\theta \in \lim_\theta(\cT_\theta)$ and 
$p_\theta = \{b_{\theta,\eta \rest i} < x < b_{\theta,\eta \rest i}:
i < \theta\}$ be such that $p_\theta$ is omitted by $N$, exists by
cardinality consideration; and so $p_\theta$ has unique extension 
$p^+_\theta$ in $\bold S(N)$ and let $e^0_\theta \in N^+$ 
realize it.  For $\theta \in \Reg 
\cap \lambda^+$ let $e^1_\theta \in N^+$ be such that $\alpha < \theta
\Rightarrow a^1_{\theta,\alpha} < e^1_\alpha$ and $(\forall a \in
N)(\bigwedge\limits_{\alpha < \theta} a^1_{\theta,\alpha} < a \Rightarrow
e^1_\alpha < a)$.  Let $e^2_\alpha \in N^+$ be such that $\alpha <
\theta \Rightarrow e^2_\theta < a^2_{\theta,\alpha}$ and $(\forall a
\in N)[\bigwedge\limits_{\alpha < \theta} a^2_\theta <
a^2_{\theta,\alpha} \rightarrow a < e^2_\theta]$.  (So the most
``economical" way is to have $a^1_{\theta,\alpha} = b_{\theta,\eta
\rest \alpha},a^2_{\theta,\alpha} = c_{\theta,\eta \rest \alpha}$ and
$e^0_\theta = e^1_\theta = e^2_\theta,\theta \in \Theta
\Rightarrow e^0_\theta = e^1_\theta$.

Now we prove the four needed inequalities
\mn
\begin{enumerate}
\item[$\boxplus_1$]  $f^{\vwt}_T(\lambda) \ge |\Theta|+1$.
\end{enumerate}
\mn
Why?  It suffices to prove that for any $\bold f:\bold S(M)
\rightarrow \bold S(N)$ such that $p \in \bold S(M) \Rightarrow 
(f(p)) \rest M = p$ we have
$|\{\bold f(p)/ \equiv_{\aut}: p \in \bold S(M)\}| \ge
|\Theta|+1$.  The types $p_0 = \tp(d_0,M)$ and $p_\theta =
\tp(e^0_\theta,M)$ for $\theta \in \Theta$ have unique
extensions in $\bold S(N)$ and clearly $\bold f(p_0),\bold
f(p_\theta),\theta \in \Theta$ are pairwise non-conjugate.
\mn
\begin{enumerate}
\item[$\boxplus_2$]  $f^{\wat}_T(\lambda) \ge |\Reg \cap \lambda^+| + 5$.
\end{enumerate}
\mn
Why?  It suffices to prove that $\bold S(N)/ \equiv_{\aut}$ has
cardinality $\ge |\Reg \cap \lambda^+|+5$.  
Now the types $\tp(d_0,N),\tp(d_1,N),\tp(d_2,N),\tp(d_3,N),
\tp(d_4,N)$ and $\tp(e^0_\theta,N)$ for $\theta \in 
\Reg \cap \lambda^+$ are pairwise non-conjugate.
\mn
\begin{enumerate}
\item[$\boxplus_3$]  $f^{\vwt}_T(\lambda) \le |\Theta| + 1$.
\end{enumerate}
\mn
Why?  It suffices to show that we can choose a model $N_*$ such that 
$M \prec N_* \in \EC_{\lambda,1}(T)$ and a function 
$\bold f:\bold S(M) \rightarrow \bold S(N_*)$ 
such that $p \in \bold S(M) \Rightarrow \bold f(p) \rest M = p$ and
$\{\bold f(p)/ \equiv_{\aut}: p \in \bold S(M)\}$ has
cardinality $\le |\Theta|+1$.  Note that $\sigma := \min(\Theta)$ is 
equal to min$\{\partial:\lambda^\partial > \lambda\}$.  Now choose 
$N_*$ such that
\mn
\begin{enumerate}
\item[$(*)$]  $(a) \quad N \prec N_* \in \EC_{\lambda,1}(T)$
\sn
\item[${{}}$]  $(b) \quad$ if $d \in \gC  \backslash M$ and
$(\theta^-_{M,d},\theta^+_{M,d}) := (\cf\{a \in M:a < d\},<_M)$,

\hskip25pt  $\cf(\{a \in M:d \in a\},>_M)) 
\notin \{(\theta,\theta):\theta \in \Theta\}$ \then

\hskip25pt the type $\tp(d,M)$ is realized in $N_*$
\sn
\item[${{}}$]  $(c) \quad N_*$ is $\sigma$-saturated
\sn
\item[${{}}$]  $(c)^+ \quad$ moreover $N_*$ is strongly
$\sigma$-saturated (i.e. every partial 

\hskip25pt  automorphism of cardinality $< \sigma$ can be extended to an

\hskip25pt automorphism)
\sn
\item[${{}}$]  $(d) \quad \cT \subseteq {}^{\lambda >} 2$ is a tree
with $\le \lambda$ nodes (and $\le \lambda$ levels) and $\bar a =
\langle a_\eta$:

\hskip25pt $\eta \in \cT\rangle$ list the members of $M$ with no
repetitions such that 

\hskip25pt for $\eta \in \cT$ we have $\alpha <
\beta < \ell g(\eta) \Rightarrow (a_{\eta \rest \alpha} < a_\eta \equiv$

\hskip25pt $a_{\eta \rest \alpha} < \bar a_{\eta \rest \beta})$ and
$\alpha < \ell g(\eta) \rightarrow (a_{\eta \rest \alpha} < a_\eta
\equiv \eta(\alpha) = 1)$
\sn
\item[${{}}$]  $(e) \quad$ if $\eta \in \cT$ \then \, for some
$e^0_\eta,e^1_\eta \in N_* \backslash M$ we have
$\{a \in M:a < e^0_\eta\} =$

\hskip25pt $\{a \in M:(\exists \alpha < \ell
g(\eta)[\eta(\alpha) = 1 \wedge \alpha < \ell g(\eta),
a \le a_{\eta \rest \alpha}]\}$

\hskip25pt  and $\{a \in M:e^1_\eta < a\} = \{a \in
M:(\exists \alpha < \ell g(\eta)[\eta(\alpha) = 0$

\hskip25pt $\wedge a_{\eta \rest \alpha} \le a]\}$
\sn
\item[${{}}$]  $(f) \quad$ if $\eta_\ell \in \cT,\delta_\ell = \ell
g(\eta_\ell)$ is a limit ordinal, $\delta_\ell = \sup\{\alpha <
\delta:\eta(\alpha) = 1\}$

\hskip25pt  for $\ell=1,2$ and $\cf(\delta_1) =
\cf(\delta_2)$ \then \, there is an automorphism 

\hskip25pt $g$ of $\gC$ mapping $N_*$ onto $N_*,e^1_{\eta_1}$ to 
$e^1_{\eta_2}$ and hence mapping 

\hskip25pt $\{a \in M:(\exists
\alpha < \delta_0)(\eta_1(\alpha) = 1 \wedge a < a_{\eta \rest
\alpha})\}$ onto 

\hskip25pt $\{a \in M:(\exists \alpha < \delta_2),\eta_2(\alpha)
= 1 \wedge a < a_{\eta_2 \rest \alpha})\}$.
\end{enumerate}
\mn
Why is this possible: for (c) as $\lambda = \lambda^{< \sigma}$, for $(b)$
as $\{\tp(d,M):d \in \gC$ and $\theta^-_{M,d} \ne \theta^+_{M,d}$ are
infinite$\}$ has $\le \lambda$ members by \cite[Ch.VIII,\S0]{Sh:c} and
$\{\tp(d,M):d \in \gC$ and 
$\aleph_0 \le \theta^-_{M,d} = \theta^+_{M,d} \notin \Theta\}$ 
has $\le \lambda$ members by the definition of
$\Theta$ (and the well known old equivalence of trees and number of
cuts); lastly $\{\tp(d,M):d \in \gC$ and $\theta^-_{M,d}
\in\{0,1\}$ or $\theta^-_{M,d} \in \{0,1\}$ has $\le \lambda$ members
trivially. Also clauses (d),(e),(f) are straight.

Now we define $\bold f$, so let $p \in \bold S(M)$.  First, if some $d
= d_p \in N_*$ realize $p$,  then let $\bold f(p) = 
\tp(d_p,N_*)$ so by clause (c)$^+$ clearly $\bold f(p),p_0 
= \tp(d_0,N_*)$ are conjugate.  Second, if $p \in \bold S(M)$ 
is not realized in $N_*$ then by clauses (d),(e),(f) 
there are $\theta \in
\Theta$ and $<_N$-increasing $\langle d^-_{p,i}:i < \theta\rangle$ and
$<_M$-decreasing $\langle d^+_{p,i}:i < \theta\}$ such that $d^-_{p,i}
<_N d^+_{p,i}$ for $i < \theta$ and $p$ include $p' = \{d^-_{p,i} < x <
d^+_{p,i}:i < \theta\}$ which $N_*$ omits hence $p$ has unique extension
$\bold f(p)$ in $\bold S(N_*)$, but for each limit $\delta < \theta$ the
types $\{d^-_{p,i} < x < a:i < \delta$ and $a \in M$ and $j < \delta
\Rightarrow d^-_{p,j} < a\},\{a < x < d^+_{p,j}:i < \delta,a \in M$
and $j < \delta \Rightarrow a < d^+_{p,j}\}$ are realized by
clauses (d),(e).  Now easily $\bold f(p)$,
tp$(e^0_\theta,N)$ are conjugate by some $g \in \aut(N)$ such
that $g(b_{\theta,i}) = d^-_{p,i},g(c_{\theta,i}) = d^+_{p,i}$,
because we can choose it in each relevant convex set by clause $(c)^+$.
\mn
\begin{enumerate}
\item[$\boxplus_4$]  $f^{\wat}_T(\lambda) \le |\alpha(*) + 6|$.
\end{enumerate}
\mn
It is simpler when $\alpha(*) \ge \omega$ and the proof is similar to
the proof of $\boxplus_3$ but use $\prec$-increasing continuous
$\langle N^*_\varepsilon:\varepsilon \le \sigma\rangle,N^*_0 = N$, etc.
\end{PROOF}

\begin{question}
\label{a21} 
1) For $T$ countable, dependent and unstable, is
 $f^{\vwt}_\lambda(T)$ essentially equal to
 $f^{\vwt}_{\Th(\bbQ,<)}(\lambda)$? at least
can we understand it (and $f^{\wat}_T(\lambda))$?  

\noindent
2) What can we say on $f^{vwt}_T(\lambda),f^{\wat}_T(\lambda)$ for 
independent $T$?, see below.
\end{question}

\begin{discussion}
\label{a23}  
1) Concerning Part (2) of \ref{a21}, it is easy to note:  if $T$ 
is independent and $|T| \le \mu < \lambda \le 2^\mu < 2^\lambda$ and
$\cf([2^\mu]^\lambda,\subseteq) > 2^\mu$ hold, e.g. if $\cf(2^\mu) \le
\mu$, then $f^{\wat}_T(\lambda) = 2^\lambda$; see more in
Kojman-Shelah \cite{KjSh:409}, \cite[4.7]{Sh:E46}.

\noindent
2) For independent $T$ the situation concerning
 $f^{\vwt}_{T,\theta}(-)$ is very different than for 
$f^{\aut}_{T,\theta}(-)$.  Why?  By the following.
\end{discussion}

\begin{claim}
\label{a24}
1) If $\lambda = \lambda^{< \lambda}  > \theta + |T| + \aleph_0$ 
and $T$ is a complete first order theory, \then \,
$f^{\vwa}_{T,\theta}(\lambda) \le 2^\theta$.

\noindent
2) Moreover for every $M \in \EC_{\lambda,1}(T)$
 there is an elementary extension $M^+ \in \EC_{\lambda,1}(T)$
 such that 
\mn
\begin{enumerate}
\item[$(*)_{M_1,M_2}$]  if $p \in \bold S^\theta(M)$ then for
 some $q = q_p \in \bold S^\theta(M^+)$ extending $p$ the 
model $M^+_{[q]}$ is saturated, see Definition \ref{z5}.
\end{enumerate}
\end{claim}

\begin{PROOF}{\ref{a24}}
1) By (2).

\noindent
2) Let $N$ be such that $M \prec N$ and every $p \in \bold
S^\theta(M)$ is realized by $\bar a_p \in {}^\theta N$.

For $\alpha < \lambda$ let $D_\alpha$ be a regular ultrafilter on
$I_\alpha = |\alpha| + \aleph_0$.  Now we choose
$(N_\alpha,M_\alpha)$ by induction on $\alpha \le \lambda$ such that 
\mn
\begin{enumerate}
\item[$(a)$]  $(N_\beta,M_\beta)$ is elementarily equivalent to
$(N,M)$ (where $(N_\beta,M_\beta)$ is the $(\tau_T \cup \{P\})$-model 
expanding $N_\beta$ by $P^{(N_\beta,M_\beta)} = |M_\beta|$, so 
$P$ a new unary predicate)
\sn
\item[$(b)$]  $(N_0,M_0) = (N,M)$
\sn
\item[$(c)$]  the sequence $\langle (N_\beta,M_\beta):\beta \le
\alpha\rangle$ is $\prec$-increasing continuous 
\sn
\item[$(d)$]  if $\alpha = \beta +1$ then there is an isomorphism
$\bold j^+_\beta$ from $(N_\alpha,M_\alpha)$ onto
$(N_\beta,M_\beta)^{I_\beta}/D_\beta$ extending the canonical embedding $\bold
j_\beta$ from$(N_\beta,M_\beta)$ into
$(N_\beta,M_\beta)^{\lambda_\beta}/D_\beta$, 
i.e. for $a \in N_\beta,\bold j_\beta(a) =
f_{a,\beta}/D_\beta$ where $f_{a,\beta}:\lambda_\beta \rightarrow
N_\beta$ is constantly $a$.
\end{enumerate}
\mn
There is no problem to carry the definition and $M^+ := M_\lambda$ is
as required.  That is, we can prove by induction on $\alpha$ that
$\|M_\alpha\| = \lambda$: if $\alpha = 0$ by clause (b) if $\alpha =
\beta +1$ as $\lambda \le \lambda^{I_\beta}/D_\alpha \le
\lambda^{|I_\beta|} \le \lambda^{< \lambda} = \lambda$ and for
$\alpha$ limit by the induction hypothesis.  Also, as $D_\beta$ is a
regular ultrafilter, clearly $M_\lambda$ is saturated hence $M_\lambda
\in \EC_{\lambda,\lambda}(T)$.  Similarly $(N_\lambda,M_\lambda)$ is
$\lambda$-saturated hence if $\bar a \in {}^{\lambda >}(N_\lambda)$
then $(M_\lambda)_{[\bar a]}$ is saturated.  We choose 
$M^+ = M_\lambda$ so indeed $M \prec M^+ \in \EC_{\lambda,\lambda}(T)$. 

Now for every $p \in \bold S^\theta(M)$ recall that $\bar a_p \in {}^\theta
N \subseteq {}^\theta(N_\lambda)$ realizes $p$, so let $q_p(\bar
x_{[\bar a]}) = \tp(\bar a_p,M_\lambda)$, so we are done.
\end{PROOF}

\begin{claim}
\label{a25}
1) Assume $T$ is a complete first order theory and $\lambda$ is strong
   limit singular of cofinality $\kappa,\lambda > \theta,\lambda > |T|
   + \aleph_0$.  \Then \, $f^{\vwa}_{T,\theta}(\lambda) \le 2^\theta$.

\noindent
2) Like \ref{a24}(2) replacing ``saturated" by ``special", see \cite{CK73}.
\end{claim}

\begin{PROOF}{\ref{a25}}
1) By part (2).

\noindent
2) Similar to the proof of \ref{a24}, but we elaborate.   Now the
 definition of ``special" says that there is $\bar M = \langle M^*_i:i <
\kappa\rangle$ which is a $\prec$-increasing continuous sequence of
 models (of $T$) with union $M$ such 
that $M_{i+1}$ is $\|M_i\|^+$-saturated and
 $i < \kappa \Rightarrow \|M^*_i\| < \lambda$.
Let $\langle \lambda_i:i < \kappa\rangle$ be an increasing sequence of
regular cardinals with limit $\lambda$.  We choose $N,\langle \bar
a_p:p \in \bold S^\theta(M)\rangle$ and $\langle D_\alpha:\alpha <
\lambda\rangle$ as in the proof of \ref{a24}.  We now choose
$(N_\alpha,\bar M_\alpha)$ by induction on $\alpha < \lambda$ such
that:
\mn
\begin{enumerate}
\item[$\boxplus$]  $(a)(\alpha) \quad \bar M_\alpha = \langle
M_{\alpha,i}:i \le \kappa\rangle$ is $\prec$-increasing continuous
\sn
\item[${{}}$]  $\qquad (\beta) \quad (N_\alpha,M_{\alpha,i})$ is
elementarily equivalent to $(N,M^*_i)$ for $i < \kappa$ 

\hskip35pt such that $\lambda_i \ge \alpha$ so $M_{\alpha,i} \prec N_\alpha$
\sn
\item[${{}}$]  $(b)(\alpha) \quad N_0 = N$
\sn
\item[${{}}$]  $\qquad (\beta) \quad M_{0,i} = M^*_i$ for $i < \kappa$
\sn
\item[${{}}$]  $\qquad (\gamma) \quad M_{0,\kappa} = M$
\sn
\item[${{}}$]  $(c) \quad \langle (N_\beta,M_{\beta,i}):\beta \le
\alpha \rangle$ is $\prec$-increasing continuous
\sn
\item[${{}}$]  $(d) \quad$ if $\alpha = \beta +1$ then
\sn
\item[${{}}$]  $\qquad (\alpha) \quad$ there is an isomorphism $\bold
j^+_\beta$ from $N_\alpha$ onto $N^{I_\beta}_\beta/D_\beta$ extending

\hskip35pt  the canonical embedding of 
$N_\beta$ into $N^{I_\beta}_\beta/D_\beta$
\sn
\item[${{}}$]  $\qquad (\beta) \quad$ if $\beta < \lambda_i$ then
$\bold j^+_\beta$ maps $M_{\alpha,i}$ onto $M^{I_\beta}_{\beta,i}/D_\beta$
\sn
\item[${{}}$]  $\qquad (\gamma) \quad$ if $\beta \ge \lambda_i$ then
$M_{\alpha,i} = M_{\beta,i}$.
\end{enumerate}
\mn
In the end $\langle M_{\lambda_i,\lambda_i}:i < \kappa\rangle$ witness
that $M^+ := \cup\{M_{\lambda_i,\lambda_1}:i < \lambda\}$ is special;
moreover, if $\bar a \in {}^\theta N$ then $q_\beta := \tp(\bar
a_p,M^+,N_\lambda)$ is as promised.
\end{PROOF}

If you do not like the use of instances of GCH, i.e. $\kappa =
\kappa^{< \kappa}$, but like to stick to essentially the same property,
we can reformulate it.

\begin{definition}
\label{a26}   
Let $f^{\aut,*}_{T,\theta}(\lambda)$,
for $\lambda$ regular be the minimal $\mu$ such that for any
$\lambda$-saturated $M \prec {\gC}$, e.g. of cardinality $2^{<
\lambda}$ we can find a subset
$\bold P$ of $\bold S^\theta(M)$ of cardinality $\le \mu$ satisfying
that:
\mn
\begin{enumerate}
\item[$(*)$]   for any $p_1(\bar x_{[\theta]}) \in 
\bold S^\theta(M)$ there is $p_2(\bar x_{[\theta]}) \in 
\bold P$ such that letting $\bar a_\ell = \langle
a_{\ell,i}:i < \theta\rangle$ realizes $p_\ell(\bar x_{[\theta]})$ in $\gC$ for
$\ell=1,2$ we have
\begin{enumerate}
\item[$\odot$]  in the E.F. (i.e. Ehrenfeucht-Fr\"aiss\'e) 
game of length $\lambda$ 
for the pair $(M_{[\bar a_1]},M_{[\bar a_2]})$ the 
ISO player has a winning strategy.
\end{enumerate}
\end{enumerate}
\end{definition}

\begin{discussion}
\label{a28}  
Concerning $f^{\aut,*}_{T,\theta}(-)$.

\noindent
1) The positive result, i.e. upper bound 
for dependent $T$ (see end of \S4) still holds as well as the negative ones.

\noindent
2) The negative results for independent $T$ holds.

\noindent
3) The question is closed to the one on ``what occurs in 
$\bold V^{\Levy(\lambda,\chi)}$ for some $\chi$".
\end{discussion}

\begin{question}
\label{a29}
Generalize to any dependent $T$ the theorem: a linear order of
cardinality $\lambda$ has $\le \lambda$ cuts of different
lower cofinality and upper cofinality.
\end{question}
\bigskip

\subsection{On the outside definable sets and uf$(p)$} \
\bigskip

\begin{definition}
\label{a31}  
1) Let $\Deef^\alpha_\Delta(M) = \{\varphi(M,\bar c):
\varphi(\bar x,\bar y) \in \Delta,\ell g(\bar x) =
\alpha$ and $\bar c \in {}^{\ell g(\bar y)}{\gC}\}$ and $\Deef_\alpha(M) =
\Deef^\alpha_{\bbL(\tau_T)}(M)$ where see below;
of course instead ${\gC}$ we can use any 
$\|M\|^+$-saturated elementary extension of $M$.

\noindent
2) $\varphi(M,\bar c) = \{\bar b:\bar b \in {}^{\ell g(\bar x)} M$ and
${\gC} \models \varphi[\bar b,\bar c]\}$ where $\varphi =
\varphi(\bar x,\bar y)$.

\noindent
2A) We say $\bold I \subseteq {}^\alpha M$ is outside definable when
it belongs to $\in \Deef_\alpha(M)$.

\noindent
3) If $p(\bar x) \in \bold S^\alpha(M)$ let $\uf(p) 
= \{D:D$ an ultrafilter on the Boolean Algebra 
$\Deef_\alpha(M)$ containing $\{\varphi(M,\bar a):
\varphi(\bar x,\bar a) \in p\}\}$.

\noindent
3A) If $p \in \bold S^\alpha(M)$ and $\Delta \subseteq \{\varphi(\bar
x_\alpha,\bar y):\bar y \in \{\bar y_{[n]}:n < \omega\},\varphi \in 
\bbL(\tau_T)\}$ \then \, let uf$_\Delta(p) = \{D \cap 
\Deef^\alpha_\Delta(M):D \in \uf(p)\}$.  If $\Delta =
\{\varphi\}$ we may write $\varphi$.

\noindent
4) We say $p$ has super multiplicity 1 \when \, $|\uf(p)| =1$.

\noindent
5) If $q(\bar x,\bar y) = \tp(\bar a \char 94 \bar b,M)$ and
p$(\bar x) = \tp(\bar a,M)$ then $\pi = \pi_{p(\bar x),
q_(\bar x,\bar y)}$ is the function from uf$(q)$ onto $\uf(p)$, we call it
the projection, such that if $D \in \uf(q)$ and $M \subseteq A 
\subset {\gC}$ and $\bar a' \char 94 \bar b'$ realizes 
$\Av(D,A)$ then $\bar a'$ realizes $\Av(\pi(D),A)$, see \ref{a34}(1) below.

\noindent
6) We say $\bold I = \langle \bar a_t:t \in I\rangle$ is an indiscernible
sequence based on $p \in \bold S^\alpha(M)$ \when \, ($I$ is a
linear order and) $\bold I$ is based on some $D \in 
\uf(p)$ which means that: for each $t \in I,\tp(\bar a_t,M 
\cup \{\bar a_s:s$ satisfies $t <_I s\} \cup M)$ is 
$\Av(D,\{\bar a_s:s$ satisfies $t <_I s\} \cup M$).  
Similarly for $p \in \bold S^m(A)$ which is finitely satisfiable in
$M$ and $\bold I$ is based on $(D,A)$.

\noindent
7) Assume $p \in \bold S^\alpha(M)$ and $D \in \uf(p)$, let
$\Dom(D) = |M|$, (we can replace it by an set).  We say $\bar a$
realizes $\tp(D,A)$ \when \, 
there is $\langle \bar a_n:n < \omega\rangle$, as in part (6), i.e.
such that $\bar a_n$ realizes $\Av(D,\{\bar a_\ell:\ell \in
(n,\omega)\} \cup \bar a \cup \Dom(D))$ for 
$n < \omega$ and $\langle \bar a\rangle
\char 94 \langle \bar a_n:n < \omega\rangle$ is an indiscernible
sequence over $A$.

\noindent
8) Above we say ``realizes 
$\tp_\Delta(D,A)$", when in the end $\langle \bar a \rangle \char 94
\langle \bar a_n:n < \omega\rangle$ is demanded only to be an
$\Delta$-indiscernible over $A$.

\noindent
9) For $D$ as above let $\lambda(D) = \min\{\sigma$: for some $A
\subseteq \Dom(D)$ of cardinality $\sigma$, no 
$\bar a \in {}^\gamma(\Dom(D))$ realizes $(A,D)\}$.

\noindent
10) $\lambda_\Delta(D)$ is defined similarly restricting ourselves to
$\Delta$.

\noindent
11) $\lambda_{\loc}(D) = \sup\{\lambda_\Delta(D):\Delta \subseteq
\bbL(\tau_T)$ finite$\}$.
\end{definition}

\begin{claim}
\label{a34}  
1) For $M,\bar a,\bar b,p(\bar x),q(\bar x,\bar y)$ as in 
Definition \ref{a31}(5), the function
$\pi_{p(\bar x),q(\bar x,\bar y))}:\text{\rm uf}(q) \rightarrow 
\uf(p)$ is well defined.

\noindent
2) Moreover it is onto.
\end{claim}

\begin{PROOF}{\ref{a34}}
  1) Should be clear.

\noindent
2) So assume $D_1 \in \uf(p(\bar x))$.  It suffices to prove that
\mn
\begin{enumerate}
\item[$(*)$]   the family ${\cX}_1 \cup {\cX}_2$
can be extended to an ultrafilter on ${}^\alpha M$ where
\begin{enumerate}
\item[$(a)$]   ${\cX}_1 := \{X'_1$: for some $X_1 \in D$ we have $X'_1 =
\{\bar a' \char 94 \bar b':\bar a' \in {}^{\ell g(\bar x)}M,\bar b'
\in {}^{\ell g(\bar y)}M$ and $\bar a' \in X_1\}\}$ and
\sn
\item[$(b)$]   ${\cX}_2 := \{\{\bar a' \char 94 \bar b':
\bar a' \in {}^{\ell g(\bar x)}M,\bar b' \in
{}^{\ell g(\bar y)}M$ and $M \models \varphi[\bar a',\bar b',\bar c]\}:
{\gC} \models \varphi[\bar a,\bar b,\bar c]$ and $\bar c
\subseteq M,\varphi = \varphi(\bar x,\bar y,\bar z) \in \bbL(\tau_T)\}$.
\end{enumerate}
\end{enumerate}
\mn
As each of the two families in the union is closed under (finite)
intersection, it suffices to prove:
\mn
\begin{enumerate}
\item[$\odot$]  assume $\varphi = \varphi(\bar x,\bar z_1),\bar c_1
\in {}^{\ell g(\bar z_1)}{\gC}$ and $X_1 := 
\varphi(M,\bar c_1) \in D_1$ define $X'_1 \in {\cX}_1$ as in $(*)(a)$ and
$\psi = \psi(\bar x,\bar y,\bar z_2),\bar c_2 \in 
{}^{\ell g(\bar z_2)}M$ such that ${\frak C} \models 
\psi[\bar a,\bar b,\bar c]$, defines $X_2 \in {\cX}_2$ as in $(*)(b)$
\then \, we can find $\bar a',\bar b'$ in $M$ such that ${\gC}
 \models \varphi[\bar a',\bar c_1] \wedge \psi[\bar a',\bar b',\bar c_2]$.
\end{enumerate}
\mn
To prove $\odot$ note that the set 
$Y_1 := \{\bar a' \in {}^{\ell g(\bar x)}M:M \models 
(\exists \bar y)\psi(\bar a',\bar y,\bar c_2)\}$ belongs to $D_1$ because $D_1
\in \uf(p),p = \tp(\bar a,M)$ and ${\gC} 
\models (\exists \bar y)\psi[\bar a,y,\bar c_2]$.  Hence $X_1
   \cap Y_1 \in D_1$ and choose $\bar a' \in X_1 \cap Y_1$.  As $\bar
   a' \in Y_1$ there is $\bar b' \in {}^{\ell g(\bar y)}M$ such that
   $M \models \psi[\bar a',\bar b',\bar c_2]$ and as $\bar a' \in X_1$
   we have $M \models \varphi[\bar a',\bar c_2]$.  Together $\bar a'
   \char 94 \bar b'$ is as required in $\odot$.
\end{PROOF}

\begin{claim}
\label{a37} 
We assume (needed really just in parts (0),(2),(4), that
$T$ is dependent.

\noindent
0) If $\bold I$ is an infinite indiscernible set, \then \, 
$\bold I$ sits stably, see \ref{a67}(2), 
(so every $p \in \bold S^{< \omega}(\cup \bold I)$ is definable).

\noindent
1) If $D \in \uf(p(\bar x)),p(\bar x) \in \bold S^\alpha(M)$ 
and $I$ is a linear order
\then \, there is an indiscernible sequence $\langle \bar a_t:t
 \in I\rangle$ over $M$ based on $D$, see Definition \ref{a31}(6).  We
can replace $M$ by a set $A$.

\noindent
2) In part (1) if $\bold I_\ell = \langle \bar a^\ell_t:t 
\in I\rangle$ is an indiscernible sequence based on $D,I$ is a linear
order with no first element and 
$\bar a^\ell_t$ realizes $\Av(D,\cup\{a^k_s:s$
satisfies $t <_I s$ and $k = 1,2\})$ \then \,  $\bold I^*_1,\bold
I^*_2$, i.e. $\bold I_1,\bold I_2$ inverted are equivalent, see \ref{a67}(5).

\noindent
3) In Definition \ref{a31}(7), it is equivalent ``for every infinite
linear order $I$ there is an
indiscernible sequence $\langle \bar a_t:t \in I \rangle$ over $M$
 based on $D$".

\noindent
4) Assume $D_\ell \in \uf_\gamma(M)$ and $\langle \bar
   a^\ell_n:n < \omega\rangle$ is an indiscernible sequence based on
$D_\ell$, see Definition \ref{a31}(6) for $\ell=1,2$, \then \,
   $D_1 = D_2$ iff $\tp(\langle \bar a^1_n:n < \omega\rangle,M) =
\tp(\langle a^2_n:n < \omega\rangle,M)$.

\noindent
4A) Assume for transparency $\gamma < \omega$ and $\Delta \subseteq
\bbL(\tau_T)$ is finite.  \Then \, for some $n_\Delta < \omega$ for
every $D_1,D_2,M,\bar a^\ell_n$ as in part (4) we have: $D_1 \cap
\Deef^\gamma_\Delta(M) = D_2 \cap \Deef^\gamma_\Delta(M)$
\Iff \, $\tp_\Delta(\langle \bar a^1_n:n < n_\Delta\rangle,M) = 
\tp_\Delta(\langle \bar a^2_n:n < n_\Delta\rangle,M)$.

\noindent
5) If $\zeta < \theta^+,M$ is $\kappa$-saturated, $\cf(\kappa) >
2^\theta$ and $p \in \bold S^\zeta(M)$ \then \,  for some
 $\bold u,A_*$ we have (we write $\bold u = \bold u_\kappa(p))$:
\mn
\begin{enumerate}
\item[$(a)$]   $\bold u$ is a non-empty subset of $\uf(p)$, see
\ref{a31}(8) 
\sn
\item[$(b)$]   if $D \in \bold u$ and $A \subseteq M$ has
cardinality $< \kappa$ \then \, some $\bar a \in {}^\zeta M$
 realizes $\tp(D,A)$
\sn
\item[$(c)$]   $A_* \subseteq M$ and $|A_*| < \kappa$
\sn
\item[$(d)$]   if $\bar a \in {}^\zeta M$ realizes $p \rest A$ then
for some $D \in \bold u,\bar a$ realizes $\tp(D,A_*)$, see
Definition \ref{a31}(8)
\sn
\item[$(e)$]   if $D \in \uf(p) \backslash \bold u$,
\then \,  no $\bar a \in {}^\zeta M$ realizes $\tp(D,A)$. 
\end{enumerate}
\end{claim}

\begin{PROOF}{\ref{a37}}
 Parts (0),(2),(4),(4A) and (5) by \cite{Sh:715}, the others are
 obvious.  
\end{PROOF}

\begin{observation}
\label{a40}  
Assume $p \in \bold S^\gamma(M)$ and $|T| + \gamma < \theta^+$.  
If $\uf(p)$ has cardinality $> 2^\theta$ \then \, for some 
$\varphi = \varphi(\bar x_\gamma,\bar y)$, also
$\uf_{\varphi(\bar x_\gamma,\bar y)}(p)$ has cardinality $>
2^\theta$.  In fact $|\uf(p)| \le \Pi\{\uf_{\{\varphi\}}(p):
\varphi = \varphi(\bar x_{[\gamma]}) \in \bbL(\tau_T)\}$. 
\end{observation}

\begin{PROOF}{\ref{a40}}
Obvious. 
\end{PROOF}
\bn
Recall question \cite[6.1]{Sh:715}.

\begin{definition}
\label{a41} 
1) $T$ has bounded directionality \when \,: if $p \in 
\bold S^\alpha(M)$, then the set $\uf(p) = \{D:D$ an ultrafilter on 
$\Deef^\alpha(M)$ such that $\Av(D,{}^\alpha M) = p\}$ has 
cardinality $\le 2^{|T|+|\alpha|}$.

\noindent
1A) We define ``finite directionality" similarly when we consider only
$p \in \bold S^{< \omega}(M)$.

\noindent
1B) We define ``unary directionality" similarly when we consider only
$p \in \bold S(M)$.

\noindent
2) We say $T$ has medium directionality \when \,  for every $p \in \bold
S^\alpha(M)$, the set uf$(p)$ has cardinality 
$\le \|M\|^{|\alpha| + |T|}$, but $T$ does not have bounded
directionality.

\noindent
3) We say that $T$ has large directionality \when \, it neither has
bounded directionality nor medium directionality.
\end{definition}

\begin{claim}
\label{a42}  
1) $T$ has bounded directionality \Iff \, $\uf_\Delta(p)$ 
is finite whenever $p \in \bold S^\varepsilon(M),
\Delta \subseteq \Gamma_\varepsilon$ is finite
\Iff \, for some $\lambda \ge |T|$ we have $M \in \EC_{\lambda,1}(T) 
\wedge \Delta \subseteq \bbL(\tau_T)$ finite $\wedge 
p \in \bold S^{<\omega}(M) \Rightarrow |\uf_\Delta(p)| < \lambda$.

\noindent
2) If $T$ has medium directionality \Iff \, for every $\lambda \ge
   |T|$ we have $\lambda = \sup\{|\uf_\Delta(p)|:p \in \bold
S^{< \omega}(M)$ and $M \in \EC_{\lambda,1}(T)$ and $\Delta
\subseteq \bbL(\tau_T)$ is finite$\}$.

\noindent
3) If $T$ has large directionality \Iff \, for every $\lambda >
|T|$ we have $\sup\{\lambda^{<\theta>_{\tr}}:\theta \le \lambda$
regular$\} = \sup\{|\uf_\Delta(p)|:p \in \bold S^{< \omega}(M),
M \in \EC_{\lambda,1}(T)$ and $\Delta \subseteq \bbL(\tau_T)$ is finite$\}$.

\noindent
4) If $T$ has medium or bounded directionality, $M \prec {\gC}, p
 = \tp(\bar a,M) \in \bold S^\varepsilon(M)$ and $D \in
\uf(p)$ and $\langle \bar a_n:n < \omega\rangle$ is an
indiscernible sequence based on $D$ \then \,  for every $\varphi =
 \varphi(\bar x^0_{[\varepsilon]},\dotsc,\bar x^{n-1}_{[\varepsilon]};\bar y)$
the type $p = \tp_\varphi(\bar a_0 \char 94 \ldots 
\char 94 \bar a_{n-1},M)$ is definable with parameters in the model
 $M_{[\bar a]}$.

\noindent
5) If $T$ has bounded directionality \then \, in part (4) the type
$p_\varphi$ is definable almost on $\emptyset$ in the model
$M_{[\bar a]}$, i.e. in $M^{\eq}_{[\bar a]}$ it is definable with parameters.
\end{claim}

\begin{PROOF}{\ref{a42}}
1) Clearly the second implies the first which implies the third.  

Lastly, assume the second fails and we shall prove the third fails.
So we are assuming $p \in \bold S^m(M)$ and $\uf_\Delta(p)$ infinite,
$\Delta$ finite.  let $\langle D_n:n < \omega\rangle$ be a member of
$\uf(p)$ such that $\langle \Av_\Delta(D_n,{\gC}):n <
\omega\rangle$ are pairwise distinct.  For $n < m < \omega$ as
$\Av(D_n,{\gC}) \ne \Av(D_m,{\gC})$ we can choose let
$\varphi_{n,m}(\bar x,\bar b_{n,m}) \in \Av(D_n,\bar b_{n,m})$
such that $\neg \varphi_{n,m}(\bar x,\bar b_{n,m}) \in \Av(D_m,\bar b_{n,n})$.

Let $M = M_0 \prec M_1 \prec M_2 \cup \{\bar b_{m,n}:m < n\} \subseteq
N_1,\bar c_n \subseteq M_2$ realizes $\Av(D_n,N)$ hence it
realizes $p \in \bold S^m(M)$.  Let $M^+ =
(N,P^{N^+}_1,Q^{N^+}_1,<^{N^+}),P^{N^+}_1 = |M|,Q^N_1 = \{\bar
b_{m,m}:m < n < \omega\},P^{N^+}_2 = |N_1|,Q_2 = \{\bar c_n:n <
\omega\}$.

Let $N^{++}_0$ be a $\mu^+$-saturated model of $\Th(M^+)$.  Without
loss of generality $N_2 = N^+ \rest \tau_T \prec {\gC}$, let $N_0
= {\gC} \rest P^{N^{++}}_1,N_1 = {\gC} \rest P^{N^{++}}_2$.  

Let $p' = \tp(\bar c,M')$ for every $\bar c \in
Q^{N^{++}}_2$.  For every $\bar c \in Q^{N^{++}}_2$ let $q_{\bar c} =
\tp(\bar c,N'_1)$, this type is finitely satisfiable in $M'$
(by $N^{++} \equiv M^+$) hence for some ultrafilter $D$ on $M'$ we
have $q_{\bar c} = \Av(D_{\bar c},N'_1)$.  Now for any $\bar
c_1 \ne \bar c_2$ in $Q^{N^{++}}_2$ we have $D_{\bar c_1} \ne D_{\bar
c_2}$ so for some $\bar b \in Q^{\text{pos}}_1,\varphi(\bar x_1,\bar b)
\equiv \neg \varphi(\bar c_2,\bar b)$ hence for some $\bold t$ we hvae
$\varphi(\bar x,\bar b)^{\bold t} \in q_{\bar c_1},\neg \varphi(\bar
x,\bar b)^{\bold t} \in q_{\bar c_2}$.

So $\langle D_{\bar c} \cap \deef_\Delta(M_0):\bar c \in
Q^{N^{++}}_0\rangle$ is a sequence of pairwise distinct members of
{\rm uf}$_\Delta(p')$.  As $|Q^{N^{++}}_2| \ge \mu^+$ we are done. 

\noindent
2),4) See Kaplan-Shelah \cite{KpSh:946} using \cite{Sh:a} and \ref{a37}(4).

\noindent
5) Obvious. 
\end{PROOF}

\begin{remark}   Can define $\uf_\Delta(p) = \{D \cap 
\deef_\Delta(M):D \in \uf(M)$ and $\Av(D,M) \supseteq p \rest 
\Delta\}$, no difference in the proof.
\end{remark}

\begin{question}
\label{a43}  
If $M_{[p(\bar x)]} \prec N_{[q(\bar x)]}$
how are $\uf(p(\bar x)),\uf(q(\bar x))$ related?
\end{question}

\begin{question}
\label{a47}  
Can we prove a substitute?  We do not deal
with it presently.  E.g. we may consider $\uf(\bold I),\bold I$ a
$k$-end-homogeneous sequence, see below and \S(5B).
\end{question}

\begin{definition}
\label{a49}  For a sequence $\bold I = \langle
\bar a_s:s \in I\rangle$ of $\zeta$-tuples, ($I$ a linear order) let

\noindent
1) $\bbB(\bold I) = \{J \subseteq I$: for some 
$\varphi(\bar x_{[\zeta]},\bar b)$ for every $t \in I$ we have
$\varphi[\bar a_t,\bar b] \Leftrightarrow t \in J\}$.

\noindent
2) $\uf(\bold I)$ is the set of ultrafilters $D$ on $\bbB(\bold I)$
   containing all co-bounded subsets, so interesting only when $\bold
   I$ has no last element.

\noindent
3) For $\Delta \subseteq \Gamma_\zeta = \{\varphi:\varphi =
\varphi(\bar x_{[\zeta]},\bar y) \in \bbL(T)\}$ let 
$\bbB_\Delta(\bold I) = \{J \subseteq I$: for some $\varphi
(\bar x_{[\zeta]},\bar y) \in \Delta$ and $\bar b \in 
{}^{\ell g(\bar y)}{\gC}$ for every $t \in I$ we have $t \in J \Rightarrow
   \varphi[\bar a_t,\bar b]\}$.

\noindent
4) $\uf_\Delta(\bold I) = \{D \cap \bbB_\Delta(\bold I):D \in 
\uf(\bold I)\}$. 
\end{definition}
\bn
Probably we may do better.

\begin{question}
\label{a51}  
Does the directionality of $T$ essentially
determine when $\lambda \rightarrow_T (\kappa)_\sigma$?  See on the
directionality, see \ref{a42} and on the arrow see \ref{a76} and on
$\lambda \rightarrow_T(\kappa)_\sigma$ for dependent $T$ see \cite{KpSh:946}.

We have divided the family of dependent unstable $T$'s to three.
\end{question}

\begin{claim}
\label{a53}  
1) Every dependent $T$ satisfies exactly
one of the following possibilities: stable, unstable (dependent with)
bounded directionality, unstable dependent with medium directionality and
unstable dependent with large directionality.

\noindent
2) Each of those classes is non-empty.
\end{claim}

\begin{remark}
\label{a55}  1) Delon, 
see Poizat \cite{Poz00}, gives an example of a
dependent $T$ with $|\uf(p)| > \|M\|$, in the present
terminology this means a dependent $T$ with large directionality.
\end{remark}

\begin{PROOF}{\ref{a53}}
1) By \ref{a42}.

\noindent
2) See Kaplan-Shelah \cite{KpSh:946}.  
\end{PROOF}

\begin{question}
\label{a57}   In the definition of medium/large directionality, can
we use $p \in \bold S(M)$?
\end{question}
\bigskip

\subsection{Indiscernibles} \
\bigskip

\begin{definition}
\label{a61}  
1) For an index model $I$ and model $M$ we say $\bold I = 
\langle \bar a_\eta:\eta \in I\rangle$ is
$(\Delta,n)$-indiscernible in $M$ when: $\bar a_\eta$ is a sequence from $M$
of length depending only $\tp_{\qf}(\eta,\emptyset,I)$ and such
that if the sequences $\bar \eta_1 = \langle \eta^1_\ell:\ell <
n\rangle,\bar \eta_2 = \langle \eta^2_\ell:\ell < n\rangle$ 
realize the same quantifier free type
in $I$ then $\bar a_{\bar\eta_1},\bar a_{\bar\eta_2}$ realize the
same $\Delta$-type in $M$ where:

\noindent
1A) For $\bar\eta = \langle \eta_\ell:\ell < n\rangle$ we let $\bar
a_{\bar\eta} := \bar a_{\eta_0} \char 94 \ldots \char 94
\bar a_{\eta_{n-1}}$.

\noindent
2) If $\Delta = \bbL(\tau_M)$ we may omit it; if $M \prec {\gC}
 = {\gC}_T$ we may omit $M$, we may write ``$<n$"
instead $n$ and omit $n$ meaning all $n$'s.

\noindent
3) \underline{Note}: saying $\bold I$ is $\{\varphi(\bar x_0,\dotsc,\bar
   x_{n-1})\}$-indiscernible in $M$ mean that we consider only $\bar
   a_{\eta_0} \char 94 \ldots \char 94 \bar a_{\eta_{n-1}},\ell
  g(\eta_\ell) = \ell g(\bar x_\ell)$, so do not allow to divide the
variables differently.

\noindent
4) $\langle \bar a_\eta:\eta \in I\rangle$ is continuously indiscernible in $M$
\when \,, say\footnote{for transparency} 
$\ell g(\bar a_\eta) = \zeta$ for every $\eta \in
 I$ and for any formula $\varphi(\bar x_0,\dotsc,\bar x_{n-1}) \in
\Gamma_{(\zeta)_n}$, with $\ell g(\bar x_\ell) = \zeta$
 for $\ell < n$ see \ref{z23}, there is a quantifier
free formula $\vartheta(y_0,\dotsc,y_{n-1}) \in \bbL(\tau_I)$
such that for every $\eta_0,\dotsc,\eta_{n-1} \in I$ we have $M
\models \varphi[\bar a_{\eta_0},\dotsc,\bar a_{\eta_{n-1}}]$
\Iff  \, $I \models \vartheta[\eta_0,\dotsc,\eta_{n-1}]$.

\noindent
5) We add ``over $B$" when we use the expansion $(M,b)_{b \in B}$
\end{definition}

\begin{claim}
\label{a64}   Let $T$ be dependent.

\noindent
1) Assume
\mn
\begin{enumerate}
\item[$(a)$]   $\Delta$ is a finite set of formulas
\sn
\item[$(b)$]  $M$ a model of $T$ and $A \subseteq M$
\sn
\item[$(c)$]  $I$ is a linear order
\sn
\item[$(d)$]  $\bold I = \langle a_{u,k,\ell}:\ell < n,k < k_\ell,
u \in [I]^\ell \rangle$ is indiscernible\footnote{so $[I]^{\le n}$ is
 defined as an index model naturally} over $A$
\sn
\item[$(e)$]  $\bar c \in {}^{\omega >}M$.
\end{enumerate}
\mn
\underline{Then} there is a finite subset $J$ of $I$ or of the completion
$\comp(I)$ of $I$ such that $\langle a_{u,k,\ell}:\ell < n,
k < k_\ell,u \in [I]^\ell
\rangle$ is $\Delta$-indiscernible over $A \cup \bar d$ above $J$. 

\noindent
2) Moreover, there is a bound on $|J|$ which depend just on $\Delta,
\langle k_\ell:\ell < n \rangle$ (and $T$), and so it is enough that
$\bold I$ is $\Delta_1$-indiscernible for appropriate finite
$\Delta_1$. 

\noindent
3) So for every $C \subseteq {\gC}$ there is $J \subseteq \comp(I)$ of 
cardinality $\le |C| + |T|$ such that 
$\bold I$ is indiscernible above $J$ over $A \cup C$.

\noindent
4) Let $I \in K_{p,\sigma}$, see Definition \ref{a73}(1) below.  If
   $\sigma$ is finite then parts (1),(2) holds.  Part (3) holds when
   we demand $J$ to be just of cardinality $\le |C| + |T| + \sigma$.
\end{claim}

\begin{PROOF}{\ref{a64}}
See \cite[\S3]{Sh:715}.
\end{PROOF}
\bn
More generally

\begin{definition}
\label{a66} 
Let ${\gk} = (K,\le_{\gk}) = (K_{\gk},\le_{\gk})$ be an a.e.c.
class of index models; normally $\le_{\gk}$ is $\subseteq$,
then we may write $K$.

\noindent
1) We say that the theory $T$ has the 
${\gk}-\theta$-indiscernibility property \when \,: if $I \in K$ (see
below) and the sequence $\bold I = \langle \bar a_\eta:\eta \in
I\rangle$ is indiscernible over $A$ in ${\gC}_T$ and $\bar b
\in {}^{\omega >}{\gC}$ \then \, there is $I_* \in K \,
   \le_{\gk}$-extending $I$ and subset $J$ of $I_*$ of cardinality
$< \theta$ such that: if $\eta_\ell = \langle \eta_{\ell,m}:m <
   n\rangle \in {}^n I$ for $\ell=1,2$ realizes the same quantifier
free type over $J$ in $I_*$ \then , the sequences 
$\bar a_{\bar\eta_\ell} := \bar a_{\eta_{\ell,0}} \char 94 \ldots \char
   94 \bar a_{\eta_{\ell,n-1}}$ for $\ell=1,2$ realize the same type
   over $A$.

\noindent
2) Writing ``${\gk}-(< \theta,n_*)$-indiscernible property" means
that above $n  \le n_*$.

\noindent
3) Writing ``${\gk}-(<\theta,\Delta,n_*)$-indiscernible property"
   means that above we restrict ourselves to the $\Delta$-type,
   i.e. which means that $\Delta \subseteq \{\varphi(\bar x_0;\bar
   x_1;\ldots;\bar x_{n-1};\bar y):\varphi(\bar x_0;\dotsc;\bar
   x_{n-1};\bar y ) \in \bbL(\tau_T);\bar y$ finite$\}$ and
   $\bar\eta_1,\bar\eta_2 \in {}^({\gC}_T)$ we use only
   $\varphi(\bar x_0,\dotsc,\bar x_{n-1};\bar y) \in \Delta$ such that
   $\ell g(\bar y) = \ell g(\bar b)$ and $\ell g(\bar x_m) = \ell
   g(\bar a_{\eta_{1,m}}) = \ell g(\bar a_{\eta_{2,m}})$.

\noindent
4) Writing ``${\gk}-\text{local}-(<\theta)$-indiscernible property"
means that ``${\gk}-(\theta,\Delta,n)$-indiscernible property",
and for every finite $\Delta$.

\noindent
5) Omitting $\theta$ means $\aleph_0$ for the local case, $|T|^+$ for
   the other case; and instead ``$<\theta^+$" we may write $\theta$.

\noindent
6) We say $I \in K$ is full\footnote{it is many reasonable to
   restrict ourselves to full $I$} \when \, for every $J \in K$ which
$\le_{\gk}$-extends $I$, every quantifier free type (in
finitely many variables) realized in $J$ is realized in $I$.

\noindent
7) We say $I \in K$ is locally full \when \, we replace above type
 by a formula.
\end{definition}

\begin{definition}
\label{a67}   1) An indiscernible sequence
$\langle \bar a_s:s \in I\rangle$ in ${\gC}_T$ is dependent
(in ${\gC}_T$) \when \, for every $\bar b \in {}^{\omega >}{\gC}$ it
satisfies the conclusion of \ref{a66} for $\kappa = |T|^+ +$ (the
number of quantifier free $(< \omega,\tau_I)$-types realized in $I$).

\noindent
1A) Above ``$\kappa$-dependent" means we use $\kappa$.

\noindent
2) If $I \in K_{\set}$, see \ref{a73} below, an 
indiscernible set $\bold I = \{\bar a_t:t \in I\}$ 
in ${\gC}$ is stable or sit stably \when \, it satisfies the
conclusion of \ref{a66}(1).

\noindent
2A) Above we say $\kappa$-stably when we use $\kappa$, superstably
when $\kappa = \aleph_0$.

\noindent
2B) An infinite indiscernible sequence $\bold I = \langle \bar a_t:t \in
I\rangle$ of $\zeta$-tuples is dependent when for every $\varphi =
\varphi(\bar x_{[\zeta]},\bar y)$ and $\bar b \in {}^{\ell g(\bar y)}
\gC$ there is a convex equivalent relation $E$ on $I$ with finitely
many equivalence classes such that $s E t \Rightarrow \gC \models
\varphi[\bar a_s,\bar b] \equiv \varphi[\bar a_t,\bar b]$.

\noindent
3) For indiscernible $\bold I = \langle \bar a_t:t \in I\rangle
\subseteq {}^\zeta {\gC}$ as in part (2) and $A_* \subseteq
{\gC}$ let $\Av(\bold I,A) = \{\varphi(\bar x_{[\zeta]},\bar a):
\bar b \in {}^{\omega >}A$ and ${\gC} \models \varphi[\bar a_t,\bar
b]$ holds for all but $< \aleph_0$ elements $t \in I\}$.

\noindent
4) For endless $I \in K_{\lin}$, see \ref{a73}, indiscernible
   sequence $\bold I = \langle \bar a_t:t \in 
I\rangle \subseteq {}^\zeta {\gC}$ and set $A$ 
let $\Av(\bold I,A) = \{\varphi(\bar x_{[\zeta]},
\bar b):\bar b \in {}^{\omega >}M$ and ${\gC} \models 
\varphi[\bar a_t,\bar b]$ for every $<_I$-large enough $t \in I\}$.

\noindent
5) We call the infinite indiscernible sequences
$\bold I,\bold J$ equivalent when $\Av(\bold I,A) = 
\Av(\bold J,A)$ for every $A$.

\noindent
6) Given endless indiscernible sequences $\bold I_\ell = \langle
   \bar a^\ell_t:t \in I_\ell\rangle$ for $\ell=1,2$, we say $\bold
   I_1,\bold I_2$ are immediate neighbours when $\bold I_\ell + \bold
   I_{3-\ell}$ naturally defined is an indiscernible sequence for some
   $\ell \in \{1,2\}$.  They are $n$-neighbours when there are $\bold
   J^*_0,\dotsc,\bold J^*_n$ such that $\bold J_k,\bold J_{k+1}$ are
endless indiscernible sequences which are immediate neighbours for
   $k=0,\dotsc,n-1$ and $\bold I_1 = \bold J_0,\bold I_2 = \bold
   J_n$.  Let being neighbours mean $n$-neighbours for some $n$.  
\end{definition}

\begin{discussion}
\label{a69}
\underline{Historical review for \S(1C)}:

Of course, Eherenfeucht-Mostowski \cite{EhMo56} use indiscernibles,
i.e. their models were generated by a sequence of indiscernibles.
Morley \cite{Mo65} prove that for $\aleph_0$-stable $T$: when $\lambda
= \mu$ is regular $\lambda \rightarrow_T (\lambda)_1$
which mean for any $a_\alpha \in {\gC}_T(\alpha < \lambda)$ for
some ${\cU} \subseteq \lambda$, of cardinality $\mu$ the sequence
$\langle a_\alpha:\alpha \in {\cU}\rangle$ is an indiscernible set,
using $\varphi(x,\bar b)$ of minimal rank such that $(\exists^\lambda
\alpha)(\varphi[a_\alpha,\bar b])$, see Definition \ref{a76}.  
The author \cite{Sh:2},\cite[III]{Sh:c},
got a parallel result for stable theory using e.g. Fodor lemma, as
 minimality does not work, when e.g. $\alpha < \lambda
\Rightarrow |\alpha|^{|T|} < \lambda$.  

Also for stable $T$:
\mn
\begin{enumerate}
\item[$(a)$]   if $\langle \bar a_\alpha:\alpha \in I\rangle$ is
an infinite indiscernible set, $I$ a set, i.e. with equality only
\begin{enumerate}
\item[$(\alpha)$]   $\varphi(\bar x,\bar b)$ can divide it only to
finite/co-finite sets, so we have average
\sn
\item[$(\beta)$]   for some $u \subseteq \lambda,|u| < \lambda,\langle
a_\alpha:\alpha \in I \backslash u\rangle$ is indiscernible over
$\bar b \cup\{\bar a_\alpha:\alpha \in u\} \cup b$.
\end{enumerate}
\end{enumerate}
\mn
On general models see \cite[\S5]{Sh:300a}.
Grossberg and the author suggest to classify first order $T$ by
$\lambda \rightarrow_T (\mu)_1$, see \ref{a76}(2)
this remains untraceable, see \cite[\S2]{Sh:702}.  
We can consider parallel to Erd\"os-Rado, see
Definition \ref{a76}(3).  This is proved for stable $T$ (and more
general context) in \cite[\S1]{Sh:300f}, e.g.
\mn
\begin{enumerate}
\item[$(b)$]   $[\lambda]^{\le 2} \rightarrow_T ([\mu]^{\le
2}])_\theta$ when $\lambda = (2^\mu)^+$ and $\mu \ge 2^{|T|}$, 
see \ref{a76}(4).
\end{enumerate}
\mn
For dependent $T$, the parallel to $(a)(\alpha)$ above is
in \cite[Ch.II,4.13,pg.77]{Sh:c} or \cite[3.2(1)]{Sh:715}
the parallel to $(\beta)$, is in Baldwin-Bendikt \cite{BlBn00} (not
seeing it is also \cite[3.2(3)]{Sh:715}).  For $I$ being essentially $I^{\le
n}$ the parallel of $(a)(\beta)$ in \cite[3.4=\ref{a64}]{Sh:715} here. 
Here we state also another generalization using end-homogeneity.

In \cite[\S3]{Sh:863} some advance was made for strongly stable
theories, $\beth_\delta \rightarrow_T (\lambda^+)_1$ 
when $\delta = (2^{|T|+\lambda})^+$,
also in \cite{Sh:863} it was suggested to look at infinite sequences
having better prospects for dichotomies and ``$T$ is $n$-dependent",
see more \cite[\S2]{Sh:886}.
\end{discussion}

\begin{question}  Is the combination reasonable?
\end{question}

\begin{definition}
\label{a73}
1) Let $K_{p,\sigma}$ be the class of $I = (I,<_I,P^I_i)_{i < \sigma}$
 where $<_I$ is a linear order of $I$ and 
$\langle P^I_i:i < \sigma\rangle$ a partition of $I$, (as
in many cases we disabuse our notation not distinguished the (index) model and
its universe).

\noindent
1A) $K_{q,\sigma}$ is the class of $I = (I,<_I,P^I_i)_{i < \sigma}$
where $<_I$ is a linear order of $I$ and $P_i$ a unary predicate.  If
$\sigma = 0$ we may omit it and so if $I$ is endless this means 
$I \in K_{\lin}$.

\noindent
2) For $I \in K_{p,\sigma}$, let $E_I = \{(s,t):s,t \in P^I_i$ for 
some $i < \sigma\}$; so $E_I$ an equivalence relation on $I$
with $\le\sigma$ equivalence classes.   

\noindent
3) Let ${\gk}_e$ be the class of $(I,<,E)$ where $<$ is a linear
order on $J$ and $E$ an equivalence relation on $I$.

\noindent
4) $K_{\set,\sigma}$ is the class of $(I,P^I_i)_{i <\sigma}$
where $\langle P^I_i:i < \sigma\rangle$ a partition, if $\sigma = 1$
we may omit $\sigma$ and $P^I_0$.
\end{definition}

\begin{remark}
\label{a74}   So by \ref{a64}(1) this case is covered, i.e. 
if $T$ is dependent \then \, it has the
$K_{p,\sigma}$-indiscernibility property.
\end{remark}

\begin{observation}
\label{a75} 
1) If $T$ is independent \then \, the conclusion of \ref{a74} fails.

\noindent
2) But there is $T$ which is unstable, but have the
$K_{\set}-\aleph_0$-indiscernible property, e.g. any
   expansion of the theory of linear order.

\noindent
3) If $T$ is a dependent theory, \then \,  it has the
   $K_{\set}$-indiscernible property (see \cite[\S1]{Sh:715}).

\noindent
4) Trivially $T$ has the $K_{\set}$-indiscernible property
\Iff \, for every $n$, every infinite indiscernible set $\bold I =
\{\bar a_\alpha:\alpha < \lambda\}$ of $n$-tuples in ${\gC}_T$
is stable (in ${\gC}_T$, see Definition \ref{a67}(2)).
\end{observation}

\begin{definition}
\label{a76}  
1) For a linear order $I$, we say that 
$\langle \bar a_t:t \in I\rangle$ is an $n$-end-homogeneous
   over $A$ \when \, if $m \le n$ and $t(0,\ell) <_I t(1,\ell) <_I
\ldots <_I t(m-1,\ell)$ for $\ell=1,2$ then the sequence $\bar
   a_{t(0,1)} \char 94 \ldots \char 94 \bar a_{t(m-1,1)}$ and $\bar
   a_{t(0,2)} \char 94 \ldots \char 94 \bar a_{t(m-1,2)}$ realize
   the same type over $\cup\{\bar a_t:t <_I t(0,1)$ and $t <_I
   t(0,2)\} \cup A$.

\noindent
1A) Replacing $n$ by ``$<n$" has the obvious meaning (and allow
$m=\omega$), $u \in [\lambda]^{<n}$.

\noindent
2) Let $\lambda \rightarrow_T (\gamma)_\sigma$ means that:
   if $\bar a_\alpha \in {}^\sigma{\gC}$ for $\alpha < \lambda$
   and $\langle \bar a_\alpha:\alpha < \lambda\rangle$ is
end-homogeneous \then \, for some $u \subseteq \lambda$, of
   order type $\gamma$ the sequence $\langle \bar a_\alpha:\alpha \in
   u\rangle$ is indiscernible.

\noindent
3) Let $\lambda \rightarrow_T (\gamma)^n_\sigma$ when: if $\bar a_u
\in {}^\sigma{\gC}$ for $u \in [\lambda]^n$ \then \, for
   some ${\cU} \subseteq \lambda$ of order type $\gamma$ the
   sequence $\langle \bar a_n:u \in [{\cU}]^n\rangle$ is 
$(<n)$-indiscernible.  Similarly with $<n$ instead of $n$.

\noindent
4) Fix ${\gk} = (K_{\gk},\le_{\gk})$ an a.e.c. of index
   models.  Then $I \rightarrow_{\gk,T}(J)_\theta$ for $I,J \in K_{\gk}$
   is defined naturally.
\end{definition}

\begin{question}
\label{a77}  Find reasonable sufficient conditions on $T$ for the
following:  for every $\sigma,\gamma$ the cardinality 
min$\{\lambda:\lambda \rightarrow_T (\gamma)_\sigma\}$ is quite small 
or at least $< \min\{\lambda:\lambda \rightarrow 
(\gamma)^{< \omega}_{\sigma_1}$ where 
$\sigma_1 = 2^{|T|+\sigma}\}$.  (Of course, Erd\"os-Rado theorem gives
lower bounds, see \cite{EHMR}.)
\end{question}
\bn
We may consider
\begin{question}
\label{a78}
\underline{The Strong Indiscernibility Question}

\noindent
1) Give sufficient conditions on $T$ for the following; where  
$|T| \le \theta$ and $\kappa = \cf(\kappa) > 2^\theta$ (or just
large enough).  For some $k_1 < \omega,T$ has the strong $k_1$-indiscernibility
existence property for $(\kappa,\theta)$, meaning: 
if $\gamma(*) < \theta^+$ and $\bar a_\alpha \in 
{}^{\gamma(*)}{\gC}$ for $\alpha < \kappa$ and $\bold I 
= \langle \bar a_\alpha:\alpha <\kappa\rangle$ is $k$-end-homogeneous
\then \, for some unbounded ${\cU} \subseteq \kappa$ the sequence
$\langle \bar a_\alpha:\alpha \in {\cU}\rangle$ is indiscernible.

\noindent
2) Similarly for ``$T$ has the $k$-strong$^+$ indiscernibility 
existence property for $\kappa$" which means that above $\bold I$ 
is mod clubs locally indiscernible.
\end{question}

\begin{discussion}
\label{a79}  
1) We will be glad even for weaker versions, anything 
better than Erd\"os cardinal.

\noindent
2) If $T$ is $\omega$-independent we are no better off than in set
   theory (because we allow $\omega$-tuples).

\noindent
3) Independent theories can satisfy strong versions of \ref{a78},
see example below.
\end{discussion}

\begin{definition}
\label{a80}  Assume $\kappa$ is regular uncountable.

We say $\bold I = \langle \bar a_\alpha:\alpha < \kappa\rangle$ is mod
clubs locally indiscernible \when \, for some club $E$ of $\kappa$
and $I \in K_{q,|T|}$ expanding $(\kappa,<)$ the sequence $\langle
\bar a_\alpha:\alpha \in I \rest E \rangle$ is locally indiscernible, see
\ref{a61}(4), this means that for every finite $\Delta$ there is a
finite $\tau_\Delta \subseteq \tau_I$ such that $\langle \bar
a_\alpha:\alpha \in (I \rest E \rest \tau_\Delta)\rangle$ is
$\Delta$-indiscernible.

Similarly $n$-indiscernible, $n$-end-homogeneous.
\end{definition}

\noindent
Recall (\cite[\S2]{Sh:886})
\begin{definition}
\label{a81}
1) We say $T$ is $2$-independent or 2 $\times$ independent
\when \ , we can find an independent sequence of formulas of the form
$\langle \varphi(\bar x,\bar b_n,\bar c_m):n,
m < \omega \rangle$ in ${\gC} = {\gC}_T$ or just in some model of $T$.

\noindent
2) ``$T$ is $2$-dependent" (or dependent/2)
means the negation of $2$-independent (see \cite[\S5 (H)]{Sh:863}). 

\noindent
3) We say $\varphi(\bar x,\bar y_0,\dotsc,\bar y_{n-1})$ is
$n$-independent (for $T$) \when \, in ${\gC}_T$ we can, for each
$\lambda < \bar \kappa$, find $\bar a^\ell_\alpha \in {}^{\ell g(\bar
y_\ell)}({\gC}_T)$ for $\alpha < \lambda,\ell < n$ such that the
sequence $\langle \varphi(\bar x,\bar a^0_{\eta(0)},\dotsc,\bar
a^{n-1}_{\eta(n-1)}):\eta \in {}^n \lambda\rangle$ is an independent
sequence of formulas.

\noindent
4) $T$ is $n$-independent \when \, some formula $\varphi(\bar x,\bar
   y_0,\bar y_1,\dotsc,\bar y_{n-1})$ is $n$-independent.

\noindent
5) $T$ is $n$-dependent (or dependent/$n$) \when \, it is not $n$-independent.
\end{definition}

\begin{example}
\label{a83} 
1) For a first order $T$ which is 3-independent assuming 
$\lambda = (2^\mu)^+$ we can find $n < \omega$ and 
$\bar d_\alpha \in {}^n{\gC}$ for $\alpha \le \lambda$ such that
$\langle \bar d_\alpha:\alpha \le \lambda\rangle$ is 
one-end-homogeneous, equivalently $\tp(\bar d_\alpha,
\cup\{\bar d_\beta:\beta < \alpha\})$ is
increasing, but for no unbounded ${\cU} \subseteq \lambda$
   and even no ${\cU}$ of cardinality $\mu^+$ is
$\langle \bar d_\alpha:\alpha \in {\cU}\rangle$ an indiscernible
sequence.

\noindent
2) For a first order $T$ which is $(k+2)$-independent and $\lambda =
   (2^\mu)^+$ we can find $n < \omega$ and $\bar d_\alpha \in
   {}^n{\gC}$ for $\alpha \le \lambda$ such that $\langle \bar
   d_\alpha:\alpha \le n\rangle$ is $k$-end-homogeneous for no $u
   \subseteq \lambda$ of cardinality $\mu^+$ is $\langle \bar
   d_\alpha:\alpha \in {\cU}\rangle$ an indiscernible sequence.
\end{example}

\begin{example}
\label{a85}  $T_{\rd}$, the theory of random graphs
has the strongly one-indiscernibility property.
\end{example}

\begin{definition}
\label{a86}  We say $T$ has bounded/medium/large
$k$-directionality \when \,: if $\bold I = \langle 
\bar a_\alpha:\alpha < \delta\rangle$ has a $k$-type-increasing ($=
k$-end-homogeneous) \then \, $\uf(\bold I)$ is defined as in Definition
\ref{a41}, replacing $\|M\|$ by $|\delta|$. 
\end{definition}

\begin{remark}  We may consider replacing well orderings by other
classes of index models.
\end{remark}

\begin{question}
\label{a88}  
Can we answer \ref{a77} or \ref{a78} when $T$ has bounded or at least medium
$k$-directionality for some $k$.
\end{question}

\begin{question}
\label{a89} 
1) Can we characterize $\bold U_{T,m,\Delta}
= \{(n_1,n_2):n_1 \rightarrow (n_2)_{T,\Delta,m}\}$? for finite
$\Delta,m$ when $T$ dependent?

\noindent
2) Similarly for $T$ \, $k$-dependent?

\noindent
3) If $T$ is $k$-dependent is there $k_1$ such that: if for $T$ is
$k$-dependent, $m < \omega,\Delta \subseteq \bbL(\tau_T)$ is
finite, then $\beth_{k_1}(n) \rightarrow (n)_{T,\Delta,m}$ for
every $n$ large enough?

\noindent
4) As in (2) for $k=1$? (i.e. $T$ dependent).
\end{question}

\begin{question}
\label{a90}   Assume $\Delta = \bbL(\tau_T)$ or
$\Delta$-finite, $p = p(\bar x)$ a $(\Delta,m)$-type over $A,\ell
g(\bar x) < \theta$ and every subset of $p$ of cardinality 
$< \kappa$ is realized in $M$.  Can we find $q
\in \bold S^m_\Delta(A)$ extending $p(\bar x)$ such that every subset
of $q$ of cardinality $< \kappa$ is realized in $M$?
\end{question}

\begin{conjecture}
\label{a91}  
Assume $M$ is a saturated model of a cardinality
$\kappa > |T|$ of a dependent complete $T$.

\noindent
1) If $p \in \bold S(M)$ \then  \, there is an indiscernible sequence
 $\bold I = \langle a_\alpha:\alpha < \kappa\rangle$ in $M$ such
that $p = \Av(I,M)$.

\noindent
2) Similarly for $p \in \bold S^\theta(M)$ where $\theta < \kappa$. 

See more on this in \S6.
\end{conjecture}
\bigskip

\subsection{Limit Models and Generic Pairs} \
\bigskip

\begin{conjecture}
\label{a93}  
We can characterize ``$M$ is $\kappa$-saturated" parallely to stable $T$, e.g. 
$M$ is a $\kappa$-saturated model of $T$
iff it is $|T|^+$-saturated and every indiscernible sequence $\langle
a_\alpha:\alpha < \delta\rangle$ of elements in $M$ of length $\delta
< \kappa$ can be continued and similarly for cuts.
\end{conjecture}

\begin{conjecture}
\label{a95}
\underline{The Generic Pair Conjecture}  

\noindent
Assume \footnote{the ``$2^\lambda = \lambda^+$" is 
just for making the formulation more transparent}
 $\lambda = \lambda^{<\lambda} > |T|,2^\lambda = \lambda^+,M_\alpha 
\in \EC_{\lambda,1}(T)$ is $\prec$-increasing continuous for $\alpha <
\lambda^+$ with $\cup\{M_\alpha:\alpha < \lambda^+\} \in 
\EC_{\lambda^+,\lambda^+}(T)$ being saturated.  \Then \, $T$ is 
dependent iff for some club $E$ of $\lambda^+$ for all 
pairs $\alpha < \beta < \lambda^+$ from $E$
of cofinality $\lambda^+,(M_\beta,M_\alpha)$ has the same isomorphism type.
\end{conjecture}

\begin{remark}   We proved in \cite{Sh:900} the ``structure" side,
i.e. the implication $\Rightarrow$ in \ref{a95} when 
$\lambda = \kappa$ is measurable, on the non-structure side of
\ref{a95}, \ref{a96}, see \cite{Sh:877}, \cite{Sh:906}.
It seemed natural to assume that the first order theories of
such pair is complicated if $T$ is independent and 
``understandable" for dependent of $T$, but this is not so, see
Kaplan-Shelah \cite{KpSh:946}.
\end{remark}

\begin{conjecture}
\label{a96}
\underline{The Unique Limit Model Conjecture}  Assume if $T$ is
dependent, $|T| < \lambda = \lambda^{< \lambda}$ and $\lambda^+ =
2^\lambda,\sigma = \cf(\sigma) < \lambda$.  If $\langle
M_\alpha:\alpha < \lambda^+\rangle$ is an increasing continuous of
models of cardinality $\lambda$ with $\lambda^+$-saturated union
\then \, for some club $E$ of $\lambda^+$, all the models in
$\{M_\alpha:\alpha$ is from $E$ and has cofinality
$\sigma\}$ are pairwise isomorphic.
\end{conjecture}
\bigskip

\centerline {$* \qquad * \qquad *$}
\bigskip

\noindent
\underline{Completions}:

For linear order the notion of completion is very important, so it is
natural to try to generalize it to dependent theories (if we accept
thesis \ref{y11}).  Note that for stable theories as every
type $p \in \bold S(A)$ is definable by formulas with parameters from
$A$, this is not so necessary (and ${\gC}^{\eq}$ is a much
less radical extension).

\begin{definition}
\label{a97.2}  
1) $\ai({\gC})$ is the set $\{\Av(\bold I,{\gC}):\bold I$ 
an infinite indiscernible sequence of finite tuples in ${\gC}\}$.

\noindent
2) $\nsp_\mu(\gC)$ is the set of $p \in \bold S^m(\gC)$ which does not
split over some $A \subseteq \gC$ of cardinality $< \mu$.

\noindent
3) $\fs_\mu({\gC})$ is the set of $p \in \bold S^m(\gC)$ 
which is $\Av(D,{\gC})$ where $D$ is an ultrafilter on ${}^m A$
   for some $m < \omega$ (or more) and $A \subseteq {\gC}$ a set
of cardinality $< \mu$.  If $\mu = \infty$, (i.e. $\|{\gC}\|$) we may omit it.
\end{definition}

\begin{thesis}
\label{a97.4}  
So the types we considered as understandable, a base
for analysis are $\fs_\mu({\gC})$ or $\nsp_\mu(\gC),\mu$ for small enough 
(hopefully $|T|^+$) and $\ai({\gC})$.
\end{thesis}

\begin{question}
\label{a97.6}  Is it reasonable to add in the completion of 
${\gC},\nsp_\infty({\gC})$ or just $\nsp_\mu({\gC}) \cup \ai({\gC})$.
\end{question}

\begin{discussion}
\label{a97.8}  
0) So our main theorems say that any $p \in \bold S(M),M \in
   \EC_{\lambda,\lambda}(T)$ is definable over $\le \beth_\omega +
   |T|$ elements from $\ai(\gC) \cup \fs_{< \beth_\omega}(\gC)$.

\noindent
1) We may prefer not to analyze complete types but
ultrafilters, i.e. the $\bar d_{\bold x}$ and $\bar c_{\bold x,i}$ are
in the full completion!  But there is no parallel to the ``recounting
of types" as there are dependent $T$ with large directionality.
However, given $D \in \uf({}^m M)$ we may choose an $\|N\|^+$-saturated
elementary extension $N$ of $M$ and let $p = \Av(D,N)$, so
analyzing $p$ is very close.

It is still reasonable that in view of later developments we may
prefer to use the ultrafilter version.

Recall that if we succeed to use $\mu = |T|^+$ for countable $T$, then
we can always use eventually indiscernible sequences, see below.  
This may be not just aesthetically nicer but helpful.
Anyhow allowing constant though not so small $\mu$, will give us the
asymptotic behaviour.

\noindent
2) To clarify our intension let us consider the class of linear
   orders.  We like to deal with the class of complete linear orders;
   or at least $(< \kappa)$-complete.  If $I$ is the completion of a
   $\kappa$-saturated dense linear order, then it is natural to add predicate
$P_{x,\theta}(x \in \{\leeft,\riight\},\theta = \cf(\theta) 
< \kappa)$ such that
\mn
\begin{enumerate}
\item[$(a)$]  $(\alpha) \quad I \models P_{\leeft,\theta}(a)$
\Iff \, $\cf(I_{< a}) = \theta$
\sn
\item[${{}}$]  $(\beta) \quad I \models P_{\riight,\theta}[a]$
\Iff \, the inverse of $I_{> a}$ has cofinality $\theta$
\sn
\item[$(b)$]  $I_1 \le_{\gk} I_2$ iff
\sn
\item[${{}}$]  $(\alpha) \quad I_1 \subseteq I_2$
\sn
\item[${{}}$]  $(\beta)_{\leeft} \quad$ if $I_1 \models 
P_{\leeft,\theta}(a)$ then $(I_1)_{<a}$ is cofinal in $(I_2)_{<a}$
\sn
\item[${{}}$]  $(\beta)_{\riight} \quad$ if $I_1 \models
P_{\riight,\theta}[a]$ then $(I_1)_{>a}$ is cofinal in the inverse of
$(I_2)_{>a}$.
\end{enumerate}
\end{discussion}

\begin{definition}
\label{a97.11} 
1) We say that $\bold I = \langle \bar a_t:t \in I\rangle$ 
is an eventually indiscernible sequence \when \,: 
$I$ is an endless linear order, 
$\ell g(\bar a_t)$ for $t \in I$ is constant and finite for
transparency, and for every finite set $\Delta \subseteq 
\bbL(\tau_T)$ there is $t(\Delta) \in I$ such that $\langle \bar a_t:t
\in I_{\ge t(\Delta)}\rangle$ is a $\Delta$-indiscernible sequence
(over $A$).
\end{definition}
\newpage

\section {Decompositions of types} 

We define $\bold x \in \pK_{\kappa,\bar\mu,\theta}$, which is a
partial analysis of $\tp(\bar d_{\bold x},M_{\bold x})$, it is related
to $K_{\kappa,\theta}$ from \cite{Sh:900} but $B_{\bold x}$, which has
cardinality $< \kappa$ there corresponds to
$B^+_{\bold x} = B_{\bold x} \cup\{\bold I_{\bold x,i}:i \in 
u_{\bold x}\}$ here.  Moreover, the set
$B_{\bold x}$ is of cardinality $< \mu_0$ rather than $<
\kappa$ but in this section 
we really do not use this.  We define
``$\bold x$ is $\sigma$-active in $i < i_{\bold x}$", which cannot occur too
many times.  We define $\qK'_{\kappa,\bar\mu,\theta}$, those for which
we ``exhaust the possible activities", this set is dense; and the related
$\qK_{\kappa,\bar\mu,\theta}$ is 
suppose to be the class of such $\bold x$'s in which we have 
fuller analysis.  For the case $\kappa = \mu$ we have
$\qK'_{\kappa,\mu,\theta} = \qK_{\kappa,\mu,\theta}$ so in this case
$\qK_{\kappa,\mu,\theta}$ is dense and we define solvability, all are 
related to \cite{Sh:900}.  But not so dealing with $(\bar
\mu,\theta)$-sets, over which the situation is similar to the one for
stable $T$ and any set $\subseteq \gC_T$; note that $B^+_{\bold x}$ is a
$(\bar\mu,\theta)$-set when $\bold x \in \pK_{\kappa,\bar\mu,\theta}$
is  smooth.  Central here are the definitions of 
similarity of decompositions and their smoothness (points which are
meaningless in \cite{Sh:900}) and we point out their basic
properties.  Those later ones indicate the possible advantages of Definition
\ref{b05}, i.e. the use of indiscernibles.  Generally, we shall concentrate on
the case $\kappa = \mu_2 \gg \mu_1 = \mu_0 > \theta \ge |T|$ so may not state
claims in full generality concerning this point.
\bigskip

\subsection{Decompositions - the basics} \
\bigskip

\begin{convention}
\label{b02}  
1) In clause (i) of Definition \ref{b05} below 
we have three options, the choice is $\iota_{\bold x} \in 
\{0,1,2\}$, usually the choice does not matter and in those cases we suppress
$\iota$; so far we can use only $\iota_{\bold x} = 2$.  Usually the
set $v$ can be a well ordering and even an ordinal but in
 disjoint amalgmation in $\sK^\oplus_{\kappa,\bar\mu,\theta}$
we shall need anti-well orderings whereas in proving density for $\qK'$ it is
natural to use just well-ordering.

\noindent
2)  Also $\bar c_{\bold x}$ consists of finite sequences and
sometimes we use $\bar c \subseteq \bar d$, normality, see Definition
\ref{b5}(7); we may demand that always $\bar c \subseteq \bar d$.  
We can work in $\gC^{\eq}$ hence use $\bar c$ a
sequence of singletons but this is immaterial in Definition \ref{b3}.  

\noindent
3) The notation is sometimes best understood as in the case when $v$ is a
set of ordinals, as the case ``$v$ an ordinal" is our prototype
so abusing notation we let, e.g. $v \cap i = \{j \in v:j <_v i\}$.
\end{convention}

\begin{definition}
\label{b05} 
Assume $\bar \mu = (\mu_2,\mu_1,\mu_0)$ and $\lambda \ge \kappa \ge \mu_0,
\lambda \ge \mu_2 \ge \mu_1 \ge \mu_0 > \theta$
but if not said otherwise in addition $\cf(\mu_2) > \theta,
\kappa \ge \mu_2$ and even $\kappa = \mu_2,\mu_1 = \mu_0$; 
usually $\theta \ge |T|$.  

We let $\pK_{\lambda,\kappa,\bar\mu,\theta}$
 be the class of objects $\bold x$ consisting of:
\mn
\begin{enumerate}
\item[$(a)$]   $M \prec {\gC}$ which is
$\kappa$-saturated of cardinality $\lambda$
\sn
\item[$(b)$]   $B = \cup\{B_i:i \in v \backslash u\}$ and each
$B_i \subseteq M$ is of cardinality $< \mu_0$
\sn
\item[$(c)$]   $\bar d \in {}^{\theta^+>}({}^{\omega >}{\gC})$ 
or even $\bar d \in {}^w({}^{\omega >}{\gC})$ 
where\footnote{This is useful when we
like to amalgamate such objects, but usually we may ignore this.  We
may work in $\gC^{\eq}$ and then use $d_i$ instead of $d_i$.
Similarly for the $\bar c_i$'s.} $w$
 is a linear order (e.g. a set of ordinals) of cardinality 
$\le \theta$; we may write $w$ as $\ell g(\bar d)$ or 
$\Dom(\bar d)$ but we usually write $i<j$ instead of $i <_w j$
and $w \cap j$ or $w_{<j}$ for $\{i \in w:i <_w j\}$, similarly for $v,u$ below
\sn
\item[$(d)$]   $\bar c = \langle \bar c_i:i \in v\rangle  \in
{}^v({}^{\omega >} \gC)$ which sometimes is treated as 
$(\ldots \char 94 \bar c_i \char
94 \ldots)_{i \in v}$ 
so\footnote{could demand $\bar c_i \in {}^\omega {\gC}$ or even
$\bar c_i \in {}^{\theta^+>}{\gC}$, in this work usually it does not
matter but not always; if we do this in \ref{b3} we can make $\bar c_{\bold
x,i} = (\ldots \char 94 \bar b_{i,\ell} \char 94 \ldots)_{\ell < \sigma}$
in Case A and a parallel demand in Case B} $\bar c_i \in
{}^{\omega >}{\gC}$; where $v$ is a linear order of cardinality
$\le \theta$; we may write $v = \ell g(\bar c)$ or $v = \Dom(\bar c)$
\sn
\item[$(e)$]  $u \subseteq v$
\sn
\item[$(f)$]   $\bar\kappa = \langle \kappa_i:i \in u \rangle$ such
that $\kappa_i = \cf(\kappa_i) \in [\mu_1,\mu_2)$
\sn
\item[$(g)$]   $\bar{\bold I} = \langle \bold I_i:i \in u\rangle$
\sn
\item[$(h)$]   $\bold I_i = \langle \bar a_{i,\alpha}:\alpha <
\kappa_i\rangle$ is an indiscernible sequence in $M$ for $i \in u$.
\sn
\item[$(i)$]   $\iota_{\bold x} \le 2$ and for $i \in v \backslash u$,
\newline

\underline{Case 0}:  $\iota_{\bold x} = 0,\tp(\bar c_i,M + 
\Sigma\{\bar c_j:j <i\})$ does not split over $B_i$
\smallskip

\underline{Case 1}: $\iota_{\bold x} = 1:\tp(\bar c_i,M + \Sigma\{\bar
c_j:j <i\})$ does not locally split over
\newline

\hskip30pt  $B_i + \Sigma\{\bar c_j:j < i\}$, see Definition
\ref{b1} below
\smallskip

\underline{Case 2}:  $\iota_{\bold x} = 2$, the type above is
finitely satisfiable in $B_i$.  

In short we may say $\tp(\bar c_i,M +
\Sigma\{\bar c_j:j < i\})$ does not $\iota_{\bold x}$-split over $B_i$
\sn
\item[$(j)$]   for $i \in u$, the type $\tp(\bar c_i,M + 
\Sigma\{\bar c_j:j < i\})$ is $\Av(\bold I_i,M_{\bold x} + \Sigma\{\bar
c_j:j< i\})$ hence $\ell g(\bar a_{i,\alpha}) = \ell g(\bar c_i)$ for
$\alpha < \kappa_i$.
\end{enumerate}
\end{definition}

\begin{definition}
\label{b1}  
1) We say that the type $p(\bar x)$
locally splits over $A$ \when \, there is $\varphi = \varphi(\bar
x,\bar y) \in \bbL(\tau_t)$ such that for every finite $\Delta
\subseteq \{\varphi:\varphi = \varphi(\bar y,\bar z) \in
\bbL(\tau_t)\}$ there are formulas $\varphi(\bar x,\bar b),\neg
\varphi(\bar x,\bar c) \in p$ where $\bar b,\bar c$ realize the
same $\Delta$-type over $A$.

\noindent
2) If $p \in \bold S^\varepsilon(B)$ does not split over $A \subseteq
B$ let the scheme of $p$ be the function $H$ defined by: if
$\varphi(\bar x_{[\varepsilon]},\bar y) \in \bbL(\tau_T)$ and $\bar b
\in {}^{\ell g(\bar y)}B$ then $H(\varphi(\bar x,\bar y)$, 
$\tp(\bar b,A))$ is the truth value of $\varphi(\bar x,\bar b) \in p$.
\end{definition}

\begin{definition}
\label{b3}  
In Definition \ref{b05} we say $i \in v_{\bold x}$ is 
$\sigma$-active (in $\bold x$) \when \, $1 \le \sigma \le
\aleph_0$ and $\sigma = 1 \Rightarrow i \notin u_{\bold x}$
 and (using notation of \ref{b5}(1) below; the default
value for $\sigma$ is 1):
\bigskip

\noindent
\underline{Case 1}:  $\sigma = 1$

We can find $\bar b_{i,0},\bar b_{i,1}$ such that
\mn
\begin{enumerate}
\item[$(a)$]  $\bar b_{i,0},\bar b_{i,1}$ realize the same
type over $\bar c_{\bold x,<i} + M_{\bold x}$, see \ref{b5}(1)
\sn
\item[$(b)$]  $\bar b_{i,0},\bar b_{i,1}$ realizes distinct
types over $\bar d_{\bold x} + \bar c_{\bold x,<i} + M_{\bold x}$
\sn
\item[$(c)$]  $\bar c_{\bold x,i} = \bar b_{i,0} \char 94 \bar b_{i,1}$.
\end{enumerate}
\bn
\underline{Case 2}:  $\sigma \ge 2$ and $i \notin u_{\bold x}$

We can find $\langle \bar b_{i,n}:n < \sigma\rangle$ such that:
\mn
\begin{enumerate}
\item[$(a)$]   $\tp(\bar b_{i,n},\cup\{\bar b_{i,m}:m \in
(n,\sigma)\} + \bar c_{\bold x,<i} + M_{\bold x}\})$ is
$\subseteq$-decreasing with $n$ and\footnote{in the other cases the
parallel statement follows} 
does not $\iota_{\bold x}$-split over 
$B_{\bold x,i}$; (i.e. does not split over $B_{\bold x,i}$ 
if $\iota_{\bold x} = 0$, does not 
locally split over $B_{\bold x,i}$
if $\iota_{\bold x} = 1$ and is finitely satisfiable in $B_{\bold x,i}$ if
$\iota_{\bold x} = 2$) 
\sn
\item[$(b)$]  $\tp(\bar b_{i,\ell},M_{\bold x} + \bar c_{\bold x,<i} + \bar
d_{\bold x})$ for $\ell=0,1$ are distinct
\sn
\item[$(c)$]   $\bar c_{\bold x,i} = \bar b_{i,0} \char 94 \bar b_{i,1}$.
\end{enumerate}
\bn
\underline{Case 3}:  $\sigma \ge 2$ and $i \in u_{\bold x}$

We can find $\bar b_{i,\alpha,n}(\alpha \le \kappa_{\bold x,i},n < 2)$
such that:
\mn
\begin{enumerate}
\item[$(a)$]   $\tp(\bar b_{i,\kappa_{\bold x,i},n},\cup\{\bar
a_{i,\kappa_{\bold x,i},m}:m \in (n,2)\} + \bar c_{\bold x,<i} + 
M_{\bold x})$ does not split over $\cup\{\bar a_{\bold x,i,\alpha}:
\alpha < \kappa_{\bold x,i}\}$
\sn
\item[$(b)$]   $\tp(\bar b_{i,\kappa_{\bold x,i},\ell},
M_{\bold x} + \bar c_{\bold x,<i} + \bar d_{\bold x})$ 
for $\ell=0,1$ are distinct
\sn
\item[$(c)$]   $\bar c_i = \bar b_{i,\kappa_{\bold x,i},0} \char 94
\bar b_{i,\kappa_{\bold x},1}$ and 
$\bar a_{\bold x,i,\alpha} = \bar b_{i,\alpha,0} 
\char 94 \bar b_{i,\alpha,1}$ for $\alpha < \kappa_{\bold x,i}$
\sn
\item[$(d)$]   $\langle \bar b_{i,\alpha,n}:(\alpha,n) \in 
\kappa_{\bold x,i} \times 2 \rangle$
is an indiscernible sequence where $\kappa_{\bold x,i} \times 2$ is
ordered lexicographically.
\end{enumerate}
\end{definition}

\begin{remark}
\label{b4}  
1) We shall return to this and to \ref{b20} and in \ref{q3}.

\noindent
2) We can use $\bar c_{\bold x,\ne i} = \langle \bar c_{\bold x,j}:j
   \in v_{\bold x} \backslash \{i\}\rangle$ instead $\bar c_{\bold
   x,<i}$ in Definition \ref{b3}, as in \cite{Sh:900}.

\noindent
3) In Case (3) of Definition \ref{b3} note that it follows that for
   every $\bar e \in {}^{\omega >}(M_{\bold x})$ for some $\beta <
   \kappa,\langle \bar b_{i,\alpha,n}:(\alpha,n) \in
   [\beta,\kappa_{\bold x,i}) \times 2\rangle$ is an indiscernible
   sequence over $\bar e$.
\end{remark}

\begin{definition}
\label{b5}  
1) For $\bold x \in \pK_{\lambda,\kappa,\bar\mu,\theta}$ let 
$\bold x= (M_{\bold x},\ldots)$ so $w_{\bold x} = w,v_{\bold x} =
v,u_{\bold x} = u,\bar c_{\bold x} = \bar c[\bold x] =
\bar c,\bar d_{\bold x} = \bar d[\bold x] = \bar d,
\kappa_i = \kappa_{\bold x,i} = \kappa(\bold x,i),
B_{\bold x,i} = B_i,B_{\bold x} = B$, etc., and let 
$B^+_{\bold x} = \cup\{\bold I_{\bold x,i}:i \in u_{\bold x}\} 
\cup B_{\bold x}$.
Let $\bar c_{<i} = \bar c[<i] = \bar c_{\bold x,<i} 
= (\ldots \char 94 \bar c_{\bold x,j} \char 94 \ldots)_{j<i},\bar c[u']
= \bar c_{\bold x}[u'] = (\ldots \char 94 \bar c_{\bold x,j} \char 94
\ldots)_{j \in u'}$ and
$\bar c_{\ne i} = \bar c_{\bold x,\ne i} := (\ldots \char 94 \bar
c_{\bold x,j} \char 94 \ldots)_{j \in v_{\bold x} \backslash \{i\}}$
for $i \in u_{\bold x}$ and for $u' \subseteq v_{\bold x}$.
We may write $\bar c,\bar d$ instead of $\bar
c_{\bold x},\bar d_{\bold x}$ when confusion is unlikely as there is only one
$\bold x$ around, in particular avoiding using, e.g. $\bar x_{\bar d_{\bold
x}}$; also we may write $\bar c[\bold x],\bar c[\bold x,<i]$, etc.
Let $\ga_{\bold x} = \{\kappa_{\bold x,i}:i \in u_{\bold x}\}$.

\noindent
1A) For $i \in v_{\bold x}$ let $D_i = 
D_{\bold x,i}$ be such that tp$(\bar c_{\bold x,i},\bar c_{\bold x,<i} +
M_{\bold x}) = \Av(D_i,\bar c_{\bold x,<i} + M_{\bold x})$ 
where $D_i$ is an ultrafilter on
${}^{\ell g(\bar c(\bold x,i))}(B_{\bold x,i})$ if $\iota_{\bold x} =
2 \wedge i \in v_{\bold x} \backslash u_{\bold x}$ and on 
$\bold I_{\bold x,i}$ if $i \in
u_{\bold x}$ such that $\alpha < \kappa_{\bold x,i} \Rightarrow \{\bar
a_{\bold x,i,\beta}:\beta \in (\alpha,\kappa_{\bold x,i})\} \in D_i$;
but only $D'_i := D_i \cap \Deef_{\ell g(\bar c(\bold x,i))}
(\Dom(D_i))$ matters so we normally use it.

\noindent
2) Concerning $\pK_{\lambda,\kappa,\bar\mu,\theta}$, 
omitting $\lambda$ means ``for some $\lambda$"; omitting
$\mu_0$ means $\mu_0 = \mu_1$, omitting also
$\mu_2$ means $\mu_2 = \kappa$, then we may write $\mu$ instead of
$\mu_1$ and of $\bar \mu$;
writing $*$ instead of $\mu_1$ means $\mu_1 =
(\theta^+ + |T|^+ + \beth_\omega)$; omitting
$\lambda,\kappa,\mu_0,\mu_1,\mu_2,\theta$ means for some such cardinals.  
Similarly in parallel definitions later.

\noindent
3) We say $i \in v_{\bold x}$ is active in $\bold x$ \when \, it is
 $\sigma$-active for some $\sigma$, equivalently for $\sigma = 1$.

\noindent
3A) We say $i \in v_{\bold x}$ is active in $\bold x$ over $u$ \when
\, $i \in v_{\bold x},u \subseteq v_{\bold x}(<i)$ and in Definition
\ref{b1} we replace $\bar c_{\bold x,<i}$ by $\bar c_{\bold x}[u]$.
Similarly in the other versions.

\noindent
4) We say $i \in v_{\bold x}$ is 
\underline{strongly active} in $\bold x$ when it is
$\aleph_0$-active.

\noindent
5) We say that $i \in v_{\bold x}$ is $(\sigma,\Delta)$-active in $\bold
   x$ \when \, $1 \le \sigma \le \aleph_0$ and $\Delta \subseteq
\Gamma^1_{\bold x,i} := \{\varphi:\varphi = 
\varphi(\bar x_{\bar d},\bar x_{\bar c[<i]},\bar y,\bar z) \in 
\bbL(\tau_T)$ and $\bar y,\bar z$ are finite$\}$ and in Definition 
\ref{b3} we replace clause (b), in all cases by
\mn
\begin{enumerate}
\item[$(b)'$]   for\footnote{No real loss if we replace $\bar y$ by
$x_{\bar c_{\bold x,i}}$.  Also no real loss if we omit $\bar z$,
absorbing $\bar a$ into $\bar c_i$ by cosmetic manipulations}
 some $\varphi(\bar x_{\bar d},\bar x_{\bar c[<i]},
\bar y,\bar z) \in \Delta$ we have
\newline

${\gC} \models \varphi[\bar d_{\bold x},\bar c_{\bold x,<i},
\bar b_{i,\kappa_{\bold x,i},0},\bar a] \wedge \neg \varphi
[\bar d_{\bold x},\bar c_{\bold x,<i},
\bar b_{i,\kappa_{\bold x,i,1}},\bar a]$ for some 
$\bar a \in {}^{\ell g(\bar z)}(M_{\bold x})$.
\end{enumerate}
\mn
5A) For $v \subseteq v_{\bold x}$ and $\Delta \subseteq 
\bbL(\tau_T)$ let $\Gamma^1_{\bold x,\Delta,i,v} = \{\varphi:\varphi \in
\Gamma^1_{\bold x,i}$ and $\varphi \in \Delta$ but 
$\bar x_{\bar c[\bold x,<i]}$ is replaced by 
$\bar x_{\bar c \rest (v \cap i)}\}$.

\noindent 
6) Let $M_{[\bold x]}$ be $(M_{\bold x})_{[B^+_{\bold x} + \bar c_{\bold x} 
+ \bar d_{\bold x}]}$, see Definition \ref{z5}.

\noindent
7) We say $\bold x$ is normal when Rang$(\bar c_{\bold x}) \subseteq
\text{ Rang}(\bar d_{\bold x})$, pedantically $\cup\{\Rang(\bar
c_{\bold x,i}:i \in v_{\bold x}\} \subseteq \cup\{\Rang(\bar d_{\bold
x,i}:i \in w_{\bold x}\}$; note that usually there is no loss in
assuming it.

\noindent
8) Let $\Gamma^1_{\bold x}$ be $\{\varphi:\varphi =
\varphi(\bar x_{\bar d},\bar x_{\bar c},\bar y) \in \bbL(\tau_T)$
for some $\bar y$, finite if not said otherwise$\}$.

\noindent
9) $\Gamma^0_{\bold x} = \{\varphi:\varphi = 
\varphi(\bar x_{\bar d,\eta},\bar y)$ and $\eta \in {}^n 
\ell g(\bar d_{\bold x})$ for some $n\}$, used in particular
when $\bold x$ is normal, see \ref{b7}(4) below.

\noindent
10) For $v \subseteq v_{\bold x},w \subseteq w_{\bold x}$ let 
$\bold x_{<v,w>} = (M_{\bold x},\bar B_{\bold x} \rest (v \backslash
u_{\bold x}),\bar c_{\bold x}
\rest v,\bar d_{\bold x} \rest w,\bar{\bold I} \rest (v \cap
u_{\bold x}))$, but if $w = \ell g(\bar d_{\bold x})$ we may omit it.

\noindent
11) We say $\bold x$ is essentially well ordered \when \, $\{i \in
u_{\bold x}:\kappa_i = \kappa_j\}$ is well ordered (by $<_v$) for
    each $j \in u_{\bold x}$.
\end{definition}

\begin{notation}
\label{b7}  
0) We may write $\bar d,\bar c$ instead of $\bar d_{\bold x},\bar
   c_{\bold x}$ when $\bold x$ is clear from the context, (usually in
   subscripts).

\noindent
1) $u,v,w$ are linear orders, members are $i,j$ but, e.g. $w \cap j =
   w_{<j} := \{i \in w:i < j\}$.

\noindent
2) If $v_1 \subseteq v_2$ let $[v_1,v_2) = \{i \in v_2:j <_{v_2} i$ for 
every $j \in v_1\}$.

\noindent
3) $\ell g(\bar d) = \Dom(\bar d)$, etc.

\noindent
4) $\bar d_\eta = \langle \bar d_{\eta(\ell)}:\ell < \ell g(\eta)\rangle$ if
$\eta $ is a function from $\ell g(\eta)$ to $\ell g(\bar d)$.

\noindent
5) $\bar d_{\bold x,\eta} = \langle \bar d_{\bold x,\eta(\ell)}:\ell < \ell
   g(\eta)\rangle$, see Definition \ref{b5}(1).

\noindent
6) $\bar x_{\bar d,\eta} = \bar x_{\bold x,\eta} = \langle
\bar x_{\bar d_{\eta(i)}}:i < \ell g(\eta)\rangle$ when $\bar d
= \bar d_{\bold x}$.

\noindent
7) $\bar x_{\bar c,\eta}$ is defined similarly.

\noindent
8) $\tp_\varphi(\bar d_{\bold x},\bar c_{\bold x}  
\dotplus A) := \{\varphi(x_{\bar d},\bar c_{\bold x},\bar b):
\bar b \in {}^{\ell g(\bar y)}A$ and ${\gC} \models
\varphi(\bar d_{\bold x},\bar c_{\bold x},\bar b)\}$ when $\varphi =
\varphi(\bar x_{\bar d[\bold x]},\bar x_{\bar c[\bold x]},\bar y)$.

\noindent
8A) We may above use $\pm \varphi,\Delta$ and/or $\varphi = 
\varphi(x_{\bar d,\eta},\bar x_{\bar c,\eta},\bar y)$.

\noindent
9) $v_2 = v_1 +1$ is defined naturally.
\end{notation}

\begin{definition}
\label{b9}  
Let $\bold x \in \pK_{\kappa,\bar\mu,\theta}$.

\noindent
1) We say that $\bar e$ solves
$(\bold x,\bar\psi,A)$ or $\bar\psi$-solves $(\bold x,A)$,
or $\bar\psi$-solves $\bold x$ over $A$, (pedantically we should add $\theta$) 
\underline{when}: 
\mn
\begin{enumerate}
\item[$(a)$]   $\bar e \in {}^\theta(M_{\bold x})$
\sn
\item[$(b)$]  $A \subseteq M_{\bold x}$
\sn
\item[$(c)$]  $\bar\psi = \langle \psi_\varphi = 
\psi_\varphi(\bar x_{\bar d},\bar x_{\bar c},\bar y_{[\theta]}):\varphi \in 
\Gamma_{\bar\psi}\rangle$ where $\Gamma_{\bar\psi} \subseteq 
\Gamma^1_{\bold x}$, recalling \ref{b5}(8)
\sn
\item[$(d)$]   $\psi_\varphi(\bar x_{\bar d},\bar c_{\bold x},\bar e) 
\vdash \{\varphi(\bar x_{\bar d},\bar c_{\bold x},\bar a):
\bar a \in {}^{\ell g(\bar y_\varphi)}A$ and ${\gC} \models 
\varphi[\bar d_{\bold x},\bar c_{\bold x},\bar a]\}$ for $\varphi \in 
\Gamma_{\bar \psi}$
\sn
\item[$(e)$]   ${\gC} \models \psi_\varphi[\bar d_{\bold x},
\bar c_{\bold x},\bar e]$.
\end{enumerate}
\mn
1A) We say that $\psi(x_{\bar d},\bar c_{\bold x},\bar e)$ solves
$(\bold x,A,\varphi)$ \when \, $\varphi = \varphi(\bar x_{\bar
d},\bar x_{\bar c},\bar y) \in \Gamma^1_{\bold x}$ and $\bar e \subseteq
M_{\bold x}$ and $\psi(\bar x_{\bar d},\bar c_{\bold x},\bar e) \vdash
\{\varphi(\bar x_{\bar d},\bar c_{\bold x},\bar b):b \in {}^{\ell
g(\bar y_\varphi)}A$ and ${\gC} \models \varphi[\bar d_{\bold x},\bar
c_{\bold x},\bar b]\}$ and $\gC \models \psi_\varphi[\bar d_{\bold
x},\bar c_{\bold x},\bar e]$.

\noindent
1B) We say that $\psi(\bar x_{\bar d},\bar x_{\bar c},\bar z)$ solves
$(\bold x,A,\varphi)$ \when \, $\psi(\bar x_{\bar d},\bar c_{\bold
x},\bar e)$ solves $(\bold x,A,\varphi)$ for some $\bar e \in {}^{\ell
g(\bar z)}M$. 

\noindent
1C) We let $\vartheta_{\bold x,\varphi}(\bar x_{\bar c},\bar z,\bar y)
= \vartheta_{\bold x,\varphi,\psi}(\bar x_{\bar c},\bar
z,\bar y)$ where $\varphi = \varphi(\bar x_{\bar d},\bar x_{\bar
c},\bar y) \in \Gamma^1_{\bold x}$ and $\psi = \psi(\bar x_{\bar d},\bar
x_{\bar c},\bar z)$  be $(\forall \bar x_{\bar d})(\psi(\bar
x_{\bar d},\bar x_{\bar c},\bar z) \rightarrow \varphi(\bar x_{\bar
d},\bar x_{\bar c},\bar y))$; we usually omit $\bold x$, being clear
from the context and similarly $\varphi$.

\noindent
1D) We say $\bar\psi = \langle \psi_\varphi(x_{\bar d},\bar x_{\bar
c},\bar e):\varphi \in \Gamma^1_{\bar\psi}\rangle$ solves $(\bold x,A)$
when $\Gamma^1_{\bar\psi} \subseteq \Gamma^1_{\bold x}$, so we usually
write $\Gamma^1_{\bar\psi}$ instead of $\Gamma_{\bar\psi}$ to stress
this (similarly in other cases), and for every $\varphi \in
\Gamma^1_{\bar\psi},\psi_\varphi(\bar x_{\bar d},\bar x_{\bar c},\bar
e)$ solves $(\bold x,A,\varphi)$.
We say $\bar\psi$ solves $(\bold x,A)$ when $\bar\psi = \langle
\psi_\varphi:\varphi \in \Gamma^1_{\bar\psi}\rangle,\psi_\varphi =
\psi_\varphi(\bar x_{\bar d},\bar x_{\bar c},\bar y_{[\theta]})$ 
and some $\bar e$ solves $(\bold x,\bar\psi,A)$.

\noindent
2) We say $\bar\psi$ is full for $\bold x$ when $\Gamma^1_{\bar\psi} =
\Gamma^1_{\bold x}$.
Omitting $\bar\psi$ means for some $\bar\psi$ full for $\bold x$.

\noindent
2A) Let ``$\tp(\bar x_{\bar d},\bar c_{\bold x},\bar e) \vdash 
\tp(\bar d,\bar c_{\bold x} \dotplus M_{\bold x})$
 according to $\bar\psi$" mean clause (d) of part (1) with
$\Gamma^1_{\bold x}$ instead of $\Gamma_{\bar\psi}$.

\noindent
3A) We say $\bar\psi$ illuminates $\bold x \in 
\pK_{\kappa,\bar\mu,\theta}$ \when \, $\bar\psi$ is as in
clause (c) of part (1) and for every
$A \subseteq M_{\bold x}$ of cardinality $< \mu_2$ some $\bar e$ does
solve $(\bold x,\bar\psi,A)$.

\noindent
3B) We say $\psi(\bar x_{\bar d},\bar x_{\bar c},\bar x_{[\theta]})$
illuminates $(\bold x,\varphi)$ \when \, 
$\varphi = \varphi(\bar x_{\bar d},\bar x_{\bar c},\bar y_\varphi) \in
\Gamma^1_{\bold x}$ and the above holds with $(\varphi,\psi)$
standing for $(\varphi,\psi_\varphi)$.

\noindent
4) We say $\bar e$ solves $(\bold x,A)$ \when \, for some $\bar\psi$ which
is full for $\bold x,\bar e$ solves $(\bold x,\bar\psi,A)$.
\end{definition}

\begin{remark}
\label{b11}  
0) Note that we use ``illuminate" rather than ``solve" when we
quantify on $A$.

\noindent
1) For the case $\iota_{\bold x} = 2,\mu_1 = \aleph_1,\theta = \aleph_0$,
i.e. for countable $T$, we can replace ``$\tp(\bar c_{\bold x,i},
\bar c_{\bold x,<i} + M_{\bold x})$ is
finitely satisfiable in some countable $B_{\bold x,i} \subseteq
M_{\bold x}$" by: $p$
is the average of an eventually indiscernible sequence $\bold I = \langle
\bar a_n:n < \omega\rangle$ from $M_{\bold x}$ 
which means that for every finite $\Delta$,
some end-segment is $\Delta$-indiscernible, see Definition
\ref{a97.11}.   Also $\Av(\bold I,A)$ is well defined.

\noindent
2) However, we cannot replace eventually 
indiscernible by indiscernible,
e.g. for $= \Th(\bbR),\bbR$ the real field, there is an eventual
indiscernible $\bold I = \langle a_n:n < \omega\rangle$ such that $a_n
\ge n$; the cut it defines cannot be defined by a really indiscernible
sequence, (well of length less than the saturation). 

\noindent
3) We can characterize when an eventually indiscernible sequence is
equivalent, (see Definition \ref{b30}(2)) to an indiscernible
sequence, but this does not always occur, by the example above.

\noindent
4) Being equivalent is well defined for eventually indiscernible
   sequences as their averages are well defined.

\noindent
5) Usually no harm is done when below in \ref{b14}(1)(b) we add
$``w_{\bold x}$ is an initial segment of $w_{\bold y}$".  Similarly in
\ref{b14}(1)(d). 
\end{remark}

\begin{definition}
\label{b14}  
1) We define a two-place relation
$\le_1$ on $\pK:\bold x \le_1 \bold y$ \when \,:
\mn
\begin{enumerate}
\item[$(a)$]  $M_{\bold x} = M_{\bold y}$
\sn
\item[$(b)$]  $\bar d_{\bold x} = \bar d_{\bold y} \rest w_{\bold x}$
and $w_{\bold x} \subseteq w_{\bold y}$ as linear orders
\sn
\item[$(c)$]   $u_{\bold x} = u_{\bold y} \cap v_{\bold x}$ and
$u_{\bold x} \subseteq u_{\bold y}$ as linear orders
\sn
\item[$(d)$]   $\bar c_{\bold x} = \bar c_{\bold y} \rest v_{\bold x}$
and $v_{\bold x} \subseteq v_{\bold y}$ as linear orders.
\end{enumerate}
\mn
2) We define $\le_2$ similarly strengthening clause (b) to
\mn
\begin{enumerate}
\item[$(b)^+$]   $\bar d_{\bold x} = \bar d_{\bold y}$.
\end{enumerate}
\mn
3) If $\bold x \le_1 \bold y$ and $\varphi = 
\varphi(\bar x_{\bar d[\bold x]},\bar x_{\bar c[\bold x]},\bar y) 
\in \Gamma^1_{\bold x}$ \then \, we may identify it with the 
$\varphi(\bar x_{\bar d[\bold y]},\bar x_{\bar c[\bold y]},\bar y) \in
\Gamma^1_{\bold y}$ naturally.

\noindent
4) For $\bold x \in \pK$ and $\varphi = 
\varphi(\bar x_{\bar d},\bar x_{\bar c},\bar y) \in \Gamma^1_{\bold x}$ let
$\supp(\varphi)$ be the pair $(w,v)$ such that $w \subseteq \ell
   g(\bar d_{\bold x}),v \subseteq \ell g(\bar c_{\bold x})$ are
   minimal (so finite) such that $\varphi \equiv \varphi(\bar x_{\bar
   d \rest u},\bar x_{\bar c \rest v},\bar y)$, moreover the omitted
   variables are dummy (= does not appear in $\varphi$, not
   just ``immaterial for satisfaction").

\noindent
4A) Similarly for $\varphi = \varphi(\bar x_{\bar d},\bar x_{\bar c},
\bar x'_{\bar d},\bar x'_{\bar c},\bar y)$, we define
$\supp(\varphi) = (w_1,v_1,w_3,v_3)$; used in \ref{c5}.

\noindent
5) For $\Delta \subseteq \Gamma^1_{\bold x}$ let $\supp_{\bold
  x}(\Delta)$ be the pair $(\cup\{w:(w,v) = \supp(u)$ for some $\varphi \in
\Delta\},\cup\{v:(w,v) = \supp(\varphi)$ for some $\varphi \in \Delta\})$. 
\end{definition}

\begin{definition}
\label{b16}  
0)
\mn
\begin{enumerate}
\item[$(A)$]  For $\bold x \in \pK_{\lambda,\kappa,\bar\mu,\theta}$
let $\Gamma^2_{\bold x} = \{\varphi:\varphi = \varphi(\bar x_{\bar
d},\bar x_{\bar c},x'_{\bar d},\bar x'_{\bar c},\bar y) \in
\bbL(\tau_T)$ so $\ell g(\bar x'_{\bar c}) = \ell g(\bar x_{\bar
c}),\ell g(x'_{\bar d}) = \ell g(\bar x_{\bar d}),\bar y$
finite$\}$;used in \ref{c5}(1).
\sn
\item[$(B)$]  For $\bold x \in \pK_{\lambda,\kappa,\bar\mu,\theta}$
let $\Gamma^3_{\bold x} := \{\varphi:\varphi 
= \varphi(\bar x_{\bar d},\bar x_{\bar c},\bar y_{[\theta]},\bar z)
\in \bbL(\tau_T),\bar z$ finite$\}$, (used in
(3A), close to $\Gamma^1_{\bold x,i}$, see Definition \ref{b5}(5),(5A)).
\end{enumerate}
\mn
1) Let $\qK'_{\lambda,\kappa,\bar\mu,\theta}
[\Delta]$ be the class of $\bold x 
\in \pK_{\lambda,\kappa,\bar\mu,\theta}$ such that for no 
$\bold y \in \pK_{\lambda,\kappa,\bar\mu,\theta}$ do
 we have $\bold x \le_2 \bold y$ and $\bold y$ is $\Delta$-active in some $i
\in v_{\bold x} \backslash v_{\bold y}$ over $v_{\bold x,<i}$,
i.e. $v_{\bold x} \cap v_{\bold y,<i}$, see\footnote{in many places
it suffices to use [$v_{\bold y},v_{\bold x}]$} Definition \ref{b3},
 \ref{b5}(3A),(5); if $\Delta = \Gamma^1_{\bold x}$ we may write
$\qK'_{\lambda,\kappa,\bar\mu,\theta}$; similarly below.

\noindent
1A) We define $\qK''_{\lambda,\kappa,\bar\mu,\theta}[\Delta]$ similarly but
restricting  ourselves to the case $v_{\bold y} = v_{\bold x}+1$.

\noindent
2) Let $\qK_{\lambda,\kappa,\bar\mu,\theta}$ be the class of
$\bold x \in \pK_{\lambda,\kappa,\bar\mu,\theta}$ such that for every
$A \in [M_{\bold x}]^{< \kappa}$ some $\bar e$ solves $\bold x$, 
see Definition \ref{b9}(4).

\noindent
3) Let $\qK^\odot_{\lambda,\kappa,\bar\mu,\theta}$ be the class of
triples $\bold n = (\bold x,\bar\psi,r)$ such that $\bold x \in
\qK_{\lambda,\kappa,\bar\mu,\theta}$ and $\bar\psi = \langle
\psi_\varphi:\varphi \in \Gamma_{\bar\psi} \subseteq 
\Gamma^1_{\bold x}\rangle$ illuminates
$\bold x$ and $r = r(\bar x_{\bar c_{\bold x}},\bar x_{\bar d_{\bold x}},
\bar y_{[\theta]})$ is a type over $M_{\bold x}$ such that: 
for every $A \subseteq M_{\bold x}$ of cardinality $< \kappa$
there is a tuple $\bar e$ from $M_{\bold x}$ such that:

$\bar e$ solves $(\bold x,\bar\psi,A)$ and 
the sequence $\bar c_{\bold x} \char 94 \bar d_{\bold x} \char 94 
\bar e$ realizes $r$.  

In this case we may say
$\bar e \subseteq M_{\bold x}$ solves $(\bold x,A,\bar\psi,r)$ or
solves $(\bold n,A)$.
If not said otherwise, $r$ is a type over $\emptyset$; in this case
we say $\bold n$ is pure.  

\noindent
3A) Let $\qK^\oplus$ be the class of $\bold n = (\bold x,\bar\psi,r)
\in \qK^\odot$ such that $\Gamma^1_{\bar\psi} = \Gamma^1_{\bold x}$.

\noindent
3B) Let $\qK^\otimes_{\kappa,\bar\mu,\theta}$ be the class of
triples $(\bold x,\bar\psi,r)$ such that\footnote{What
is the difference with part (3)?  Here in the end, $\bar e$ appears in
$\varphi$}:
\mn
\begin{enumerate}
\item[$(a)$]  $\bold x \in \qK_{\lambda,\kappa,\bar\mu,\theta}$
\sn
\item[$(b)$]  $r = r(\bar x_{\bar c},\bar x_{\bar d},
\bar y_{[\theta]}),\bar\psi = \langle \psi_\varphi:\varphi \in
\Gamma^3_{\bold x}\rangle$ satisfy
\sn
\begin{enumerate}
\item[$\bullet_1$]  $\psi_\varphi = \psi(\bar x_{\bar d},\bar x_{\bar
c},\bar y_{[\theta]})$
\sn
\item[$\bullet_2$]  for every $A \subseteq M_{\bold x}$ of cardinality
$< \kappa$ for some $\bar e \in {}^\theta(M_{\bold x})$ we have: if
$\varphi \in \Gamma^3_{\bold x}$ then $\gC \models \psi_\varphi[\bar
d_{\bold x},\bar c_{\bold x},\bar e]$ and $\psi_\varphi(\bar x_{\bar
d},\bar c_{\bold x},\bar y_{[\theta]}) \vdash \{\varphi(\bar x_{\bar
d},\bar c_{\bold x},\bar y_{[\theta]},\bar b):\gC \models ``\varphi[\bar
d_{\bold x},\bar c_{\bold x},\bar e,\bar b]"$ and $\bar b \in
{}^{\ell g(\bar y)}A\}$.
\end{enumerate}
\end{enumerate}
\mn
4) We say $\bar e$ universally solves the triple 
$(\bold x,\bar\psi,r) \in \qK^\odot_{\lambda,\kappa,
\bar\mu,\theta}$ \when \,
for every $A \in [M_{\bold x}]^{< \kappa}$ there is $\bar e'$ 
as in part (3) such that $\bar e,\bar e'$ realizes the
 same type over $\bar c_{\bold x} + \bar d_{\bold x} + A$, see
\ref{b9}(4) and Theorem \ref{d33}.

\noindent
4A) Similarly for $(\bold x,\bar\psi,r) \in
\qK^\oplus_{\lambda,\kappa,\bar\mu,\theta}$ but ``$\bar e'$ is as in
part (3A)".

\noindent
5) We define the partial orders $\le_1,\le_2$ on $\qK',\qK,
\qK^\odot,\qK^\oplus$ naturally.
\end{definition}

\begin{remark}
\label{b17f}
Concerning Definition \ref{b16}(3),(4) note that
$\qK^\odot$ is used in the end and in $(*)_5$ of the proof of
\ref{d33} and in the proof of \ref{d39} only.
\end{remark}

\begin{observation}
\label{b18}  
Let $\ell=0,1,2$.

\noindent
0) If the type $p(\bar x)$ is finitely satisfiable in $A$ \then \,
 $p(\bar x)$ does not locally split over $A$ 
and this in turn \underline{implies} that $p(\bar x)$ does not split 
over $A$ (hence the corresponding implications hold 
for the variants of Definition \ref{b05}).

\noindent
1) $\le_\ell$ is a partial order on $\pK$.

\noindent
2) If $\bar{\bold x} = \langle \bold x_\varepsilon:\varepsilon <
\delta\rangle$ is $\le_\ell$-increasing sequence of members of 
$\pK_{\lambda,\kappa,\bar\mu,\theta}$ and 
$\delta < \theta^+$ is a limit ordinal 
\then \, $\bar{\bold x}$ has a $\le_\ell$-lub, essentially the union,
naturally defined and it belongs to $\pK_{\lambda,\kappa,\bar\mu,\theta}$. 

\noindent
3) If $\bold x \in \pK_{\kappa,\bar\mu,\theta}$ and we define $\bold
   y$ like $\bold x$ replacing $\bar d_{\bold x}$ by $\bar d_{\bold x}
   \char 94 \bar c_{\bold x}$, \then \, $\bold y \in 
\pK_{\kappa,\bar\mu,\theta}$ is normal and $\bold x \in 
\qK_{\kappa,\bar\mu,\theta} \Leftrightarrow \bold y \in 
\qK_{\kappa,\bar\mu,\theta}$ and $\bold x \in
\qK'_{\kappa,\bar\mu,\theta} \Leftrightarrow 
\bold y \in \qK'_{\kappa,\mu,\theta}$ and $\bold x \in
\qK''_{\kappa,\bar\mu,\theta} \Leftrightarrow \bold y \in
   \qK''_{\kappa,\bar\mu,\theta}$ and $\bold x \le_1 \bold y$ and no loss if
systematically we use only normal $\bold x$ (except when we like e.g. $\ell
 g(\bar d_{\bold x})$ to be finite).

\noindent
4) Parts (1),(2) apply\footnote{But for
 $\qK^\oplus_{\kappa,\bar\mu,\theta},\qK^\otimes_{\kappa,\bar\mu,\theta},
\qK^\odot_{\kappa,\bar\mu,\theta}$ parts (1),(3) are O.K. but part (2)
 is a different, harder matter; for $\qK''_{\kappa,\bar\mu,\theta}$ all
 are not clear.} also to $\qK,\qK'$.

\noindent
5) Assume $\kappa > \theta \ge |T|$ and $\bar \mu$ are as in 
\ref{b05}.  If $M$ is $\kappa$-saturated, $w$ a linear order of
cardinality $< \theta^+$ and $\bar d \in {}^w({}^{\omega >}{\gC})$, 
\then \, for one and only one $\bold x \in \pK_{\kappa,\bar\mu,\theta}$ we 
have $M_{\bold x} =M,
\bar d_{\bold x} = \bar d,v_{\bold x} = \emptyset$ hence 
$\bar c_{\bold x} = \langle \rangle,B_{\bold x} = \emptyset$.

\noindent
6) Let $\bold x \in \pK_{\lambda,\kappa,\bar\mu,\theta}$.
\mn
\begin{enumerate}  
\item[$(a)$]   If $\cf(\mu_0) > \theta$ \then \, 
$B_{\bold x} = \cup\{B_{\bold x,i}:i \in v_{\bold x} 
\backslash u_{\bold x}\} \subseteq M_{\bold x}$ has cardinality $<
\mu_0$.
\sn
\item[$(b)$]    If $\cf(\mu_2) > \theta$ \then \, $\cup\{\bold
   I_{\bold x,i}:i \in u_{\bold x}\}$ has cardinality $< \mu_2$.
\sn
\item[$(c)$]   If $\cf(\mu_2) > \theta$, (hence $\mu_0 = 
\mu_2 \Rightarrow \cf(\mu_0) > \theta)$ \then \, 
also $|B^+_{\bold x}| < \mu_2$.
\sn
\item[$(d)$]   Always $|B_{\bold x}| \le \mu_0,|B^+_{\bold x}| \le \mu_2$.
\end{enumerate}
\end{observation}

\begin{PROOF}{\ref{b18}}  
Easy, concerning part (2) for $\ell=1,2$ note that the union, it is
not uniquely defined as if $i \in v_{\bold x_\varepsilon} \backslash
u_{\bold x_\varepsilon},\varepsilon < \delta$ then 
$\langle B_{\bold x_\zeta,i}:\zeta \in [\varepsilon,\delta)\rangle$ is not necessarily
constant, but we can use any one of them.  Similarly for $i \in
u_{\bold x_\varepsilon}$.
\end{PROOF}

\begin{claim}
\label{b20}  
1) If $\theta \ge |T|$ \then \, in $\pK_{\kappa,\bar\mu,\theta}$
there is no $\le_2$-increasing sequence $\langle 
\bold x_\varepsilon:\varepsilon < \theta^+\rangle$ such that: 
if $\varepsilon < \theta^+$ then $\bold x_{\varepsilon +1}$ is active 
in some $i \in v(\bold x_{\varepsilon +1}) 
\backslash v(\bold x_\varepsilon)$.

\noindent
2) For finite\footnote{if we restrict ourselves to 
$\bar d_{\bold x} \restriction u$ for some finite 
$u \subseteq \ell g(\bar d_{\bold x})$ then any finite $\Delta
\subseteq \bbL(\tau_T)$ is O.K.}
 $\Delta \subseteq \Gamma^1_{\bold x}$, \underline{there is} 
$n_\Delta = n_{\Delta,T} < \omega$ such that 
there is no $\le_2$-increasing chain $\langle
   \bold x_\ell:\ell \le n_\Delta\rangle$ of members of
$\pK_{\kappa,\bar\mu,\theta}$ such that $x_{\ell +1}$ is
   $\Delta$-active in some $i \in [v(\bold x_\ell),v(\bold x_{\ell +1}))$.

\noindent
3) In part (1), the sequence may be just $\le_1$-increasing if
   $\{\varepsilon < \theta^+:\bar d_{\bold x_\varepsilon} = 
\cup\{\bar d_{\bold x_\zeta}:\zeta < \varepsilon\}$ is a stationary subset of
$\theta^+$.
\end{claim}

\begin{PROOF}{\ref{b20}}
A similar proof appears in Case 1 of the proof of
\ref{k4} or see \cite[2.8=tp25.33]{Sh:900} recalling Definition
\cite[2.6=tp25.32]{Sh:900}.   
\end{PROOF}

\begin{claim}
\label{b22}  
1) If $\bold x \in \pK_{\kappa,\bar\mu,\theta}$ \then \, there is
$\bold y \in \qK'_{\kappa,\bar\mu,\theta}$ such that $\bold x \le_2 \bold y$.

\noindent
2) If the finite $\Delta$ is as in \ref{b20}(2) and $\bold x \in
\pK_{\kappa,\bar\mu,\theta}$ \then \, there is $\bold y \in
\pK_{\kappa,\bar\mu,\theta}$ such that $\bold x \le_2 \bold y$ 
and $v_{\bold y} \backslash v_{\bold x}$ is finite and there is no $\bold z
\in \pK_{\kappa,\bar\mu,\theta}$ such that $\bold y \le_2
\bold z$ and some $i \in [v_{\bold y},v_{\bold z})$ or just $i \in
v_{\bold z} \backslash v_{\bold x}$ is $\Delta$-active in $\bold z$.

\noindent
2A) If above we restrict $\bold z$ to the case $v_{\bold z} = v_{\bold y}
+1$, \then \, we can demand $v_{\bold y} \subseteq v_{\bold x} +
n_\Delta$ when $n_\Delta$ is\footnote{and see ind$(\Delta)$ in \S3}
from \ref{b20}(2).

\noindent
2B) In part (2), if we restrict the assumption to the case
$v_{\bold y} < v_{\bold x} + \omega$, i.e. $v_{\bold y} =
v_{\bold x} + n$ for some $n$ \then \, this is O.K. provided
that we restrict the conclusion to the case 
$v_{\bold y} \trianglelefteq v_{\bold z}$ (actually just $v_{\bold y}
\subseteq v_{\bold z} \wedge v_{\bold x} \trianglelefteq v_{\bold z}$).

\noindent
3) If $\bold x \in \qK'_{\kappa,\kappa,\theta}$ or just
 $\bold x \in \qK'_{\kappa,\bar\mu,\theta}$ and $\mu_0 = \kappa$
\then \, \footnote{the ``$\mu_0 = \kappa$" is of course undesirable,
 but eliminating it is the reason of much of the work here.} $\bold x
 \in \qK_{\kappa,\bar\mu,\theta}$, that is, 
for every $A \in [M_{\bold x}]^{< \kappa}$ some $\bar\psi$ solves 
$(\bold x,A)$, see Definition \ref{b9}(1D).

\noindent
4) [Local version\footnote{we may use \ref{b22}(3),(4) replacing $A
\in [M]^{< \mu_0}$ by $\in [M]^{< \kappa}$ as the definition of $\qK$.}]; 
if $\varphi \in \Gamma^1_{\bold x}$ and $\bold x \in \qK'_{\kappa,\bar
\mu,\theta}$ or just is as $\bold y$ in \ref{b22}(2) or
   just \ref{b22}(2A), \then \, for 
every $A \in [M_{\bold x}]^{< \mu_0}$ 
there is $\psi = \psi(\bar x_{\bar d},\bar x_{\bar c},\bar z) 
\in \bbL(\tau_T),\bar z$ finite and $\bar e \in {}^{\ell g(\bar z)}M$
such that $\psi(\bar x_{\bar d},\bar c_{\bold x},\bar e)$ solves
$(\bold x,A,\varphi)$.

\noindent
5) Assume $\iota_{\bold x} = 2$.  
The inverse of part (3) holds, i.e. if $\bold x \in 
\pK_{\kappa,\bar\mu,\theta}$ and $\bar{\bold x} \in
\qK_{\kappa,\bar\mu,\theta}$, i.e. for every 
$A \subseteq M_{\bold x}$ of cardinality $< \kappa$ there is a 
solution \then \, $\bold x \in \qK'_{\kappa,\bar\mu,\theta}$ 
(and see \ref{c27}(2)).

\noindent
6) Assume $\mu_0 = \kappa$.  Assume $\bold x \in 
\pK_{\kappa,\bar\mu,\theta}$ but $\bold x \notin \qK'_{\kappa,\bar\mu,\theta}$
\then \, there is a pair $(\bold y,\varphi)$ such that:
\mn
\begin{enumerate}
\item[$(a)$]   $\bold x \le_1 \bold y \in \pK_{\kappa,\bar\mu,\theta}$
\sn
\item[$(b)$]  $\varphi = \varphi(\bar x_{\bar d},\bar x_{\bar
c},\bar y,\bar z) \in \Gamma^1_{\bold x}$
\sn
\item[$(c)$]   $\bold y$ is $\{\varphi\}$-active in some $i \in
[v_{\bold x},v_{\bold y})$, so $\bar c_i = \bar c_{i,0} \char 94 \bar
c_{i,1},\bar c \subseteq \Rang(\bar c_{\bold y,<i}),\ell
g(\bar c_{i,0}) = \ell g(\bar y) = \ell g(\bar c_{i,1})$ and ${\gC} 
\models \varphi[\bar d_{\bold x},\bar c_{\bold x},\bar c_{i,\ell},
\bar c]$ iff $\ell=1$, etc., see Definition \ref{b3}.
\end{enumerate}
\end{claim}

\begin{remark}
\label{b23}
Note that in part (6), if $\bigwedge\limits_{\ell} \mu_\ell = \mu$ for
 transparency then we allow $\mu > \theta + |T| \ge 
\cf(\mu)$ and $|B_{\bold x}| = \mu$; also note that $B^+_{\bold x} =
B_{\bold x}$.
\end{remark}

\begin{PROOF}{\ref{b22}}
By \cite[2.10=tp25.36,2.11=tp25.38]{Sh:900} this should
be clear, still:

\noindent
5) Toward a contradiction assume that $\bold y,i \in v_{\bold y}
\backslash v_{\bold x},\varphi,\bar b_{i,0},\bar b_{i,1}$ exemplify $\bold
  x \notin \qK'_{\kappa,\bar\mu,\theta}$ so $\bar b_{i,0},\bar b_{i,1}$
are as in Definition \ref{b3} in particular there are $\bar b^*$ from
$M_{\bold x}$ and $\varphi$ such that $\gC \models ``\varphi[\bar d_{\bold
x},\bar c_{\bold x,<i},\bar b_{i,0},\bar b^*] \wedge 
\neg \varphi[\bar d_{\bold x},\bar c_{\bold x,<i},\bar b_{i,1},
\bar b^*]"$ and $\bar c_{\bold y,i} = \bar
b_{i,0} \char 94 \bar b_{i,1}$ and $\bar b_{i,0},\bar b_{i,1}$ realize the
same type over $\bar c_{\bold x,<i} + M_{\bold x}$ which is finitely
satisfiable in $B_{\bold y,i}$. 
Let $A$ be $B_{\bold y,i} + \bar b^*$ if $i \notin u_{\bold y}$ and be 
$\cup \bold I_{\bold y,i} + \bar b^*$ if $i \in u_{\bold x}$;
 so $A \subseteq M_{\bold x}$ has cardinality $< \kappa$ hence
for some $\psi = \psi(\bar x_{\bar d},\bar c_{\bold x},\bar z_\psi)$ and
$\bar e \in {}^{\ell g(\bar z_\psi)}(M_{\bold x})$ we have
\mn
\begin{enumerate}
\item[$(*)$]   $\psi(\bar x_{\bar d},\bar c_{\bold x},\bar e) \in
\tp(\bar d_{\bold x},\bar c_{\bold x} \dotplus M_{\bold
x})$ satisfies $\psi(\bar x_{\bar d},\bar c_{\bold x},\bar e) \vdash
\tp_\varphi(\bar d_{\bold x},\bar c_{\bold x} \dotplus A)$.
\end{enumerate}
\mn
Hence
\mn
\begin{enumerate}
\item[$(*)'$]   $(a) \quad \gC \models \psi[\bar d_{\bold x},\bar
c_{\bold x},\bar e]$
\sn
\item[${{}}$]  $(b) \quad$ for every $\bar b \subseteq A^{\ell g(\bar
b_{i,\ell})}$ and for some truth value $\bold t$ we have

\hskip25pt  $\gC \models
(\forall \bar x_{\bar d})[\psi(\bar x_{\bar d},\bar c_{\bold x},e)
\rightarrow \varphi(\bar x_{\bar d},\bar b)^{\bold t}]$.
\end{enumerate}
\mn
Now for $\ell=0,1$ we know 
that $\tp(\bar b_{i,\ell},M_{\bold x} + \bar
c_{\bold x,<i})$ is finitely satisfiable in $A$ and does
not depend on $\ell$, easy contradiction.
\end{PROOF}

\begin{observation}
\label{b24}  
1) Assume $\mu_0 = \kappa$ and $\cf(\kappa) > \theta + |T|$.
If $\bold x \in \qK'_{\kappa,\bar\mu,\theta}$ \then \, for 
some full $\bar\psi$ we have 
$\bold n := (\bold x,\bar\psi,\emptyset) \in 
\qK^\oplus_{\kappa,\bar\mu,\theta}$, see \ref{b16}(3A).  Moreover there is
$\bold n = (\bold x,\bar\psi',\emptyset) \in 
\qK^\otimes_{\kappa,\bar\mu,\theta}$.

\noindent
2) If $(\bold x,\bar\psi,r) \in \qK^\oplus_{\mu,\mu,\theta}$ and
the model $M$ is $\kappa$-saturated and $\kappa > \mu > \cf(\mu)$
\then \, for some $\bar\psi'$ we have $(\bold x,\bar\psi',r) \in 
\qK^\oplus_{\kappa,\mu^+,\theta}$.
\end{observation}

\begin{PROOF}{\ref{b24}}
1) First, for each $\varphi = \varphi(\bar x_{\bar d},\bar
x_{\bar c},\bar y) \in \Gamma^1_{\bold x}$ there 
is $\psi_\varphi = \psi_\varphi(\bar x_{\bar
d},\bar x_{\bar c},\bar z_\varphi)$ illuminating $(\bold
x,\varphi)$.  Why?  for every $A \subseteq M_{\bold x}$ of cardinality
$< \kappa$ by \ref{b22}(4) there is 
$\psi(\bar x_{\bar d},\bar c_{\bold x},\bar e)$ solving
$(\bold x,\varphi,A)$.  The set $\Lambda$ of candidates $\psi
=\psi(\bar x_{\bar d},\bar x_{\bar c},\bar z_\varphi)$ has cardinality $\le
\theta + |T|$ and if $\psi \in \Lambda$ fails there is a set $A_{\varphi,\psi}$
exemplifying it.  As $\cf(\kappa) > \theta + |T|$ the set $A_\varphi =
\cup\{A_{\varphi,\psi}:\psi \in \Lambda_\varphi\}$ has cardinality $<
\kappa$, so there are $\psi$ and $\bar e$ such that $\psi(\bar x_{\bar
d},\bar c_{\bold x},\bar e)$ solves $(\bold x,\varphi,A_\varphi)$, hence it
contradicts the choice of $A_{\varphi,\psi}$.  So $\psi_\varphi$
exists.

Renaming the $\bar z_\varphi$'s we have $\langle \psi_\varphi(\bar
x_{\bar d},\bar x_{\bar c},\bar z_{[\theta]}):\varphi \in 
\Gamma^1_{\bold x}\rangle$ as required for $\bold n := (\bold
x,\bar\psi,\emptyset) \in \qK^\oplus_{\kappa,\bar\mu,\theta}$.

Second, to get the ``moreover", let $\bar\varphi = 
\langle \varphi_i(\bar x_{\bar d},
\bar x_{\bar c},\bar y_{[\theta]},\bar z_i):i < \theta \rangle$
list the formulas of this form.  For $i < \theta$ let
$u_i \subseteq \theta$ be finite such that $\varphi_i = \varphi_i(\bar
x_{\bar d},\bar x_{\bar c},\bar y_{[u_i]},\bar z_i)$ and \wilog \, we
choose the sequence $\bar\varphi$ such that $u_i \subseteq i$.  Let $\psi_i =
\psi_i(\bar x_{\bar d},\bar x_{\bar c},\bar y^*_i)$ be as above for
$\varphi = \varphi_i(\bar x_{\bar d},\bar x_{\bar c};\bar y_{[u_i]}
\char 94 \bar z_i)$, so $\ell g(\bar y^*_i) < \omega$, let $\alpha_i =
\Sigma\{\ell g(\bar y^*_j):j < i\}$ and let $\bar y_i = \langle
y_{\alpha_i + \ell}:\ell < \ell g(\bar y^*_i)\rangle$ and now let
$\bar\psi^* = \langle \psi^*_{\varphi_i}: \psi_i(x_{\bar d},
\bar x_{\bar c},\bar y^*_i):i < \theta\rangle$.

Now given $A \subseteq M_{\bold x}$ of cardinality $< \kappa$ we
choose $\bar e_i = \langle e_{\alpha_i + \ell}:\ell < \ell g(\bar
y^*_i)\rangle$ by induction on $i < \theta$ such that $\psi_i(\bar
x_{\bar d},\bar x_{\bar c},\bar e_i)$ solves $(\bold x,A \cup
\{e_\alpha:\alpha < \alpha_i\},\varphi_i(\bar x_{\bar d},\bar x_{\bar
c},\bar y_{[u_i]} \char 94 \bar z_i))$.  So $\bar e \in
{}^\theta(M_{\bold x})$ is well defined and satisfies the requirements. 

\noindent
2) Let $\bar\psi = \langle \psi_\varphi:\varphi \in \Gamma^1_{\bold
  x}\rangle$. 
  For $\varphi = \varphi(\bar x_{\bar d},\bar x_{\bar c},\bar y) \in
\Gamma^1_{\bold x}$ let $\vartheta_{0,\varphi} =
\varphi,\vartheta_{1,\varphi} = \psi_{\vartheta_{0,1}}$ and recalling
  \ref{b9}(1C) let $\vartheta_{2,\varphi} =  
\psi_{\vartheta_1,\varphi}$.  Lastly, $\bar\psi := \langle
\psi_{\vartheta_{2,\varphi}}:\varphi \in \Gamma^1_{\bold x}\rangle$ 
satisfies $(\bold x,\bar\psi',r) \in
  \qK^\oplus_{\kappa,\mu^+,\theta}$; compare with the proof of \ref{b47}.
\end{PROOF}

\noindent
Note that we shall use \ref{b24} in \ref{d39}.
\bigskip

\subsection{Smoothness and $(\bar \mu,\theta)$-sets} \
\bigskip

We like to show that in some sense there are few decompositions, so 
toward this we define smooth ones, show that for a saturated model, 
the smooth decompositions are few up to being conjugate and every
$\bold x \in \pK_{\kappa,\bar\mu,\theta}$ is equivalent to a smooth
one modulo the relevant equivalence relation; this suffices.

\begin{definition}
\label{b27}  1) The decomposition $\bold x \in 
\pK_{\kappa,\bar\mu,\theta}$ is called smooth \when \,: if $\kappa \in
\ga_{\bold x}$, see end of \ref{b5}(1), then $\bold I^+_{\bold
x,\kappa}$ is an indiscernible sequence over $\cup\{\bold
I^+_{\bold x,\kappa_1}:\kappa_1 \in {\ga}_{\bold x} \backslash
\{\kappa\}\} \cup B_{\bold x} = \cup\{\bold I_{\bold x,i,\alpha}:i \in
u_{\bold x}$ and $\alpha < \kappa_i$ but $\kappa_{\bold x,i} \ne \kappa\} 
\cup B_{\bold x}$ in the sense of Definitions \ref{a61}(1), \ref{a73},
\underline{where} 

\noindent
1A) We define 

\[
\bold I^+_{\bold x,\kappa} =
\langle \bar a_{\bold x,\kappa,\alpha}:\alpha \in I_{\kappa,u(\bold
x,\kappa)} = I_{\bold x,\kappa}\rangle
\]

\mn
for $\kappa \in {\ga}_{\bold x} := 
\{\kappa_{\bold x,i}:i \in u_{\bold x}\}$ where 
$u_{\bold x,\kappa} = u(\bold x,\kappa) := 
\{i \in u_{\bold x}:\kappa_i = \kappa\}$ and $I_{\kappa,u(\bold
x,\kappa)} = (\kappa \times u_{\bold x,\kappa},<,
P_\varepsilon)_{\varepsilon \in (u(\bold x,\kappa))}$, where $<$
ordered $\kappa \times u_{\bold x,\kappa}$ lexicographically (if
$u_{\bold x,\kappa}$ is well ordered we can use $\kappa$) where
$\langle P_\varepsilon:\varepsilon < \otp(u_{\bold x,\kappa})\rangle$
is a partition to unbounded subsets, in fact, $P_\varepsilon = \kappa
\times \{\varepsilon\}$
and $\bar a_{\bold x,\kappa,\beta}$ is $\bar a_{\bold x,i,\alpha}$ when
$\beta = (\alpha,\varepsilon)$.

\noindent
2) For $\bold x \in \pK_{\kappa,\mu,\theta}$ and $h \in
 \Pi{\ga}_{\bold x}$ let $\bold x_{[h]}$ be defined like $\bold x$ but
 $\bar{\bold I}_{\bold x}$ is replaced by $\bar{\bold I}_{\bold x,h} = \langle
 \bold I_{\bold x,i,h(\kappa(\bold x,i))}:i \in u_{\bold x}\rangle$
 where $\bold I_{\bold x,i,\alpha} = \langle \bar a_{\bold x,i,\beta}:\beta 
\in [\alpha,\kappa_{\bold x,i})\rangle$ and $\bold I^+_{\bold
 x_{[h]},\kappa} = \langle a_{\bold x,i,\alpha}:\alpha \in
[\alpha,\kappa),i \in u_{\bold x,\kappa}\rangle$.

\noindent
3) We say $\bar b_1,\bar b_2$ are $\bold x$-similar \when \, for every
$n,(\forall i_0 \in u_{\bold x})(\forall^{\kappa_{i(0)}} \alpha_0 <
\kappa_i)(\forall i_1 \in u_{\bold x}) \ldots (\forall i_{n-1} \in
u_{\bold x})(\forall^{\kappa_{i(n-1)}} \alpha_{n-1} <
\kappa_{i(n-1)})[\tp(\bar b_1,
\cup\{\bar a_{\kappa_{i(\ell)},\alpha_\ell}:\ell < n\} \cup B_{\bold x}) =
 \tp(\bar b_2,\cup\{\bar a_{\kappa_{i(\ell)},\alpha_\ell}:\ell < n\}
 \cup B_{\bold x})]$, where we stipulate $i_\ell = i(\ell)$.
\end{definition}

\begin{definition}
\label{b30}  
1) We say the decompositions 
$\bold x,\bold y \in \pK_{\kappa,\bar\mu,\theta}$ are very 
similar \when \,:
\mn
\begin{enumerate}
\item[$(a)$]   $M_{\bold x} = M_{\bold y},w_{\bold x} = w_{\bold y},
\bar d_{\bold x} = \bar d_{\bold y},v_{\bold x} = 
v_{\bold y},u_{\bold x} = u_{\bold y},\bar
c_{\bold x} = \bar c_{\bold y}$ (so $\bar c_{\bold x,i} = \bar c_{\bold
y,i}$ for every $i$) and\footnote{We may consider weakening it.}
 $B_{\bold x,i} = B_{\bold y,i}$ for $i \in v_{\bold x} \backslash
u_{\bold x}$
\sn 
\item[$(b)$]   for $i \in u_{\bold x}$, the indiscernible sequences 
$\bold I_{\bold x,i},{\bold I}_{\bold y,i}$ are equivalent, (i.e. have the
same average over $M_{\bold x}$, equivalently over ${\gC}$)
and\footnote{usually this follows, but  not for stable indiscernible sets}
$\kappa_{\bold x,i} = \kappa_{\bold y,i}$.
\end{enumerate}
\mn
2) We say $\bold x,\bold y \in \pK_{\kappa,\bar\mu,\theta}$ are
similar when $v_{\bold x} = v_{\bold y},u_{\bold x} = u_{\bold y}$
and there is an elementary mapping $g$ of ${\gC}$ witnessing it which
means:
\mn
\begin{enumerate}
\item[$(a)$]   $g(B_{\bold x}) = g(B_{\bold y}),
g(\bar c_{\bold x}) = \bar c_{\bold y},g(\bar d_{\bold x})=
\bar d_{\bold y}$ and $g(B_{\bold x,i}) = B_{\bold y,i}$ for $i \in
v_{\bold x} \backslash u_{\bold x}$
\sn
\item[$(b)$]   for $i \in v_{\bold x} \backslash u_{\bold x},
g(\bar c_{\bold x,i}) = \bar c_{\bold y,i}$
and the scheme defining tp$(\bar c_{\bold x,i},M_{\bold x})$
(equivalently $\tp(\bar c_{\bold x,i},\bar c_{\bold x,<i} + M_{\bold
x}))$ is mapped to the scheme defining $\tp(\bar c_{\bold y,i},M_{\bold
y})$; so if $\iota_{\bold x} = 2$ this means
$g(D_{\bold x,i}) = D_{\bold y,i}$, i.e. 
$g(D'_{\bold x,i}) = D'_{\bold y,i}$, pedantically $g(D_{\bold x} \cap 
\Deef_{\ell g(\bar c_{\bold x,i})},B_{\bold x,i}) = D_{\bold y,i} \cap
\Deef_{\bold y,i} \cap \Deef_{\ell g(\bar c_{\bold x,i})}(B_{\bold y,i})$
\sn
\item[$(c)$]   $g(\bold I_{\bold x,i}),\bold I_{\bold y,i}$ 
are equivalent indiscernible sequences and $\kappa_{\bold x,i} =
\kappa_{\bold y,i}$ for $i \in u_{\bold x}$.
\end{enumerate}
\mn
3) Above we say weakly similar \when \, 
(so possible ${\ga}_{\bold x}
\ne {\ga}_{\bold y}$) as in part (2) but for each 
$i \in u_{\bold x}$ we replace the ``are equivalent" in clause (c), 
by the indiscernible sequences $h(\bold I_{\bold x,i}),
\bold I_{\bold y,i}$ being neighbors, (see here \ref{a67}(6)) and
$\kappa_{\bold x,i} = \kappa_{\bold x,j} \Leftrightarrow \kappa_{\bold
y,i} = \kappa_{\bold y,j}$.

\noindent
4) If $\bold x,\bold y$ are 
smooth we say they are smoothly immediately weakly similar 
\when \, in part (2) we replace clause (c) by
\mn
\begin{enumerate}
\item[$(c)'$]    there is a one-to-one function $h$ from 
${\ga}_{\bold x}$ onto ${\ga}_{\bold y}$ 
such that $\kappa_{\bold x,i} = \kappa \Rightarrow 
\kappa_{\bold y,i} = h(\kappa)$ and for some one-to-one order
preserving function from some infinite $u \subseteq \kappa$ into
$h(\kappa)$, we have $\alpha \in u \wedge \kappa_{\bold x,i} = \kappa
\Rightarrow a_{\bold x,\kappa_{\bold x,i},\alpha} = \ga_{\bold
y,\kappa_{\bold y,i},\alpha}$.
\end{enumerate}
\mn
5) For $\bold x,\bold y \in \pK$ we say they are essentially similar
\when \, there are smooth $\bar{\bold x}',\bold y' \in 
\pK_{\kappa,\mu,\theta}$ which are very similar to $\bold x,\bold y$  
respectively and are similar (by the definition in part (2);
note that e.g. $B_{\bold x,i},B_{\bold y,i}$ for $i \in v_{\bold x}
\backslash u_{\bold x}$ may be different). 
\end{definition}

\begin{claim}
\label{b33}  
Let $\kappa,\bar\mu$ and $\theta \ge |T|$ be as in \ref{b05} and
we let $\mu'_0$ be $\mu_0$ if $\cf(\mu_0) > \theta$ 
and $\mu'_0 = \mu^+_0$ otherwise, similarly for $\mu'_2$.

\noindent
1) Being similar, very similar, essentially similar and also 
weakly similar are equivalence relations.

\noindent
1A) Being very similar implies being similar which implies being
weakly similar which implies being smoothly immediately weakly similar.

2) For $\kappa$-saturated $M \prec {\gC}$, the number of $\bold x
\in \pK_{\kappa,\bar\mu,\theta}$ up to weak similarity is $\le 2^{< \mu'_0}$.

\noindent
3) For $\kappa$-saturated $M \prec {\gC}$ if $\mu_2 = 
\mu^{+ \alpha}_1$ and $\mu > |T| + \theta$, \then \, 
the number of $\bold x \in \pK_{\kappa,\bar\mu,\theta}$ up
to similarity is $\le 2^{< \mu_0} + |\alpha|^\theta$. 

\noindent
4) For $\bold x,\bold y \in \pK_{\kappa,\bar\mu,\theta}$ 
we have: $\bold x,\bold y$
are very similar \Iff \, $\bold x,\bold y$ are $\le_1$-equivalent,
i.e. $\bold x \le_1 \bold y \le_1 \bold x$.

\noindent
5) If $\bold x,\bold y$ are very similar and $\bar b_1,\bar b_2 \in
{}^\zeta{\gC}$ for some $\zeta < \mu_1$, \then \, $\bar b_1,\bar b_2$ are
$\bold x$-similar iff $\bar b_1,\bar b_2$ are $\bold y$-similar;
 see Definition \ref{b30}.
\end{claim}

\begin{PROOF}{\ref{b33}}
Easy (for essentially similar use \ref{b35}(1) below).   
\end{PROOF}

\begin{claim}
\label{b35}  
1) For every $\bold x \in \pK_{\kappa,\bar\mu,\theta}$ there
is a smooth $\bold y \in \pK_{\kappa,\bar\mu,\theta}$ very similar to
$\bold x$. 

\noindent
2) If $\bold x \in \pK_{\kappa,\mu,\theta}$ and $h \in \prod
{\ga}_{\bold x}$ \then \, $\bold x_{[h]} \in 
\pK_{\kappa,\bar\mu,\theta}$ is very similar to $\bold x$; see \ref{b27}(2).

\noindent
3) In part (2), if $\bold x$ is smooth \then \, so is $\bold x_{[h]}$.

\noindent
4) If $\bold x \in \qK_{\kappa,\bar\mu,\theta}$ and 
$\bold y \in \pK_{\kappa,\bar\mu,\theta}$ is very similar to $\bold
x$ \then \, $\bold y \in \qK_{\kappa,\bar\mu,\theta}$.

\noindent
4A) Similarly for $\qK'_{\kappa,\bar\mu,\theta}$ and
$\qK''_{\kappa,\bar\mu,\theta}$ .

\noindent
5) If $\bold x \in \pK_{\kappa,\bar\mu,\theta}$ is smooth 
\then \, for every $\bar a \in {}^{\mu_0 >}(M_{\bold x})$ for some $h \in
\Pi{\ga}_{\bold x}$ also $((M_{\bold x})_{[\bar a]},
\bar B_{\bold x},\bar c_{\bold x},\bar d_{\bold x},
\bar{\bold I}_{\bold x,h}) \in \pK_{\kappa,\bar\mu,\theta}$ is 
smooth, even replacing ${\gC}$ by
${\gC}_{[\bar a]}$; also if $\bar a = (\ldots \char 94
\bar a_i \char 94 \ldots)_{i \in v(\bold x) \backslash u(\bold x)}$
where $\bar a_i \in {}^{\omega >}(M_{\bold x})$ or just $\bar a_i \in
{}^{\mu_0 >}(M_{\bold x})$ for every $i \in v_{\bold x} 
\backslash u_{\bold x}$ \then \, for some $h \in \Pi \ga_{\bold x}$
the tuple $(M_{\bold x},\langle B_{\bold x,i} + \bar a_i:i \in u_{\bold
x}\rangle,\bar c_{\bold x},\bar d_{\bold x},
\bar{\bold I}_{\bold x,h}) \in \pK_{\kappa,\bar\mu,\theta}$
is smooth, see \ref{b27}(2).

\noindent
6) $\bold S^{< \theta^+}(B^+_{\bold x})$ has cardinality 
$\le |B^+_{\bold x}|^\theta$ \when \, $\bold x \in 
\pK_{\kappa,\bar\mu,\theta}$ is smooth.
\end{claim}

\begin{PROOF}{\ref{b35}}
E.g. for parts (5),(6) use \ref{a64}(1).
\end{PROOF}
\bn
We may formalize how ``small" is $B^+_{\bold x}$ for smooth $\bold x \in
\pK_{\kappa,\bar\mu,\theta}$.

\begin{definition}
\label{b38}  We say that $\bold f = (\bar B,\bar{\bold I})$ is a 
$(\bar\mu,\theta)$-set or a $(\bar\mu,\theta)$-smooth set
\when \, $\bar \mu = (\mu_2,\mu_1,\mu_0)$ and for some $u,v$ we have:
\mn
\begin{enumerate}
\item[$(a)$]   $v$ is a linear order of cardinality 
$< \theta^+$ and $u \subseteq v$
\sn
\item[$(b)$]   $\bar B = \langle B_i:i \in v \backslash u
\rangle$, we let $B = \cup\{B_i:i \in v \backslash u\}$ and
each $B_i$ is of cardinality $< \mu_0$; 
but $\bold f = (B,\bar{\bold I})$ means $i \in v
\backslash u \Rightarrow B_i = B$ so in this case $|B| < \mu_0$
\sn
\item[$(c)$]  $\bar{\bold I} = \langle \bold I_i:i \in u\rangle$
\sn
\item[$(d)$]   $\bold I_i = \langle \bar a_{i,\alpha}:\alpha <
\kappa_i\rangle$ is an indiscernible sequence of finite tuples,
$\kappa_i \in \Reg \cap \mu_2 \backslash \mu_1$
\sn
\item[$(e)$]   $(B,\bar{\bold I})$ satisfies the smoothness demand,
clause (e) in Definition \ref{b27} and $\ga,I_\kappa,\bold I^+_\kappa,
\bar a_{\kappa,\alpha}$ (for $\kappa \in \ga,
\alpha \in I^+_\kappa)$ are defined as there.
\end{enumerate}
\end{definition}

\begin{definition}
\label{b41}  
1) For $\bold f$ as in \ref{b38} we let\footnote{This is an abuse of
our notation as $\bold f$ does not determine $\mu_\ell$ in Definition
\ref{b38}, pedantically we can expand $\bold f$ to have this
information.}: $\mu_{\bold f,\ell} =
\mu_\ell$ for $\ell=0,1,2,v_{\bold f} = v,
u_{\bold f} = u,B_{\bold f,i} = B_i,
B_{\bold f} = \cup\{B_{\bold f,i}:i \in v_{\bold f} \backslash
u_{\bold f}\},\bar{\bold I}_f = \bar{\bold I},\bold I_{\bold f,\kappa}= 
\bold I_\kappa,\bar a_{\bold f,i,\alpha} = \bar a_{i,\alpha},
\ga_{\bold f} = \ga = \{\kappa_i:i \in u\}$, 
etc.; for $u \subseteq u_{\bold f}$ let
${\frak a}_{\bold f,u}  = \{\kappa_i:i \in u\}$.

\noindent
1A) If $\bold x \in \pK$ is smooth \then \, $\bold f = \bold f_{\bold x}$ is
defined by $v_{\bold f} = v_{\bold x},u_{\bold f} = u_{\bold
x},B_{\bold f,i}= B_{\bold x,i},\bold I_{\bold f,j} = \bold I_{\bold
x,j}$ hence $\bar a_{\bold f,i,\alpha} = \bar a_{\bold x,i,\alpha}$
for $i \in v_{\bold x} \backslash u_{\bold x}$ and $j \in u_{\bold
x},\alpha < \kappa_{\bold x,j}$.

\noindent
2) For $u \subseteq u_{\bold f}$ let 
$B^+_{\bold f,u} = \cup\{\bar a_{\bold f,i,\alpha}:\alpha <
 \kappa_i$ and $i \in u\} \cup B_{\bold f}$, if we omit $u$ 
we mean $u = u_{\bold f}$.

\noindent
2A) For $v \subseteq v_{\bold f}$ let $B^\pm_{\bold f,v} =
\cup\{B_{\bold f,i}:i \in v \backslash u_{\bold f}\} \cup \bigcup\{\bar
a_{\bold f,i,\alpha}:i \in v \cap u_{\bold f}\}$.

\noindent
3) We say $\bold f$ is an infinitary $(\bar\mu,\theta)$-set \when \,
$\ell g(\bar a_{\bold f,i,\alpha})$ is just $< \theta^+$ 
for every $i \in u_{\bold f},\alpha < \kappa_{\bold f,i}$ instead of
being finite.

\noindent
4) Let $\bar{\bold I}_{\bold f,h} =  \langle 
\bold I_{\bold f,i,h(\kappa_{\bold f,i})}:i \in u_{\bold f}\rangle 
= \langle \langle \bar a_{\bold f,i,\alpha}:
\alpha < \kappa_i$ and $i \in \Dom(h)
   \Rightarrow h(\kappa_{\bold f,i}) \le \alpha \rangle:
i \in u_{\bold f}\rangle$ so $h \in \Pi{\ga}_{\bold f}$.
Let $B^+_{\bold f,u,h}$ be defined as in part (2) using $\bar{\bold I}_{\bold
   f,h}$.  Let $\bold f_{[h]} = (\bar B_{\bold f},\bar{\bold I}_{\bold f,h})$.

\noindent
5) For $g \in \Pi \ga_{\bold f}$ let
$\cF_{\bold f,u,g} = \{h:h \in \prod\limits_{i \in u}
(\kappa_{\bold f,i} \backslash g(\kappa_{\bold f,i}))$ and 
if $i<j$ are from $u$ and
$\{i_1 \in u:\kappa_{i_1} = \kappa_i\}$ is well ordered and
$\kappa_{\bold f,i} = \kappa_{\bold f,j}$ then $h(i) < h(j)\}$; also
 for $h \in \cF_{\bold f,u,g}$ let $\bar a_{\bold f,u,h} := 
\langle \bar a_{\bold f,i,h(i)}:i \in u\rangle$.  

\noindent
6) We say that $\bold f$ is essentially well ordered \when \, for each
   $\kappa$ the set $\{i \in u_{\bold f}:\kappa_{\bold f,i}=\kappa\}$
   is well ordered by $\le_v$; compare with Definition \ref{b5}(11).
\end{definition}

\begin{claim}
\label{b45}  
1) If $\bold f$ is a $(\bar\mu,\theta)$-set, $|B^+_{\bold f}| \ge 2$
   for simplicity and $\varepsilon < \theta^+$
\then \, $\bold S^\varepsilon(B^+_{\bold f})$ has cardinality 
$\le |B^+_{\bold f}|^\theta$.

\noindent
1A) If $m < \omega$ and $\Delta \subseteq \bbL(\tau_T)$ is finite, 
\then \, for some $k$ we have 
$|\bold S^m_\Delta(B^+_{\bold f})| \le |B^+_{\bold f}|^k$ whenever
$\bold f$ is a $(\bar\mu,\theta)$-set, $u_{\bold f}$ is finite.

\noindent
2) If $\bold x \in \pK_{\kappa,\bar\mu,\theta}$ is smooth
\then \, $(\bar B_{\bold x},\bar{\bold I}_{\bold x})$ is a
$(\bar\mu,\theta)$-set.

\noindent
3) If $\bold f$ is an essentially well ordered 
$(\bar\mu,\theta)$-set and $\bar e \in
{}^{\mu_1 >}{\gC}$ \then \, for some $h \in \Pi{\ga}_{\bold f}$ 
for some type $q$ we have: $g \in \prod\limits_{i \in
u_{\bold f}} \kappa_{\bold f,i} \wedge \bigwedge\limits_{i \in
u_{\bold f}} h(\kappa_i) \le g(i) \wedge \bigwedge\limits_{\kappa_i =
\kappa_j,i <_v j} g(i) \le g(j) \Rightarrow 
\tp((\ldots \char 94 \bar a_{\bold f,i,g(\alpha)} \char 94
\ldots),B_{\bold f} + \bar e) = q$.

\noindent
4) If $\bold f = (B_{\bold f},\bar{\bold I})$ is 
a $(\bar\mu,\theta)$-set and $C \subseteq {\gC}$ has cardinality $<
\mu_0$ \then \, $(B_{\bold f} + C,\bar{\bold I}_{\bold f,h})$ 
is a $(\bar\mu,\theta)$-set for some $h \in \Pi{\ga}_{\bold f}$. 

\noindent
5) If $\bold x \in \pK_{\kappa,\bar \mu,\theta}$ is smooth \then \,
   $\bold f_{\bold x}$ is a $(\bar \mu,\theta)$-set, see Definition
   \ref{b41}(1A).
\end{claim}

\begin{PROOF}{\ref{b45}}
E.g. part (1) by \ref{a64}(4) using $\{(\kappa,\alpha):\kappa \in
\ga_f$ and $\alpha < \kappa\}$ ordered lexicographically;
part (4) by \ref{a64}(1) as in \ref{b35}(6), for part (3)
  recall the smoothness demand.
\end{PROOF}
\bigskip

\subsection{Measuring non-solvability and reducts} \
\bigskip

The following is needed in \S4, \S5, it measures how far solutions are
missing.
\begin{definition}
\label{b43}  
1) For $\bold x \in \pK_{\kappa,\bar\mu,\theta}$ let $\ntr(\bold x)$, the
non-transitivity of $\bold x$ be the minimal cardinal $\lambda$ such
that for some $A \subseteq M_{\bold x}$ of cardinality $\lambda$ for no
$\bar e \in {}^\theta(M_{\bold x})$ do we have $\tp(\bar d_{\bold
x},\bar c_{\bold x} + \bar e) \vdash \tp(\bar d_{\bold x},
\bar c_{\bold x} + A)$.

\noindent
2) For $\bold x \in \pK_{\kappa,\bar\mu,\theta}$ let 
$\ntr_{\lc}(\bold x)$ be the minimal cardinal $\lambda$
such that for some $\varphi = \varphi(\bar x_{\bar d},\bar x_{\bar c},
\bar y) \in \Gamma^1_{\bold x}$, we have $\ntr_\varphi(\bold x) =
\lambda$, see below.

\noindent
3) For $\varphi \in \Gamma^1_\varphi$, let 
$\ntr_\varphi(\bold x)$ be the minimal $\lambda$ such that no
   $\psi$ does $\lambda^+$-illuminates $(\bold x,\varphi)$, i.e. there
is $A \subseteq M_{\bold x}$ of cardinality $\lambda$ such that
for no $\psi(\bar x_{\bar d},\bar c_{\bold x},
\bar e) \in \tp(\bar d_{\bold x},\bar c_{\bold x}
\dotplus M_{\bold x})$ do we have $\psi(\bar x_{\bar d},\bar
 c_{\bold x},\bar d) \vdash \tp_\varphi(\bar d_{\bold x},
\bar c_{\bold x} \dotplus A)$.

\noindent
4) Let $\ntr_{\varphi,\psi}(\bold x)$ be defined naturally.

\noindent
5) We say that $\bar\psi$ does $\lambda$-illuminate $\bold x \in
\pK_{\kappa,\bar\mu,\theta}$ \when \, $\Gamma_{\bar\psi} =
\Gamma^1_{\bold x}$ and for every $A \subseteq M_{\bold x}$ 
of cardinality $< \lambda$, some $\bar e
\in {}^\theta(M_{\bold x})$ solves $(\bold x,\bar\psi,A)$, 
(see \ref{b9}(1),(3A)).

\noindent
6) We say $\bar\psi$ does $\lambda$-illuminate $(\bold x,\Gamma)$ or
$\bar\psi$ illuminate $(\bold x,\lambda,\Gamma)$ \when \, $\Gamma \subseteq 
\Gamma^1_{\bold x},\bar \psi =
\langle \psi_\varphi:\varphi \in \Gamma\rangle$ and
   for every $A \subseteq M_{\bold x}$ of cardinality $< \lambda$ for
   some $\bar e \in {}^\theta(M_{\bold x})$ the sequence $\bar e$
   solves $(\bold x,\bar\psi,A)$.

\noindent
7) Similarly when $\Gamma_{\bar\psi} \subseteq \Gamma^3_{\bold x}$.

As in \ref{b24}.
\end{definition}

\begin{observation}
\label{b47}  
1) If $\bold x \in \pK_{\kappa,\bar\mu,\theta}$ 
and $\lambda = \ntr(\bold x) > \theta(\ge |T|)$ \then \,
$\ntr(\bold x)$ is a regular cardinal.

\noindent
2) If $x \in \pK_{\kappa,\bar\mu,\theta}$ and $\lambda = 
\ntr_{\lc}(\bold x)$ is singular \then  \, $\cf(\lambda) \le \theta  + |T|$.

\noindent
3) If $\theta < \cf(\lambda) \le \lambda \le \ntr(\bold x)$ \then \,
some $\bar\psi = \langle \psi_\varphi:\varphi \in 
\Gamma^1_{\bold x}\rangle$ does $\lambda$-illuminate $\bold x$.
\end{observation}

\begin{PROOF}{\ref{b43}}  
1) Why is $\lambda$ regular?  If $\lambda > \cf(\lambda)$, let 
$A \subseteq M_{\bold x}$ exemplify the choice of
$\lambda$, let $\langle A_\alpha:\alpha < \cf(\lambda)\rangle$ 
be $\subseteq$-increasing, each $A_\alpha$ being of
cardinality $< \lambda$ such that $A = \cup\{A_\alpha:\alpha < \lambda\}$.  For
each $\alpha < \lambda$ by the choice of $\lambda$ there is $\bar
e_\alpha \in {}^\theta(M_{\bold x})$ such that $\tp(\bar d_{\bold
x},\bar c_{\bold x} + \bar e_\alpha) \vdash 
\tp(\bar d_{\bold x},\bar c_{\bold x} + A_\alpha)$.

Let $A_* = \cup\{\bar e_\alpha:\alpha < \cf(\lambda)\}$ so
$|A_*| \le \theta + \cf(\lambda) < \lambda$ hence for some
$\bar e \in {}^\theta(M_{\bold x})$ we have $\tp(\bar d_{\bold x},\bar
c_{\bold x} + \bar e) \vdash \tp(\bar d_{\bold x},\bar c_{\bold
x} + A_*)$.  Clearly $\bar e$  contradicts the choice of $A$.

\noindent
2) Similarly (as in \ref{b24}(2)), changing the $\psi_\varphi$'s.

\noindent
3) Similarly, as in the proof of \ref{b24}. 
\end{PROOF}

\begin{definition}
\label{b49}  
1) If $\tau \subseteq \tau_T$ and $\bold x \in 
\pK_{\kappa,\bar\mu,\theta}$ let $\bold x \rest
\tau$ be defined like $\bold x$ but ${\gC}$ is replaced by ${\gC} 
\rest \tau$ and $M_{\bold x}$ is replaced by $M_{\bold x} \rest \tau$.

\noindent
2)  If $\bold x \in \pK_{\kappa,\bar\mu,\theta}$ and $v
 \subseteq v_{\bold x},w \subseteq w_{\bold x}$ \then \, $\bold y =
   \bold x \rest (v,w)$ is defined by $M_{\bold y} = M_{\bold
   x},w_{\bold y} = w,\bar d_{\bold y} = \bar d_{\bold x} \rest
 w,v_{\bold y} = v,\bar c_{\bold y} = \bar c_{\bold x}
   \rest v,u_{\bold y} = u_{\bold x} \cap v,B_{\bold y,i} = B_{\bold
   x,i}$ for $i \in v_{\bold y} \backslash u_{\bold y}$ and 
$\bold I_{\bold y,i} = \bold I_{\bold x,i}$ for $i \in u_{\bold y}$. 
\end{definition}

\begin{observation}
\label{b51}
Membership in $\pK_{\kappa,\bar\mu,\theta}$ is preserved under 
reducts, i.e. if $\tau \subseteq \tau(T)$ then
$\bold x \rest \tau \in \pK_{\kappa,\bar\mu,\theta}[{\gC} \rest \tau]$; also 
and $\bold x \rest (v,w) \in \pK_{\kappa,\bar\mu,\theta}$ in 
the cases above.  Also smoothness, ``very similar", etc. are preserved.
\end{observation}

\begin{PROOF}{\ref{b51}}
Straightforward.
\end{PROOF}

\begin{claim}
\label{b68}  
1) If $\bold x \in \pK_{\kappa,\bar\mu,\theta}$ 
and $\iota_{\bold x} = 2$ \then \, $\tp(\bar c_{\bold x},M_{\bold x})$ is 
finitely satisfiable in $B^+_{\bold x}$
hence for some ultrafilter $D_{\bold x}$ on $\Deef_{\ell g(\bar c_{\bold x})}
(B^+_{\bold x})$, we have $\tp(\bar c_{\bold x},M^+_{\bold x}) =
\Av(D_{\bold x},M_{\bold x})$ in fact $D_{\bold x}$
is unique.

\noindent
2) If $\bold x \in \pK_{\kappa,\bar\mu,\theta}$ and
$\iota_{\bold x} = 0$ \then \, $\tp(\bar c_{\bold x},
M_{\bold x})$ does not split over $B^+_{\bold x}$.

\noindent
3) If $\bold x \in \pK_{\kappa,\bar\mu,\theta}$ and
$\iota_{\bold x} = 1$ \then \, $\tp(\bar c_{\bold x},M_{\bold x})$
does not locally split over $B^+_{\bold x}$.
\end{claim}

\begin{PROOF}{\ref{b68}}
Straightforward.
\end{PROOF}

\noindent
We can elaborate \ref{b68}(1)
\begin{definition}
\label{b72}  
1) Let $D_\ell$ be an ultrafilter
on $\Deef_\varepsilon(A_\ell)$ for $\ell=1,2$.  We say $D_1,D_2$ are
equivalent when $\Av(D_1,C) = \Av(D_2,C)$ for every set $C \subseteq {\gC}$.

\noindent
2) We say an ultrafilter $D$ on $\Deef_\varepsilon(A)$ is
$(\bar\mu,\theta)$-smooth \when \, ${}^\varepsilon(B^+_{\bold f}) 
\in D$ for some $(\bar\mu,\theta)$-set $\bold f$.
\end{definition}

\begin{definition}
\label{b74} 
1) For $\bold x \in \pK$ such that $\iota_{\bold x} = 2$ let
$D_{\bold x}$ be the following ultrafilter:
\mn
\begin{enumerate}
\item[$(a)$]  $D_{\bold x}$ is an ultrafilter 
on ${\cC}_{\bold x} = \{\langle \bar c'_i:i 
\in v_{\bold x}\rangle:\bar c'_i \in {}^{\ell g(\bar c_i)}(B_i)$ 
if $i \in v_{\bold x} \backslash u_{\bold x}\}$ and
$\bar c'_i \in \{\bar a_{\bold x,i,\alpha}:\alpha < \kappa_{\bold
x,i}\}$ if $i \in u_{\bold x}\}$ 
\sn
\item[$(b)$]   $\{\bar c' \in \cC_{\bold x}:\gC
\models \varphi(\bar c',\bar b]\} \in \Av(D,{\gC})$ 
\Iff \, letting $\varphi$ depend just on
$(\bar x_{\bar c_{i(0)}},\dotsc,\bar x_{\bar c_{i(n-1)}}),i(0) > i(1) >
\ldots > i(n-1)$ and the formula $\varphi(\bar x_{\bar
c_{i(0)}},\dotsc,\bar x_{\bar c_{i(n-1)}},\bar b)$ 
belongs to $\Av(D_{\bold x,i(0)} \times D_{\bold x,i(1)} 
\times \ldots \times D_{\bold x,i(n-1)},\bar b)$.
\end{enumerate}
\mn
2) Above $D_{\bold x,i}$ is the natural ultrafilter.  
\end{definition}

\begin{definition}
\label{b77}  
For a $(\bar \mu,\theta)$-set $\bold f$, set $v$ 
 and $A \subseteq \gC$ of cardinality $< \mu_1$, we define an
   equivalence relation ${\cE}^v_{\Delta,A}$ on $\bold S^v_\Delta(A
+ B^+_{\bold f,u_n})$ as follows: (if $\Delta = \bbL(\tau_T)$ we
may omit $\Delta$)  

$\tp_\Delta(\bar b_1,A + B^+_{\bold f,u_n}) 
{\cE}^v_\Delta \tp_\Delta(\bar b_2,A + B^+_{\bold f,u_n})$ \Iff \,
($\ell g(\bar b_1) = v = \ell g(b_2)$ and) for some $h \in 
\Pi{\ga}_{\bold f}$, the types $\tp_\Delta(\bar b_2,A + 
B^+_{\bold f,u_n,h}),\tp_\Delta(\bar b_2,A + B^+_{\bold f,h})$ are equal.
\end{definition}

\begin{observation}
\label{b79}
1) On $\bold S^v_{A,\bold f,h} := \{\tp
(\bar e \rest v,A + B^+_{\bold f,u,h}):(B_{\bold f,u_n} + 
\bar e,\bar{\bold I}_{\bold f,h})$ is a
$(\bar\mu,\theta)$-set$\}$ the equivalence relation ${\cE}^v_A$ is
the equality.

\noindent
2) For $\bold x \in \pK_{\kappa,\bar\mu,\theta}$ such that $(\forall
\alpha < \kappa)(|\alpha|^\theta < \kappa)$ and
$(\forall \mu < \mu_0)(2^\mu < \kappa)$ and $\mu_2 = \kappa$
letting $\bold f = \bold f_{\bold x}$ see Definition \ref{b33}(1A) and 
$n < \omega$ and $A \subseteq M_{\bold x}$ of cardinality $< \mu_0$ 
the equivalence relation ${\cE}^n_{\Delta,A}$
has $< \kappa$ equivalence classes.
\end{observation}

\begin{PROOF}{\ref{b79}}
Straightforward.
\end{PROOF}
\newpage

\section {Strong analysis} 

In \S2 we have dealt with $\pK$ and $\qK$, here we use $\tK$ and
$\vK$.  
Now $\tK$ is the ``really analyzed" case, one essentially
with ``$\bar d_{\bold x}$ universally solve itself" so it
is a central notion here.  But we have problems in proving its
density in enough cases (i.e. cardinals), so we use also a relative $\vK$, 
weak enough for the density proof,
strong enough for the main desired consequence.  We do not forget $\tK$ as
it is more transparent and says more.
We give some consequences of $\bold x \in 
\tK_{\kappa,\bar\mu,\theta}$ or $\bold x \in 
\vK_{\kappa,\bar\mu,\theta}$.  
First, $M_{[\bold x]} = M_{[B^+_{\bold x} +
\bar c_{\bold x} + \bar d_{\bold x}]}$ is $(\bold D_{\bold x},\kappa)$-
sequence-homogenous (see \ref{z25}(1); so pcf, see \cite{Sh:g} 
appears naturally when we try to analyze
$\bold D_{\bold x}$ but this is not really used here).  
This implies uniqueness, so indirectly few
types up to conjugacy; this will solve the recounting 
problems from \S(1A) but only 
when we shall prove density of $\tK$ or $\vK$ 
in $(\pK_{\kappa,\bar\mu,\theta},\le_1)$.  We give sufficient
condition for existence, using existence of universal solutions and
prove it for $\kappa$ weakly compact when $\|M_{\bold x}\| = \kappa$.
We end with criterions for indiscernibility related to $\tK$.
\newline

Note that $\tK$ is better than $\qK$, but the relevant density result is for
$\le_1$ rather than $\le_2$, i.e. you may say that we add more
variables to the type analyzed.  
\bigskip

\subsection{Introducing $\rK,\tK,\vK$} \
\bigskip

So a central definition is
\begin{definition}
\label{c1}  
Let $\tK_{\lambda,\kappa,\bar\mu,\theta}$ be the class of $\bold x \in
\pK_{\lambda,\kappa,\bar\mu,\theta}$ such that: for every
$A \subseteq M_{\bold x}$ of cardinality $< \kappa$ there is
$(\bar c_*,\bar d_*)$ which strongly solves $(\bold x,A)$ 
which means: $\bar c_* \char 94 \bar d_*$ is from
$M_{\bold x}$ and it realizes $\tp(\bar c_{\bold x} \char 94 \bar
d_{\bold x},A)$, of course $\ell g(\bar c_*) = \ell g(\bar c_{\bold
x}),\ell g(\bar d_*) = \ell g(\bar d_{\bold x})$ and $\tp(\bar d_{\bold
x},\bar c_{\bold x} + \bar d_* + \bar c_*) \vdash 
\tp(\bar d_{\bold x},\bar c_{\bold x} + \bar d_* + \bar c_* + A)$ by
some $\bar\psi$. 
\end{definition}

\begin{remark}
\label{c2} 
1) For $\le_1$-increasing chains in $\tK_{\kappa,\bar\mu,\theta}$ the union is 
naturally defined (essentially see in \ref{b18}(2))
but it is not a priori clear it belongs to 
$\tK_{\kappa,\bar\mu,\theta}$, i.e. if $\langle \bold x_\alpha:
\alpha < \delta\rangle$ is $\le_1$-increasing in 
$\tK_{\kappa,\bar\mu,\theta}$ and $\delta < \theta^+$ \then \, does the 
union belongs to $\tK_{\kappa,\bar\mu,\theta}$?

\noindent
2) To have enough cases when this holds we define a relative of $\pK$
   which carries more information.

\noindent
3) Note that below
\mn
\begin{enumerate}
\item[$(a)$]   $\rK,\tK,\uK,\vK$ are subsets of $\pK$
\sn
\item[$(b)$]   $\rK^\oplus,\tK^\oplus,\vK^\oplus,\vK^\odot$ 
has the form $(\bold x,\bar\psi,r)$ with $\psi_\varphi$'s only for
some $\varphi \in \Gamma^2_{\bold x}$
\sn
\item[$(c)$]   $\sK^\oplus,\uK^\otimes,\uK^\otimes$ are
similar but with $\psi_\varphi$'s only for (some) $\varphi 
\in \Gamma^3_{\bold x}$
\sn
\item[$(d)$]   the $\vK$'s and $\uK$'s use so-called duplicates
  (defined below)
\sn
\item[$(e)$]   $\uK_{\kappa,\bar\mu,\theta}$ is a parallel of
$\qK_{\kappa,\bar\mu,\theta}$ when we allow duplication, see
\ref{b16}(1), \ref{c23}(3c).
\end{enumerate}
\end{remark}

\begin{definition}
\label{c5}   
1) Let $\rK^\oplus_{\lambda,\kappa,\bar\mu,\theta}$ be the class of
triples $\bold n = (\bold x,\bar\psi,r)$ such that
\mn
\begin{enumerate}
\item[$(a)$]   $\bold x \in \pK_{\lambda,\kappa,\bar\mu,\theta}$
\sn
\item[$(b)$]   $r$ is a type in the variables $\bar x_{\bar d} \char
94 \bar x_{\bar c} \char 94 \bar x'_{\bar d} \char 94 \bar x'_{\bar c}$, 
over $\emptyset$ if not said otherwise
\sn
\item[$(c)$]  $\bar\psi = \langle \psi_\varphi(\bar x_{\bar d},
\bar x_{\bar c},\bar x'_{\bar d},\bar x'_{\bar c}):\varphi \in
\Gamma^2_{\bar\psi}\rangle$ recalling\footnote{but we may allow one
$\varphi$ to appear more than once} \ref{b16}(0)(A)
\sn
\item[$(d)$]  $\Gamma_{\bar\psi} = \Gamma^2_{\bar\psi} 
\subseteq \Gamma^2_{\bold x} :=
\{\varphi:\varphi = \varphi(\bar x_{\bar d},\bar x_{\bar c},\bar
x'_{\bar d},\bar x'_{\bar c},\bar y) \in \bbL(\tau_T)\}$
\sn
\item[$(e)$]  $\psi_\varphi(\bar x_{\bar d},\bar x_{\bar c},
\bar x'_{\bar d},\bar x'_{\bar c}) \in r$ for 
every $\varphi \in \Gamma^2_{\bar\psi}$
\sn
\item[$(f)$]   if $A \subseteq M_{\bold x}$ has cardinality $<
\kappa$ \then \, some $(\bar c',\bar d')$ solves 
$(\bold x,\bar\psi,r,A)$ or solves $(\bold n,A)$; we may write $\bar c' 
\char 94 \bar d'$ instead $(\bar c',\bar d')$; which means:
\sn
\begin{enumerate}
\item[$(\alpha)$]  $\bar c' \char 94 \bar d'$ is from 
$M_{\bold x}$ and realizes $\tp(\bar c_{\bold x} \char 94 \bar d_{\bold x},A)$
\sn
\item[$(\beta)$]  $\bar c_{\bold x} \char 94 \bar d_{\bold x} 
\char 94 \bar c' \char 94 \bar d'$ realizes $r$, 
of course, $\ell g(\bar c') = \ell g(\bar c_{\bold x}),\ell g(\bar d') 
= \ell g(\bar d_{\bold x})$
\sn
\item[$(\gamma)$]  if $\varphi \in \Gamma^2_{\bold x}$ then
$\psi_\varphi(\bar x_{\bar d},\bar c_{\bold
x},\bar d',\bar c') \vdash \tp_\varphi(\bar d_{\bold x},\bar
c_{\bold x} \char 94 \bar d' \char 94 \bar c' \dotplus A)$ 
recalling the latter means $\{\varphi(\bar x_{\bar d},\bar c_{\bold
x},\bar d',\bar c',\bar b):\bar b \in {}^{\ell g(\bar y_\varphi)}A$
and ${\gC} \models \varphi(\bar d_{\bold x},\bar c_{\bold x},\bar
d',\bar c',\bar b)\}$. 
\end{enumerate}
\end{enumerate}
\mn
2) Let $\sK^\oplus_{\lambda,\kappa,\bar\mu,\theta}$ be the class of
tuples $\bold m = (\bold x,\bar\psi,r)$ such that:
\mn
\begin{enumerate}
\item[$(a)$]  $\bold x \in \pK_{\lambda,\kappa,\bar\mu,\theta}$
\sn
\item[$(b)$]  $r$ is a type in the variables $\bar x_{\bar c} \char 94
\bar x_{\bar d} \char 94 \bar y_{[\theta]}$, over $\emptyset$ if not
said otherwise
\sn
\item[$(c)$]  $\bar\psi = \langle \psi_\varphi(\bar x_{\bar d},\bar
x_{\bar c},\bar y_{[\theta]}):\varphi \in \Gamma^3_{\bar\psi}\rangle$

where recalling \ref{b16}(1B)
\sn
\item[$(d)$]  $\Gamma_{\bar\psi} = \Gamma^3_{\bar\psi} \subseteq
\Gamma^3_{\bold x} = \{\varphi:\varphi = \varphi(\bar x_{\bar d},\bar
x_{\bar c},\bar y_{[\theta]},\bar z) \in \bbL(\tau_T)\}$
\sn
\item[$(e)$]  $\psi_\varphi(\bar x_{\bar d},\bar x_{\bar c},\bar
y_{[\theta]}) \in r$ for every $\varphi \in \Gamma^3_{\bar\psi}$
\sn
\item[$(f)$]  if $A \subseteq M_{\bold x}$ has cardinality $< \kappa$
\then \, some $\bar e$ solves $(\bold x,\bar\psi,r,A)$ or
solves $(\bold m,A)$ which means:
\sn
\begin{enumerate}
\item[$(\alpha)$]  $\bar c \char 94 \bar d \char 94 \bar e$ realizes
$r$
\sn
\item[$(\beta)$]  $\psi_\varphi(\bar x_{\bar d},\bar c_{\bold x},\bar
e) \vdash \tp_\varphi(\bar d_{\bold x},\bar c_{\bold x} \char
94 \bar e \dotplus A)$ for every $\varphi \in \Gamma^3_{\bar\varphi}$.
\end{enumerate}
\end{enumerate}
\mn
3) We define ``very similar/similar/weakly similar" on $\rK^\oplus$
and $\sK^\oplus$ naturally, (and they are equivalence relations).
\end{definition}

\begin{remark}
\label{c21}  
1) So arbitrary $\bar b \subseteq \Rang(\bar c_{\bold x})$ is 
not allowed in clauses $(f)(\gamma)$ of \ref{c5}(1) and 
$(f)(\beta)$ of \ref{c5}(2).  The reason is in the
proof of \ref{b22}(3),(4), i.e. \cite[2.10,2.11]{Sh:900}.  
We can partially allow it, see \ref{b22}(4), the ``moreover", 
but not needed now.

\noindent
2) Note that for singular $\mu_2$ we get a
better result for free (as in the case $\kappa = \mu^+,\mu$ strong
limit singular of cofinality $> \theta$ is easier, see \ref{b24}(2)
and the proof of \ref{d39}.  

\noindent
3) In Definition \ref{c23} below note that $\vK$ is a weak form of
   $\tK$ and $\uK$ a weak form of $\qK$.
\end{remark}

\begin{discussion}
\label{c22}
Concerning Definition \ref{c1}, \ref{c5} and \ref{c23} below:

\noindent
0) For the $\vK$'s, $\uK$'s instead of dealing with some $\varphi
\in \Gamma^2_{\bold x}/\Gamma^3_{\bold x}$ allows us to deal with a
so called duplicate.

\noindent 
1) Note that $\tK^\oplus,\rK^\oplus,\vK^\oplus,\vK^\odot$ deals with
   $\Gamma^2_{\bold x}$ while $\sK^\oplus,\uK^\oplus,\uK^\otimes$ deals
with $\Gamma^3_{\bold x}$.

\noindent
2) Note that $\uK^\otimes,\vK^\odot$ have the witness $\bar{\bold w}$ as part
   of the $\bold m$ while $\uK^\oplus,\vK^\oplus$ do not.

\noindent
3) $\tK,\uK^\oplus,\uK^\otimes,\vK^\oplus,\vK^\odot$ deal with all
   formulas unlike $\rK^\oplus,\sK^\oplus$.

\noindent
4) $\vK,\uK$ is the projection of $\vK^\oplus,\uK^\oplus$ respectively 
to $\pK$.

\noindent
5) What is the point of $\vK^\oplus_{\lambda,\kappa,\bar\mu,\theta}$?
   We do not deal with every $\varphi = \varphi_0 \in \Gamma^2_{\bold
   x}$, we ``translate" the problem of $\varphi_0$ to a ``duplicate"
   $\varphi_2$ similar enough which is in $\Gamma^2_{\bar\psi}$.

\noindent
6) We can formulate \ref{c23}(3A)(d), more like
\ref{c23}(3C)$(\alpha)$.

\noindent
7) $\rK$ is intended as a step toward $\tK$ (or $\vK$).

\noindent
8) Note $\uK/\uK^\odot/\uK^\oplus/\uK^\otimes$
 relate like $\qK/\qK^\odot/\qK^\oplus/\qK^\otimes$.

\noindent
9) $\tK^\oplus$ is parallel to $\vK^\odot$ just as it is parallel to
   $\vK^\oplus$.

\noindent
10) Note that $\vK$ is defined as the projection of $\vK^\oplus$
    whereas $\tK$ is only provably the projection of $\tK^\oplus$ when
    $\cf(\kappa) > 2^\theta$.

\noindent
11) So $\vK^\oplus/\vK^\odot$ is not parallel to
$\uK^\oplus/\uK^\otimes$ but the latter is parallel to 
$\qK^\oplus/\qK^\otimes$.  
\end{discussion}

\begin{definition}
\label{c23}  
1) In \ref{c1}, \ref{c5} we adopt the conventions of \ref{b5}(2)
concerning the cardinals.

\noindent
1A) If $\bold m$ belongs to $\rK^\oplus$, let $\bold m =
 (\bold x_{\bold m},\bar\psi_{\bold m},r_{\bold m}) = 
(\bold x[\bold m],\bar\psi[\bold m],r[\bold m])$ and $M_{\bold m} =
 M_{\bold x[\bold m]}$, etc. and $\Gamma^2_{\bold m} = 
\Gamma^2_{\bold x[\bold m]}$, 
this may well be $\ne \Gamma^2_{\bar\psi[\bold m]}$, see \ref{c5}(1)(d).

\noindent
2) We define a two-place relation $\le_1$ on
$\rK^\oplus_{\kappa,\bar\mu,\theta}: (\bold x_1,\bar\psi_1,r_1) 
\le_1 (\bold x_2,\bar\psi_2,r_2)$ \when \,
 $\bold x_1 \le_1 \bold x_2$ (in $\pK_{\kappa,\bar\mu,\theta}$),
$\bar\psi_1 = \bar\psi_2 \rest \Gamma^2_{\bar\psi_1}$ 
(but dummy variables may be added) and $r_1 \subseteq r_2$.

\noindent
2A) Similarly for $\sK^\oplus_{\kappa,\bar\mu,\theta}$.

\noindent
3) Let $\tK^\oplus_{\lambda,\kappa,\bar\mu,\theta}$ be the class of 
$(\bold x,\bar\psi,r) \in \rK^\oplus_{\lambda,\kappa,\bar\mu,\theta}$
such that $\Gamma^2_{\bar\psi} = \Gamma^2_{\bold x}$ and 
$r$ is a complete type, over $\emptyset$ if not said otherwise.

\noindent
3A) Let $\vK^\oplus_{\lambda,\kappa,\bar\mu,\theta}$ be the class of
$\bold m = (\bold x,\bar\psi,r) \in 
\rK^\oplus_{\lambda,\kappa,\bar\mu,\theta}$ such that:
for every $\varphi \in \Gamma^2_{\bold x}$ there is an
$(\bold m,\varphi)$-duplicate $\bold w = (\eta_0,\nu_0,\eta_1,\nu_1,
\eta_2,\nu_2,\eta_3,\nu_3,\varphi_0,\varphi_1,\varphi_2)$ which
means\footnote{may use a $\bar d_\eta \rest u,\bar x_{\bar d_\eta
\rest u} = \bar x_u$ instead of $\bar d_\eta,\bar x_{\bar d,\eta}$; as we
may omit $\bar x'_{\bar c}$, no real change, in particular for normal
$\bold x$ it is the same}
\mn
\begin{enumerate}
\item[$(a)$]  $\varphi = \varphi_0 = \varphi_0(\bar x_{\bar d},\bar x_{\bar
c},\bar x'_{\bar d},\bar x'_{\bar c},\bar y) \in \Gamma^2_{\bold x}$ (as in
Definition \ref{b16}(0)(A), \ref{c5}(1)(d))
\sn
\item[$(b)$]   $\eta_0,\eta_1,\eta_2,\eta_3 
\in {}^{\omega >}\ell g(\bar d_{\bold
x})$ and $\ell g(\eta_0) = \ell g(\eta_1),\ell g(\eta_2) = \ell g(\eta_3)$
\sn
\item[$(c)$]  $\nu_0,\nu_1,\nu_2,\nu_3 \in {}^{\omega >}\ell g(\bar c_{\bold
x})$ and $\ell g(\nu_0) = \ell g(\nu_1),\ell g(\nu_2) = \ell g(\nu_3)$
\sn
\item[$(d)$]  $\varphi_1 = \varphi_1(\bar x_{\bar d,\eta_1},\bar
x_{\bar c,\nu_1},\bar x'_{\bar d,\eta_3},\bar x'_{\bar c,\nu_3},\bar
y) \equiv \varphi_0$
\sn
\item[$(e)$]  $\varphi_1(\bar x_{\bar d,\eta_0},\bar x_{\bar
c,\nu_0},\bar x'_{\bar d,\eta_2},\bar x'_{\bar c,\nu_2},\bar y) \equiv
\varphi_2 = \varphi_2(\bar x_{\bar d},
\bar x_{\bar c},\bar x'_{\bar d},\bar x'_{\bar c},\bar y) \in 
\Gamma^2_{\bold x}$
\sn
\item[$(g)$]  ${\gC} \models \varphi_1[\bar d_{\bold x,\eta_0},
\bar c_{\bold x,\nu_0},\bar d_{\eta_2},\bar c_{\nu_2},\bar b] 
\equiv \varphi_1 [\bar d_{\bold x,\eta_1},\bar c_{\bold x,\nu_1},
\bar d_{\eta_3},\bar c_{\eta_3},\bar b]$ for every 
$\bar b \in {}^{\ell g(\bar y)}(M_{\bold x})$
\sn
\item[$(h)$]  $\varphi_2 \in \Gamma^2_{\bar\psi}$.
\end{enumerate}
\mn
3B) For $\bold x \in \pK_{\kappa,\bar\mu,\theta}$ we say
$\Gamma$ is $\bold x-\vK$-large where $\Gamma \subseteq \Gamma^2_{\bold
x}$ \when \, for every $\varphi \in 
\Gamma^2_{\bold x}$ there is $\bold w$
satisfying clause (a)-(g) of part (3A) and $\varphi_2 \in
\Gamma$.

\noindent
3C) Let $\uK_{\lambda,\kappa,\bar \mu,\theta}$ be the class of $\bold x
\in \pK_{\lambda,\kappa,\bar\mu,\theta}$
such that: for every $\varphi(\bar x_{\bar d},\bar x_{\bar c},\bar y)
\in \Gamma^1_{\bold x}$ there is a weak $(\bold x,\varphi)$-duplicate
$\bold w = (\eta_0,\nu_0,\eta_1,\nu_1,\varphi_0,\varphi_1,\varphi_2)$
meaning $\varphi_0 = \varphi$ and\footnote{We may demand
  $\nu_1=\nu_0$; it seems there is no serious diference.}:
\mn
\begin{enumerate}
\item[$(a)$]  $(\eta_1,\nu_1) = \supp(\varphi)$, (see
\ref{b14}(4)), i.e. $\eta_1,\nu_1$ list $w,v$ respectively for some
$(w,v) \in \supp(\varphi)$
\sn
\item[$(b)$]  $\eta_0,\eta_1 \in {}^{\omega >}\ell g(\bar d_{\bold
x})$ and $\ell g(\eta_0) = \ell g(\eta_1)$
\sn
\item[$(c)$]  $\nu_0,\nu_1 \in {}^{\omega >} \ell g(\bar c_{\bold x})$ and
$\ell g(\nu_0) = \ell g(\nu_1)$
\sn
\item[$(d)$]  $\varphi_1 = \varphi_1(\bar x_{\bar d,\eta_1},
\bar x_{\bar c,\nu_1},\bar y) \equiv \varphi_0 = \varphi_0(x_{\bar
d},\bar x_{\bar c},\bar y) \in \Gamma^1_{\bold x}$
\sn
\item[$(e)$]  $\varphi_1(\bar x_{\bar d,\eta_0},
\bar x_{\bar c,\nu_0},\bar y) \equiv \varphi_2 = 
\varphi_2(\bar x_{\bar d},\bar x_{\bar c},\bar y) \in \Gamma^1_{\bold x}$
\sn
\item[$(f)$]  ${\gC} \models ``\varphi_1[\bar d_{\bold x,\eta_0},
\bar c_{\bold x,\nu_0},\bar b] \equiv \varphi_2
[\bar d_{\bold x,\eta_1},\bar c_{\bold x,\nu_1},\bar b]"$ for
every $\bar b \in {}^{\ell g(\bar y)}(M_{\bold x})$
\sn
\item[$(g)$]  some $\psi$ illuminates $(\bold x,\varphi_2)$.
\end{enumerate}
\mn
3D) For $\bold x \in \pK_{\kappa,\bar\mu,\theta}$, we say
$\Gamma \subseteq \Gamma^1_{\bold x}$ is $\bold x-\uK$-large 
\when \, for every $\varphi \in \Gamma^1_{\bold x}$ there is a weak $(\bold
x,\varphi)$-duplicate $\bold w$, see part (3C).

\noindent
3E) Let $\vK_{\lambda,\kappa,\bar\mu,\theta}$ be the class of $\bold x
\in \pK_{\lambda,\kappa,\bar\mu,\theta}$ such that for every 
$A \subseteq M_{\bold x}$ of cardinality $< \kappa$ 
we can find $\bar\psi$ and $(\bar c',\bar d')$ such that:
\mn
\begin{enumerate}
\item[$(\alpha)$]  as in \ref{c5}(1)(f)
\sn
\item[$(\beta)$]  if $\varphi \in \Gamma^2_{\bar\psi}$ then
  $\psi_\varphi(\bar x_{\bar d},\bar c_{\bold x},\bar d',\bar c')
  \vdash \tp_\varphi(\bar d_{\bold x},\bar c_{\bold x} \char 94 \bar
  d' \char 94 \bar c' \dotplus A)$
\sn
\item[$(\gamma)$]  $\Gamma^2_{\bar\psi}$ is $\bold x-\vK$ large, see
  part (3B).
\end{enumerate}
\mn
3F) We define $\uK^\oplus_{\lambda,\kappa,\bar\mu,\theta}$ as the 
class of triples $\bold n = (\bold x,\bar\psi,r) \in 
\sK^\oplus_{\lambda,\kappa,\bar\mu,\theta}$ such that 
$\Gamma^3_{\bar\psi} \subseteq \Gamma^1_{\bold x}$ is $\bold
x-\uK$-large. We define $\vK^\oplus_{\lambda,\kappa,\bar\mu,\theta}$
as the class of $\bold n = (\bold x,\bar\psi,r)$ which belongs to
$\rK^\oplus_{\lambda,\kappa,\bar\mu,\theta}$ and $\Gamma^2_{\bar\psi}
\subseteq \Gamma^2_{\bold x}$ is $\bold x-\vK$-large.

\noindent
3G) We define $\uK^\otimes_{\kappa,\bar\mu,\theta}$ as the class of
$\bold n = (\bold x,\bar\psi,r,\bar{\bold w})$ such that $(\bold
x,\bar\psi,r) \in \sK^\oplus_{\kappa,\bar\mu,\theta}$ so
$\Gamma^3_{\bar\psi} \subseteq \Gamma^3_{\bold x},\bar\psi 
= \langle \psi_\varphi(\bar x_{\bar d},\bar x_{\bar
c},\bar y_{[\theta]}):\varphi = \varphi(\bar x_{\bar d},\bar x_{\bar
c},\bar y_{[\theta]},\bar z) \in 
\Gamma^3_{\bar\psi}\rangle,\bar{\bold w} = \langle \bold
w_\varphi:\varphi \in \Gamma^3_{\bold x}\rangle$ recalling
\ref{b16}(0)(B) and: for every
$\varphi = \varphi(\bar x_{\bar d},\bar x_{\bar c},\bar
y_{[\theta]},\bar z) \in \Gamma^3_{\bold x}$, $\bold w_\varphi$ is a
witness of the form
$(\eta_0,\nu_0,\eta_1,\nu_1,\varphi_0,\varphi_1,\varphi_2)$ (see below) and for
every $A \subseteq M_{\bold x}$ of cardinality $< \mu_2$ there is a
solution $\bar e$, i.e. $\bar e$ solves $(\bold n,A)$ 
which means: $\bar e \in {}^\theta(M_{\bold x})$
solves $(\bold x,\bar\psi,r)$ recalling \ref{c5}(2)(f) and for 
every $\varphi \in \Gamma^3_{\bold x}$ 
the witness $\bold w_\varphi$ satisfies:
\mn
\begin{enumerate}
\item[$(a)$]  $(\eta_1,\nu_1) = \supp(\varphi(\bar x_{\bar
d},\bar x_{\bar c};\bar y_{[\theta]},\bar z))$
\sn
\item[$(b)$]  $\eta_0,\eta_1 \in {}^{\omega >}\ell g(\bar d_{\bold
x})$ and $\ell g(\eta_0) = \ell g(\eta_1)$
\sn
\item[$(c)$]  $\nu_0,\nu_1 \in {}^{\omega >}\ell g(\bar c_{\bold x})$
and $\ell g(\nu_0) = \ell g(\nu_1)$
\sn
\item[$(d)$]  $\varphi_2 = \varphi_2(\bar x_{\bar d,\eta_1},
\bar x_{\bar c,\nu_0},\bar y_{[\theta]},\bar z) \equiv \varphi_0 =
\varphi(\bar x_{\bar d},\bar x_{\bar c},\bar y_{[\theta]},\bar z) \in
\Gamma^3_{\bold x}$
\sn
\item[$(e)$]  $\varphi_2(\bar x_{\bar d,\eta_0},
\bar x_{\bar c,\nu_0},\bar y_{[\theta]},\bar z) \equiv \varphi_2 =
\varphi_2(\bar x_{\bar d},\bar x_{\bar c},\bar y_{[\theta]},\bar z) \in
\Gamma^3_{\bold x}$
\sn
\item[$(f)$]  $\gC \models ``\varphi_1[\bar d_{\bold x,\eta_0},\bar
c_{\bold x,\nu_0},\bar e,\bar b] \equiv \varphi_2[\bar d_{\bold
x,\eta_1},\bar c_{\bold x,\nu_2},\bar e,\bar b]"$ for every $\bar b
\in {}^{\ell g(\bar z)}A$
\sn
\item[$(g)$]  $\varphi_2 \in \Gamma^3_{\bar\psi}$.
\end{enumerate}
\mn
4) Let $\le^+_1$ be the following two-place relation on $\rK^\oplus$:

\noindent
$(\bold x_1,\bar\psi_1,r_1) \le^+_1 (\bold x_2,\bar\psi_2,r_2)$ 
\Iff \, $(\bold x_1,\bar\psi_1,r_1) \le_1 (\bold x_1,\bar\psi_2,r_2)$ and
$\Gamma^2_{\bold x_1} \subseteq \Gamma^2_{\bar\psi_2}$.

\noindent
4A) Let $\le^\odot_1$ be the following two-place relation on
$\rK^\oplus$:

$\bold m \le^\odot_1 \bold n$ \Iff \, $\bold m \le_1 \bold n$ and
if $\varphi = \varphi(x_{\bar d[\bold m]},\bar x_{\bar c[\bold m]},\bar
x'_{\bar d[\bold m]},\bar x'_{\bar c[\bold m]},\bar y) \in
\Gamma^2_{\bold x[\bold m]}$ \then \, some $\bold w$ is an $(\bold
n,\varphi)$-duplicate, see part (3A).

\noindent
4B) We define $\le^+_{1,\Delta},<^\odot_{1,\Delta}$ similarly where
$\Delta \subseteq \Gamma^2_{\bold x_1}$ and we demand to deal only with
$\varphi \in \Delta$.

\noindent
4C) We define $\vK^\otimes_{\lambda,\kappa,\bar\mu,\theta}$ as the class of
$\bold n = (\bold x,\bar\psi,r,\bar{\bold w})$ such that:
\mn
\begin{enumerate}
\item[$(a)$]  $(\bold x,\bar\psi,r) \in 
\vK^\oplus_{\lambda,\kappa,\bar\mu,\theta}$
\sn
\item[$(b)$]  $\bar{\bold w} = \langle \bold w_\varphi:\varphi \in
\Gamma^2_{\bold x}\rangle$
\sn
\item[$(c)$]  for $\varphi = \varphi(\bar x_{\bar d},\bar x_{\bar
c},\bar x'_{\bar d},\bar x'_{\bar c},\bar y) \in \Gamma^2_{\bold x}$
we have\footnote{may use $\bold u(\varphi,u_0)$ but this can be 
absorbed as we consider $u_1 = \{0\}$
 for $\varphi = (x_{d_0} = x_{d_0})$}  $\bold w_\varphi$ is  
a\footnote{we may restrict ourselves to normal $\bold x$ (and
 $\bold m$) and then demand $v_1 = w_1 \cap \ell g(\bar c_{\bold x})$} 
$((\bold x,\bar\psi,r),\varphi)$-duplicate, see part (3A).
\end{enumerate}
\mn
5) If $\langle \bold m_\varepsilon:\varepsilon < \delta\rangle$ is
$\le_1$-increasing in $\rK^\oplus$ then we let $\bold m_\delta :=
\cup\{\bold m_\varepsilon:\varepsilon < \delta\}$ be naturally
defined (uniquely up to ``very similar") but it is not clear that
$\bold m_\delta \in \rK^\oplus$, the problem is with \ref{c5}(1)(f).  
Similarly in the other cases.

\noindent
6) We define reducts, $\bold m \rest \tau$ for $\tau \subseteq
\tau(T)$ naturally.
\end{definition}
\bn
Note
\begin{observation}
\label{c27}  
Let $\kappa > \theta$ and $\bar \mu$ be as in Definition \ref{b05}.

\noindent
1) If $\iota_{\bold x}=2$ and $\bold x \in \qK_{\kappa,\bar\mu,\theta}$ and 
$\Rang(\bar d_{\bold x}) \nsubseteq M_{\bold x}$ \then \, 
$\tp(\bar d_{\bold x},\bar c_{\bold x} + B^+_{\bold x})$ is 
not realized in $M_{\bold x}$, 
moreover if $d'\in \Rang(\bar d_{\bold x}) \backslash
M_{\bold x}$ then $\tp(d',B^+_{\bold x} + \bar c_{\bold x})$ is not
realized in $M_{\bold x}$.

\noindent
1A) In part (1), even if $\iota_{\bold x} \ne 2$, it suffices to
assume ``$\bold x \in \qK'_{\kappa,\bar\mu,\theta}$ and $d' \in
\Rang(\bar d_{\bold x}) \backslash M_{\bold x}$. 

\noindent
2) $\bold x \in \tK_{\lambda,\kappa,\bar\mu,\theta} \Rightarrow \bold x \in 
\qK_{\lambda,\kappa,\bar\mu,\theta}$ and $\tK_{\lambda,\kappa,\bar\mu,\theta}
\subseteq \uK_{\lambda,\kappa,\bar\mu,\theta}$ and
$\tK^\oplus_{\lambda,\kappa,\bar\mu,\theta} \subseteq 
\vK^\oplus_{\lambda,\kappa,\bar\mu,\theta} \subseteq
\vK^\otimes_{\lambda,\kappa,\bar\mu,\theta}$. 

\noindent
3) For every $\kappa$-saturated $M$ there is $\bold x \in 
\tK_{\lambda,\kappa,\bar\mu,\theta}$ with 
$M_{\bold x} = M,\bar d_{\bold x} =
\langle \rangle = \bar c_{\bold x}$ hence $w_{\bold x} = \emptyset =
v_{\bold x},B_{\bold x} = \emptyset$.

\noindent
4) Assume $\cf(\kappa) > 2^\theta + |T|$.  \Then \, 
$\bold x \in \tK_{\lambda,\kappa,\bar\mu,\theta}$ iff for
some $\bar\psi,r$ we have $(\bold x,\bar\psi,r) \in 
\tK^\oplus_{\lambda,\kappa,\bar\mu,\theta}$ with $r$ a complete
type over $\emptyset$.

\noindent
4A) Similarly for $\vK_{\lambda,\kappa,\bar\mu,\theta},
\vK^\oplus_{\lambda,\kappa,\bar\mu,\theta}$.

\noindent
4B) If $\bold m = (\bold x,\bar\psi,r) \in 
\vK^\oplus_{\lambda,\kappa,\bar\mu,\theta}$ \then \, for some $\bold w$ we
have $(\bold x,\bar\psi,r,\bold w) \in 
\vK^\odot_{\lambda,\kappa,\bar\mu,\theta}$.

\noindent
5) If $\bold m \in \tK^\oplus_{\lambda,\kappa,\bar\mu,\theta}$ \then \,
$\bold x_{\bold m} \in \tK_{\lambda,\kappa,\bar\mu,\theta}$.

\noindent
5A)  $\bold m \in \vK^\oplus_{\kappa,\bar\mu,\theta}$ \then \,
$\bold x_{\bold m} \in \vK_{\kappa,\mu,\theta}$.

\noindent
6) If $\bold x \in \tK_{\lambda,\kappa,\bar\mu,\theta}$ and
$\bold y$ is defined like $\bold x$ replacing $\bar d_{\bold x}$
 by $\bar d_{\bold x} \char 94 \bar c_{\bold x}$ \then \,
  $\bold y \in \tK_{\lambda,\kappa,\bar\mu,\theta}$ and is normal.
   Similarly for $\vK_{\lambda,\kappa,\bar\mu,\theta},
\uK_{\lambda,\kappa,\bar\mu,\theta}$.
\end{observation}

\begin{PROOF}{\ref{c27}}
1) It is enough to prove the ``moreover".  Let $d'' \in
M_{\bold x}$ realize $\tp(d',B^+_{\bold x}+ \bar c_{\bold x})$, let
$t_* \notin v_{\bold x},v' = v + \{t_*\}$ and
$\bar c'' = \bar c_{\bold x} \char 94  \langle \bar c''_{t_*}
\rangle$, where $\bar c''_{t_*} = \langle d'' \rangle$ 
and let $\bold y$ be like $\bold x$
replacing $\bar c_{\bold x}$ by $\bar c'',v_{\bold x}$ by $v'$ and letting
$B_{\bold y,t_*} = \{d''\}$.  Clearly $\tp(\langle d'' \rangle,\bar c_{\bold
x} + M_{\bold x})$ is satisfiable in $B_{\bold y,t_*}$ hence by
\ref{b22}(5) does not $\iota_{\bold x}$-split over it.

So $\bold x \le_1 \bold y$ and clearly $\bold y$ is active in $t_*$,
hence ``$\bold x \notin \qK'_{\kappa,\bar\mu,\theta}"$, 
see Definition \ref{b16}(1).  Now if $\iota_{\bold x} = 2$ by claim
\ref{b22}(5) we have $\bold x \notin \qK_{\kappa,\bar\mu,\theta}$.

\noindent
1A) By the proof of part (1).

\noindent
2) For the first statement recall 
(from Definition \ref{c1}) that a consequence of $\bold x \in 
\tK_{\lambda,\kappa,\bar\mu,\theta}$ is the existence of solutions, 
 but this consequence for $\bold x \in 
\pK_{\lambda,\kappa,\bar\mu,\theta}$ implies $\bold x \in 
\qK_{\lambda,\kappa,\bar\mu,\theta}$ by Definition \ref{b16}(2)
so indeed $\bold x \in 
\tK_{\lambda,\kappa,\bar\mu,\theta} \Rightarrow \bold x \in
\qK_{\lambda,\kappa,\bar\mu,\theta}$.  Also for the other statements see
the definitions.

\noindent
3) Obvious (and see \ref{b18}(5)).

\noindent
4),4A),4B),5), 5A)  Easy. 
Read the definitions for $\Leftarrow$ and immitate \ref{b24} for
 $\Rightarrow$. 

\noindent
6) Straight (as in \ref{b18}(3)).  
\end{PROOF}

\begin{observation}
\label{c29}  
0) If $\bold m \in \rK^\oplus_{\kappa,\bar\mu,\theta}$ 
\then \, $\bold m \in \tK^\oplus_{\kappa,\bar\mu,\theta} 
\Leftrightarrow \bold m \le^+_1 \bold m$ 
and $\bold m \in \vK^\oplus_{\kappa,\bar\mu,\theta}
  \Leftrightarrow \bold m \le^\odot_1 \bold m$.

\noindent
1) $\le^+_1$ partially ordered $\rK^\oplus_{\kappa,\bar\mu,\theta}$ except
   that possibly $\neg(\bold m \le^+_1 \bold m)$.

\noindent
1A) Similarly $\le^\odot_1$.

\noindent
2) Also on $\rK^\oplus_{\kappa,\bar\mu,\theta}$ we have
$\le^+_1 \subseteq \le^\odot_1 \subseteq
\le_1$ and $\bold m_1 \le_1 \bold m_2 \le^+_1 \bold m_3 \le_1 \bold
m_4$ implies $\bold m_1 \le^+_1 \bold m_4$ and $\bold m_1 \le_1 \bold
m_2 \le^\odot_1 \bold m_3 \le_1 \bold m_4$ implies $\bold m_1
\le^\odot_1 \bold m_4$.

\noindent
3) If $\bold x \in \pK_{\kappa,\bar\mu,\theta}$ \then \, $(\bold
x,\langle \rangle,\emptyset) \in \rK^\oplus_{\kappa,\bar\mu,\theta}$.

\noindent
4) If $\bold m = (\bold x,\bar\psi,r) \in 
\rK^\oplus_{\kappa,\bar\mu,\theta}$ and $\bold x,\bold y \in 
\pK_{\kappa,\bar\mu,\theta}$ are very similar \then \,
$\bold n := (\bold y,\bar\psi,r) \in 
\rK^\oplus_{\kappa,\bar\mu,\theta}$ and $\bold m \in 
\tK^\oplus_{\kappa,\bar\mu,\theta} \Leftrightarrow \bold n
   \in \tK^\oplus_{\kappa,\bar\mu,\theta}$ and $\bold m \in
   \vK^\oplus_{\kappa,\bar\mu,\theta} \Leftrightarrow \bold n
 \in \vK^\oplus_{\kappa,\mu,\theta}$.
\end{observation}

\begin{observation}
\label{c31} 
1) If $\bold x,\bold y \in \pK_{\kappa,\bar\mu,\theta}$ are
very similar \then \, $\bold x \in 
\qK_{\kappa,\bar\mu,\theta}$ iff $\bold y \in 
\qK_{\kappa,\bar\mu,\theta}$.

\noindent
2) Similarly for $\qK'_{\kappa,\bar\mu,\theta},
\rK^\oplus_{\kappa,\bar\mu,\theta},
\tK^\oplus_{\kappa,\bar\mu,\theta},\vK^\oplus_{\kappa,\bar\mu,\theta},
\uK_{\kappa,\bar\mu,\theta}$ and $\uK^\oplus_{\kappa,\bar\mu,\theta}$.
E.g., for such $\bold x,\bold y$: for any $\bar\psi,r$ we have
$(\bold x,\bar\psi,r) \in \rK^\oplus_{\kappa,\bar\mu,\theta}$ iff
$(\bold y,\bar\psi,r) \in \rK^\oplus_{\kappa,\bar\mu,\theta}$. 
\end{observation}
\bigskip

\subsection{Sequence homogeneity and indiscernibles} \
\bigskip

We now try to prove that decompositions from $\tK$ and $\vK$ are ``good" and
``helpful".  We prove for $\bold x \in \tK_{\kappa,\bar\mu,\theta}$
that $M_{[\bold x]} = M_{[B^+_{\bold x} + \bar c_{\bold x} + \bar
d_{\bold x}]}$ defined in \ref{b5}(6), is
$\kappa$-sequence-homogeneous, see \ref{z5}, this is nice,
and help to prove that there are few types up to conjugacy 
because if $M,N$ are $(\bold D,\kappa)$-sequence homogeneous models 
of cardinality $\kappa$ then they are isomorphic.
\bigskip

\begin{theorem}
\label{c32}   
\underline{The sequence homogeneous Theorem}
1) If $\bold x \in \tK_{\kappa,\bar\mu,\theta}$ 
\then \, $M_{[\bold x]}$ is a
$\kappa$-sequence-homogeneous model for the finite diagram which
we call $\bold D_{\bold x}$; see Definition \ref{z25}(1), \ref{b5}(6).

\noindent
1A) Similarly for $\bold x \in \vK_{\kappa,\bar\mu,\theta}$.

\noindent
2) Moreover, if $(\bold x,\bar\psi,r) \in
\tK^\oplus_{\kappa,\bar\mu,\theta}$ or just 
$(\bold x,\bar\psi,r,\bar{\bold w}) \in 
\vK^\otimes_{\kappa,\bar\mu,\theta}, \bold x$ is smooth  
and $r$ is a complete type \then \, $\bold D_{\bold x}$ depends just on
$T,\bar\psi,r,\bar{\bold w},B^+_{\bold x},\langle D_{\bold x,i}:i \in
v_{\bold x} \backslash u_{\bold x}\rangle$ and 
$\tp(\bar d_{\bold x} \char 94 
\bar c_{\bold x},B^+_{\bold x})$.  That is, if $\bold m _\ell = (\bold
x_\ell,\bar\psi,r,\bar{\bold w}) \in \vK^\otimes_{\kappa,\mu,\theta}$
for $\ell=1,2$ and $\bold x_1,\bold x_2$ are smooth and similar as witnessed by
$g$, see Definition \ref{b30}(2) \then \, 
$g$ maps $\bold D_{\bold x_1}$ onto $\bold D_{\bold x_2}$.
\end{theorem}

\begin{remark}
\label{c34}  We use a little less than the requirements in the definitions
of $\tK_{\kappa,\bar\mu,\theta},\vK_{\kappa,\bar\mu,\theta}$; see the
proof, i.e. in $(*)_1$ below there is $\psi(\bar x_{\bar d},\bar
c,\bar d_*) \in \tp(\bar d,\bar c \char 94 \bar d_*)$ such that
$\psi(\bar x_{\bar d},\bar c \char 94 \bar d_*) \vdash \varphi(\bar x_{\bar
d},\bar c,\bar b,a_1)$ \underline{but} $\psi$ may depend on $b_1,a_1$.
\end{remark}

\begin{PROOF}{\ref{c32}}
1) Let $B = B^+_{\bold x}$
and as usual let $\bar c = \bar c_{\bold x},\bar d = \bar d_{\bold x}$.
So it suffices to prove that $M^+ := M_{[\bold x]} = 
M_{[B + \bar c + \bar d]}$ is a $\kappa$-sequence-homogeneous model.

Let $f$ be an elementary mapping from $A_1 \subseteq  M^+$ onto $A_2
\subseteq M^+$ in the sense of $M^+$ and $|A_1| < \kappa$ 
and $b_1 \in M$ and we should find
such $g \supseteq f$ for which $b_1 \in \Dom(g)$, this
suffices.  Let $A = B + A_1 + A_2 + b_1$.
Let $f_0 = f,f_1 = f \cup \id_B$ and $f_2 = f_1
\cup \id_{\bar c + \bar d}$.
By the definition of $M^+$ the mappings $f_1,f_2$ are
elementary (in the sense of ${\gC}$, the 
default value).  As $A \subseteq M$ has cardinality $< \kappa$,
recalling $\bold x \in \tK_{\kappa,\bar\mu,\theta}$
there is $\bar c_* \char 94 \bar d_*$ in
$M_{\bold x}$ realizing $\tp(\bar c \char 94 \bar d,A)$ such that:
\mn
\begin{enumerate}
\item[$\odot_0$]   $\tp(\bar d,\bar c + \bar d_* + \bar c_*) \vdash 
\tp(\bar d,\bar c + \bar d_* + \bar c_* +A)$.  
\end{enumerate}
\mn
But actually we need just
\mn
\begin{enumerate}
\item[$\odot'_0$]   $\tp(\bar d,\bar c + \bar d_* + \bar c_*) \vdash 
\tp(\bar d,\bar c + A)$.  
\end{enumerate}
\mn
By the choice of $(\bar c_*,\bar d_*)$, clearly 
the following function $h$ is elementary for ${\gC}$:
\mn
\begin{enumerate}
\item[$\odot''_0$]  $\Dom(h) = A + \bar c + \bar d$ and $h \rest A$ 
is the identity, $h(\bar c \char 94 \bar d) = \bar c_* \char 94 \bar d_*$.
\end{enumerate}
\mn
Let $f'_2 = f_1\cup \id_{\bar c_* + \bar d_*}$, so $f'_2 = h
\circ f_2 \circ h^{-1}$ but $h$ and $f_2$ are elementary so $f'_2$
is elementary too.  Clearly $\Dom(f'_2) = \Dom(f_1) + \bar c_* + \bar
d_*= B+A_1 + \bar c_* + \bar d_* \subseteq M_{\bold x}$ and
$\Rang(f'_2) = \Rang(f_1) + \bar c_* + \bar d_* = B + A_2 + \bar c_* +
\bar d_* \subseteq M_*$. 
Hence there is an
elementary mapping $g_1$ such that $g_1 \supseteq f'_2$ and $\Dom(g_1) =
B + A_1 + \bar c_* + \bar d_* + b_1$ and \wilog \, 
$b_2 := g_1(b_1)$ belongs to $M_{\bold x}$ recalling
$\Rang(f'_2) \subseteq M_{\bold x}$ and $M_{\bold x}$ is
$\kappa$-saturated.

Let $g_2 = g_1 \cup \id_{\bar c}$, next:
\mn
\begin{enumerate}
\item[$\odot_1$]  $(a) \quad g_2$ is an elementary mapping
\sn
\item[${{}}$]  $(b) \quad g_2$ is with domain $B + A_1 + \bar c_* +
\bar d_* + b_1 + \bar c$.
\end{enumerate}
\mn
[Why?  Clause (a) as $\tp(\bar c,M_{\bold x})$ does not split over $B$ and $g_2
\supseteq g_1 \supseteq f'_2 \supseteq f_1 \supseteq \id_B$.
Clause (b) holds as $\Dom(g_2) = \Dom(g_1) \cup \bar c$ by the choice
of $g_1$.] 

Now assume for awhile:
\mn
\begin{enumerate}
\item[$\odot_2$]   $\bar a_1 \in {}^{\omega >}(B + A_1)$ 
and\footnote{we can strengthen the demand on $\bar a_1$ to
$\bar a_1 \in {}^{\omega >}(A_1 + B + \bar d_* + \bar c_*)$ and change
according in later cases}
${\gC} \models \varphi[\bar d,\bar c,b_1,\bar a_1]$; let
$\bar a_2 = f_1(\bar a_1)$.
\end{enumerate}
\mn
Now
\mn
\begin{enumerate}
\item[$(*)_0$]  $(a) \quad f_2 \supseteq f_1$ and $g_2 \supseteq g_1
\supseteq f'_2 \supseteq f_1$
\sn
\item[${{}}$]  $(b) \quad \bar a_1 \subseteq (B + A_1) = 
\Dom(f_1)$, hence
\sn
\item[${{}}$]  $(c) \quad g_2(\bar a_1) = \bar a_2$; also
\sn
\item[${{}}$]  $(d) \quad g_2(b_1) = g_1(b_1) = b_2$
\sn
\item[${{}}$]  $(e) \quad g_2$ is the identity on $B + \bar c_* + \bar
d_* + \bar c$.
\end{enumerate}
\mn
[E.g. why clause (e)?  By their choice, $f_1$ is the identity on
$B,f'_2$ is the identity on $\bar c_* + \bar d_*$ and $g_2$ is the
identity on $\bar c$ hence by clause (a) we are done.]

We know that $\tp(\bar d,\bar c + \bar d_* + \bar c_*) 
\vdash \tp(\bar d,\bar c + A)$ by $\odot_0$ or $\odot'_0$
hence, (recalling $\bar a_1 \subseteq B + A_1 \subseteq A$, see
$(*)_0(b)$ + the choice of $A$ and $b_1 \in A$ by the choice of $A$):
\mn
\begin{enumerate}
\item[$(*)_1$]  $\tp(\bar d,\bar c + \bar d_* + \bar c_*) 
\vdash \varphi(\bar x_{\bar d},\bar c,b_1,\bar a_1)$.
\end{enumerate}
\mn
So applying $g_2$ recalling $(*)_0(e)$
\mn
\begin{enumerate}
\item[$(*)_2$]  $g_2$ maps $\tp(\bar d,\bar c + \bar d_* + \bar c_*)$ to itself
\end{enumerate}
\mn
As $(*)_1 + (*)_2$ hold and $g_2(b_1) =
b_2,g_2(\bar a_1) = \bar a_2$ and $g_2(\bar c) = \bar c$ (by
$(*)_0(d),(c),(e)$ respectively) recalling $g_2$ is an elementary
mapping by $\odot_1(a)$ we get
\mn
\begin{enumerate}
\item[$(*)_3$]  $\tp(\bar d,\bar c + \bar d_* + \bar c_*) 
\vdash \varphi(\bar x_{\bar d},\bar c,b_2,\bar a_2)$.
\end{enumerate}
\mn
So it follows that:
\mn
\begin{enumerate}
\item[$\odot_3$]  ${\gC} \models \varphi[\bar d,\bar c,b_2,\bar a_2]$.
\end{enumerate}
\mn
We have proved $\odot_2 \Rightarrow \odot_3$ when $\bar a_1$ was any
finite sequence from $B + A_1$.
Recalling $(*)_0(e)$ and $g_2(b_1) = b_2,g_2(\bar a_1) =
\bar a_2$ and $g_2(\bar c) = \bar c$, this means that 
$g_3 := (g_2 \rest (B + A_1 + b_1 + \bar c)) \cup \id_{\bar d}$ 
is an elementary mapping, so the function $g_3 =  
g_3 \rest (B + A_1 + b_1 + \bar c + \bar
d)$ is an elementary mapping of ${\gC}$, so as $g_3 \rest (B + \bar c +
\bar d)$ is the identity clearly $g := g_3 \rest (A_1 + b_1)$ is an
elementary mapping in the sense of $M^+$, so $g$ is as required.

\noindent
1A) The proof above works now, too, except that not necessarily $\odot_0$
holds (and so $\odot'_0$, too) which was used 
only in proving $(*)_1$ so in proving $\odot_2
\Rightarrow \odot_3$.  
Again it suffices to prove $\odot_3$ assuming $\odot_2$.  Let
$\varphi = \varphi(\bar x_{\bar d},\bar x_{\bar c};z,\bar y)$ so there
is $\varphi_0 = \varphi_0(\bar x_{\bar d,\eta_1},\bar x_{\bar
c,\nu_1};z,\bar y)$ equivalent to $\varphi$ for some $\eta_1 \in {}^{\omega
>} \ell g(\bar d_{\bold x}),\nu_1 \in {}^{\omega >}\ell g(\bar
c_{\bold x})$, hence by the Definition
\ref{c23}(3A), a degenerated case\footnote{this is
a weak version of \ref{c23}(3A) as $\eta_2,\nu_2$ disappear}, 
there are $\eta_0,\nu_0$ such that
\mn

\begin{enumerate}
\item[$\oplus$]  $(a) \quad \eta_0 \in {}^{\omega >} \ell g(\bar d_{\bold
x}),\nu_0 \in {}^{\omega >} \ell g(\bar c_{\bold x})$
\sn
\item[${{}}$]  $(b) \quad \ell g(\eta_0) = \ell g(\eta_1)$ and
$\ell g(\nu_0) = \ell g(\nu_1)$, all finite
\sn
\item[${{}}$]  $(c) \quad {\gC} \models ``\varphi_0[\bar d_{\bold x,\eta_1},
\bar c_{\bold x,\nu_1},b',\bar a] \equiv \varphi_0[\bar d_{\bold x,\eta_0},
\bar c_{\bold x,\nu_0},b',\bar a]"$

\hskip25pt  for every $b' \in M_{\bold x},\bar a' 
\in {}^{\ell g(\bar y)}(M_{\bold x})$
\sn 
\item[$\oplus'$]  there is $\bar d_* \char 94 \bar c_*$ from $M_{\bold
x}$ realizing $\tp(\bar d_{\bold x} \char 94 \bar c_{\bold x},A)$ such
that
\sn
\item[${{}}$]   $(d) \quad \tp(\bar d_{\bold x,\eta_0},
\bar c_{\bold x,\nu_0} + \bar d_* + \bar c_*) 
\vdash \{\varphi_0(\bar x_{\bar d,\eta_0},
\bar c_{\bold x,\nu_0},b',\bar a'):b' \in A,\bar a' \in
{}^{\ell g(\bar a_1)}A$

\hskip25pt and ${\gC} \models \varphi_0[\bar d_{\bold x,\eta_0},
\bar c_{\bold x,\nu_0},b',\bar a']\}$.
\end{enumerate}
\mn
Now as in the proof of part (1)  above we assume
\mn
\begin{enumerate}
\item[$\odot_2$]  $\bar a_1 \in {}^{\omega >}(B + A_1)$ and
${\gC} \models \varphi[\bar d,\bar c,b_1,\bar a_1]$.
\end{enumerate}
\mn
By the choice of $\varphi_0$ this means that 
${\gC} \models \varphi_0[\bar d_{\bold x,\eta_1},
\bar c_{\bold x,\nu_1},b_1,\bar a_1]$ hence by $\oplus(c)$ we have:
\mn
\begin{enumerate}
\item[$\odot'_2$]  ${\gC} \models \varphi_0[\bar d_{\bold
x,\eta_0},\bar c_{\bold x,\nu_0},b_1,\bar a_2]$.
\end{enumerate}
\mn
Recalling $\oplus'(d)$, for this formula the proof of $\odot_2 \Rightarrow
\odot_3$ in part (1) works so ${\gC} \models \varphi_0
[\bar d_{\bold x,\eta_0},\bar c_{\bold x,\nu_0},b_2,\bar a_2]$.  Using
$\oplus(c)$ again this implies 
$\gC \models \varphi_0[\bar d_{\bold x,\eta_1},\bar
c_{\bold x,\nu_1},b_2,\bar a_2]$.  As this holds for any $\bar a_1 \in
{}^{\omega >}(B + A_1)$ we finish as in part (1).

\noindent
2) Assume
\mn
\begin{enumerate}
\item[$\boxplus_{\tK}$]   $\bold m_\ell = (\bold x_\ell,\bar\psi,r) \in 
\tK^\oplus_{\kappa,\bar\mu,\theta}$ for $\ell=1,2$

\underline{or}
\sn
\item[$\boxplus_{\vK}$]   $\bold m_\ell =
(\bar{\bold x}_\ell,\psi,r,\bar{\bold w}) \in 
\vK^\otimes_{\kappa,\bar\mu,\theta}$ for $\ell=1,2$.
\end{enumerate}
\mn
Assume further that $g$ witnesses $\bold x_1,\bold x_2$ are similar;
the proof is like the proof of parts (1),(1A) but we give some
details.  \Wilog \, $g$ witnesses $\bold x_1,\bold x_2$ are very
similar, so $g_0 = g \rest B^+_{\bold x_1}$ is an elementary mapping
from $B^+_{\bold x_1}$ onto $B^+_{\bold x_2}$.

Let $\bar c_\ell = \bar c_{\bold x_\ell},\bar d_\ell = \bar d_{\bold
x_\ell}$.
Assume $\bar a_\ell \in {}^{\omega >}(M_{\bold x_\ell})$ for
   $\ell=1,2$.  Let $\ell \in \{1,2\}$ choose $\bar c^\ell_* \char 94
   \bar d^\ell_*$ as in Definition for $A_\ell := B^+_{\bold x_\ell} \cup
   \bar a_\ell$, so
\mn
\begin{enumerate}
\item[$(*)_{1,\ell}$]  $(a) \quad \bar d^\ell_*,\bar c^\ell_*$ are
from $M_{\bold x_\ell}$
\sn
\item[${{}}$]  $(b) \quad \bar d^\ell_*,\bar c^\ell_*$ realize
$\tp(\bar d_\ell \char 94 \bar c_\ell,A_\ell)$
\sn
\item[${{}}$]  $(c) \quad \bar d_\ell \char 94 \bar c_\ell \char 94
\bar d^\ell_*  \char 94 \bar c^\ell_*$ realizes $r_\ell := r[\bold m_\ell]$
\sn
\item[${{}}$]  $(d)(\alpha) \quad$ if $\boxplus_{\tK}$ then
$\tp(\bar d_\ell,\bar c + \bar d^\ell_* +
\bar c^\ell_*) \vdash \tp(\bar d_\ell,\bar c + A_\ell)$
\sn
\item[${{}}$]  $\quad (\beta) \quad$ if $\boxplus_{\vK}$ then $(\bar
c^\ell_*,\bar d^\ell_*)$ solves $(\bold x_\ell,\bar\psi,r,A)$, see
\ref{c5}(1)(f). 
\end{enumerate}
\mn
Let $h_\ell$ be the elementary mapping with domain $A_\ell + \bar
c_\ell + \bar d_\ell,h_\ell \rest A_\ell = \id_{A_\ell},h_\ell(\bar
c_\ell \char 94 \bar d_\ell) = \bar c^\ell_* \char 94 \bar d^\ell_*$,
it is well defined and elementary by the choice of $\bar d^\ell_*$.

Now $g^*_1 := h_2 \circ g \circ h^{-1}_1$ is an elementary mapping by
$(*)_{1,\ell}(b)$, for $\ell=1,2$.  As $g^*_1$'s domain is $\subseteq M_{\bold
x_1}$ and its range is $\subseteq M_{\bold x_2}$ there is an extension
$g^*_2$ of $g_2$ to an elementary mapping with domain $\supseteq A_1$
but $\subseteq M_{\bold x_1}$ and range $\supseteq A_2$ but $\subseteq
M_{\bold x_2}$. Next extend $g_2$ to the mapping $g_3$ by letting $g_3(\bar
c) = \bar c$, easily also $g_3$ is an elementary mapping.  Let $g_4$
be the mapping with domain Rang$(\bar d_1 \char 94 \bar c_1 \char 94
\bar d^1_* \char 94 \bar c^1_*)$ mapping $\bar d_1,\bar c_1,\bar
d^1_*,\bar c^1_*$ to $\bar d_2,\bar c_2,\bar d^2_*,\bar c^2_*$
respectively, now $g_4$ is an elementary mapping by the assumption on
$r$.  Easily $g_3(\bar a_1)$ witness $g$ maps $\tp(\bar
a_1,\emptyset,M_{\bold x_1}) \in \bold D_{\bold x_1}$ to the member
$\tp(g_3(\bar a_1),\emptyset,M_{[\bold x_1]})$ of $\bold D_{\bold
x_2}$.  Similarly for $g^{-1}_3,\bar a_2,M_{[\bold x_2]},M_{[\bold
x_1]}$.

So we are done.
\end{PROOF}

\begin{discussion}
\label{c37}
1) Now we can 
start to see the relevance of $\tK,\vK$ to the recounting of types; of
course, the following conclusion will be helpful only if we prove the density
of $\tK$ (or of $\vK$).

\noindent
2)  Note that if we below like to use \ref{c32}(1),(1A)
rather than \ref{c32}(2), we lose little using
$\Sigma\{2^{2^\partial}:\partial < \mu_0\}$ instead $2^{< \mu_0}$.
\end{discussion}

\begin{conclusion}
\label{c39}  
1) If $\kappa,\bar\mu,\theta$ are as in Definition \ref{b05}, $\kappa =
\kappa^{< \kappa} = \mu^{+\alpha}$ and\footnote{without assuming it
we have just to replace $\tK_{\kappa,\bar\mu,\theta}/
\vK_{\kappa,\bar\mu,\theta}$ by $\tK^\oplus_{\kappa,\bar\mu,\theta}/
\vK^\otimes_{\kappa,\bar\mu,\theta}$.}  $2^\theta < \kappa$ and
 $M \in \EC_{\kappa,\kappa}(T)$ \then \, the number of 
$\{\tp(\bar d,M)$: for some $\bold x \in \tK_{\kappa,\bar\mu,\theta}$ we have 
$\bar d \trianglelefteq \bar d_{\bold x},M_{\bold x} = M\}$ up to conjugacy is 
$\le 2^{< \mu_0} + |\alpha|^\theta$. 

\noindent
2) Similarly for $\bold x \in \vK_{\kappa,\bar\mu,\theta}$.
\end{conclusion}

\begin{PROOF}{\ref{c39}}
  1) By part (2) recalling \ref{c27}(2).

\noindent
2) By \ref{b35}(1) + \ref{c29}(4) we can restrict ourselves to smooth
$\bold x \in \vK_{\kappa,\bar\mu,\theta}$.  
By \ref{c27}(4A),(4B) we can deal with $\{\tp(\bar d,M_{\bold
m}):\bold m \in \vK^\otimes_{\kappa,\bar\mu,\theta}$ satisfies
$\bar d \trianglelefteq \bar d_{\bold m},M_{\bold m} = M$ and, as said
above, $\bold x_{\bold m}$ is smooth$\}$.

Now if $(\bold x_\ell,\bar\psi,r,\bar{\bold w}) \in 
\vK^\otimes_{\kappa,\bar\mu,\theta}$, see \ref{c23}(4C) are smooth
for $\ell=1,2$ and $\bold x_1,\bold x_2$ are similar as witnessed by $g$
then $g$ maps $\bold D_{\bold x_1}$ onto $\bold D_{\bold x_2}$, see
\ref{c32}(2) hence by the uniqueness of the $(\bold
D,\kappa)$-sequence-homogeneous model of cardinality $\kappa$ there is
an automorphism $f$ of $M$ such that $g \cup f$ is an elementary
mapping.  Hence $\tp(\bar d_{\bold x_1},M),\tp(\bar d_{\bold x_2},M)$
are conjugate.  We are done as: the number of relevant triples
$(\bar\psi,r,\bar{\bold w})$ is $\le 2^{|T|}$ and
the number of $\bold m \in \vK^\otimes_{\kappa,\bar\mu,\theta}$ 
with $M_{\bold m} = M,(\bar\psi_{\bold m},r_{\bold m},\bar{\bold w}_{\bold m})
= (\bar\psi,r,\bold w)$ up to similarly
is $\le 2^{< \mu_0} + |\alpha|^\theta$ if $\cf(\mu_0) > \theta$ and
$2^{\mu_0} + |\alpha|^\theta$ if $\cf(\mu_0) < \theta$.  The
$2^{\mu_0}/2^{< \mu_0}$ comes from the type of $\bar b_{\bold x}$ 
where $\bar b_{\bold x}$ consists of: $\bar b_{\bold x,i}$
listing $B_i$ for $i \in v_{\bold m} \backslash u_{\bold m}$ and 
$(\bar a_{\bold x,i,0} \char 94 \ldots \char 94 \bar a_{\bold x,i,n}
\char 94 \ldots)_{n < \omega}$ (and of course the respective lengths,
etc.); the $|\alpha|^\theta$ is for the choice of the $\langle
\kappa_{\bold x,i}:i \in u_{\bold x}\rangle$.

Now for each $\bold x \in \vK_{\kappa,\bar\mu,\theta}$ 
the set $\{\bar d:\bar d \trianglelefteq \bar
d_{\bold x}\}$ is $\le \theta$ (even allowing $\{\bar d:\bar d$ a
sub-sequence of $\bar d_{\bold x}\}$ gives $2^\theta \le 2^{< \mu_0}$).   
\end{PROOF}

Now we turn to proving sufficient conditions for (some versions of)
indiscernibility, they are naturally related to $\tK$ and $\vK$.
\begin{claim}
\label{c42}  
$\langle \bar c_s \char 94 \bar d_s:s \in I\rangle$ is an
indiscernible sequence over $B$ \when \,:
\mn
\begin{enumerate}
\item[$(a)$]  $I \in K_{p,\sigma}$, see Definition \ref{a73}
\sn
\item[$(b)$]  if $s <_I t$ are $E_I$-equivalent then $\tp(\bar c_s
\char 94 \bar d_s,B_s) \subseteq \tp(\bar c_t \char 94 \bar d_t,B_t)$ where
\begin{enumerate}
\item[$(\alpha)$]  $E_I = \{(s,t):(\exists i < \sigma)(s,t \in P^I_i)\}$ 
\sn
\item[$(\beta)$]  $B_t = \cup\{\bar c_s \char 94 \bar b_s:s <_I t\}
\cup B$
\sn
\item[$(\gamma)$]  $\ell g(\bar d_s),\ell g(\bar c_s)$ for 
$s \in I$ depend just on $s/E_I$.
\end{enumerate}
\item[$(c)$]   $\tp(\bar c_s,B_s)$ does not split over $B$
\sn
\item[$(d)$]  if $s \in P^I_i$ and $t \in P^I_j$ and $s <_I t$ then
$r_{i,j} = \tp(\bar c_s \char 94 \bar d_s \char 94 \bar c_t \char 94 
\bar d_t,\emptyset)$, i.e. depend only on $(i,j)$
\sn
\item[$(e)$]  $\tp(\bar d_t,\bar c_t + \bar c_s + \bar d_s) 
\vdash \tp(\bar d_t,\bar c_t + \bar d_s + \bar c_s + B_s)$ 
when $s <_I t$

\underline{or} just 
\sn
\item[$(e)'$]  if $s <_I t,\eta_1 \in {}^{\omega >}\ell g(\bar
d_t),\nu_1 \in {}^{\omega >}\ell g(\bar c_t),\eta_3 \in {}^{\omega
>}\ell g(\bar d_s),\nu_3 \in {}^{\omega >}\ell g(\bar c_s)$ and
$\varphi = \varphi(\bar x_{\bar d_t,\eta_1},\bar x_{\bar
c_t,\nu_1},\bar x_{\bar d_s,\eta_3},\bar x_{\bar d_s,\nu_3},\bar y)$
\then \, for some $\eta_0,\nu_0,\eta_2,\nu_2$ (depending on
$(s/E_I,t/E_I,\eta_1,\nu_1,\eta_3,\nu_3)$ but not on $(s,t))$
we have
\mn
\begin{enumerate}
\item[$(\alpha)$]  $\eta_0 \in {}^{\ell g(\eta_1)}(\ell g(\bar
d_t))$ and $\nu_0 \in {}^{\ell g(\nu_1)}(\ell g(\bar c_t))$ and
$\eta_2 \in {}^{\ell g(\eta_3)}(\ell g(\bar d_s))$ and $\nu_2 \in
{}^{\ell g(\nu_3)}(\ell g(\bar c_s))$
\sn
\item[$(\beta)$]  if $\bar b \in {}^{\ell
g(\bar y)}(B_s)$ then 
\newline
${\gC} \models ``\varphi[\bar d_{t,\eta_0},
\bar c_{t,\nu_0},\bar d_{s,\eta_2},\bar c_{s,\nu_2},\bar b] 
\equiv \varphi[\bar d_{t,\eta_1},\bar c_{t,\nu_1},
\bar d_{s,\eta_3},\bar c_{s,\nu_3},\bar b]"$
\sn
\item[$(\gamma)$]   $\tp(\bar d_{t,\eta_0},\bar c_t 
+ \bar d_s + \bar c_s) \vdash \{\varphi(\bar x_{\bar d_{t,\eta_0}},
\bar c_{t,\nu_0},\bar d_{s,\eta_2},\bar c_{s,\nu_2},\bar b):
\bar b \in {}^{\ell g(\bar y)}(B_s)$ and ${\gC}
\models \varphi[\bar d_{t,\eta_0},\bar c_{t,\nu_0},\bar d_{s,\eta_2},
\bar c_{s,\nu_2},\bar b]\}$. 
\end{enumerate}
\end{enumerate}
\end{claim}

\begin{PROOF}{\ref{c42}}
Recall $E = \{(s,t):s,t \in P^I_i$ for some $i < \sigma\}$.
We prove by induction on $n$ that
\mn
\begin{enumerate}
\item[$(*)_n$]   if $s_0 < \ldots < s_{n-1}$ and $t_0 < \ldots
< t_{n-1}$ and $\ell < n \Rightarrow s_\ell E t_\ell$ 
\then \, the sequence 
$\bar c_{s_0} \char 94 \bar d_{s_0} \char 94 \ldots
\char 94 \bar c_{s_{n-1}} \char 94 \bar d_{s_{n-1}}$ 
and the sequence $\bar c_{t_0} \char 94 \bar d_{t_0} 
\char 94 \ldots \char 94 \bar c_{t_{n-1}} \char 94 
\bar d_{t_{n-1}}$, realize the same type over $B_{\min\{s_0,t_0\}}$.
\end{enumerate}
\bigskip

\noindent
\underline{The case $n=0$}:  The desired conclusion is trivial.
\bigskip

\noindent
\underline{The case $n=1$}:  By clause (b) of the assumption, i.e. for
any $t_* \in I$ and $i < \sigma$, 
the type $p_t = \tp(\bar c_t \char 94 \bar d_t,B_t)$ is 
increasing with $t \in \{s:s \in P^I_i$ and $t_* \le_I s\}$.
\bigskip

\noindent
\underline{The case $n=m+1,m \ne 0$}:

By clause (b) of the claim assumption, \wilog \, $s_m = t_m$ call it
$t(*)$ and let $s(*) = \min\{s_0,t_0\}$.

Let $f_0 = \id_{B_{s(*)}}$, let $f_1$ be the function with
domain $B_{s(*)} + \bar c_{s_0} \char 94 \bar d_{s_0} + \ldots + \bar
c_{s_{m-1}} \char 94 \bar d_{s_{m-1}}$ such that $f_1 \supseteq f_0$
and $f_1(\bar c_{s_\ell} \char 94 \bar d_{s_\ell}) = \bar c_{t_\ell} \char
94 \bar d_{t_\ell}$ for $\ell < m$, it is an elementary mapping by the
induction hypothesis.  Let $f_2 = f_1 \cup \id_{\bar c_{t(*)}}$, it 
is an elementary mapping as $\tp(\bar c_{t(*)},
B_{t(*)})$ does not split over $B$ by clause (c) of the claim assumption.

Let $f_3$ be an elementary mapping (in ${\gC}$) extending $f_2$ with domain
$\Dom(f_2) + \bar d_{t(*)} = \Dom(f_1) + \bar c_{s_m} 
\char 94 \bar d_{s_m} = B_{s(*)} 
+ \bar c_{s_0} \char 94 \bar d_{s_0} \char 94 \ldots \char
94 \bar c_{s_m} \char 94 \bar d_{s_m}$ and let $\bar d'_{t_m} =
f_3(\bar d_{s_m}) = f_3(\bar d_{t(*)}) = 
f_3(\bar d_{t_m}) = f_3(\bar d_{t(*)})$.  
Let $i,j < \sigma$ be such that $s_{m-1} \in P^I_i,s_m \in P^I_j$.

So 
\mn
\begin{enumerate}
\item[$\boxplus_1$]  $\tp(\bar d'_{t_m},\bar c_{t_m} + \bar
d_{t_{m-1}} + \bar c_{t_{m-1}}) = \tp(\bar d_{t_m},\bar c_{t_m} + \bar
d_{t_{m-1}} + \bar c_{t_{m-1}})$.
\end{enumerate}
\mn
[Why?  By clause (d) of the claim assumption as $s_\ell E t_{\ell}$ for
$\ell=m-1,m$ and $s_{m-1} <_I s_m,t_{m-1} <_I t_m$ we have $\tp(\bar
d_{t_m} \char 94 \bar c_{_m} \char 94 \bar d_{t_{m-1}} \char 94 \bar
c_{t_{m-1}},\emptyset,\gC) = \tp(\bar d_{s_m} \char 94 \bar c_{s_m}
\char 94 \bar d_{s_{m-1}} \char 94 \bar c_{s_{m-1}},\emptyset,\gC)$.  By
recalling the choice of $f_2,f_3$ and $\bar d'_{t_m}$ we have
$\tp(\bar d'_{t_m} \char 94 \bar c_{t_m} \char 94 \bar d_{t_{m-1}}
\char 94 \bar c_{t_{m-1}},\emptyset,\gC) = \tp(\bar d_{s_m} \char 94
\bar c_{s_m} \char 94 \bar d_{s_{m-1}} \char 94 \bar
c_{s_{m-1}},\emptyset,\gC)$.

Together we are done.]

We shall be done proving for $n$ (hence finish the proof):
\bigskip

\noindent
\underline{Case 1}:  Clause (e) of the claim assumption holds.
\mn
\begin{enumerate}
\item[$\boxplus_2$]  $\bar d_{t_m},\bar d'_{t_m}$ realize the same
type over $\bar c_{t_m} +C$ where $C = B_{s(*)} + \bar c_{t_0} \char
94 \bar d_{t_0} + \ldots + \bar c_{t_{m-1}} \char 94 \bar
d_{t_{m-1}}$ that is $\Rang(f_1)$.
\end{enumerate}
\mn
Why $\boxplus_2$ holds? 

Clearly $\tp(\bar d_{t_m},\bar c_{t_m} + \bar d_{t_{m-1}} + \bar c_{t_{m-1}})
\vdash \tp(\bar d_{t_m},\bar c_{t_m} + C)$ and together
with $\boxplus_1$ we are done.
\bigskip

\noindent
\underline{Case 2}:  Clause (e)$'$ of the claim assumption holds.

So assume
\mn
\begin{enumerate}
\item[$\odot_1$]  ${\gC} \models \varphi[\bar d_{t_m,\eta_1},
\bar c_{t_m,\nu_1},\bar d_{t_{m-1},\eta_3},\bar c_{t_{m-1},\nu_3},\bar
b]$ where $\bar b \in {}^{\ell g(\bar y)}C$ and 
$\eta_1 \in {}^{\omega >}\ell g(d_{t_m}),\nu_1 \in {}^{\omega >} \ell g(\bar
c_{t_m}),\eta_3 \in {}^{\omega >} \ell g(\bar d_{t_{m-1}}),\nu_3 \in
{}^{\omega >} \ell g(\bar c_{t_{m-1}})$ all finite.
\end{enumerate}
\mn
By clause (e)$'$ we can find $\eta_0,\nu_0,\eta_2,\nu_2$ as there.

By $\odot_1$ and subclause $(\beta)$ of $(e)'$ we have
\mn
\begin{enumerate}
\item[$\odot_2$]  ${\gC} \models \varphi[\bar d_{t_m,\eta_0},
\bar c_{t_m,\nu_0},\bar d_{t_{m-1},\eta_2},\bar c_{t_{m-1},\nu_2},\bar b]$.
\end{enumerate}
\mn
By $\odot_2$ and subclause $(\gamma)$ of (e)$'$ we have:
\mn
\begin{enumerate}
\item[$\odot_3$]  $\tp(\bar d_{t_m,\eta_0},\bar c_{t_m} + 
\bar d_{t_{m-1}} + \bar c_{t_{m-1}}) \vdash \varphi(\bar x_{\bar d_m,\eta_0},
\bar c_{t_m,\nu_0},\bar d_{t_{m-1},\eta_2},\bar c_{t_{m-1},\nu_2},\bar b)$.
\end{enumerate}
\mn
But $\tp(\bar d'_{t_m,\eta_0},\bar c_{t_m} + 
\bar d_{t_{m-1}} + \bar c_{t_{m-1}}) = \tp(\bar d_{t_m,\eta_0},
\bar c_{t_{m-1}} + \bar d_{t_m} + \bar c_{t_{m-1}})$ by $\boxplus_1$ so 
by $\odot_3$ we have
\mn
\begin{enumerate}
\item[$\odot_4$]  $\tp(\bar d'_{t_m,\eta_0},\bar c_{t_m} + 
\bar d_{t_{m-1}} + \bar c_{t_{m-1}}) \vdash 
\varphi(\bar x_{\bar d_{t(n)},\eta_0},
\bar c_{t_m,\nu_0},\bar d_{t_{m-1},\eta_2},\bar c_{t_{m-1},\nu_2},\bar b]$
\end{enumerate}
hence
\mn
\begin{enumerate}
\item[$\odot_5$]  ${\gC} \models \varphi[\bar d'_{t_m,\eta_0},
\bar c_{t_m,\nu_0},\bar d_{t_{m-1},\eta_2},\bar c_{t_{m-1},\nu_2},\bar b]$.
\end{enumerate}
\mn
We can apply the elementary mapping $f^{-1}_3$ whose range include all
of the elements of $\gC$ appearing in $\odot_5$ hence we get
\mn
\begin{enumerate}
\item[$\odot_6$]  $\gC \models \varphi[d_{s_m,\eta_0},\bar
c_{s_m,\nu_0},\bar d_{s_{m-1},\eta_2},\bar
c_{s_{m-1},\nu_2},f^{-1}_s(b)]$.
\end{enumerate}
\mn
By subclause $(\beta)$ of $(e)'$ and $\odot_6$ we have (recalling from
$(e)'$ depending on $s/E_J$, but not on $(s,t))$
\mn
\begin{enumerate}
\item[$\odot_7$]   ${\gC} \models \varphi[\bar d_{s_m,\eta_1},
\bar c_{s_m,\nu_1},\bar d_{t_{m-1},\eta_2},\bar c_{s_{m-1},\nu_2},\bar b]$.
\end{enumerate}
\mn
As this holds for any such $\varphi$ we have finished proving
$(*)_n$ also in Case 2, so we are done.
\end{PROOF}

\begin{discussion}
\label{c44} 
1) Naturally we can prove finitary versions of \ref{c42} in some
senses.  Below we deal with $\bold k$-indiscernibility; another 
variant deals with $\Delta$-indiscernible.

\noindent
2) See \ref{c49}.
\end{discussion}

\begin{claim}
\label{c47}  
The sequence $\langle \bar c_{s,0} \char 94 \bar d_{s,0}:s \in
I\rangle$ is $\bold k$-indiscernible over $B_0$ \underline{when} the sequences
$\langle(\bar d_{s,\ell},\bar c_{s,\ell}):\ell \le \bold k,s \in
I\rangle,\langle B_\ell:\ell \le \bold k\rangle$ satisfy:
\mn
\begin{enumerate}
\item[$(a)$]   $I \in K_{p,\sigma}$ and $\ell g(\bar c_{s,\ell})$ and $\ell
g(\bar d_{s,\ell})$ depend just on $\ell$ and $\bold i(s) :=$ the
unique $i$ such that $s \in P^I_i$; also $\bar c_{s,\ell}
\trianglelefteq \bar c_{s,\ell +1},\bar d_{s,\ell} \trianglelefteq
\bar d_{s,\ell +1}$ for $\ell < \bold k,s \in I$
\sn 
\item[$(b)$]  $\tp(\bar d_{s,\bold k} \char 94 \bar c_{s,\bold k},B_s)$ 
is $\subseteq \tp(\bar c_{t,\bold k} \char 94 \bar
d_{t,\bold k},B_t)$ when $s <_I t \wedge \bold i(s) = \bold i(t)$ 
(hence this holds for $\ell < \bold k$, too) where $B_t = B_{t,\bold k}$ and 

$B_{t,k} = \cup\{\bar d_{s,k} \char 94
\bar c_{s,k}:s <_I t\} \cup B_k$ for $k \le \bold k$
\sn
\item[$(c)$]   $\tp(\bar c_{s,k},B_{s,k})$ does not split over
$B_k$ and $B_0 \subseteq B_1 \subseteq \ldots \subseteq B_{\bold k}$
\sn
\item[$(d)$]    if $s \in P^I_i,t \in P^I_j$ 
and $s <_I t$ then $r_{i,j} = \tp(\bar c_{s,\bold k} 
\char 94 \bar d_{s,\bold k} \char 94 \bar c_{t,\bold k} \char 94
\bar d_{t,\bold k},\emptyset)$
\sn
\item[$(e)$]   $\tp(\bar d_{t,\ell},\bar c_{t,\ell +1} + 
\bar d_{s,\ell +1} + \bar c_{s,\ell +1}) \vdash 
\tp(\bar d_{t,\ell},\bar c_{t,\ell} +
\bar d_{s,\ell} + \bar c_{s,\ell} + B_{s,\ell})$ when $s <_I t$
\newline
\underline{or just}
\sn
\item[$(e)'$]   if $s <_I t,\ell < \bold k,\eta_1 \in {}^{\omega >} \ell g(\bar
d_{t,\ell}),\nu_1 \in {}^{\omega >} \ell g(\bar c_{t,\ell}),\eta_3 \in
{}^{\omega >} \ell g(\bar d_{s,\ell}),\nu_3 \in {}^{\omega >} 
\ell g(\bar c_{s,\ell})$ are all
finite and 
\newline
$\varphi'(\bar x_{\bar d_{t,\ell}},\bar x_{\bar
c_{t,\ell}},\bar x'_{\bar d_{t,\ell}},x'_{\bar c_{t,\ell}},\bar y) =
\varphi = \varphi(\bar x_{\bar d_{t,\ell},\eta_1},
\bar x_{\bar c_{t,\ell},\nu_1},\bar x'_{\bar d_{s,\ell},\eta_3},
\bar x'_{\bar c_{s,\ell},\nu_3},\bar y) \in \bbL(\tau_T)$,
\then \, we can find $\eta_0 \in {}^{\omega >} 
\ell g(\bar d_{t,\ell}),\nu_0 \in {}^{\omega >} 
\ell g(\bar c_{t,\ell +1}),\eta_2 \in
{}^{\omega >} \ell g(\bar d_{s,\ell +1}),\nu_2 \in {}^{\omega >} 
\ell g(\bar c_{s,\ell +1})$ (depending on
$s/E_I,t/E_I,\ell,\eta_1,\nu_1,\eta_3,\nu_3$ but not on $(s,t))$ such that:
\sn
\begin{enumerate}
\item[$(\alpha)$]   $\ell g(\eta_0) = \ell g(\eta_1)$ and $\ell
g(\nu_0) = \ell g(\eta_1)$ and $\ell g(\eta_2) = \ell g(\eta_3)$ and
$\ell g(\nu_2) = \ell g(\nu_3)$
\sn
\item[$(\beta)$]  if $\bar b \in {}^{\ell g(\bar y)}(B_t)$ then 

${\gC} \models \varphi[\bar d_{t,\ell,\eta_0},\bar c_{t,\ell +1,\nu_0},
\bar d_{s,\ell +1,\eta_2},\bar c_{s,\ell +1,\nu_2},\bar b] 
\equiv \varphi[\bar d_{t,\ell,\eta_1},
\bar c_{t,\ell,\nu_1},\bar d_{s,\ell,\eta_3},\bar c_{s,\ell,\nu_3},\bar b]$
\sn
\item[$(\gamma)$]    $\tp(\bar d_{t,\ell,\eta_0},
\bar c_{t,\ell +1} + \bar d_{s,\ell +1} + \bar c_{s,\ell +1}) 
\vdash \tp_\varphi(\bar d_{t,\ell +1,\eta_0},
(\bar c_{t,\ell +1,\nu_0} + \bar d_{s,\ell +1,\eta_2} + 
\bar c_{s,\ell +1,\nu_2}) \dotplus B_{s,\ell})$.
\end{enumerate}
\end{enumerate}
\end{claim}

\begin{remark}
\label{c45}
1) Note that $(e) \Rightarrow (e)'$.

\noindent
2) In clause $(e)'$ we may use ``$\eta_0 \in 
{}^{\omega >}\ell g(\bar d_{t,\ell})$" 
rather than ``$\eta_0 \in {}^{\omega >}\ell g(\bar d_{t,\ell +1})$".  
We may add $\bar d^*_{s,\ell}$ such that $\bar
   d_{s,\ell} \trianglelefteq \bar d^*_{s,\ell} \trianglelefteq \bar
   d_{s,\ell +1}$ and use $\bar d^*_{s,\ell}$ instead $\bar
   d_{s,\ell}$ in clause (d) and above.
\end{remark}

\begin{PROOF}{\ref{c47}}
We prove by induction on $k < \bold k$ that:
\mn
\begin{enumerate}
\item[$(*)$]   if $\varrho_1,\varrho_2 \in {}^{k+1}I$ are
$<_I$-decreasing, $(\forall \ell \le k)(\exists i <
 \sigma)[\varrho_1(\ell),\varrho_2(\ell) \in P^I_i]$ 
and $s \le_I \varrho_1(k),\varrho_2(k)$ \then \, the sequences
$\bar d_{\varrho_\ell(0),\bold k-k} \char 94
\bar c_{\varrho_\ell(0),\bold k-k} \char 94 \ldots \char 94
\bar d_{\varrho_\ell(k),\bold k-k} \char 94
\bar c_{\varrho_\ell(k),\bold k-k}$ for $\ell =1,2$ 
realize the same type over $B_{s,\bold k-k}$.
\end{enumerate}
\bigskip

\noindent
\underline{Case $k=0$}:  This holds by clause (b) of the claim
as $\ell g(\varrho_1) = 1 = \ell g(\varrho_2)$.
\bigskip

\noindent
\underline{Case $k > 0$}:   For $\iota = 1,2$, let $\rho_\iota =
\langle \varrho_\iota(1+m):m < k \rangle$ and for $i \le \bold k$ let
$\bar c_{\rho_\iota,i} = \bar c_{\rho_\iota(0),i} \char 94 \ldots 
\char 94 \bar c_{\rho_\iota(k-1),i}$
and $\bar d_{\rho_\iota,i} = \bar d_{\rho_\iota(0),i} \char 94 \ldots
\char 94 \bar d_{\rho_\iota(k-1),i}$ so the induction hypothesis applies
and if $k=1$ then $\bar c_{\rho_\iota,i} = \bar c_{\eta_\iota(1),i},\bar
d_{\rho_\iota,i} = \bar d_{\eta_\iota(1)}$ for $i=0$.

Note that $\bar c_{\rho_\iota,i}$ is a subsequence of 
$\bar c_{\rho_\iota,i+1}$ and $\bar d_{\rho_\iota,i}$ is a 
subsequence of $\bar d_{\rho_\iota,i+1}$.

By the case $k=0$, by clause (b) \wilog \, $\varrho_1(0) = \varrho_2(0)$
call it $t$.  So assume
\mn
\begin{enumerate}
\item[$(*)_1$]  $\varphi = \varphi(\bar x_{\bar d_{t,\bold n-k}},
\bar x_{\bar c_{t,\bold n-k}},\bar x'_{\bar d_{\rho_1,\bold n-k}},
\bar x'_{\bar c_{\rho_1,\bold n-k}},\bar z)$ and

$\bar b \in {}^{\iota g(\bar z)}
(B_{s,\bold k-k})$ and ${\gC} \models \varphi[\bar d_{t(1),\bold n-k},
\bar c_{t(1),\bold n-k},\bar d_{\rho_1,\bold n-k},
\bar c_{\rho_1,\bold n-k},\bar b]$.
\end{enumerate}
\mn
We should prove the parallel statement for $\varrho_2$, 
i.e. for $t$ and $\varrho_2$.
\bigskip

\noindent
\underline{Subcase 1}:  Clause (e) of the assumption.

Hence by clause (e) there is a formula $\psi_\varphi =
\psi_\varphi(\bar x_{\bar d_{s,\bold k-k}},
\bar x_{\bar c_{s,\bold k-k+1}},\bar x'_{\bar d_{\rho_1(0),\bold k-k+1}},
\bar x'_{\bar c_{\rho_1(0),\bold k-k+1}})$ such that
\mn
\begin{enumerate}
\item[$(*)_2$]  $(a) \quad {\gC} \models \psi_\varphi
[\bar d_{t,\bold k-k},\bar c_{t,\bold k-k+1},\bar d_{\rho_1(0),\bold
k -k+1},\bar c_{\rho_1(0),\bold k-k+1}]$
\sn
\item[${{}}$]   $(b) \quad \psi_\varphi(\bar x_{\bar d_{t,\bold k-k}},
\bar c_{t,\bold k-k+1},\bar d_{\rho_1(0),\bold k-k+1},\bar
c_{\rho_1(0),\bold k -k+1}) \vdash$

\hskip25pt $\varphi(\bar x_{\bar d_{t,\bold k-k}},
\bar c_{t,\bold k-k},\bar d_{\rho_1,\bold k-k},\bar c_{\rho_1,\bold k-k},
\bar b)$.
\end{enumerate}
\mn
Hence
\mn
\begin{enumerate}
\item[$(*)_3$]   $(a) \quad {\gC} \models \vartheta_\varphi
[\bar c_{t,\bold k-k+1},\bar d_{\rho_1,\bold k -k+1},
\bar c_{\rho_1,\bold k-k+1},\bar b]$ where
\sn
\item[${{}}$]  $(b) \quad \vartheta_\varphi
(\bar x_{\bar c_{t,\bold k-k+1}},\bar x'_{\bar d_{\rho_1,\bold k-k+1}},
\bar x'_{\bar c_{\rho_1,\bold k-k+1}},\bar z) :=$

\hskip25pt $(\forall \bar x_{\bar d_{s,\bold k-k}})[\psi_\varphi
(\bar x_{\bar d_{t,\bold k -k}},\bar x_{\bar c_{t,\bold k-k+1}},
\bar x'_{\bar d_{\rho_1,\bold k-k+1}},\bar x'_{\bar c_{\rho_1,,\bold
k-k+1}})$

\hskip25pt $\rightarrow \varphi(\bar x_{\bar d_{t,\bold k-k}},
\bar x_{\bar c_{t,\bold k-k}},\bar x'_{\bar d_{\rho_1,\bold k-k}},
\bar x'_{\bar c_{\rho_1,\bold k-k}},\bar z)]$.
\end{enumerate}
\mn
Now
\mn
\begin{enumerate}
\item[$(*)_4$]  $\bar d_{\rho_\iota,\bold k-k+1} \char 94
\bar c_{\rho_\iota,\bold k-k+1} \char 94 \bar b$ realize the same type
over $B_{s,\bold k-k+1}$ for $\iota=1,2$.
\end{enumerate}
\mn
[Why?  By the induction hypothesis as $\bar b$ is from 
$B_{s,\bold k -k} \subseteq B_{s,\bold k-k+1}$.]
\mn
\begin{enumerate}
\item[$(*)_5$]  $\bar c_{t,\bold n-k+1} \char 94
\bar d_{\rho_\iota,\bold k-k+1} \char 94 \bar c_{\rho_\iota,\bold k-k+1}
\char 94 \bar b$ realize the same type over $B_{s,\bold k-k+1}$ 
for $\iota=1,2$.
\end{enumerate}
\mn
[Why?  As first, $\tp(\bar c_{t,\bold k-k+1},B_{t,k+1})$ does not split
over $B_{s,k+1}$ by clause (c) of the assumption, second
$\bar d_{\rho_\ell,\bold k-k+1},\bar c_{\rho_\ell,\bold k-k+1},\bar b$
are included in $B_{t,k+1}$ and third $(*)_4$.]
\mn
\begin{enumerate}
\item[$(*)_6$]   in $(*)_3(a)$ we can replace $\rho_1$ by $\rho_2$,
i.e. ${\gC} \models \vartheta_\varphi[\bar c_{t,\bold k -k+1},
\bar d_{\rho_2,\bold k-k+1},\bar c_{\rho_2,\bold k-k+1},\bar b]$.
\end{enumerate}
\mn
[Why?  By $(*)_5$ and $(*)_3(a)$.]
\mn
\begin{enumerate}
\item[$(*)_7$]  ${\gC} \models \psi_\varphi[\bar d_{t,\bold k-k},
\bar c_{t,\bold k-k+1},\bar d_{\rho_2(0),\bold k-k+1},
\bar c_{\rho_2(0),\bold k-k+1}]$.
\end{enumerate}
\mn
[Why?  By clause (d) of the hypothesis of the claim and $(*)_2(a)$.]
\mn
\begin{enumerate}
\item[$(*)_8$]   ${\gC} \models \varphi
[\bar d_{t,\bold k-k},\bar c_{t,\bold k-k},\bar d_{\rho_2,\bold k -k},
\bar c_{\rho_2,\bold k -k},\bar b]$.
\end{enumerate}
\mn
[Why?  By $(*)_6 + (*)_7$ and the definition of $\vartheta$ in
$(*)_3(b)$.]

So we are done.
\bigskip

\noindent
\underline{Subcase 2}:  Clause (e)$'$ of the assumption holds.

Similarly as in the proof of \ref{c42} and see the proof of \ref{c56}.
\end{PROOF}

\begin{claim}
\label{c49}  
The conclusions of \ref{c42}, \ref{c47} and \ref{c56} below still hold
(and even $(*)$ from its proof holds) even under the 
following weaker assumptions
\mn
\begin{enumerate}
\item[$(a)$]    we add $I_\iota \subseteq I$ for $\iota = 1,2$ and $t \in I
\Rightarrow t \in I_1 \vee t \in I_2$
\sn
\item[$(b),(c),(d)$]  the same
\sn
\item[$(e),(e)'$]  the same but only for $I_1$ and for $I_2$
(but not for $I_1 \cup I_2$)
\sn
\item[$(f)$]  if $t_1 <_I t_2,\iota \in \{1,2\}$ 
and $t_1 \in I_\iota \backslash I_{3-\iota}$
and $t_2 \in I_{3-\iota} \backslash I_\iota$ and $\eta \in {}^{\bold k}
\sigma$ (for \ref{c47}, \ref{c56} 
$\bold k$ is given, otherwise any $\bold k < \omega$)
\then \, we can find $s_0 < \ldots < s_{\bold k-1}$ from $(t_1,t_2)_{I_1
\cap I_2}$ such that $s_\iota \in P^I_{\eta(\iota)}$ for $\ell < \bold k$.
\end{enumerate}
\end{claim}

\begin{remark}
\label{c51}
1) The case which suffice in \ref{c70} below is 
$I_1 = [0,\omega + \omega),I_2 = [0,\omega + \omega
+ 1) \backslash \{\omega\}$ which is somewhat easier.

\noindent
2)  In the natural case, for decreasing $\varrho \in {}^{n+1}I$ we have
$\tp(\bar d_{\varrho(0)} \char 94 \ldots \char 94 \bar
d_{\varrho(n-1)},\bar c_{\varrho(0)} \char 94 \ldots \char 94 
\bar c_{\varrho(n)} + \bar d_{\varrho(n)}) \vdash
\tp(\bar d_{\varrho(0)} \char 94 \ldots \char 94 
\bar d_{\varrho(n-1)-1},\bar c_{\varrho(0)} \char 94 \ldots \char 94 
\bar c_{\varrho(n)} + d_{\varrho(n)} + B_{\varrho(n)})$ and it is
quite natural to use this.

\noindent
3) A variant is: e.g.  $\langle \bar c_s \char 94 \bar d_s:s \in
   I\rangle$ is an indiscernible sequence over $B$ when we assume
$(a) + (b)$ of \ref{c42} and (a),(f) of \ref{c49} and
\mn
\begin{enumerate}
\item[$(g)$]  if $s \in I$ and $\iota \in \{1,2\}$ then $\langle \bar
  c_t \char 94 \bar d_t:t \in I_1$ and $t \ge s\rangle$ is an
  indiscernible sequence over $B$.
\end{enumerate}
\end{remark}

\begin{PROOF}{\ref{c49}}
It is enough to prove this when $I_1 \backslash I_2,I_2
\backslash I_1$ is finite (and has $\le \bold k$ members for
\ref{c47}, \ref{c56}) by induction on $|I_1 \backslash
I_2| + |I_2 \backslash I_1|$ (probably losing appropriately in $\ell$
for \ref{c47}).  So
\wilog \, this number is 2.  
By symmetry \wilog \,
$I_1 \backslash I_2 = \{t_1\},I_2 \backslash I_1 = \{t_2\}$ and
$t_1 <_I t_2$.  The rest should be clear by the transitivity of the equality of
types.
I.e. for notational simplicity concerning \ref{c42}, by it we know
\mn
\begin{enumerate}
\item[$\boxplus$]  if $\iota \in \{1,2\},t \in I$, then $\langle \bar
c_s \char 94 \bar d_s:s \in (I_\iota)_{\ge t}\rangle$ is an
indiscernible sequence over $\cup\{\bar c_s \char 94 \bar d_s:s \in
I_{<t}\} \cup B$.
\end{enumerate}
\mn
It suffices to prove
\mn
\begin{enumerate}
\item[$\oplus$]  if $s_0 <_I \ldots <_I s_{n-1}$ \then \, for some
$r_0 <_{I_1} \ldots <_{I_1} r_{n-1}$, (so all from $I_1$) the sequence
$\bar c_{s_0} \char 94 \bar d_{s_0} \char 94 \ldots \char 94 \bar
c_{s_{n-1}} \char 94 \bar d_{s_{n-1}}$ realizes the same type over $B$
as 
\newline
$\bar c_{r_0} \char 94 \bar d_{r_0} \char 94 \ldots \char 94 \bar c_{r_{n-1}}
\char 94 \bar d_{r_{n-1}}$.
\end{enumerate}
\mn
Why $\oplus$ holds?  
Now if $t_2 \notin \{s_0,\dotsc,s_{n-1}\}$ this is obvious, so assume
$t_2 = s_{k(2)}$, and let $k(1)$ be minimal such that $t_1 <_I
s_{k(1)}$, so $k(1) \le k(2)$; we can even demand $t_1 = s_{k(1)-1}$,
but not used).  By clause (f) there are $r_{k(1)} <_I
\ldots <_I r_{k(2)}$ from $(t_1,t_2)_{I_1 \cap I_2}$ such that $k \in
[k(1),k(1)] \wedge i < \sigma \Rightarrow r_k \in P^I_i \Leftrightarrow
s_k \in P^I_i$ and let $r_k = s_k$ for $k < n$ such that $k \notin
[k(1),k(2)]$. 

So applying $\boxplus$ for $\iota=2,t=t_{k(1)}$ we know that $\bar
c_{r_{k(1)}} \char 94 \bar d_{r_{k(1)}} \char 94 \ldots \char 94 \bar
c_{r_{n-1}} \char 94 \bar d_{r_{n-1}}$ realizes over
$B_{\min\{s_{k(1)},r_{k(1)}\}}$ the same type 
as $\bar c_{s_{k(1)}} \char 94 \bar d_{s_{k(1)}} \char 94 
\ldots \char 94 \bar c_{s_{n-1}} \char 94 \bar d_{r_{n-1}}$.  
As $\bar c_{s_k} \char 94 \bar d_{s_k} = \bar c_{r_k}
\char 94 \bar d_{r_k}$ is from $B_{\min\{s_{k(1)},r_{k(1)}\}}$ for $k<
k(1)$ and $\{r_k:k<n\} \subseteq I_1$ we are done proving $\oplus$
hence the claim.
\end{PROOF}

\noindent
The following may be used in \ref{c47}, \ref{c56}.
\begin{definition}
\label{c53}
Assume $\bold b = \langle (\bar c_s,\bar d_s\rangle:s \in
   I)\rangle$ where $I \in K_{p,\sigma}$ but below we omit $\bold b$
   if clear from the context, and if we have $\langle(\bar
   c_{s,\ell},\bar d_{s,\ell}):s \in I\rangle$ for $\ell \le \bold k$
   we shall write $\ell$ instead of $\bold k$.

\noindent
1) For $k < \omega$ and $\varrho \in {}^k I$ let

\[
\bar d_{\varrho,\bold b} = \bar d_{\varrho(0),\bold b} \char 94 \ldots
\char 94 \bar d_{\varrho(k-1),\bold b}
\]

\[
\bar c_{\varrho,\bold b} = \bar c_{\varrho(0),\bold b} \char 94 \ldots
\char 94 \bar c_{\varrho(k-1),\bold b}.
\]
\mn
2) The sequences 
$\eta_1,\eta_2 \in {}^k I$ are called similar when they realize
the same quantifier-free types in $I$.
\end{definition}

\noindent
The following generalizes \ref{c47}: using only formulas for some
$\Delta$'s following the quantifier-free types in $I$ and using 
a parallel of $\vK$ rather than of $\tK$.
\begin{claim}
\label{c56}  
The sequence $\langle \bar d_{s,0} \char 94 \bar c_{s,0}:s \in
I\rangle$ is $(\Delta_{\bold k},\bold k)$-indiscernible over $B_0$ 
\when \, the sequence
$\langle(\bar d_{s,\ell},\bar c_{s,\ell}):\ell \le \bold k,s \in I\rangle$
satisfies 
\mn
\begin{enumerate}
\item[$(a)$]  $(\alpha) \quad I \in K_{p,\sigma}$ and\footnote{if we
  assume (e) then \wilog \, $\bar d^*_{s,\ell} = \bar d_{s,\ell}$}
 $\ell g(\bar d_{s,\ell}),\ell g(\bar c^*_{s,\ell}),
\ell g(\bar d^*_{s,\ell})$ depend just on $\ell$ and $\bold i(s) :=$ the
unique $i < \sigma$ such that $s \in P^I_i$; also $\bar c_{s,\ell}
\trianglelefteq \bar c_{s,\ell +1},\bar d_{s,\ell} \trianglelefteq
\bar d_{s,\ell +1}$ for $\ell < \bold k,s \in I$
and $B_0 \subseteq B_1 \subseteq \ldots \subseteq B_n$; and $\bar
  c_{s,\ell},\bar d_{s,\ell}$ are finite\footnote{this helps in
  phrasing the demands on the $\Delta_\ell$'s}
\sn
\item[${{}}$]  $(\beta) \quad$ for $i < \sigma$ and $k \le \bold k$
let $w_{i,k} = \ell g(\bar d_{s,k}),v_{i,k} = \ell g(\bar c_{s,k})$
for any 

\hskip25pt $s \in P^I_i$
\sn
\item[${{}}$]   $(\gamma) \quad$ for $k \le \bold k$ let $R_k = 
\{\tp_{\qf}(\varrho,\emptyset,I):\varrho \in {}^{\bold k-k}I$ 
is $<_I$-decreasing$\}$ and let 

\hskip25pt $\bold i(r,\ell) = \bold i(\eta(\ell)) \Leftrightarrow \varrho(\ell)
\in P^I_{\bold i(t,\ell)},\ell g(r) = \ell g(\varrho)$ so $=k$ when 

\hskip5pt  $r = \tp_{\qf}(\varrho,\emptyset,I) \in R_k$
\sn
\item[${{}}$]   $(\delta) \quad$ for $\ell \le \bold k,\Delta_\ell$
  and also $\Delta^*_\ell$ is the closure of a finite set of formulas
  (each with finite set of variables) under permuting the variables,
  negation and adding dummy variables
\sn
\item[${{}}$]   $(\varepsilon) \quad \Delta_i \subseteq \Delta^*_\ell
  \subseteq \Delta_{\ell +1} \subseteq \Delta_{\ell +1}$
\sn
\item[$(b)$]  $\tp_{\Delta_{\bold k}}(\bar d_{s,\bold k} 
\char 94 \bar c_{s,\bold k},B_s)$ 
is $\subseteq \tp_{\Delta_k}(\bar c_{t,\bold k} \char 94 
\bar d_{t,\bold k},B_t)$ when $s <_I t \wedge \bold i(s) = \bold i(t)$ 
(hence this holds for $k \le \bold k$, too) where $B_t = B_{t,\bold k}$ and 
$B_{t,k} = \cup\{\bar d_{s,k} \char 94
\bar c_{s,k}:s <_I t\} \cup B_k$ for $k \le \bold k$
\sn
\item[$(c)$]   $\tp_{\Delta^*_\ell}(\bar c_{s,\ell},B_{s,\ell +1})$ does not 
$\Delta_{\ell +1}$-split over $B_s$ 
\sn
\item[$(d)$]    if $s \in P^I_i,t \in P^I_j$ 
and $s <_I t$ then $r_{i,j} = \tp_{\Delta_k}(\bar c_{s,\bold k} 
\char 94 \bar d_{s,\bold k} \char 94 \bar c_{t,\bold k} \char 94
\bar d_{t,\bold k},\emptyset)$
\sn
\item[$(e)$]   $\tp(\bar d_{t,\ell},\bar c_{t,\ell +1} + 
\bar d_{s,\ell} + \bar c_{s,\ell +1}) \vdash \tp
(\bar d_{t,\ell},\bar c_{t,\ell} +
\bar d_{s,\ell} + \bar c_{s,\ell},B_\ell)$ when $s <_I t$
\newline
\underline{or just}
\sn
\item[$(e)'$]   if $s <_I t,\ell < \bold k,\eta_1 \in {}^{\omega >} \ell g(\bar
d_{t,\ell}),\nu_1 \in {}^{\omega >} \ell g(\bar c_{t,\ell}),\eta_3 \in
{}^{\omega >} \ell g(\bar d_{s,\ell}),\nu_3 \in {}^{\omega >} 
\ell g(\bar c_{s,\ell})$ are all
finite and $\varphi' = \varphi'(\bar x_{\bar d_{t,\ell},\eta_1},
\bar x_{\bar c_{t,\ell},\nu_1},\bar x'_{\bar d_{s,\ell},\eta_3},
\bar x'_{\bar c_{s,\ell},\nu_3},\bar y) \in \bbL(\tau_T)$,
\then \, we can find $\eta_0 \in {}^{\omega >} \ell g(\bar d_{t,\ell
+1}),\nu_0 \in {}^{\omega >} \ell g(\bar c_{t,\ell +1}),\eta_2 \in
{}^{\omega >} \ell g(\bar d_{s,\ell +1}),\nu_2 \in {}^{\omega >} 
\ell g(\bar c_{s,\ell +1})$, (depending on
$\ell,s/E_1,t/E_I,\eta_1,\nu_1,\eta_3,\nu_3$ and $\varphi$ but not on
$(s,t))$ such that:
\sn
\begin{enumerate}
\item[$(\alpha)$]   $\ell g(\eta_0) = \ell g(\eta_1)$ and $\ell
g(\nu_0) = \ell g(\eta_1)$ and $\ell g(\eta_2) = \ell g(\eta_3)$ and
$\ell g(\nu_2) = \ell g(\nu_3)$
\sn
\item[$(\beta)$]  if $\bar b \in {}^{\ell g(\bar y)}(B_t)$ then 

${\gC} \models \varphi'[\bar d_{t,\ell +1,\eta_0},
\bar c_{t,\ell +1,\nu_0},\bar d_{s,\ell,\eta_2},\bar
c_{s,\ell,\nu_2},\bar b] \equiv \varphi'[\bar d_{t,\ell,\eta_1},
\bar c_{t,\ell,\nu_1},\bar d_{s,\ell,\eta_3},\bar c_{s,\ell,\nu_3},\bar b]$
\sn
\item[$(\gamma)$]    $\tp(\bar d_{t,\ell +1,\eta_0},
\bar c_{t,\ell +1} + \bar d_{s,\ell +1},\bar c_{s,\ell +1}) 
\vdash \tp_\varphi(\bar d_{t,\ell +1,\eta_0},
(\bar c_{t,\ell +1,\nu_0} + \bar d_{s,\ell +1,\eta_2} + 
\bar c_{s,\ell +1,\nu_2}) + B_s)$.
\end{enumerate}
\end{enumerate}
\mn
But by the assumptions on the $\Delta_\ell$'s, $(e)'$ is equivalent to:
\mn
\begin{enumerate}
\item[$(e)''$]  if $s <_I t,\ell < \bold k$ and $\varphi = \varphi_0 =
  \varphi_0(\bar x_{\bar d_{t,\ell}},\bar x_{c_{t,\ell}},\bar x_{\bar
  d_{s,\ell}},\bar x_{\bar c_{s,\ell}},\bar z) \in \Delta_\ell$ \then
  \, we can find $\varphi_1 = \varphi_1(\bar x_{\bar
d^*_{t,\ell}},\bar x_{\bar c^*_{t,\ell}},\bar x'_{\bar d_{s,\ell +1}},
\bar x'_{\bar c_{s,\ell +1}},\bar z)$ and $\psi = \psi_\varphi =
\psi_\varphi(\bar x_{\bar d^*_{t,\ell}},\bar x_{\bar c_{t,\ell +1}},
\bar x'_{\bar d_{s,\ell +1}},\bar x'_{\bar c_{s,\ell +1}})$ depending
  only if $\bold i(s),\bold i(t),\ell$ and $\varphi$ such that:
\sn
\begin{enumerate}
\item[$(a)$]  $\gC \models \psi[\bar d^*_{t,\ell},\bar c_{t,\ell
    +1},\bar d_{s,\ell +1},\bar c_{s,\ell +1}]$
\sn
\item[$(b)$]   for every $\bar b \in {}^{\ell g(\bar z)}(B_{s,\ell})$
  we have $\gC \models \varphi_0[\bar d_{t,\ell},\bar c_{t,\ell},\bar
  d_{s,\ell},\bar c_{s,\ell},\bar b]$ iff $\gC \models \varphi_1[\bar
  d^*_{t,\ell},\bar c^*_{t,\ell},\bar d_{s,\ell +1},\bar c_{s,\ell
  +1}]$
\sn 
\item[$(c)$]  $\psi(\bar x_{\bar d^*_{t,\ell}},\bar c_{t,\ell +1},
\bar d_{s,\ell +1},\bar c_{s,\ell +1}) \vdash 
\{\varphi_1(\bar x_{\bar d^*_{t,\ell}},
\bar c^*_{t,\ell},\bar d_{s,\ell+1},\bar c_{s,\ell +1},
\bar b):\bar b \in {}^{\ell g(\bar z)}(B_{s,\ell})$ and $\gC
  \models \varphi_1[\bar d^*_{t,\ell},\bar c^*_{t,\ell},\bar d_{s,\ell
  +1},\bar e_{s,\ell +1};b]\}$
\sn 
\item[$(d)$]  $\vartheta_\varphi = \vartheta_\varphi(\bar x_{\bar
  c_{t,\ell +1}},\bar x'_{\bar d_{s,\ell +1}},\bar x'_{\bar c_{s,\ell
  +1}},\bar z) \in \Delta^*_\ell$ where $\vartheta_\varphi$ is
\newline
$(\forall \bar x_{\bar d^*_{t,\ell}})[\psi_\varphi(\bar x_{\bar
  d^*_{t,\ell}},\bar x_{\bar c_{t,\ell+1}},\bar x_{\bar
  d_{s,\ell+1}},\bar x_{\bar c_{s,\ell+1}}) \rightarrow \varphi_i(\bar
  x_{\bar d^*_{t,\ell}},\bar x_{\bar c^*_{t,\ell}},\bar x_{\bar
  d_{s,\ell+1}},\bar x_{\bar c_{s,\ell+1}},\bar z)]$.
\end{enumerate}
\end{enumerate}
\end{claim}

\begin{remark}
\label{c58}
Used in \ref{c70}, \ref{c72} below (for the case we use $(c)',(c)$;
respectively).
\end{remark}

\begin{PROOF}{\ref{c58}}
We prove by induction on $k < \bold k$ that:
\mn
\begin{enumerate}
\item[$(*)$]   if $\varrho_1,\varrho_2 \in {}^{k+1}I$ are
$<_I$-decreasing, $(\forall \ell \le k)(\exists i <
 \sigma)[\varrho_1(\ell),\varrho_2(\ell) \in P^I_i]$ 
and $s \le_I \varrho_1(k),\varrho_2(k)$ \then \, the sequences
$\bar d_{\varrho_\ell(0),\bold k-k} \char 94
\bar c_{\varrho_\ell(0),\bold k-k} \char 94 \ldots \char 94
\bar d_{\varrho_\ell(k),\bold k-k} \char 94
\bar c_{\varrho_\ell(k),\bold k-k}$ for $\ell =1,2$ 
realize the same $\Delta_\ell$-type over $B_{s,\bold k-k}$.
\end{enumerate}
\bigskip

\noindent
\underline{Case $k=0$}:  This holds by clause (b) of the claim
as $\ell g(\varrho_1) = 1 = \ell g(\varrho_2)$.
\bigskip

\noindent
\underline{Case $k > 0$}:   For $\iota = 1,2$, let $\rho_\iota =
\langle \varrho_\iota(1+m):m < k \rangle$ and for $i \le \bold k$ let
$\bar c_{\rho_\iota,i} = \bar c_{\rho_\iota(0),i} \char 94 \ldots 
\char 94 \bar c_{\rho_\iota(k-1),i}$
and $\bar d_{\rho_\iota,i} = \bar d_{\rho_\iota(0),i} \char 94 \ldots
\char 94 \bar d_{\rho_\iota(k-1),i}$ so the induction hypothesis applies
and if $k=1$ then $\bar c_{\rho_\iota,i} = \bar c_{\eta_\iota(1),i},\bar
d_{\rho_\iota,i} = \bar d_{\eta_\iota(1)}$ for $i=0$.

Note that $\bar c_{\rho_\iota,i}$ is a subsequence of 
$\bar c_{\rho_\iota,i+1}$ and $\bar d_{\rho_\iota,i}$ is a 
subsequence of $\bar d_{\rho_\iota,i+1}$.

By the case $k=0$, i.e. by clause (b) \wilog \, $\varrho_1(0) = \varrho_2(0)$
call it $t$.  So assume
\mn
\begin{enumerate}
\item[$(*)_1$]  $\varphi = \varphi(\bar x_{\bar d_{t,\bold n-k}},
\bar x_{\bar c_{t,\bold n-k}},\bar x'_{\bar d_{\rho_1,\bold n-k}},
\bar x'_{\bar c_{\rho_2,\bold n-k}},\bar z) \in \Delta_\ell$ and

$\bar b \in {}^{\ell g(\bar z)}
(B_{s,\bold k-k})$ and ${\gC} \models \varphi[\bar d_{t,\bold n-k},
\bar c_{t,\bold n-k},\bar d_{\rho_1,\bold n-k},
\bar c_{\rho_1,\bold n-k},\bar b]$.
\end{enumerate}
\mn
We should prove the parallel statement for $\varrho_2$, 
i.e. for $t$ and $\rho_2$; this will suffice.
\bigskip

\noindent
\underline{Subcase 1}:  Clause (e) of the assumption.

Follows by the second subcase, (and has easier proof).
\bigskip

\noindent
\underline{Subcase 2}:  Clause (e)$'$ of the assumption but we shall
use the formulation of (e)$''$.

Hence 
\mn
\begin{enumerate}
\item[$(*)_2$]  choose $\varphi_1,\psi_\varphi$ and
  $\vartheta_\varphi$ as in clause (e)$''$ for $\rho(0),t$ (chosen
  above) and $\varphi$ from $(*)_1$, hence in particular
\sn
\begin{enumerate}
\item[$(a)$]  ${\gC} \models \psi_\varphi[\bar d^*_{t,\bold k-k},
\bar c_{t,\bold k-k+1},\bar d_{\rho_1(0),\bold k -k+1},
\bar c_{\rho_1(0),\bold k-k+1}]$
\sn
\item[$(b)$]  $\psi_\varphi(\bar x_{\bar d^*_{t,\bold k-k}},
\bar c_{t,\bold k-k+1},\bar d_{\rho_1(0),\bold k-k+1},
\bar c_{\rho_1(0),\bold k -k+1}) \vdash \varphi(\bar x_{\bar d_{t,\bold k-k}},
\bar c_{t,\bold k-k},\bar d_{\nu_1,\bold k-k},\bar c_{\nu_1,\bold k-k},
\bar b)$.
\end{enumerate}
\end{enumerate}
\mn
Hence
\mn
\begin{enumerate}
\item[$(*)_3$]   $(a) \quad {\gC} \models \vartheta_\varphi
[\bar c_{t,\bold k-k+1},\bar d_{\rho_1,\bold k -k+1},
\bar c_{\rho_1,\bold k-k+1},\bar b]$ where
\sn
\item[${{}}$]  $(b) \quad \vartheta_\varphi
(\bar x_{\bar c_{t,\bold k-k+1}},\bar x'_{\bar d_{\rho_1,\bold k-k+1}},
\bar x'_{\bar c_{\rho_1,\bold k-k+1}},\bar z) :=$

\hskip25pt $(\forall \bar x_{\bar d^*_{s,\bold k-k}})[\psi_\varphi
(\bar x_{\bar d^*_{t,\bold k -k}},\bar x_{\bar c_{t,\bold k-k+1}},
\bar x'_{\bar d_{\rho_1,\bold k-k+1}},\bar x'_{\bar c_{\rho_1,,\bold
k-k+1}})$

\hskip25pt $\rightarrow \varphi(\bar x_{\bar d_{t,\bold k-k}},
\bar x_{\bar c_{t,\bold k-k}},\bar x'_{\bar d_{\rho_1,\bold k-k}},
\bar x'_{\bar c_{\rho_1,\bold k-k}},\bar z)]$.
\end{enumerate}
\mn
Now
\mn
\begin{enumerate}
\item[$(*)_4$]  $\bar d_{\rho_\iota,\bold k-k+1} \char 94
\bar c_{\rho_\iota,\bold k-k+1} \char 94 \bar b$ realize the same
$\Delta_{\ell +1}$- type over $B_{s,\bold k-k+1}$ for $\iota=1,2$.
\end{enumerate}
\mn
[Why?  By the induction hypothesis as $\bar b$ is from 
$B_{s,\bold k -k} \subseteq B_{s,\bold k-k+1}$.]
\mn
\begin{enumerate}
\item[$(*)_5$]  $\bar c_{t,\bold n-k+1} \char 94
\bar d_{\nu_\iota,\bold k-k+1} \char 94 \bar c_{\nu_\iota,\bold k-k+1}
\char 94 \bar b$ realizes the same $\Delta_\ell$-type over $B_{s,\bold k-k+1}$ 
for $\iota=1,2$.
\end{enumerate}
\mn
[Why?  As first, $\tp_{\Delta^*_\ell}
(\bar c_{t,\bold k-k+1},B_{t,k+1})$ does not $\Delta_{\ell +1}$-split
over $B_{s,k+1}$ by clause (c) of the assumption, second
$\bar d_{\rho_\ell,\bold k-k+1},\bar c_{\rho_\ell,\bold k-k+1},\bar b$
are included in $B_{t,k+1}$ and third $(*)_4$.]
\mn
\begin{enumerate}
\item[$(*)_6$]   in $(*)_3(a)$ we can replace $\rho_1$ by $\rho_2$,
i.e. ${\gC} \models \vartheta_\varphi[\bar c_{t,\bold k -k+1},
\bar d_{\rho_2,\bold k-k+1},\bar c_{\rho_2,\bold k-k+1},\bar b]$.
\end{enumerate}
\mn
[Why?  By $(*)_5$ and $(*)_3(a)$.]
\mn
\begin{enumerate}
\item[$(*)_7$]  ${\gC} \models \psi_\varphi[\bar d_{t,\bold k-k},
\bar c^*_{t,\bold k-k+1},\bar d_{\rho_2(0),\bold k-k+1},
\bar c_{\rho_2(0),\bold k-k+1}]$.
\end{enumerate}
\mn
[Why?  By clause (d) of the hypothesis of the claim and $(*)_2(a)$.]
\mn
\begin{enumerate}
\item[$(*)_8$]   ${\gC} \models \varphi
[\bar d_{t,\bold k-k},\bar c_{t,\bold k-k},\bar d_{\rho_2,\bold k -k},
\bar c_{\rho_2,\bold k -k},\bar b]$.
\end{enumerate}
\mn
[Why?  By $(*)_6 + (*)_7$ and the definition of $\vartheta$ in
$(*)_3(b)$.]

So we are done.
\end{PROOF}
\bigskip

\subsection{Toward Density of $\tK$} \
\bigskip

We first show that the existence of $\le^+_1$-extension for every
$\bold x \in \rK^\oplus_{\kappa,\mu,\theta}$ suffice
for existence (i.e. for density) for $\tK$.  The main case in \ref{c70},
\ref{c72} is $\sigma = \omega$.  Then in \ref{c78} we prove this
sufficient condition for weakly compact $\kappa$.  
Note that for $\rK^\oplus$ closure under
union is not obviously true.

\begin{claim}
\label{c70}  
We have $\bold x_\delta \in \tK_{\kappa,\bar\mu,\theta}$, 
moreover $\bold m_\delta = (\bold x_\delta,\bar\psi_\delta,r_\delta) \in 
\tK^\oplus_{\kappa,\bar\mu,\theta}$
and $\varepsilon < \delta \Rightarrow \bold x_\varepsilon \le_1 
\bold x_\delta$ \when \, ($\delta < \theta^+$ is a limit ordinal and):
\mn
\begin{enumerate}
\item[$\boxplus$]  $(a) \quad \bold m_\varepsilon = 
(\bold x_\varepsilon,\bar\psi_\varepsilon,r_\varepsilon) 
\in \rK^\oplus_{\kappa,\bar\mu,\theta}$ for $\varepsilon <
\delta$ is $\le_1$-increasing
\sn
\item[${{}}$]   $(b) \quad r_\varepsilon$ is a complete type, (over
  the empty set)
\sn
\item[${{}}$]  $(c) \quad \bold m_\varepsilon \le^+_1 
\bold m_{\varepsilon +1}$, see Definition \ref{c23}(4) or just
\sn
\item[${{}}$]  $(c)' \quad$ if
$\varepsilon < \delta$ and $\varphi \in \Gamma^2_{\bold
m_\varepsilon}$ \then \, for some $\zeta \in [\varepsilon,\delta)$ we
have $\varphi \in \Gamma^2_{\bar\psi[\bold m_\zeta]}$
\sn
\item[${{}}$]  $(d) \quad \bold m_\delta = \cup\{\bold
m_\varepsilon:\varepsilon  < \delta\}$, see \ref{c23}(6), i.e.
\begin{enumerate}
\item[${{}}$]  $(\alpha) \quad \bold x_\delta = 
\cup\{\bold x_\varepsilon:\varepsilon < \delta\}$, see \ref{b18}(2)
\sn
\item[${{}}$]   $(\beta) \quad \bar\psi_\delta$ is the
limit\footnote{But note that for $\varepsilon < \zeta$ 
the formulas in $\bar\psi_\zeta$ has more dummy variables than those
in $\bar\psi_\varepsilon$}
 of $\langle \bar\psi_\varepsilon:\varepsilon < \delta\rangle$
\sn
\item[${{}}$]  $(\gamma) \quad r_\delta = 
\cup\{r_\varepsilon:\varepsilon  < \delta\}$.
\end{enumerate}
\end{enumerate}
\end{claim}

\begin{claim}
\label{c72}  
We have $\bold x_\delta \in \vK_{\kappa,\bar\mu,\theta}$ 
moreover $\bold m_\delta = (\bold  x_\delta,
\bar\psi_\delta,r_\delta) \in 
\vK^\oplus_{\kappa,\bar\mu,\theta}$ and $\varepsilon < \delta 
\Rightarrow \bold x_\varepsilon \le_1 \bold x_\delta$ 
\when \, as in \ref{c70} except that we replace $(c),(c)'$ by
\mn
\begin{enumerate}
\item[$(c)_1$]  $\bold m_\varepsilon \le^\odot_1 \bold
m_{\varepsilon +1}$, see \ref{c23}(4A)

or at least
\item[$(c)'_1$]  for every $\varepsilon < \delta$ and $\varphi \in
\Gamma^2_{\bold m_\varepsilon}$ for some $\zeta \in
(\varepsilon,\delta)$ we have $\bold m_\varepsilon
\le^\odot_{1,\varphi} \bold m_\zeta$, see Definition \ref{c23}(4B).
\end{enumerate}
\end{claim}

\begin{remark}
\label{c73}  
1) We may weaken clause (b), i.e. $r_\varepsilon,r_\delta$ are 
not necessarily complete, still need sufficient condition for
indiscernibility in proving $\boxplus_3$ below.

\noindent
2) In \ref{c72} we can use $\vk^\otimes_{\kappa,\bar\mu,\theta}$.
\end{remark}

\begin{PROOF}{\ref{c70}}
\underline{Proof of \ref{c70}, \ref{c72}}  
For simplicity we assume
(c),(c)$_1$ in \ref{c70}, \ref{c72}, respectively, otherwise we
have to use \ref{c56} (or use compactness).
Let $\bar d_\varepsilon = \bar d_{\bold x_\varepsilon},\bar
d_\delta = \bar d_{\bold x_\delta},\bar c_\varepsilon = 
\bar c_{\bold x_\varepsilon}$ and $\bar c_\delta = 
\bar c_{\bold x_\delta}$ for $\varepsilon < \delta$.  The main point
is proving clause (f) from Definition \ref{c5}(1).

Let $A \subseteq M_{\bold x}$ be of cardinality $< \kappa$ and \wilog
\, $\varepsilon < \delta \Rightarrow B^+_{\bold x_\varepsilon} 
\subseteq A$.  We now
choose $A_\alpha,\bar d_{\alpha,\delta},\bar c_{\alpha,\delta},
\bar d_{\alpha,\varepsilon},\bar
c_{\alpha,\varepsilon}$ (for $\varepsilon < \delta$) by induction on
$\alpha < \theta^+$, really $\alpha < \delta + \delta$ suffice, such that:
\mn
\begin{enumerate}
\item[$\oplus$]  $(a) \quad \bar d_{\alpha,\delta},\bar
d_{\alpha,\varepsilon},\bar c_{\alpha,\varepsilon}$ are sequences
from $M_{\bold x}$
\sn
\item[${{}}$]  $(b) \quad \ell g(\bar d_{\alpha,\varepsilon}) = 
\ell g(\bar d_\varepsilon)$ and $\varepsilon < \zeta <
\delta \Rightarrow \bar d_{\alpha,\varepsilon} = \bar d_{\alpha,\zeta}
\rest \ell g(\bar d_\varepsilon)$ 
\sn
\item[${{}}$]  $(c) \quad \ell g(\bar c_{\alpha,\varepsilon}) 
= \ell g(\bar c_\varepsilon)$ and $\varepsilon < \zeta <
\delta \Rightarrow \bar c_{\alpha,\varepsilon} = \bar c_{\alpha,\zeta}
\rest \ell g(\bar c_\varepsilon)$
\sn
\item[${{}}$]  $(d) \quad \bar d_{\alpha,\delta} = \cup\{\bar
d_{\alpha,\varepsilon}:\varepsilon < \delta\}$ and 
$\bar c_{\alpha,\delta} = \cup\{\bar
c_{\alpha,\varepsilon}:\varepsilon < \delta\}$; actually follows
\sn
\item[${{}}$]  $(e) \quad$ if $\varepsilon \le \delta$ then 
$\bar c_{\alpha,\varepsilon} \char 94 \bar d_{\alpha,\varepsilon}$ and 
$\bar c_\varepsilon \char 94 \bar d_\varepsilon$

\hskip25pt  realize
the same type over $A_\beta := A + \Sigma\{\bar c_{\beta,\delta} 
\char 94 \bar d_{\beta,\delta}:\beta < \alpha\}$
\sn
\item[${{}}$]  $(f) \quad$ if $\alpha = \varepsilon \mod \delta$
\then \, 
the sequence $\bar c_\varepsilon \char 94 \bar d_\varepsilon 
\char 94 \bar c_{\alpha,\varepsilon}
\char 94 \bar d_{\alpha,\varepsilon}$ realizes $r_\varepsilon$
\sn
\item[${{}}$]   $(g) \quad$ if $\alpha = \varepsilon \mod \delta$ and
$\varphi = \varphi(\bar x_{\bar d_\varepsilon},\bar x_{\bar
c_\varepsilon},\bar x'_{\bar d_{\varepsilon +1}},\bar x'_{\bar
c_{\varepsilon +1}},\bar y) \in 
\Gamma^2_{\bar\psi[\bold m_\varepsilon]}$ \then \,

\hskip25pt  $\psi_\varphi(\bar x_{\bar d_\varepsilon},\bar
c_\varepsilon,\bar d_{\alpha,\varepsilon},\bar c_{\alpha,\varepsilon}) 
\vdash \tp_\varphi(\bar d_\varepsilon,(\bar c_\varepsilon \char
94 \bar d_{\alpha,\varepsilon} \char 94 \bar d_{\alpha,\varepsilon}) 
\dotplus A_\alpha)$.
\end{enumerate}
\mn
This is possible by the assumptions recalling the definitions, that
is, if $\varepsilon < \delta + \delta$ and $\alpha = \varepsilon \mod \delta$
then first we choose $\bar d_{\alpha,\varepsilon},\bar
c_{\alpha,\varepsilon}$ as required in clauses (f),(g), second we
choose $(\bar d_{\alpha,\delta},\bar c_{\alpha,\delta})$ such that
clauses (a),(e) holds and $\bar d_{\alpha,\varepsilon} = \bar
d_{\alpha,\delta} \rest \ell g(\bar d_\varepsilon),\bar
c_{\alpha,\varepsilon} = \bar c_{\alpha,\delta} \rest \ell g(\bar
c_\varepsilon)$ and third define $\bar d_{\alpha,\zeta},\bar
c_{\alpha,\zeta}$ as $\bar d_{\alpha,\delta} \rest \ell g(\bar
d_\zeta),\bar c_{\alpha,\delta} \rest \ell g(\bar c_\zeta)$
for $\zeta < \delta$ so clauses (b),(c) hold.

So let $u_\varepsilon = (\varepsilon,\delta) \cup (\delta +
\varepsilon,\delta + \delta)$; now
\mn
\begin{enumerate}
\item[$\boxplus_1$]   if $\zeta < \delta$ and 
$\alpha \in u_{\zeta +1}$, then
\begin{enumerate}
\item[$(a)$]  $\bar c_\zeta \char 94 \bar d_\zeta \char 94 
\bar c_{\alpha,\zeta} \char 94 \bar d_{\alpha,\zeta}$ realizes $r_\zeta$
\sn
\item[$(b)$]  $\bar c_\zeta \char 94 \bar d_\zeta$ and
$\bar c_{\alpha,\zeta} \char 94 \bar d_{\alpha,\zeta}$ realize the same type
over $A_\alpha$
\sn
\item[$(c)$]   \underline{for \ref{c70}}
$\tp(\bar d_{\alpha +1,\zeta},\bar c_{\alpha +1,\zeta +1} + \bar
d_{\alpha,\zeta +1} + \bar c_{\alpha,\zeta +1}) 
\vdash \tp(\bar d_{\alpha +1,\zeta},\bar c_{\alpha +1,\zeta +1} 
+ \bar d_{\alpha,\zeta +1} + \bar c_{\alpha,\zeta +1} + A_\alpha)$.

\underline{for \ref{c72}}: if $\varphi = \varphi(\bar x_{\bar d},\bar
x_{\bar c},\bar x'_{\bar d},\bar x'_{\bar c},\bar y) \in
\Gamma^2_{\bar\psi_\varepsilon}$ then $\tp(\bar d_{\alpha +1,\zeta},\bar
c_{\alpha +1,\zeta +1} + \bar d_{\alpha,\zeta +1} + \bar
c_{\alpha,\zeta +1}) \vdash \tp_\varphi(\bar d_{\alpha +1,\zeta},
(\bar c_{\alpha +1,\zeta +1} \char 94 \bar d_{\alpha,\zeta +1} 
\char 94 \bar c_{\alpha,\zeta +1}) \dotplus A_\alpha)$.
\end{enumerate}
\end{enumerate}
\mn
[Why?  Let $\alpha = \varepsilon \mod \delta$ so $\zeta < \varepsilon$; 
clause (b) holds by clause (e) of $\oplus$, for clause (a) uses clause
(f) of $\oplus$ noting that $\bold m_\zeta \le_1 \bold m_\varepsilon$ 
hence $r_\zeta \subseteq r_\varepsilon$ (so $\zeta \le \varepsilon$ 
suffices for (a),(b)).
For clause (c), first assume clause (c) \underline{of \ref{c70}}.
Note that $\bold m_\zeta \le^+_1
\bold m_{\zeta +1}$ hence (by Definition \ref{c23}(4)) we
have $\Gamma^2_{\bold x[\bold m_\varepsilon]} \subseteq 
\Gamma^2_{\bar\psi[\bold m_\zeta]}$.  Second, assume clause (c)
\underline{of \ref{c72}}: 
similarly using $\bold m_\varepsilon \le^\odot_1 \bold m_\zeta$.]

\noindent
Let $D$ be an ultrafilter on $\delta$ to which every co-bounded subset
of $\delta$ belongs and let
\mn
\begin{enumerate}
\item[$\boxplus_2$]   $q = q(\bar x_{\bar d_\delta},
\bar x_{\bar c_\delta}) := \{\vartheta(\bar x_{\bar d_\delta},
\bar x_{\bar c_\delta},\bar b):\bar b \in {}^{\delta >}(A_{\delta + \delta})$
and for some ${\cU} \in D$ we have $\alpha \in {\cU}
\Rightarrow M_{\bold x} \models \vartheta[\bar d_{\alpha,\delta},
\bar c_{\alpha,\delta},\bar b]\}$ so $q(\bar x_{\bar d_\delta},
\bar x_{\bar c_\delta})$ is a complete
type over $A_{\delta + \delta} \subseteq M_{\bold x}$.
\end{enumerate}
Let $\bar d'_\delta \char 94 \bar c'_\delta$ be a sequence 
from $M_{\bold x}$ realizing $q(\bar x_{\bar d_\delta},
\bar x_{\bar c_\delta})$.

Let $\bar d'_\varepsilon,\bar c'_\varepsilon$ be such that $\bar
d'_\varepsilon \triangleleft \bar d'_\delta,\ell g(\bar d'_\varepsilon) = 
\ell g(\bar d_\varepsilon)$ and $\bar c'_\varepsilon \triangleleft
\bar c'_\delta,\ell g(\bar c'_\varepsilon) = 
\ell g(\bar c_\varepsilon)$.
\mn
\begin{enumerate}
\item[$\boxplus_3$]  if $\varepsilon < \delta$ and $\varepsilon = n \mod
\omega$ and $\gamma \in u_\varepsilon$ then $\langle \bar
d_{\alpha,\varepsilon} \char 94 
\bar c_{\alpha,\varepsilon}:\alpha \in u_{\varepsilon + 2n} 
\backslash \gamma\rangle \char 94 
\langle \bar d_\varepsilon \char 94 \bar c_\varepsilon \rangle$ is
an $n$-indiscernible sequence over $A_\gamma$.
\end{enumerate}
\mn
[Why?  For claim \ref{c70}, by claim \ref{c47} the version with clause
(e), for claim \ref{c72} by claim \ref{c47} the version 
with clause (e)$'$.]
\mn
\begin{enumerate}
\item[$\boxplus_4$]   if $\varepsilon < \delta$ and $\varepsilon = n$ mod
$\omega$ then $\langle \bar d_{\alpha,\varepsilon} \char 94 
\bar c_{\alpha,\varepsilon}:\alpha \in u_{2 \varepsilon + 2n} \cap 
\delta\rangle \char 94 \langle \bar d'_\varepsilon \char 94 
\bar c'_\varepsilon\rangle \char 94 \langle \bar
d_{\alpha,\varepsilon} \char 94 \bar c_{\alpha,\varepsilon}:
\alpha \in u_{\varepsilon + 2n} \backslash \delta\rangle$ 
is an $n$-indiscernible sequence over $A_{\varepsilon + 2n}$.
\end{enumerate}
\mn
[Why?  By $\boxplus_3$, the choice of $q(\bar x_{\bar d_\delta},\bar
x_{\bar c_\delta})$ and the choice of $(\bar d'_\delta,\bar
c'_\delta),(\bar d'_\varepsilon,\bar c_\varepsilon)$.  
Note that $\delta \notin u_\varepsilon$ 
by the definition of $u_\varepsilon$.]
\mn
\begin{enumerate}
\item[$\boxplus_5$]   if $\varepsilon < \delta$ and $\varepsilon = n$
mod $\omega$ and $\beta \in u_{\varepsilon +2n} \backslash \delta$ and
$v \subseteq u_{\varepsilon +2n} \cap \beta$ and $|v| < n$ \then \, 
$\bar c_\varepsilon \char 94 \bar d_\varepsilon$ and 
$\bar c_{\beta +1,\varepsilon} \char 94 
\bar d_{\beta +1,\varepsilon}$ realize the same type over 
$A_{v,\varepsilon} + \bar d'_\varepsilon + \bar c'_\varepsilon$ where 
$A_{v,\varepsilon} = A_0 + \Sigma\{\bar d_{\alpha,\varepsilon} \char 94 
\bar c_{\alpha,\varepsilon}:\alpha \in v\}$.
\end{enumerate}

We elaborate the more complicated case.

\medskip

\noindent
\underline{Proof of $\boxtimes_5$ for \ref{c72}}:

Let 
$v_1=v,v_2 = v \cup \gamma$ where $\varepsilon + 2n + 3 < \gamma < \delta$.

So assume
\mn
\begin{enumerate}
\item[$\odot_1$]   $(a) \quad \varphi = 
\varphi(\bar x_{\bar d_\varepsilon},\bar x_{\bar c_\varepsilon},\bar y)$ 
\sn
\item[${{}}$]  $(b) \quad \bar b_1 \in {}^{\ell g(\bar y)}
(A_{v,\varepsilon} + \bar d'_\varepsilon + \bar c'_\varepsilon)$ 
\sn
\item[${{}}$]   $(c) \quad {\gC} \models \varphi[\bar d_{\beta +1,\varepsilon},
\bar c_{\beta+1,\varepsilon},\bar b_1]$.
\end{enumerate}
\mn
We choose $\bar b_2$ such that 
\mn
\begin{enumerate}
\item[$\odot_2$]   $\bar b_2 \in {}^{\ell g(\bar y)}(A_{v_2,\varepsilon})$
and $\bar b_1,\bar b_2$ realize the same type over $A_0 + \bar
c_{\beta +1,\varepsilon +1} + \bar d_{\beta +1,\varepsilon +1} + 
\bar c_{\beta,\varepsilon +1} + \bar d_{\beta,\varepsilon +1}$.
\end{enumerate}
\mn
[Why possible?  By $\boxplus_4$ and the choice of $\gamma$.]

So
\mn
\begin{enumerate}
\item[$\odot_3$]   $(a) \quad \gC \models \varphi[\bar d_{\beta
+1,\varepsilon},\bar c_{\beta +1,\varepsilon},\bar b_2]$
\sn
\item[${{}}$]  $(b) \quad {\gC} \models \varphi[\bar d_\varepsilon,\bar
c_\varepsilon,\bar b_2]$.
\end{enumerate}
\mn
[Why?  Clause (a) follows by clause $\odot_1(c)$ and the choice of
$\bar b_2$, i.e. $\odot_2$.  Clause (b) follows from clause (a) by
$\oplus(e)$.] 

Let $(\eta_1,\nu_1) = \supp_{\bold x}(\varphi)$ and let
$(\eta_0,\nu_0),\psi$ be as guaranteed in (a degenerated case of)
Definition \ref{c23}(3A) for $\varphi$ and $\varphi \equiv
\varphi'(\bar x_{\bar d_\varepsilon,\eta_1},
\bar x_{\bar c_\varepsilon,\nu_1},\bar y)$.

So
\mn
\begin{enumerate}
\item[$\odot_4$]   ${\gC} \models \varphi'[\bar d_{\varepsilon
+1,\eta_1},\bar c_{\varepsilon +1,\nu_1},\bar b_2]$.
\end{enumerate}
\mn
So by the choice of $(\eta_0,\nu_0)$ and $\psi$
\mn
\begin{enumerate}
\item[$\odot_5$]   $(a) \quad {\gC} \models 
\varphi'[\bar d_{\varepsilon +1,\eta_0},\bar c_{\varepsilon +1,\nu_0},
\bar b_2]$
\sn
\item[${{}}$]  $(b) \quad \psi = \psi(\bar x_{\bar d_{\varepsilon
+1,\eta_0}},x_{\bar c_{\varepsilon +1,\nu_0}},
\bar x'_{\bar d_{\varepsilon+1}},\bar x'_{\bar c_{\varepsilon+1}})$
\sn
\item[${{}}$]  $(c) \quad {\gC} \models \psi[\bar
d_{\varepsilon +1,\eta_0},\bar c_{\varepsilon +1,\nu_0},
\bar d_{\beta,\varepsilon+1},\bar c_{\beta,\varepsilon+1}]$
\sn
\item[${{}}$]  $(d) \quad {\gC} \models \vartheta[\bar
c_{\varepsilon +1,\nu_0},\bar d_{\beta,\varepsilon +1},\bar
c_{\beta,\varepsilon +1},\bar b_2]$ where

\hskip25pt  $\vartheta(\bar x_{\bar c_{\varepsilon +1,\nu_0}},
\bar x'_{\bar d_{\varepsilon +1}},
\bar x'_{\bar c_{\varepsilon+1}},\bar y) :=$

\hskip25pt $(\forall \bar x_{\bar d_{\varepsilon +1,\eta_0}})
(\psi(\bar x_{\bar d_{\varepsilon +1,\eta_0}},
\bar x_{\bar c_{\varepsilon +1,\nu_0}},\bar x'_{\bar
d_{\varepsilon +1}},\bar x'_{\bar c_{\varepsilon +1}}) \rightarrow
\varphi'(\bar x_{\bar d_{\varepsilon +1,\eta_0}},
\bar x_{\bar c_{\varepsilon +1,\nu_0}},\bar y))$.
\end{enumerate}
\mn
Next
\mn
\begin{enumerate}
\item[$\odot_6$]  ${\gC} \models 
\vartheta[\bar c_{\varepsilon +1,\nu_0},
\bar d_{\beta,\varepsilon +1},\bar c_{\beta,\varepsilon +1},\bar b_1]$.
\end{enumerate}
\mn
[Why?  By $\odot_2$ the sequences
$\bar d_{\beta,\varepsilon +1} \char 94 \bar
c_{\beta,\varepsilon +1} \char 94 \bar b_1,\bar d_{\beta,\varepsilon +1} \char
94 \bar c_{\beta,\varepsilon +1} \char 94 \bar b_2$ realize the 
same type over $A_0 \supseteq B^+_{\bold x}$
hence also over $B^+_{\bold x} + \bar c_{\varepsilon +1}$, 
so by $\odot_5(d)$ we get the statement in $\odot_6$.]
\mn
\begin{enumerate}
\item[$\odot_7$]   ${\gC} \models \varphi'[\bar d_{\varepsilon
+1,\eta_0},\bar c_{\varepsilon +1,\nu_0},\bar b_1]$.
\end{enumerate}
\mn
[Why?  Recall by $\odot_5(c)$ we have
${\gC} \models \psi[\bar d_{\varepsilon +1,\eta_0},
\bar c_{\varepsilon +1,\nu_0},\bar d_{\beta,\varepsilon +1},\bar
c_{\beta,\varepsilon +1}]$ so by $\odot_6$ and the definition of
$\vartheta$ we get $\odot_7$.]

By the choice of $(\eta_0,\nu_0)$
\mn
\begin{enumerate}
\item[$\odot_8$]   ${\gC} \models \varphi'[\bar d_{\varepsilon,\eta_1},
\bar c_{\varepsilon,\nu_1},\bar b_1]$. 
\end{enumerate}
\mn
This means
\mn
\begin{enumerate}
\item[$\odot_9$]   ${\gC} \models
\varphi[\bar d_\varepsilon,\bar c_\varepsilon,\bar b_1]$.
\end{enumerate}
\mn
As for any $\varphi = \varphi(\bar x_{\bar d_\varepsilon},
\bar x_{\bar c_\varepsilon},\bar y)$ and $\bar b_1 \in {}^{\ell g(\bar
y)}(A_{v,\varepsilon} + \bar d' + \bar c')$ for some truth value
$\bold t$ the statement $\odot_1(c)$ holds for $\varphi^{\bold t}$,
i.e. $\gC \models \varphi^{\bold t}[\bar d_{\beta +1,\varepsilon},\bar
c_{\beta +1,\varepsilon},\bar b_1]$ hence by the above, see
$\odot_9$, we get $\gC \models \varphi[\bar d_\varepsilon,\bar
c_\varepsilon,\bar b_1]$.  Hence we get the equality of types stated in
$\boxplus_5$, so $\boxplus_5$ holds.
\mn
\begin{enumerate}
\item[$\boxplus_6$]  if $\varepsilon < \delta$ and $\varepsilon =
n$ mod $\omega$ then $\langle \bar d_\alpha \char 94 \bar
c_\alpha:\alpha \in u_{2 \varepsilon + n+2} \cap \delta\rangle 
\char 94 \langle \bar d' \char 94 \bar c'\rangle \char 94 
\langle \bar d_\alpha \char 94 \bar c_\alpha:
\alpha \in u_{\varepsilon + 2n+2} \backslash \delta\rangle \char 94
\langle \bar c_\varepsilon \char 94 
\bar d_\varepsilon \rangle$ is an $n$-indiscernible sequence over $A_0$.
\end{enumerate}
\mn
[Why?  By $\boxplus_5$ recalling \ref{c49}; the main point is that
clause (b) there holds, except that for $\bar d_\varepsilon \char 94 \bar
c_\varepsilon$ we can use $\boxplus_4$ and for this case we use
$\boxplus_5$.] 

This shows that for each finite $v \subseteq \delta$ and 
$\varepsilon < \delta$, the pair $(\bar d'_\varepsilon,
\bar c'_\varepsilon)$ solves $(\bold m_\varepsilon,A_v)$, but this means 
that $(\bar d',\bar c')$ solves 
$(\bold m_\delta,A)$ which is what we need. 
\end{PROOF}

\begin{conclusion}
\label{c76}  
1) If $\delta < \theta^+$ is a limit ordinal and $\langle(\bold
x_\varepsilon,\bar\psi_\varepsilon,\bar r_\varepsilon):\varepsilon <
\delta\rangle$ is a $\le_1$-increasing sequence of members of
$\tK^\oplus_{\kappa,\mu,\theta}$, see Definition \ref{c1}, \then \, the
limit $(\bold x_\delta,\bar\psi_\delta,r_\delta)$ belongs to
$\tK^\oplus_{\kappa,\bar\mu,\theta}$ and is a $\le_1$-lub of the sequence. 

\noindent
2) Similarly for $\vK_{\kappa,\bar\mu,\theta}$.
\end{conclusion}

\begin{PROOF}{\ref{c76}}
1) By \ref{c70} as $\bold m \in \tK^\oplus_{\kappa,\mu,\theta},
\bold m \le_1 \bold n \in \rK^\oplus_{\kappa,\mu,\theta} 
\Rightarrow \bold m \le^+_1 \bold n$ so we use (c) rather than (c)$'$ there.

\noindent
2) Similarly by \ref{c72}(2), so we use (c)$_1$ rather than (c)$'_1$ there.  
\end{PROOF}

\begin{claim}
\label{c78}   
Assume $\kappa$ is weakly compact $> \theta \ge |T|$.  

\noindent
1) If $(\bold x_1,\bar\psi_1,\bar r_1) \in 
\rK_{\kappa,\kappa,\theta}$ and $M_{\bold x}$ has cardinality
$\kappa$, \then \, there is $(\bold x_2,\bar\psi_2,r_2) \in 
\rK^\oplus_{\kappa,\kappa,\theta}$ which is $\le^+_1$-above $(\bold
x_1,\bar\psi_1,r_1)$.

\noindent
2) If $M \in \EC_{\kappa,\kappa}(T)$ and $\bar d \in {}^{\theta^+ >}\gC$
   \then \, for some $\bold m \in \tK^\oplus_{\kappa,\kappa,\theta}$ we
have $M_{\bold x[\bold m]} = M$ 
and $\bar d \trianglelefteq d_{\bold x[\bold m]}$.

\noindent
3) If $M \in \EC_{\kappa,\kappa}(T)$ \then \, $\gS^\theta_{\aut}(M)$ has
   cardinality $\le \kappa$.
\end{claim}

\begin{remark}
\label{c65}
Compare with \cite[\S4]{Sh:900}.

So for $\kappa$ weakly compact we can prove
the density of $\tK^\oplus_{\kappa,\kappa,\theta}$ 
(by \ref{c70} + \ref{c78} above), hence 
using the $(\bold D_{\bold x},\kappa)$-sequence 
homogeneity (see Theorem \ref{c32}, Conclusion \ref{c39}) 
we can prove that there are few types (i.e. $\le \kappa$) 
up to conjugacy on saturated
model (the proof in the end of \S4 use only this).
To get it for some smaller cardinals we shall need a replacement of
weak compactness which is the major point of \S4 and to get it for
all large enough $\kappa$ we use $\vK_{\kappa,\bar\mu,\theta}$.
\end{remark}

\begin{PROOF}{\ref{c78}}
1) By \ref{b22}(1) there is $\bold y \in 
\qK'_{\kappa,\kappa,\theta}$ such that $\bold x_1 \le_2 \bold y$ so
$\bar d_{\bold y} = \bar d_{\bold x_1}$ and as we are using
$\rK_{\kappa,\kappa,\theta}$  \wilog \,
$u_{\bold y} = \emptyset$ .  Let
$\langle M_\alpha:\alpha < \kappa\rangle$ be $\prec$-increasing
continuous with union $M_{\bold x}$ such that $\|M_\alpha\| < \kappa$
for $\alpha <\kappa$.
As $(\bold x_1,\bar\psi_1,r_1) \in 
\rK_{\kappa,\kappa,\theta}$, for $\alpha < \kappa$ we can choose $(\bar
c_\alpha,\bar d_\alpha)$ from $M_{\bold x}$ solving $(\bold
x_1,\bar\psi_1,r_1,\bar e_\alpha + M_\alpha)$, see \ref{c5}(1)(f).  
As $\bold y \in \qK_{\kappa,\kappa,\theta}$ by \ref{b22}(3) we have
$\bold y \in \qK_{\kappa,\kappa,\theta}$ so, see Definition \ref{b16}(2) 
for each $\alpha < \kappa$ there are $\bar e_\alpha \in
{}^\theta(M_{\bold x})$ and $\bar\psi_\alpha$ such that
$\tp(\bar d_{\bold y},\bar c_{\bold y} + \bar e_\alpha) \vdash 
\tp(\bar d_{\bold y},\bar c_{\bold y} \dotplus M_\alpha))$ according to
$\bar\psi_\alpha$. 
As $\kappa$ is weakly compact we can find
$(\bar\psi_*,f)$ such that
\mn
\begin{enumerate}
\item[$(*)$]   $(a) \quad f$ is an increasing function from
$\kappa$ to $\kappa$ so $\alpha \le f(\alpha)$,
\sn
\item[${{}}$]  $(b) \quad \bar\psi_{f(\alpha)} = \bar\psi_*$ for
$\alpha < \kappa$
\sn
\item[${{}}$]  $(c) \quad \langle \tp(\bar d_{f(\alpha)}
\char 94 \bar c_{f(\alpha)} \char 94 \bar e_{f(\alpha)},\bar
d_{\bold y} + \bar c_{\bold y} + M_\alpha):\alpha < \kappa\rangle$ is
$\subseteq$-increasing.
\end{enumerate}
\mn
Let $(\bar c,\bar d,\bar e)$ from ${\gC}$ be such that $\bar c
\char 94 \bar d \char 94 \bar e$ realize 
$\cup\{\tp(\bar c_{f(\alpha)} \char 94 \bar d_{f(\alpha)} \char 94 \bar
e_{f(\alpha)},\bar c_{\bold x} + \bar d_{\bold x} + M_\alpha):
\alpha < \kappa\}$, but the pairs $(\bar c,\bar d),(\bar c_{\bold y},\bar
d_{\bold y})$ realize the same type over $M_{\bold x}$ so \wilog \,
$(\bar c,\bar d) = (\bar c_{\bold y},\bar d_{\bold y})$.

Hence
\mn
\begin{enumerate}
\item[$(*)$]  for $\alpha < \kappa$ the sequences $\bar c_{\bold y}
\char 94 \bar d_{\bold y} \char 94 \bar e$ and $\bar c_{f(\alpha)}
\char 94 \bar d_{f(\alpha)} \char 94 \bar e_{f(\alpha)}$ realize 
the same type over $M_\alpha$.
\end{enumerate}
\mn
Now we can define $(\bold x_2,\bar\psi_2,r_1)$ as follows:
\mn
\begin{enumerate}
\item[$\circledast$]   $(a) \quad M_{\bold x_2} = M_{\bold x_1}$
\sn
\item[${{}}$]  $(b) \quad \bar d_{\bold x_2} = \bar d_{\bold x_1}
\char 94 \bar e = \bar d_{\bold y} \char 94 \bar e$ and\footnote{no
  harm in demanding $u_{\bold x_1} = \emptyset$}
\sn
\item[${{}}$]  $(c) \quad \bar c_{\bold x_2} = \bar c_{\bold y}$
\sn
\item[${{}}$]  $(d) \quad \bar B_{\bold x_2} = \bar B_{\bold y}$ and
$v_{\bold x_2} = v_{\bold y}$ and $u_{\bold x_2} = u_{\bold x_1}$
\sn
\item[${{}}$]  $(e) \quad \bar\psi_2$ is just putting together
$\bar\psi_1$ and $\bar\psi_*$
\sn
\item[${{}}$]  $(f) \quad r_2$ is such that $r_2 = 
\tp(\bar c_{\bold x_2} \char 94 \bar d_{\bold x_2} \char 94
\bar e \char 94 \bar c_{f(\alpha)} \char 94 \bar d_{f(\alpha)} \char 94 
\bar e_{f(\alpha)},\emptyset)$

\hskip25pt  for unboundedly many $\alpha < \kappa$.
\end{enumerate}
\mn
Clearly $(\bold x_2,\bar\psi_2,r_2)$ is as required.

\noindent
2) By part (1) and \ref{c76}(1).

\noindent
3) By part (2) and \ref{c32}.
\end{PROOF}
\newpage

\section {Density} 

Our immediate goal is, concentrating on countable $T$, to prove 
density for $\tK_{\kappa,\mu,\theta}$ in some ZFC cases: 
$\kappa = \mu^{+n},\mu = \beth_\delta,\cf(\delta) > \theta$.  
We do it in \ref{d39} when $n=1$ and shall do it in \S(5A) when $n \ge 1$.  
The point is proving some $\bar e$ universally
solves a given $\bold x \in \qK_{\kappa,\kappa,\theta}$ done in
\S(4B) and for this we use the partition Theorem \ref{d10}.  
Theorem \ref{d21} is a partition theorem which 
is nicer per se, and is
more transparent (and stronger in some respects, see also \ref{d23}), 
but it is not enough for helping in the proofs in decompositions.  
\bigskip

\subsection {Partition theorems for Dependent $T$} \
\bigskip

The following partition theorem will be crucial (in the proof
of \ref{d33} and also will be used in \ref{p4}).  We prove a nicer one
later, but not useful here.  We can below use ``$v$
finite, $\bold k=1$" in \ref{d10}, see \ref{d15}(3).  For a case when
the conclusion of \ref{d10} can be nicely phrased, see \ref{d19}.  In
\ref{d10} we do not explicitly demand $T$ to be dependent
\underline{but} clause (i) holds if $T$ is dependent.
\begin{theorem}
\label{d10}
\underline{The partition Theorem}   

There are ${\cD}_1$-positive sets $\cS_{1,n}$ for 
$n < \bold k$ and a $\cD_2$-positive set ${\cS}_2$ 
and $h \in \Pi{\ga}_{\bold f}$  
and $q_n \in {\bold S}^{v_n+w_n}_{\Delta_n}(B^+_{\bold f,h} +
\Sigma_n C_n)$ such that 
for every $n < \bold k,\cS_{1,n+1} \subseteq \cS_{1,n}$  and\footnote{can
add $s \in \mathscr{S}_{1,n} \Rightarrow \tp
(\bar e_s,B^+_{\bold f,u_n,h}) = q_n$} for each $s \in 
\cS_{1,n}$ for ${\cD}_2$-almost every $t \in 
\cS_2$ we have\footnote{could ask just $q \rest \Delta_n= 
\tp_{\Delta_n}(\bar e_s \char 94 \bar e_i \char 94 \bar a_{\bold
f,u_n,h},B_{\bold f,u_n}+C)$ for every 
$h \in \Pi\{\kappa_{\bold f,i} \backslash
h_n(i):i \in u_n\}$ does not really matter.} $q_n = 
\tp_{\Delta_n}((\bar e^1_s \rest v_n) \char 94 (\bar e^2_t \rest
w_n),B^+_{\bold f,u_n,h} + C_n)$, see \ref{b41}(2) \when \,:
\medskip

\noindent
\begin{enumerate}
\item[$\oplus$]  $(a) \quad \bold k \le \omega$ and $k_n <
\omega,k_n \ge 1$ for $n < \bold k$ 
\smallskip

\noindent
\item[${{}}$]  $(b) \quad$ for $n < \bold k,\Delta_n \subseteq 
\bbL(\tau_T)$ is finite, each $\varphi \in \Delta_n$ has
the form $\varphi(\bar x_{[v_n]},\bar y_{[w_n]},\bar z)$
\smallskip

\noindent
\item[${{}}$]  $(c)(\alpha) \quad \bold f$ is a
$(\bar\mu,\theta)$-set, see Definition \ref{b38}, \ref{b41}, we use
  their notations
\sn
\item[${{}}$]  $\quad (\beta) \quad v_n \subseteq v$ is
finite\footnote{instead ``$v_n$ finite" we can use $v_n = \varepsilon$
but $\Delta^1_n \subseteq \Gamma_{[\varepsilon],n} \subseteq 
\bbL(\tau_T)$, see \ref{z23}(4)} and $\subseteq$-increasing with $n <
\bold k$
\sn
\item[${{}}$]  $\quad (\gamma) \quad u_n \subseteq v_{\bold f}$ for $n
< \bold k$, note $v,v_{\bold f}$ are unrelated objects

\hskip25pt and let $u_{n,2} = u_n \cap u_{\bold f}$ and $u_{n,1} = u_n
\backslash u_{\bold f}$
\sn
\item[${{}}$] $\quad (\delta) \quad \cf(\Pi\{\kappa_{\bold f,i}:
i \in u_{n,2}\}) = \cf(\Pi{\frak a}_{\bold f,u_n,2})< \kappa$
\sn
\item[${{}}$]  $\quad (\varepsilon) \quad w_n \subseteq w_{n+1} 
\subseteq w = \cup\{w_m:m < \omega\}$
\sn
\item[${{}}$]  $\quad (\zeta) \quad$ for 
simplicity\footnote{As we can work in $\gC^{\eq}$ this
is not a loss.} $\bar a_{\bold f,i,\alpha}$ is a singleton for $i \in
u_{\bold f},\alpha < \kappa_{\bold f,i}$
\sn
\item[${{}}$]  $(d) \quad \Delta^1_n \subseteq 
\bbL(\tau_T)$ and $C_n \subseteq \gC$ for $n < \bold k$
\sn
\item[${{}}$]   $(e) \quad (\alpha) \quad \kappa = \cf(\kappa)$ and
$\min\{\kappa_{\bold f,i}:i \in u_n \cap u_{\bold f}\}$ are 

\hskip25pt  $> \beth_{k_n}(|B_{\bold f,u_n \backslash u_{\bold f}}| 
+ |C_n| + \theta),\theta \ge \aleph_0 + |v_n| + |u_n| + 
|\Delta^1_n|$ 

\hskip25pt  and $C_n \subseteq \gC$

\underline{or}
\sn
\item[${{}}$]   $\qquad \quad (\beta) \quad \kappa = \cf(\kappa)$ and
$B^+_{\bold f,u_n \backslash u_{\bold f}},v_n,w_n,C_n,\Delta^1_n$ are finite
\smallskip

\noindent
\item[${{}}$]   $(g) \quad \bar e^1_s \in {}^v \gC$ for $s \in I_1$ 
\sn
\item[${{}}$]  $(h) \quad \bar e^2_t \in {}^w \gC$ for $t \in I_2$
\sn
\item[${{}}$]  $(i) \quad$ if $n < \omega,B \subseteq \gC,\bar e \in 
{}^{w_n} \gC$ and $\langle \bar a_\ell:
\ell < \omega\rangle$ is a $(\Delta^1_n,k_n)$-indiscernible 

\hskip25pt sequence over $B$ where $\ell < \omega \Rightarrow
\bar a_\ell \in {}^{v_n} \gC$ \then \, the set

\hskip25pt  $\{\ell < \omega:\tp_{\Delta_n}(\bar a_{\ell}
\char 94 \bar e,B) \ne \tp_{\Delta_n}(\bar a_{\ell +1} \char
94 \bar e,B)\}$ is finite
\smallskip

\noindent
\item[${{}}$]  $(j) \quad {\cD}_2$ is a $\kappa$-complete filter on $I_2$
\smallskip

\noindent
\item[${{}}$]  $(k)(\alpha) \quad \cD_1$ is a filter on $I_1$
\sn
\item[${{}}$]  $\quad (\beta) \quad$ if $\bold k = \omega$ then
\sn
\item[${{}}$]  \hskip35pt $\bullet \quad \cD_1$ is $\aleph_1$-complete 
\sn
\item[${{}}$]  \hskip35pt $\bullet \quad \alpha < \kappa
\Rightarrow |\alpha|^{\aleph_0} < \kappa$.
\end{enumerate}
\end{theorem}

\begin{remark}
\label{d11}
Similarly if $T$ is strongly dependent 
(hence by \cite[4.1]{Sh:863} we already 
get some existence of indiscernibles) we can get more.
\end{remark}

\begin{claim}
\label{d12}  
1) Assume $T$ is strongly dependent.
In \ref{d10} we can use: $\Delta_n = \Delta^1_n =
\{\varphi:\varphi = \{\varphi(\bar x_{[v_n]},\bar y_{[w]},\bar z) \in
\bbL(\tau_T)\}$ but demand $\kappa = \cf(\kappa) > \beth_{|T|^+}(|B_{\bold
f,u_n}| + |C_n| + \theta)$.

\noindent
2) Similarly in \ref{d21}.
\end{claim}

\begin{PROOF}{\ref{d12}}
1) Like the proof \ref{d10} but we just use 
\cite[4.1]{Sh:863} instead Erd\"os-Rado theorem.

\noindent
2) Similarly.
\end{PROOF}

\begin{remark}
\label{d15}  
1) A nice case is $\cup\{\Delta_n:n < \bold k = \omega\} =
   \{\varphi(\bar x_{[v_n]},\bar x_{[w]}):\varphi \in \bbL(\tau_T)$
   and $n < \omega\}$ and
$\cup\{v_n:n < \omega\} = v$ and $\cup\{\Delta^1_n:n < \omega\}
 = \bbL(\tau_T)$ and each $\Delta_n,\Delta^1_n$ is finite so $T$ is countable.

\noindent
2) If we first replace $(\bar c_{\bold x},
\bar d_{\bold x})$ by $(\langle \rangle,\bar c_{\bold x} \char 94 
\bar d_{\bold x})$ the way back is problematic!

\noindent
3) We could use $v$ finite.  Also we may use $I_\ell = I
= [M]^{< \kappa}$ and $\cD_\ell$ is a normal filter on $I$, it is
   natural in the application here (similarly for the definition 
of a $(\bar\mu,\theta)$-sets!)

\noindent
4) In clause (c) of \ref{d10}, by \ref{b45}(1A), it suffices to
demand
\medskip

\noindent
\begin{enumerate}
\item[$(c)'$]  $\bold f$ is a $(\bar\mu,\theta)$-set, $\Delta_n
\subseteq \bbL(\tau_T)$ for $n < \bold n$ and $i \in v_{\bold f}
\backslash u_{\bold f} \Rightarrow B_{\bold f,i} = B_{\bold f}$.
\end{enumerate}

\noindent
5) We have considerable leeway in the proof.

\noindent
6) In order to use infinite $\Delta_n$ at 
present we need a stronger assumption on $T$,
see \ref{d12}.

\noindent
7) For transparency assuming $\bold k = 1$, we can get also that for some
$q'$ for each $t \in \cS_2$ for ${\cD}$-almost $s \in \cS_{1,0}$ we 
have $q' = 
\tp_{\Delta_0}(\bar e_s \char 94 \bar e_t,B^+_{\bold f,u_n,h} + C)$.

\noindent
8) In the proof we can demand that $\bar s_{n,k}$ has length $k+1$ and so
can demand that the game is of a fix finite number of moves,
e.g. $\Sigma\{2 \times \ind(\varphi):\varphi \in \Delta_n\} + 1$, on
$\ind(\varphi)$, see \ref{q3}.

\noindent
9) We can assume $\cS^*_\ell \in \cD^+_\ell$ for $\ell=1,2$ and
demand $\cS_{1,n} \subseteq \cS^*_1,\cS_2 \subseteq \cS^*_2$ but 
this does not add anything because we
may just use $\cD'_\ell = \cD^*_\ell \rest \cS^*_\ell := 
\{\cS \cap \cS^*_\ell:\cS \in \cD_\ell\}$.

\noindent
10) There is no real harm if in \ref{d10} we assume $v=w$,
i.e. $\Dom(\bar e^1_s) = \Dom(\bar e^2_t)$ for $s \in I_1,t \in I_2$.

\noindent
11) Assume $\bold f_1$ satisfies $\cf(\Pi \ga_{\bold f_1}) < \kappa$ (or
$\kappa \notin \pcf(\ga_{\bold f_1})$; so we have $h_{\bar e_1,
\bar e_2} \in \Pi \ga_{\bold f_1}$.  Can we find $h$ satisfying $(\forall
s_0 \in \cS^1_{\alpha,0})(\forall^{\cD_2} s_1 \in
\cS_1)(h_{\bar e_{s_0},\bar e_{s_1}} \le h)$? see \ref{p2}.

\noindent
12) We could have asked $u_n \subseteq v_{\bold f}$ instead $u_n
    \subseteq u_{\bold f}$ and use $B^\pm_{\bold f,u_n}$ instead of
    $B_{\bold f,u_n}$.
\end{remark}

\begin{remark}
\label{d19}
If you are interested in weakening the generality of the theorem for
having a somewhat more transparent proof, note that the 
statement of \ref{d10} is simplified when we use a model $M$ of
cardinality $\kappa$ to which all relevant elements belong (as in the
proof).  Let $\langle M_\alpha:\alpha <
\kappa\rangle$ be $\prec$-increasing continuous with union $M$ such
that $\alpha < \kappa \Rightarrow \|M_\alpha\| < \kappa$.  So we can
decide $I_\ell = \kappa,\cD_\ell$ a normal filter on $\kappa$, e.g. the
club filter hence instead of ``for every $s \in \cS_{1,m}$ for 
$\cD_2$-almost every $t \in \cS_2$" we
have: if $s \in \cS_{1,n},t \in \cS_2$ and $s < t$ as ordinals \then
\, $q_n = \tp_{\Delta_n}((\bar e^1_s \rest v_n) 
\char 94 (\bar e^2_t \rest w_n),B^+_{\bold f,u_n,h} + C_n)$.  
This is by normality of the filter.
\end{remark}

\begin{PROOF}{\ref{d10}}  
\underline{Proof of \ref{d10}}  Let $M \prec \gC$ include $\cup\{\bar
e^\ell_s:s \in I_\ell$ and $\ell=1,2\} \cup \{C_n:n < \omega\} \cup 
B^+_{\bold f}$.  Let $\mu = \min\{\mu_\ell:\ell < 2\}$ hence
$\mu \le \min(\ga_{\bold f})$.
We can choose $\bar a_{*,k} = \langle \bar a_{*,i,k}:
i \in u_{\bold f}\rangle$ for $k < \omega$ such that
(recalling Definition \ref{b41}(5)):
\medskip

\noindent
\begin{enumerate}
\item[$\boxplus_1$]   for every $A \subseteq M$ of cardinality $<
\mu$ and $k < \omega$, for some $g \in \Pi \ga_{\bold f,u_{\bold f}}$ we have: 
$g \le h \in \Pi \ga_{\bold f,u_{\bold f}} \Rightarrow 
\text{ tp}(\bar a_{\bold f,u_{\bold f},h},A + B_{\bold f} + C_n 
+ \sum\limits_{m > k} \bar a_{*,m}) = \tp(\bar a_{*,k},A + 
B_{\bold f} + C_n + \sum\limits_{m > k} \bar a_{*,m})$.
\end{enumerate}

\noindent
Let
\medskip

\noindent
\begin{enumerate}
\item[$\boxplus_2$]   $C^+_n = \cup\{\bar a_{*,n}:n < \omega\} \cup
B_{\bold f} \cup C_n$ but if $\bold k < \omega$, i.e. the assumptions
of clause $(\beta)$ of clause (e) of $\oplus$
fails, \then \, 
$C^+_n = \cup\{\bar a_{u_n}:n < k_n\} \cup B_{\bold f} \cup C_n$
\end{enumerate}

\noindent
hence
\medskip

\noindent
\begin{enumerate}
\item[$\boxplus_3$]   if $n < \bold k,\zeta < \mu$ and
$\bar e_1,\bar e_2 \in {}^\zeta M$ realize the same type
over $C^+_n$ \then \, for some $g \in \Pi{\ga}_{\bold f,u_n}$, we
have $(C^+_n + \bar e_1 + \bar e_2,\bar{\bold I}_{\bold f,u_{n,2},g})$ 
is a $(\bar \mu,\theta)$-set and $\bar e_1,\bar e_2$ realize the same type over
$C^+_n + \bar{\bold I}_{\bold f,u_{n,2},g}$.
\end{enumerate}

\noindent
Choose ${\cF}_n$ such that
\medskip

\noindent
\begin{enumerate}
\item[$\boxplus_4$]   ${\cF}_n$ is a cofinal subset of
$\Pi{\ga}_{\bold f,u_{n,2}}$ of cardinality cf$(\Pi{\ga}_{\bold f,u_{n,2}})$, 
\end{enumerate}

\noindent
hence by clause $(c)(\delta)$ of the assumption
\medskip

\noindent
\begin{enumerate}
\item[$\boxplus_5$]  $|{\cF}_n| < \kappa$.
\end{enumerate}

\noindent
Without loss of generality
\medskip

\noindent
\begin{enumerate}
\item[$\boxplus_6$]  $(B_{\bold f,u_{n,1}} + C^+_n,
\bar{\bold I}_{\bold f,u_{n,2}})$
is a $(\bar\mu,\theta)$-set for each $n < \bold k$
\smallskip

\noindent
\item[$\boxplus_7$]  $|\tau_T| = \Sigma\{|\Delta_n| + 
|\Delta^1_n|:n <\bold k\}$ or both are $\le \aleph_0$.
\end{enumerate}

Next
\begin{enumerate}
\item[$(*)_1$]  if $\ell = 1,2$ and $n < \bold k$ and $\cS \in 
{\cD}^+_\ell$ \then \, for some $q = q_{\cS,n}$ and $h =
h_{\cS,n} \in {\cF}_n \subseteq \Pi{\ga}_{\bold f,u_n}$ 
we have $\cS_{h,q} \in {\cD}^+_\ell$ where $\cS_{h,q} := 
\{s \in \cS:(C^+_n + \bar e^\ell_s,\bar{\bold I}_{\bold f,u_n,h})$ is a
$(\bar\mu,\theta)$-set and $q = \tp_{\bbL(\tau_n)}(\bar e_s \rest 
v_n,B^+_{\bold f,u_n,h} + C^+_n)\}$.
\end{enumerate}
\medskip

\noindent
[Why?  For each $s \in \cS$ by \ref{b45}(4) there is a
function $h_s \in \Pi{\ga}_{\bold f,u_n}$ such that 
$(C^+_n + \bar e^\ell_s,\bar{\bold I}_{\bold f,h})$ is a
$(\bar\mu,\theta)$-set; \wilog \, $h_s \in {\cF}_n$.

But by clause $(j)$ of the assumption, i.e. the
$\kappa$-completeness of $\cD_2$ 
and $\boxplus_5$ there is a function $h \in {\cF}_n 
\subseteq \Pi{\ga}_{\bold f,u_n}$ such that $\cS' := \{s \in
\cS:h_s \le h$ and we can add $h_s = h\} \in {\cD}^+_\ell$.
Now $2^{|C^+_n|+|\Delta^1_n|} < \kappa$ by clauses (a),(d),(e) of
$\oplus$, hence for some $(\cS'',q)$ we
have: $\cS'' \subseteq \cS',\cS'' \in {\cD}^+_\ell$ and $s \in 
\cS'' \Rightarrow q = \tp_{\bbL(\tau_n)}
(\bar e^\ell_s,C^+_n)$.  By the
choice of $h_s=h$ and $C^+_n$ it follows that $\cS'' \subseteq
\cS_{h,q}$ so $\cS_{h,q} \in {\cD}^+_\ell$ is as required in $(*)_1$.]

We now define some games; for 
any $n < \bold k,g \in {\cF}_n,q \in {\bold S}^{v_n}_{\bbL(\tau_n)}
(B^+_{\bold f,u_{n,2},g})$ and $\cS \in {\cD}^+_1$ we define a game
$\Game_{\cS,n,g,q}$: a play last $\omega$ moves, in the $\ell$-th move
the antagonist chooses $\cX_\ell \in {\cD}_1$ and the 
protagonist chooses $s_\ell \in \cX_\ell \cap \cS$.

In the end of a play the protagonist wins the play \Iff \,:
\medskip

\noindent
\begin{enumerate}
\item[$\otimes$]   $(a) \quad (C^+_n + \Sigma_{\ell < \omega} 
\bar e^1_{s_\ell},\bar{\bold I}_{\bold f,u_n,g})$ is a $(\bar\mu,\theta)$-set
\smallskip

\noindent
\item[${{}}$]  $(b) \quad 
\langle \bar e^1_{s_\ell} \rest v_n:\ell < \omega\rangle$ is an
$(\bbL(\tau_T),k_n)$-indiscernible sequence over 

\hskip25pt $B^+_{\bold f,u_n,g} + 
C^+_n$ and\footnote{we may omit clause (c)}
\smallskip

\noindent
\item[${{}}$]  $(c) \quad q = \tp_{\bbL(\tau_n)}(\bar e^1_{s_\ell} \rest v_n,
B^+_{\bold f,u_n,g} + C^+_n)$ so is the same for every $\ell < \omega$.
\end{enumerate}
\mn
Alternatively\footnote{we use the second; presently it does not make a
difference what we use} (by $(e)(\beta)$ of $\oplus$) we can define
$\Game'_{\cS,n,g,q}$ similarly but in the end of the play the
protagonist wins the play \Iff \,:
\mn
\begin{enumerate}
\item[$\otimes$]  $(a)' \quad$ as of $\otimes$
\sn
\item[${{}}$]  $(b)' \quad \langle \bar e^1_{s_\ell} \rest v_n:\ell <
\omega\rangle$ is a $(\Delta^1_n,k_n)$-indiscernible sequence over
$B^+_{\bold f,u_n,g} + C^+_n$
\smallskip

\noindent
\item[${{}}$]  $(c)' \quad \tp_{\Delta_n}
(\bar e^1_{s_\ell} \rest v_n,C^+_n) = (q \rest \Delta_n) 
\rest C^+_n$ hence is the same for all $\ell < \omega$.
\end{enumerate}

\noindent
So only $q_n := (q \rest \Delta_n) \rest C^+_n$ matters and we may write 
$\Game'_{\cS,n,g,q_n}$. 
\medskip

\noindent
\begin{enumerate}
\item[$(*)_2$]  if $\cS \in {\cD}^+_1$ and $n < \bold k$ \then \, 
for some $h \in {\cF}_n$ and $q$,
the protagonist wins in the game $\Game'_{\cS,n,h,q}$.
\end{enumerate}

\noindent
[Why?  Let $\lambda = 2^{|C^+_n|+|\Delta^1_n|} + \aleph_1$, so by
Erd\"os-Rado theorem and $\oplus(e)$ clearly $\kappa
\rightarrow (\omega)^{k_n}_\lambda$.

For each $h \in {\cF}_n$ and $q \in \bold S^{v_n}(C^+_n)$ 
the game $\Game'_{\cS,n,h,q}$
is determined being closed for the protagonist, 
so toward contradiction let {\bf st}$_{\cS,n,h,q}$ 
be a winning strategy for the antagonist.  We choose
$s_\alpha \in \cS \subseteq I_1$ by induction on 
$\alpha < \kappa$ such that:  for any relevant $h$ and $q$ in any
finite initial segment of a play of $\Game'_{\cS,n,h,q}$
in which the antagonist uses the strategy {\bf
st}$_{\cS,n,h,q}$ and the 
protagonist chooses members of $I_1$ from $\{s_\beta:\beta <
\alpha\}$, the last move of the antagonist is a member $\cX$ of ${\cD}_1$ to
which $s_\alpha$ belongs.  So $s_\alpha$ just have to belong to
$\le |[\alpha]^{< \aleph_0}| + |{\cF}_n| + |\bold S^{v_n}_{\Delta_n}(C^+_n)| 
< \kappa$ sets\footnote{Why is it $< \kappa$?  Recall
$|{\cF}_n| < \kappa$ by $\boxplus_4$ and 
$|\bold S^{v_n}_{\Delta_n}(C^+_n)| < \kappa$ by clause (e) of the
assumption.} $\cX \in {\cD}_1$, but ${\cD}_1$ is a $\kappa$-complete filter
so this is possible.  As $\kappa = \text{ cf}(\kappa)$ is large enough
\wilog \, $\langle \text{tp}_{\Delta_n}(e^1_{s_\alpha} \rest
v_n,C^+_n):\alpha < \kappa\rangle$ is constant.  Now letting
$\lambda_n = |\bold S^{k_n \times |v_n|}(C^+_n)|$, by clause (e) of
the assumption we
have $\kappa \rightarrow (\omega)^{k_n}_{\lambda_n}$, so for some increasing
sequence $\langle \alpha(i):i <\omega\rangle$ of ordinals $< \kappa$
the sequence $\langle \bar e^1_{s_{\alpha(i)}}:i < \omega\rangle$ is
$(\Delta^1_n,k_n)$-indiscernible over $C^+_n$.  We can find $h \in 
{\cF}_n$ such that $(\cup\{\bar e^1_{t_{\alpha(i)}}:i < \omega\} + C^+_n,
\bar{\bold I}_{\bold f,u_n,h})$ is a $(\bar\mu,\theta)$-set.  By the
choice of $C^+_n$ and of $\langle
\alpha(i):i < \omega\rangle$ it follows that $\langle 
\bar e^1_{s_{\alpha(i)}}:i < \omega\rangle$ is 
$(\Delta^1_n,k_n)$-indiscernible also
over $B^+_{\bold f,u_n,h} + C^+_n$.  So for some $q$ the sequence
$\langle \bar e^1_{s_i}:i < \omega\rangle$ is a result of a play of the game 
$\Game'_{\cS,n,h,q}$ in which the protagonist 
wins, easily a contradiction.]
\medskip

\noindent
\begin{enumerate}
\item[$(*)_3$]  for any $\cS \in {\cD}^+_1$ and $n < \omega$ we choose
$({\cW},q,h$,\,{\bf st},$\bold T) = ({\cW}_{\cS,n},q_{\cS,n},h_{\cS,n}$,\,
{\bf st}$_{\cS,n},\bold T_{\cS,n})$ satisfying:
\begin{enumerate}
\item[$(a)$]  $h \in {\cF}_n$
\smallskip

\noindent
\item[$(b)$]  ${\cW} \in {\cD}^+_1$ and ${\cW} \subseteq \cS$
\smallskip

\noindent
\item[$(c)$]   if $s \in {\cW}$ then $q = \tp_{\Delta_n}
(\bar e^1_s,B^+_{\bold f,u_n,h} + C^+_n)$
\smallskip

\noindent
\item[$(d)$]   {\bf st} is a winning strategy in
$\Game'_{{\cW},n,h,q}$ for the protagonist 
\smallskip

\noindent
\item[$(e)$]   $\bold T = \{\bar s$: for some $m < \omega$ the
sequence $\bar s = \langle \cX_\ell,s_\ell:\ell < m\rangle$ is a 
finite initial segment 
of the play in which the protagonist uses {\bf st}$\}$.
\end{enumerate}
\end{enumerate}
\medskip

\noindent
[Why?  Easy by $(*)_1 + (*)_2$.]

Let $\chi$ be large enough and 
\medskip

\noindent
\begin{enumerate}
\item[$(*)_4$]   let ${\cN}$ be the set of $N \prec ({\cH}(\chi),
\in,<^*_\chi)$ such that:
\begin{enumerate}
\item[$(a)$]   ${\gC},\bold f,M,\langle C_n:n < \bold k \rangle,
\langle \bar e^\ell_t:t \in I_\ell\rangle$ for $\ell=1,2$ 
and $\langle {\cF}_n:n < \bold k\rangle$ belongs to $N$
\smallskip

\noindent
\item[$(b)$]   the following function belongs to $N$:

$\quad (\cS,n) \mapsto ({\cW}_{\cS,n},
q_{\cS,n},h_{\cS,n}$, {\bf st}$_{\cS,n},\bold T_{\cS,n})$
\smallskip

\noindent
\item[$(c)$]  $N$ has cardinality $< \kappa$ and $N \cap \kappa \in
\kappa$ (hence, e.g. $n < \bold k \Rightarrow \cF_n \subseteq N$)
\smallskip

\noindent
\item[$(d)$]  if $\bold k = \omega$ \then \, $[N]^{\aleph_0} \subseteq N$.
\end{enumerate} 
\end{enumerate}
\mn
Now
\mn
\begin{enumerate}
\item[$(*)_5$]   for every $N \in {\cN}$ we can choose $t_N \in
\cap\{\cS \in \cD_2:\cS \in N\}$.
\end{enumerate}
\medskip

\noindent
[Why?  Because $\cD_2$ is a $\kappa$-complete filter and $N$ is of
cardinality $< \kappa$.]
\medskip

\noindent
\begin{enumerate}
\item[$(*)_6$]   for every $N \in {\cN}$ choose 
$(\cS_{1,n},q_n,h_n$,\,{\bf st}$_n,\bold T_n) 
= (\cS_{N,n},q_{N,n},h_{N,n}$,\,{\bf st}$_{N,n},\bold T_{N,n})$ 
by induction on $n < \bold k$ such that:
\begin{enumerate}
\item[$(a)$]    $(\cS_{1,n},q_n,h_n$,\,{\bf st}$_n,\bold T_n) 
\in N$ is as in $(*)_3$
\smallskip

\noindent
\item[$(b)$]   $\cS_{1,n} \supseteq \cS_{1,n+1}$
\smallskip

\noindent
\item[$(c)$]  $\cS_{1,n} \in {\cD}^+_1 \cap N$
\smallskip

\noindent
\item[$(d)$]  if $s \in \cS_{1,n+1} \cap N$ then $q_n = \tp_{\Delta_n}
((\bar e^1_s \rest v_n) \char 94 (\bar e^2_{t_N} \rest w_n),
B^+_{\bold f,u_n,h_n} + C^+_n)$.
\end{enumerate} 
\end{enumerate}

\noindent
We can carry the inductive construction.

\noindent
[Why?  For $n=0$ choose $\cS_{1,n},q_n,h_n$,{\bf st}$_n,\bold T_n$ 
as in $(*)_3$ with $(I_1,n)$ here standing for $(\cS,n)$ there 
and as we are assuming $N \in \cN$ without loss of generality 
they belong to $N$.  Assume that the tuple $(\cS_{1,n},q_n,h_n$,
\,{\bf st}$_n,\bold T_n)$ was chosen and 
$n+1 < \bold k$. We try to choose $\bar s_{n,k} = \langle
\cX_\ell,s_\ell:\ell \le \ell_k\rangle$ by induction 
on $k < \omega$ such that: $\bar s_{n,k}$ is a finite 
initial segment of a play of the game $\Game'_{\cS_{1,n},n,h_n,q_n}$
in which the antagonist uses the strategy {\bf st}$_n$ 
and $\bar s_{n,k} \in N$ and if $k=m + 1$, then for
some $\ell \in [\ell_m,\ell_k-1)$ we have tp$_{\Delta_n}
((e^1_{s_\ell} \rest v_n) \char 94 (\bar e^2_{t_N} \rest w_n),
B^+_{\bold f,u_n,h_n} + 
C^+_n) \ne \tp_{\Delta_n}((\bar e^1_{s_{\ell +1}} \rest v_n) 
\char 94 (\bar e^2_{t_N} \rest w_n),B^+_{\bold f,u_n,h_n} + C^+_n)$. 

If we can choose\footnote{of course, as $\Delta^1_n$ is finite we can
 use a finite long enough game; part of our leeway}
all $\bar s_{n,k},k < \omega$ we get a contradiction to clause (i) of
 the assumption of the theorem.  
Obviously, we can choose $\bar s_{n,0}$.  So for
 some $k,\bar s_{n,k}$ is well defined but we cannot choose $\bar
 s_{n,k+1}$.

Let

\[
\cS'_{1,n} := \{s_{\ell_{k+1}}: \text{ for some } \bar s = \langle
\cX_\ell,t_\ell: \ell \le \ell_k +1\rangle \in \bold T_n \text{ we
  have } \bar s_{n,k} \triangleleft \bar s\}.
\]

\mn
Let $q_n = \tp_{\Delta_n}((\bar e^1_{s_n,\ell_k} \rest v_n) \char 94
(\bar e^2_{t_N} \rest w_n),B^+_{\bold f,u_n,h_n} + C^+_n)$.

Lastly, choose $(\cS_{N,n+1},\bold q_{n+1},h_{n-1}$,{\bf
  st}$_{n+1},\bold T_{n+1})$ as in $(*)_3$ with $(\cS'_{1,n+1},n+1)$
  here standing for $(\cS,n)$ there.  Clearly we are done proving
$(*)_6$, i.e. we can carry the induction.]

So we have chosen $\langle (h_{N,n},q_{N,n},\cS_{1,n}):n < k\rangle$
for each $N \in \cN$ and it belongs to $N$: if $\bold k < \omega$
trivially by $(*)_6$ and if $\bold k = \omega$ by clause (d) of
$(*)_4$ and let $q_N = \cup\{q_{,n}:n < \bold k\}$.  Also $\cN$ is a
stationary subset of $[\cH(\chi)]^{< \kappa}$ by $(*)_3$ and clause (k) of the
assumption.  Hence using the club filter on $\cN$:
\mn
\begin{enumerate}
\item[$(*)_7$]  for some $\cS_{1,n},h_n,r_n,q^*_n$ for $n < \bold k$
  the set $\cS^2 := \{N \in \cN:q_{N,n} = q^*_n$ and $\cS_{N,n} =
  \cS_{1,n}$ for every $n < \bold k\}$.
\end{enumerate}
\mn
Let $h \in \Pi \ga_{\bold f}$ be $\sup\{h_n:n < \omega\}$.  So clearly
any suitable $q \supseteq \cup\{q^*_n:n < k\},h,\langle \cS_{1,n}:n <
\omega\rangle$ and $\cS_2$ are as required.
\end{PROOF}

\noindent
The following is a transparent ``$n(*)$-dimensional" relative of 4.1
\begin{theorem}
\label{d21}
Assume $\kappa$ is regular uncountable, $\cD_0$ is a filter on
$I_0,\cD_\ell$ is a $\kappa$-complete filter on $I_\ell$ for non-zero
$\ell < n,\bar e^\ell_s \in {}^{m(\ell)}\gC$ for $\ell < n,s \in
I_\ell$ and $\Delta \subseteq \bbL(\tau_T),C \subseteq \gC$ are
finite.  \Then \, there are a type $q$ and $\cS_\ell \in \cD^+_\ell$
for $\ell < n$ such that $(\forall^{\cD_0} s_0 \in
\cS_0)(\forall^{\cD_1} s_1 \in \cS_1) \ldots (\forall^{\cD_{n-1}}
s_{n-1} \in \cS_{n-1})[q = \tp_\Delta(\bar e^0_{s(0)} \char 94 \bar
  e^1_{s(1)} \char 94 \ldots \char 94 \bar
  e^{n-1}_{s_{n-1}},\emptyset)]$.
\end{theorem}

\begin{PROOF}{\ref{d21}}
Let $m(< i) = \Sigma\{m(j):j<i\}$.
\medskip

\noindent
\underline{Stage A}:   We prove it by induction on $n$; for $n=0$ it
says nothing, for $n=1$ it holds by the pigeon-hull principle.  So 
assume we know it for $n \ge 1$ and we shall prove it for
$n+1$.

Let $I = \prod\limits_{\ell < n} I_\ell$ and $\bar e_{\bar s} = \bar
e^0_{s_0} \char 94 \ldots \char 94 \bar e^{n-1}_{s_{n-1}} \in
{}^{m(<n)}\gC$ for $\bar s \in I$ and let $\Gamma = \{\{\bar s \in
I:\gC \models \varphi[\bar e_{\bar s},\bar c]\}:
\varphi = \varphi(\bar x_{[m(<n)]},\bar y) \in \bbL(\tau_T)$ 
and $\bar c \in {}^{\ell g(\bar y )}\gC$ and for some finite
$C_1 \subseteq \gC$ and finite $\Delta_1 \subseteq \bbL(\tau_T)$ 
there are no $(\Delta_1,m(<n))$-type $q$ on $C_1$ and sequence
$\langle \cS_\ell:\ell < n\rangle \in 
\prod\limits_{\ell < n} \cD^+_\ell$ such that 
$(\forall^{\cD_0} s_0 \in \cS_0) \ldots
(\forall^{\cD_{n-1}} s_{n-1} \in \cS_{n-1})
[q = \tp_{\Delta_1}(\bar e_{\langle s_\ell:\ell < n\rangle},C_1)]$
and $\neg \varphi(\bar x_{[m(<n)]},\bar c) \in q\}$.

By the induction hypothesis $\Gamma$ is a filter on $I$
hence there is an ultrafilter $\cD_*$ on $I$ extending it.
\medskip

\noindent
\underline{Stage B}: 

Choose finite $\Delta^1 \subseteq \bbL(\tau_T)$ large enough, i.e. such that
\medskip

\noindent
\begin{enumerate}
\item[$(*)_1$]  if $\bar e_\ell \in {}^{m(<n)}\gC$ for $\ell < \omega$
and $\langle \bar e_\ell:\ell < \omega\rangle$ is a
$\Delta^1$-indiscernible sequence over some set $C_1 \subseteq \gC$
\then \, for no formula $\varphi(\bar x,\bar y,\bar z) \in
\Delta,\ell g(\bar x) = m(<n)$ and $\bar b \in {}^{\ell g(\bar y)}\gC$
is the set $\{\ell < \omega$: for some $\bar c \in 
{}^{\ell g(\bar z)} C_1$ we have $\gC \models \varphi(\bar e_\ell,
\bar b,\bar c) \equiv \neg
\varphi(\bar e_{\ell +1},\bar b,\bar c):\ell< \omega\}$ infinite.
\end{enumerate}
\mn
Choose $\chi$ and define $\cN$ as in $(*)_4$ the 
proof of \ref{d10}.  For $C \subseteq \gC$ define a game $\Game_C$.  A
play last $\omega$ moves (really $n_* < \omega$ large enough
suffice).  In the $\ell$-th move the $\IND$ player chooses 
$\cX_\ell \in \cD_* \cap \Deef_{m(<n)}(M)$ and 
the antogonist chooses $\bar s_\ell \in
\mathop{\Pi}\limits_{m<n} I_m$ such that $\bar e_{\bar s_\ell} 
\in \cX_\ell$.  In
the end of the play the $\IND$ player wins the play when $\langle \bar
e_\ell:\ell < \omega\rangle$ is a $\Delta^1$-indiscernible sequence over $C$.

As in the proof of \ref{d10}, see $(*)_2$ there, 
the $\IND$ player has a winning strategy {\bf st}, and let
$\bold N \in \cN$ be such that {\bf st} $\in \bold N$ and choose $t_*
\in I_n$ such that $\cX \in \cD_n \cap \bold N \Rightarrow t_* \in \cX$,
possible as $\cD_n$ is $\kappa$-complete because $n \ge 1$.  Let $\ell
= \emptyset$ and we now simulate a play of 
$\Game_C$ called $\langle (\cX_\ell,\bar s_\ell):\ell < \omega \rangle$
such that:
\mn
\begin{enumerate}
\item[$(*)_2$]   $(a) \quad$ the $\IND$ player uses 
{\bf st} to choose $\cX_\ell$
\smallskip

\noindent
\item[${{}}$]   $(b) \quad$ the antagonist chooses 
$\bar s_\ell \in I$ and in $\cX_\ell$
such that if $\ell > 0$ and 

\hskip25pt  it is possible then $\tp_{\Delta^1}(\bar e^1_{\bar s_\ell} 
\char 94 \bar e^2_{t_*},C^+) \ne \tp_{\Delta^1}
(\bar e^1_{\bar s_{\ell-1}} \char 94 \bar e^2_{t_*},C)$.
\end{enumerate}

\noindent
It follows that $(\cX_\ell,\bar s_\ell) \in \bold N$ for $\ell < \omega$ 
and that for some $\ell(*) > 0$ the demand in clause (b) of $(*)_2$ 
is not possible.  So for some $q$
\medskip

\noindent
\begin{enumerate}
\item[$(*)_3$]  $(a) \quad \cX_{\ell(*)} \in D$
\smallskip

\noindent
\item[${{}}$]  $(b) \quad$ tp$_{\Delta^1}(\bar e^1_{\bar s} \char 94 \bar
e^2_{t_*},C) = q$ for every $\bar s \in \cX_{\ell(*)}$.
\end{enumerate}

\noindent
By the definition of $\cD_*$ and of the game there is $\langle \cS_\ell:\ell
< n\rangle \in \prod\limits_{\ell < n} \cD^+_\ell$ as there, such that
$\prod\limits_{\ell < n} \cS_\ell \subseteq \cX$ and \wilog \,
$\langle \cS_\ell:\ell < n \rangle \in \bold N$. 

We finish as in the proof of \ref{d10}.   
\end{PROOF}
\bigskip

\subsection{Density of $\tK$ in ZFC occurs} \

\begin{theorem}
\label{d33}
\underline{The universal solution theorem}   Assume $T$ is
countable, $(\kappa,\bar\mu,\theta)$ as usual, $\mu_0 \ge \beth_\omega$
and $\theta = \aleph_0$ and $\cf(\kappa) > 2^\theta$.

\noindent
1) If $\bold m_1 = (\bold x_1,\bar\psi_1,r_1) \in 
\rK^\oplus_{\kappa,\bar\mu,\theta}$ and 
$\bold x_1 \le_1 \bold y \in \qK_{\kappa,\bar\mu,\theta}$ 
\then \, we can find $\bold m_2$ such that 
$\bold m_1 \le^+_1 \bold m_2 \in 
\rK^\oplus_{\kappa,\bar\mu,\theta}$ and $\bold y \le_1 \bold x_{\bold m_2}$.

\noindent
2) Similarly but in the assumption $\bold y \in 
\uK_{\kappa,\bar\mu,\theta}$ and in the conclusion 
$\bold m_1 \le^\odot_1 \bold m_2$.
\end{theorem}

\begin{remark}
\label{d34} 
0) Note that this theorem restricts the cardinals lightly, but for
density of $\tK$ we shall have heavy restrictions, still ZFC ones.

\noindent
1) Part (2) is not needed for this subsection.

\noindent
2) If $M_{\bold x} \in \EC_{\kappa,\kappa}(T)$ the proof is
   somewhat easier, similarly in \ref{d10}.

\noindent
3) There is no real difference between the two parts.  We just deal
with the set of pairs $(\varphi,\psi)$ where $\varphi \in
\Gamma_{\bold x[\bold m_1]}$ and $\psi$ illuminate $(\bold
m,\varphi)$.

\noindent
4) In \ref{d33} we use $\iota(\bold x_1) = 2$ but with minor changes
 $\iota(\bold x_1) = 1$ is O.K., too; the changes are in $\oplus_5
   - \oplus_7$.

\noindent
5) Concerning \ref{d33}(2) see \ref{c2}(3)(e).
\end{remark}

Before proving note 
\begin{observation}
\label{d36}  
1) If $(\bold x_1,\bar\psi_1,r_1) \in \rK^\oplus_{\kappa,\bar\mu,\theta}$ and
$\bold x_1 \le_1 \bold y \in \pK_{\kappa,\bar\mu,\theta}$,
\then \, $(\bold x_1,\bar\psi_1,r_1) \le_1
(\bold y,\bar\psi_1,r_1) \in \rK^\oplus_{\kappa,\bar\mu,\theta}$.

\noindent
1A) If $\bold m_1 \in \rK^\oplus_{\kappa,\bar\mu,\theta}$ and  
$\cf(\kappa) > 2^\theta,\theta \ge |T|$ \then \, for some
$r_2 = r_2(\bar x_{\bar c_{\bold y}},\bar x_{\bar d_{\bold y}},\bar x'_{\bar
c_{\bold y}},\bar x'_{\bar d_{\bold y}})$ which extends $r_{\bold m}$ and 
is complete (over $\emptyset$) we have $\bold m \le_1 (\bold x_{\bold
m},\bar\psi_{\bold m},r_2) \in \rK^\oplus_{\kappa,\bar\mu,\theta}$.

\noindent
1B) If in (1A) in addition $\cf(\kappa) > 2^{|B_{\bold x}|}$ we can demand
$r_2$ is a complete type over $B_{\bold y}$; (similarly for
$B^+_{\bold x}$ when $|\bold S^\theta(B^+_{\bold x})| < \kappa$). 

\noindent
2) If $\bold m_1 \in \rK^\oplus_{\kappa,\bar\mu,\theta}$ 
and $\bold x_{\bold m_1} \le_1 \bold y \in 
\qK_{\kappa,\bar\mu,\theta}$ and $\cf(\kappa) > \theta \ge |T|$
\then \, 
\mn
\begin{enumerate}
\item[$(a)$]  for some pair $(\bar\psi,r)$ we have
$\bar\psi_{\bold m_1} \trianglelefteq \bar\psi$ and $\Gamma_{\bar\psi}
  = \Gamma_{\bold x[\bold m_1]} \cup \Gamma^1_{\bold x[\bold m]}$ and
  $(\bold y,\bar\psi,r_{\bold m_1}) \in
  \rK^\oplus_{\kappa,\bar\mu,\theta}$
\sn
\item[$(b)$]  similarly replacing $\Gamma^1_{\bold x[\bold m_1]}$ by
$\Gamma^3_{\bold x[\bold m_1]}$.
\end{enumerate}
\mn
2A) Like part (2) replacing
$\qK_{\kappa,\bar\mu,\theta},\qK^\oplus_{\kappa,\bar\mu,\theta}$ by
$\uK^\oplus_{\kappa,\bar\mu,\theta},\vK^\oplus_{\kappa,\bar\mu,\theta}$
respectively.

\noindent
3) Assume $\cD$ is a $\kappa$-complete filter on a set $I,\bar e_t \in
{}^\zeta{\gC}$ for $t \in I,\kappa = \cf(\kappa) >
2^{|B|+|\zeta|}$ and $\bold f$ is a $(\bar\mu,\theta)$-set 
and $\kappa > \cf(\Pi\{\kappa_{\bold f,i}:i \in u_{\bold x})\}$. 
\Then \, for some $q \in \bold S^\zeta(B_{\bold f,h}),\cS \in \cD^+$
and $h \in \Pi \ga_{\bold f}$ we have: $q = \tp
(\bar e_t,B^+_{\bold f,h})$ for every $t \in \cS$.

\noindent
4) Similarly for $\sK^\oplus_{\kappa,\bar\mu,\theta}$.
\end{observation}

\begin{remark}
\label{d37}
1) Variants: for $\bar\psi$ enough if $\cf(\kappa) > \theta$, see \ref{b47}.

\noindent
2)  Compare with \ref{b24} and \ref{c27}.
\end{remark}

\begin{PROOF}{\ref{d36}}
Straightforward, e.g. part (2) as in the proof of \ref{c78}(1).
\end{PROOF}

\begin{PROOF}{\ref{d33}}
\underline{Proof of \ref{d33}}   

1) Without loss of generality $\bold y$ is smooth and for notational
simplicity $\omega$ is disjoint to $v_{\bold y},w_{\bold y}$ and let
$\bold x = \bold x_1$
\mn
\begin{enumerate}
\item[$(*)_0$]  let $\bar v,\bar w,\bar u$ be such that
\sn
\begin{enumerate}
\item[$(a)$]  $\bar v = \langle v_n:n < \omega \rangle$ and $\bar w =
\langle w_n:n < \omega \rangle$
\sn
\item[$(b)$]  $v_n \subseteq v_{n+1} \subseteq v_{\bold y}$ and $w_n
\subseteq w_{n+1} \subseteq w_{\bold y}$
\sn
\item[$(c)$]  $v_n$ is finite and $w_n$ is finite
\sn
\item[$(d)$]  $v_{\bold y} = \cup\{v_n:n < \omega\}$ and $w_{\bold y}
= \cup\{w_n:n < \omega \}$
\sn
\item[$(e)$]  $u = v_{\bold y} \cup w_{\bold y} \cup \omega$
\sn
\item[$(f)$]  let $\bar u = \langle u_n:n < \omega \rangle$ where $u_n = v_n
\cup w_n \cup \{0,\dotsc,n-1\}$
\sn
\item[$(g)$]  let $u(n) = u_n,v(n) = v_n,w(n) = w_n$.
\end{enumerate}
\end{enumerate}
\mn
We choose $\langle(\Delta_n,\Delta^1_n,k_n,m_n):n < \omega\rangle$ such that
\medskip

\noindent
\begin{enumerate}
\item[$(*)_1$]    $(a) \quad m_n < \omega$ is increasing with $n$
\sn
\item[${{}}$]  $(b) \quad \Delta_n \subseteq \{\varphi(\bar x_{[u_n]},
\bar z_{[n]}):\varphi \in \bbL(\tau_T)\}$ is finite
\sn
\item[${{}}$]  $(c) \quad \Delta_n \subseteq \Delta_{n+1}$ in the
  natural sense, i.e. up to equivalence
\smallskip

\noindent
\item[${{}}$]  $(d) \quad \Delta_\omega = \cup\{\Delta_n:n < \omega\} =
\{\varphi(\bar x_{[u]},z_{[\omega]}):\varphi \in \bbL(\tau_T)\}$
up to equivalence
\smallskip

\noindent
\item[$(*)_2$]  $(a) \quad \Delta^1_n \subseteq \{\varphi(\bar
x_{[u_n]},\bar y_{[u_n]},\bar z_{[m_n]}):\varphi \in \bbL(\tau_T)\}$ is
finite
\sn
\item[${{}}$]  $(b) \quad \Delta^1_n \subseteq \Delta^1_{n+1}$
\sn
\item[${{}}$]  $(c) \quad \Delta^1_\omega = \cup\{\Delta^1_n:n <
\omega\} = \{\varphi(\bar x_{[u]},\bar y_{[u]},\bar z):
\bar z = \bar z_{[n]}$ for some

\hskip25pt  $n,\varphi \in \bbL(\tau_T)\}$ up to equivalence,
i.e. adding dummy variables
\sn
\item[${{}}$]  $(d) \quad$ if $\varphi(\bar x_{[u_n]},\bar z_{[n]})
\in \Delta_n$ then some $\varphi'(\bar x_{[u_n]},\bar y_{[u_n]},\bar
z_{[m_n]}) \in \Delta^1_n$ is 

\hskip25pt equivalent to $\varphi(\bar x_{[u_n]},\bar
z_{[n]})$ and some $\varphi''(\bar x_{[u_n]},\bar y_{[u_n]},\bar
z_{[m_n]})$ is 

\hskip25pt equivalent to $\varphi(\bar y_{[u_n]},\bar z_{[n]})$
\sn
\item[${{}}$]  $(e) \quad$ the finite $\Delta^1_n \subseteq
\Delta^1_\omega$ and $k_n < \omega$ are such that clause (i) of $\oplus$

\hskip25pt of \ref{d10} holds for $\Delta_n$, see the way we use this
proving $(*)_6$ below.
\end{enumerate}
\mn
This is possible as $T$ is countable, and for clause (e) of $(*)_2$ as
$T$ is dependent.
\medskip

\noindent
\begin{enumerate}
\item[$(*)_3$]   $I := ([M_{\bold x}]^{< \kappa},\subseteq)$ is
$\cf(\kappa)$-directed and let $\cD_I$ be the club filter on $I$.
\end{enumerate}

\noindent
Clearly 
\medskip

\noindent
\begin{enumerate}
\item[$(*)_4$]    $(a) \quad$ if $\alpha < \kappa \Rightarrow
|\alpha|^{\aleph_0} < \kappa$ then 
$\bold S^\omega_{\bbL(\tau_T)}(B^+_{\bold y} + \bar c_{\bold y} +
\bar d_{\bold y})$ has cardinality $< \kappa$
\smallskip

\noindent
\item[${{}}$]   $(b) \quad \bold S^{u(*)}_{\Delta_n}(B^+_{\bold
y,v_n} + \bar c_{\bold y} + \bar d_{\bold y})$ has cardinality $<
\kappa$ for any finite $u(*)$.
\end{enumerate}

\noindent
[Why?  Recall that we are assuming that $\bold y$ is smooth hence
$(B^+_{\bold y} + \bar c_{\bold y} + \bar d_{\bold
y},\bar{\bold I}_{\bold y})$ is a $(\bar\mu,\theta)$-set by
\ref{b45}(2); now for clause (a) apply \ref{b45}(1) and for clause
(b) apply \ref{b45}(1A).]

So together by observation \ref{d36}(2).
\medskip

\noindent
\begin{enumerate}
\item[$(*)_5$]   there are $\bar\psi_2$ and $r_2$ such that, (recall
  \ref{b16}(3),(3A)): 
\begin{enumerate}
\item[$(a)$]    $(\bold y,\bar\psi_2,r_2)$ belongs\footnote{moreover we
can demand it belongs to $\qK^\oplus_{\kappa,\bar\mu,\theta}$.}
to $\qK^\odot_{\kappa,\bar\mu,\theta}$ for part (1) and belongs to
  $\uK^\odot_{\kappa,\bar\mu,\theta}$ for part (2)
\sn
\item[$(b)$]   $\bar\psi_2 = \langle \psi_{2,\varphi}:\varphi \in
\Gamma_1\rangle$ where $\Gamma_1 = \Gamma^1_{\bold y}$ for part (1) and 
$\Gamma_1 \subseteq \Gamma^1_{\bold y}$ is $\bold y-\uK$-large 
for part (2), see Definition \ref{c23}(3B)
\smallskip

\noindent
\item[$(c)$]   $r_2 = r_2(\bar x_{\bar d_{\bold y}},
\bar x_{\bar c_{\bold y}},\bar x_{[\omega]})$ is a type over $\emptyset$ in 
$\gC$; may add $r$ is complete but should at least contain
$\psi_{2,\varphi}(\bar x_{\bar d[\bold y]},\bar x_{\bar c[\bold y]},\bar
x_{[\omega]})$ for $\varphi \in \Gamma_1$
\smallskip

\noindent
\item[$(d)$]   for every $A \in [M_{\bold x}]^{< \kappa}$ 
some $\bar e \in {}^\omega(M_{\bold x})$ solve $(\bold
y,\bar\psi_2,r_2,A)$, see Definition \ref{b16}(3) or \ref{c5}(2)(f).
\end{enumerate}
\end{enumerate}
\mn
For $A \in I$ first choose $(\bar d_A,\bar c_A)$ 
solving $(\bold y,\bar\psi,r,A)$
and second choose $\bar e_A$ as in $(*)_5(d)$ for $A + \bar d_1 + \bar
c_A$ and let $\bar e^+_A = \bar e_A \char 94 \bar d_A \char 94 \bar c_A$.

Next
\medskip

\noindent
\begin{enumerate}
\item[$(*)_6$]   there are $q^0_n,q^1_n,h_*$ and 
$\cS_{\ell,n}$ (for $\ell < 3,n < \omega$) such that:
\begin{enumerate}
\item[$(a)$]    $\cS_{\ell,n} \in \cD^+_I$
\smallskip

\noindent
\item[$(b)$]   $\cS_{\ell,n} \subseteq \cS_{\ell,m}$ when $n = m+1$
\smallskip

\noindent
\item[$(c)$]   $q^\ell_n \in \cS^{u(n)+u(n)}_{\Delta^1_n}
(B^+_{\bold y,u_n,h_{*,n}} + \bar d_{\bold y} + \bar c_{\bold y})$
\smallskip

\noindent
\item[$(d)$]   if $s_0 \in \cS_{0,n}$ then for the $(\cD_I +
\cS_{1,n})$-majority of $s_1 \in \cS_{1,n}$ (say for
every $s_1 \in \cS_{1,n,s_0})$ we have $q^0_n = 
\tp_{\Delta^1_n}((\bar e^+_{s_0} \rest m_n) \char 94 
(\bar e^+_{s_1} \rest m_n),B^+_{\bold y,v_n,h_*} + \bar
d_{\bold y} + \bar c_{\bold y})$ so $\cS_{1,n,s_0} \subseteq
\cS_{1,n}$ belongs to $\cD_I + \cS_{1,n}$
\smallskip

\noindent
\item[$(e)$]  if $s_1 \in \cS_{1,n}$ then for the $(\cD_I +
\cS_{2,n})$-majority of $s_2 \in \cS_{2,n}$ (say for
every $s_2 \in \cS_{2,n,s_1}$) we have $q^1_n = 
\tp_{\Delta^1_n}((\bar e^+_{s_1} \rest m_n) \char 94 (\bar e^+_{s_2} 
\rest m_n),B^+_{\bold y,v_n,h_*} + \bar d_{\bold y} + \bar c_{\bold y})$.
\end{enumerate}
\end{enumerate}

\noindent
[Why?  We do it by induction on $n$ replacing $h_*$ by $h_n$.  
For $n=0$, without loss of
generality $m_n =0,h_0$ constantly zero and 
we can let $\cS_{\ell,n} = I$ for $\ell=0,1,2$.  
For $n=m+1$ we do it in two steps.  First, 
letting $\bold f_n = (\langle
B_{\bold y,i}:i \in v_n \backslash u_{\bold y}\rangle,\langle \bold
I_{\bold y,i}:i \in v_n \cap u_{\bold y}\rangle)$ and applying
Theorem \ref{d10} for $\bold k = 1$ with $(\cD_I +
\cS_{\ell,m},\cS_{\ell,m},
\langle \bar e_A \rest m_n:A \in \cS_{\ell,m}\rangle,M_{\bold x},
\bold f_n,\bar c_{\bold y} + \bar d_{\bold y},\Delta_n,
\Delta^1_n,k_n\rangle)_{\ell < 2}$ here standing
for $(\cD_\ell,\cS_\ell,\langle \bar e^\ell_A:A \in I_\ell\rangle,M,
\bold f,C,\Delta_0,\Delta^1_0,k_0)_{\ell < 2}$ there.

We get $h^0_n,\cS_{0,n},\cS'_{1,m},q^0_n \in \bold S^{2m_n}_{\Delta_n}
(B^+_{\bold y,v_n,h^0_n} + \bar c_{\bold y} + \bar d_{\bold x})$.

Second, let $\cS'_{2,m} = \cS_{2,m}$ and we apply Theorem \ref{d10} for
$\bold k =1$ with $(\cD_I + \cS'_{\ell,m},\langle \bar e_A:A \in 
\cS'_{\ell,m}\rangle,M,\bold f_n,\bar c_{\bold y} + \bar d_{\bold
y},\Delta_n,\Delta^1_n,\bold k)_{\ell=1,2}$ here standing for
$(\cD_\ell,\cS_\ell,\langle \bar e^\ell_A:A \in I_\ell\rangle,M,\bold
f,C,\Delta_0,\Delta^1_0,k_0)_{\ell=0,1}$ there.  We get
$h^1_n,\cS_{1,n},\cS_{2,m},q^1_n \in \cS^{2m_n}_{\Delta_n}
(B^+_{\bold y,v_n,h^1_n} + \bar c_{\bold y} +
\bar d_{\bold y})$. 

Let $h_* = \sup\{h^\ell_n:n < \omega$ and $\ell=0,1\}$, 
i.e. $\kappa \in \ga_{\bold y}
\Rightarrow h_*(\kappa) = \sup\{h^\ell_n(\kappa):n < \omega$ and
$\ell=0,1\}$.  Now for $\ell=0,1$ and
$A_\ell \in \cS_{\ell,n}$ let $\cS_{\ell,n,A_\ell} 
= \{A_{\ell +1} \in \cS_{\ell +1,n}$:  we have $q^\ell_n =
\tp_{\Delta_n}((\bar e_{A_\ell} \rest m_n) \char 94 
(\bar e_{A_{\ell +1}} \rest m_n),B^+_{\bold y,v_n,h_*} +
\bar c_{\bold y} + \bar d_{\bold y})\}$.]
\medskip

We choose an ultrafilter ${\cD}$ on $I$ extending
${\cD}_I + \cS_{1,n}$ for every $n < \omega$ so clearly
$A \in I \Rightarrow \{B:A \subseteq B \in I\} \in {\cD}$.  
Let $\bar e_* \in {}^\omega {\gC}$ realizes 
$p_*(\bar y_{[\omega]}) := \text{ Av}({\cD},
\langle \bar e_A:A \in I\rangle,\bar d_{\bold y} + \bar c_{\bold y}
+ M_{\bold x})$, i.e. (clause (a) by the definition of Av 
and clause (b) follows)
\medskip

\noindent
\begin{enumerate}
\item[$(*)_7$]   $(a) \quad$ for $\bar b \in {}^{\omega >}
(M_{\bold x})$ we have
${\gC} \models \varphi[\bar d_{\bold y},\bar c_{\bold y},\bar e_*,\bar
b]$ iff 

\hskip25pt $\varphi(\bar d_{\bold y},\bold c_{\bold y},\bar y_{[\omega]},\bar
b) \in p_*(\bar y_{[\omega]})$ iff

\hskip25pt  $\{A \in I:{\gC} \models \varphi[\bar d_{\bold y},
\bar c_{\bold y},\bar e_A,\bar b]\} \in {\cD}$
\smallskip

\noindent
\item[${{}}$]   $(b) \quad \bar e_*$ exemplify $\bar\psi$,
i.e. is as in $(*)_5(d)$ except that it may be $\notin M_{\bold x}$
\sn
\item[$(*)_8$]  $(a) \quad$ let $\Gamma_1$ be $\Gamma^1_{\bold y}$
\sn
\item[${{}}$]  $(b) \quad$ we define $\bar\psi^*_2 = \langle 
\psi^*_{2,\varphi}: \varphi \in \Gamma_1 \rangle$ as follows: 
letting $\varphi =
\varphi(\bar x_{\bar d[\bold y]},\bar x_{\bar c[\bold y]},\bar z) \in$

\hskip25pt $\Gamma_3$, for some $n_0(\varphi),\varphi \equiv
 \varphi_0(\bar x_{\bar d[\bold y]},\bar x_{\bar c[\bold y]},\bar
 y_{[n_0(\varphi)]},\bar z)$ and let $\varphi_1 =$

\hskip25pt $\psi_{2,\varphi_0}(\bar x_{\bar d[\bold y]},
\bar x_{\bar c[\bold y]},\bar y)$ where $\bar\psi_2$ is from $(*)_5(a)$ 

\hskip25pt so really $\varphi_1 \equiv
\varphi_2(\bar x_{\bar d[\bold y]},\bar x_{\bar c[\bold y]},
\bar y_{[n_1(\varphi)]})$ for some $n_1(\varphi) \ge n_0(\varphi)$

\hskip25pt  and lastly let $\psi^*_{2,\varphi} = \psi_{\varphi_2}
(\bar x_{\bar d[\bold y]},\bar
x_{\bar c[\bold y]},\bar y_{[\omega]})$.
\sn
\item[$(*)_9$]  \wilog \, $\tp(\bar e_* \char 94 \bar d_{\bold y}
\char 94 \bar c_{\bold y},M_{\bold x})$ is recalling $u$ is from $(*)_0(e)$
\newline
 $p_{**}(\bar x_{[u]}) = 
\Av(\cD,\langle \bar e_A \char 94 \bar d_A \char 94 c_A:A \in
I\rangle,M_{\bold x})$.
\end{enumerate}
\mn
[Why?  By the choice of $\bar e_*$ it is enough to have $\tp(\bar
d_{\bold y} \char 94 \bar c_{\bold y},M_{\bold x}) = \Av(\cD,\langle
\bar d_A \char 94 \bar c_A:A \in I\rangle,M_{\bold x})$ which is easy
by $\cD \supseteq \cD_I$ and the choice of $(\bar d_A,\bar c_A)$ for $A \in I$
after $(*)_5$.]
\mn
\begin{enumerate}
\item[$(*)_{10}$]  let $e^+_* = \bar e_* \char 94 \bar d_{\bold y}
\char 94 \bar c_{\bold y}$.
\end{enumerate}
\mn
We now consider the statement
\mn
\begin{enumerate}
\item[$\boxtimes$]  for every $A_* \subseteq M_{\bold x}$ of
cardinality $< \kappa$, i.e. $A_* \in I$ there are 
$\bar e \in {}^\omega(M_{\bold x}),\bar d \in {}^{\ell g(\bar d[\bold
y])}(M_{\bold x})$ and $\bar c \in {}^{\ell g(\bar c[\bold x])}
(M_{\bold x})$ such that
\sn
\begin{enumerate}
\item[$(a)$]  $\bar e$ solves $(\bold y,\bar\psi^*,r_2,A_*)$
\sn
\item[$(b)$]  $(\bar d,\bar c)$ solves $(\bold m_1,A_*)$ so $\ell g(\bar c)
  = \ell g(\bar c_{\bold x}),\ell g(\bar d) = \ell g(\bar d_{\bold x})$
\sn
\item[$(c)$]  $\bar e \char 94 \bar d \char 94 \bar c$ realizes
$\tp(\bar e_* \char 94 \bar d_{\bold y} \char 94 \bar c_{\bold
x},A_*)$, but we do not say ``$\tp(\bar e_* \char 94 \bar d_{\bold y}
\char 94 \bar c_{\bold y},A_*)"$.
\end{enumerate}
\end{enumerate}
\mn
\underline{Why proving $\boxtimes$ is enough?}

We define $\bold x_2$ as $(M_{\bold x_1},\bar B_{\bold y},
\bar c_{\bold y},\bar d_{\bold x_1} \char 94 \bar e_*,
\bar{\bold I}_{\bold y})$; so clearly $\bold y \le_1 \bar{\bold x}_2 \in 
\pK_{\kappa,\bar\mu,\theta}$.  We choose 
\mn
\begin{enumerate}
\item[$(*)_{11}$]  we define $r_3$ by: $r_3 = r_1
\cup\{\varphi(\bar x_{\bar d[\bold y]},\bar x_{\bar c[\bold y]},\bar
x'_{\bar e_*}):\varphi(\bar x_{\bar d[\bold y]},
\bar x_{\bar c[\bar{\bold y}]},\bar x_{[\omega]}) \in r_2\}$ recalling
$r_1 = r[\bold m_1]$ and $r$ is from $(*)_5$.
\end{enumerate}
\mn
Lastly, let
\mn
\begin{enumerate}
\item[$(*)_{12}$]  $\bar\psi_3 = \langle \psi_{3,\varphi}:\varphi \in
\Gamma_{\bar\psi_3}\rangle$, where\footnote{if $I=\kappa,\langle
M_\alpha:\alpha < \kappa\rangle$ our problem will be to choose the $\bar
e'_\alpha$ such that $\langle \tp(\bar e'_\alpha,M_\alpha + \bar c_{\bold y}
+ \bar d_{\bold y}):\alpha < \kappa\rangle$ is increasing}
 $\Gamma_{\bar\psi_3} = \Gamma_{\bar\psi_1} \cup \Gamma_1$, 
see $(*)_5(b)$ 
(here it is convenient to allow repetitions of $\varphi$'s in $\bar\psi$;
for part (2) we have to change more) where:
\sn
\begin{enumerate}
\item[$(a)$]  $\psi_{3,\varphi} = 
\psi_{1,\varphi}$ if $\varphi \in \Gamma^2_{\bar\psi_1}$, adding dummy
variables recalling $\bar\psi_1 = \bar\psi_{\bold m_1}$

\hskip25pt so $\bar\psi_1 = \langle \psi_{1,\varphi}:\varphi \in
\Gamma^2_{\bar\psi_1}\rangle$ 
\sn
\item[$(b)$]  $\psi_{3,\varphi}$ is $\psi^*_{2,\varphi}$ if
  $\varphi \in \Gamma^1_{\bold y}$
using $w_{\bold y} + \omega$ instead $w_{\bold y}$,
i.e.

\hskip25pt  $\psi_{3,\varphi}(\bar x_{\bar d[\bold y]} \char 94 \bar
x_{[\omega]},\bar x_{\bar c[\bold y]},\bar x'_{\bar d[\bold y]} \char 94
\bar x'_{[\omega]},\bar x'_{\bar c[\bold y]})$ 

\hskip25pt $= \psi^*_{2,\varphi}(\bar x_{\bar d[\bold y]},
\bar x_{\bar c[\bold y]},\bar x'_{[\omega]})$.
\end{enumerate}
\end{enumerate}
\mn
We shall show that the triple $\bold m_2 = (\bold x_2,\bar\psi_3,r_3)$ 
is as required.  First, $\bold m_2 \in \rK^\oplus_{\kappa,\bar\mu,\theta}$:
 all the requirements are obviously satisfied, e.g. for
clause (f) of Definition \ref{c5}(1), given $A \in I$ we can choose 
$(\bar e,\bar d,\bar c)$ as in $\boxtimes$
so $\ell g(\bar c) = \ell g(\bar c_{\bold x_1}),\ell g(\bar d) = \ell g(\bar
d_{\bold x_1}) = \ell g(\bar d_{\bold y})$ and let $\bar c' \in
{}^{\ell g(\bar c[\bold y])}(M_{\bold x})$ be such that $\bar c =
\bar c' \rest \ell g(\bar c)$ and $\bar c' \char 94 \bar d \char
94 \bar e$ realizes $\tp(\bar c_{\bold y} \char 94 \bar d_{\bold y}
\char 94 \bar e_*,A)$; this is possible as $\bar d \char 94 \bar c$
realizes $\tp(\bar d_{\bold x} \char 94 \bar c_{\bold x},A)$ because
by $\boxtimes(b)$ the pair $(\bar d,\bar c)$ solves $\bold m_1$ and $\bar
d,\bar c$ are from $M_{\bold x}$.  We shall now show that $(\bar d
\char 94 \bar e,\bar c)$ solves $\bold m_2$.

We have to check clauses $(\alpha),(\beta),(\gamma)$ of Definition
\ref{c5}(1)(f).

By the choice of $\bar c',\bar c' \char 94 (\bar d \char 94 \bar e)$
realizes $\tp(\bar c_{\bold x_2} \char 94 \bar d_{\bold x_2},A) =
\tp(\bar c_{\bold y} \char 94 (\bar d_{\bold y} \char 94 \bar e_*),A)$
so clause $(\alpha)$ there holds.

Second, $\bar d_{\bold x_2} \char 94 \bar c_{\bold x} \char 94 \bar d
\char 94 \bar c = \bar d_{\bold y} \char 94 \bar c_{\bold x} \char 94
\bar d \char 94 \bar c$, realizes $r_1 = r_{\bold m_1}$ by clause (b)
of $\boxtimes$ and the definition; in addition $\bar d_{\bold y} \char 94 \bar
c_{\bold y} \char 94 \bar e$ realizes $r_2$ by clause (a) of
$\boxtimes$.  Together $(\bar d_{\bold y} \char 94 \bar e^+_*) \char 94 \bar
c_{\bold y} \char 94 (\bar d \char 94 \bar e) \char 94 \bar c'$
realizes $r_3$ by the choice of $r_2$ above and the previous sentence;
so clause $(\beta)$ there holds.

For clause $(\gamma)$ there, recalling that $\Gamma^2_{\bar\psi_3} =
\Gamma^2_{\bar\psi_1} \cup \Gamma^1_{\bold y}$ and we have to check a
condition for each $\varphi \in \Gamma^2_{\bar\psi_3}$.  Now if $\varphi
\in \Gamma^2_{\bar\psi_1}$ the desired conclusion holds by clause (b) of
$\boxtimes$ and the definition of ``solve".

If $\varphi \in \Gamma^1_{\bold y}$, use clause (a) of $\boxtimes$.

So we have proved indeed that $\bold m_2 = (\bold x_2,\bar\psi_3,r_3)
\in \rK^\oplus_{\kappa,\bar\mu,\theta}$.  In addition obviously
$\bold m_1 \le_1 \bold m_2$.  Lastly, $\bold m_1 \le^\odot_1 \bold
m_2$ by the choice of $\bar\psi_3 \rest \Gamma^1_{\bold y}$, so we are done.

So we are left with 

\noindent
\underline{proving $\boxtimes$ holds}:

Let $A_* \subseteq M_{\bold x}$ be
of cardinality $< \kappa$ and we shall show that there are sequences
$\bar e,\bar d,\bar c$ as required for $A_*$ in $\boxtimes$, this
suffices.  We can choose $\langle A_{*,n}:n < \omega\rangle$ 
and $A_{**}$ such that, \wilog \, 
\mn
\begin{enumerate}
\item[$\boxplus_0$]   $A_* \cup B^+_{\bold y} \subseteq A_{*,n} 
\in \cS_{0,n}$ and let $A_{**} = \cup\{\bar e^+_{A_{*,n}} + 
A_{*,n}:n < \omega\} \in [M]^{< \kappa}$.
\end{enumerate}

\noindent
Recalling $(*)_6$ let $\cS'_{1,n} = \cap\{\cS_{1,m,A_{*,m}}:m \le n\}$ but
$\cS'_{1,n} \subseteq \cS_{1,n,A_{*,m}} \in \cD_I +
\cS_{1,m} \subseteq \cD_I + \cS_{1,n}$ for $m \le n < \omega$ 
hence $\cS'_{1,n} \in \cD_I + \cS_{1,m} $.

Recalling $(*)_8$ let\footnote{we could use parameters just from 
$\bar d_{\bold y} + \bar c_{\bold y} + \Sigma_n \bar e_{A_{*,n}}$} 
$\Lambda = \{p:p$ a finite subset of $p_{**}(\bar y_{[u]})$
with parameters from $A_{**}\}$; so
clearly $|\Lambda| < \kappa$ and let $\Lambda_{\ge n} = \{p \in
\Lambda:|p| \ge n\}$.  By
the choice of $\bar e_*,{\cD}$ we can find $\langle A(p):p \in
\Lambda \rangle$ such that:
\medskip

\noindent
\begin{enumerate}
\item[$\boxplus_1$]   $A(p) \in \cS'_{1,n} \subseteq I$ 
and $\bar e_{A(p)}$ realizes $p$ and $A_{**} \subseteq A(p)$ 
\when \, $p \in \Lambda$ and $|p| = n$.
\end{enumerate}
\mn
For $n < \omega$ let $C_n$ be a member of $\cS_{2,n}$ which includes
$\cup \{\bar e^+_{A(p)} + A(p):p \in \Lambda\} 
\cup A_{**}$ such that $p \in \Lambda \wedge |p| =
n \Rightarrow C_n \in \cS_{2,n,A(p)}$, possible by
\ref{d10} the ``${\cD}_2$-almost", i.e. recalling $\cD_I$ is from $(*)_3$
as $\cS_{2,n,A(p)} \in
\cD_I + \cS_{2,n}$ by $(*)_6(e)$ and the choice of $\cD$.  
Let $\cD_*$ be an ultrafilter
on $\Lambda$ such that $p_1 \in \Lambda \Rightarrow \{p \in \Lambda:p_1
\subseteq p\} \in \cD_*$.  Let $q_*(\bar x_{[\omega]}) = 
\Av(\cD_*,\langle \bar e_{A(p)}:p \in \Lambda\rangle,
\cup\{C_n + \bar e_{C_n}:n <\omega\})$, it
is a type in $M_{\bold x}$ of cardinality $< \kappa$.  Recalling $(*)_8$ let
$q_{**}(\bar x_{[u]}) = 
q_{**}(\bar x_{[\omega]},\bar x_{\bar d[\bold y]},\bar x_{\bar
c[\bold y]}) = \Av(\cD_*,\langle \bar e_{A(p)} \char 94 \bar d_{A(p)}
\char 94 \bar c_{A(p)} = e^+_{A(p)}:p \in \Lambda\rangle,\cup\{C_n +
\bar e_{C_n} + \bar d_{C_n} + \bar c_{C_n}:n <\omega\})$, it is a type
in $M_{\bold x}$ of cardinality $< \kappa$ and it extends $q_*(\bar
x_{[\omega]})$.  Lastly, let $\bar c \in {}^{\ell g(\bar c_{\bold
y})}M_{\bold x},\bar d \in {}^{\ell g(\bar d_{\bold y})}(M_{\bold x})$ and
$\bar e \in {}^\omega(M_{\bold x})$ be such that the sequence $e^+ =
\bar e \char 94 \bar d \char 94 \bar c$ realizes $q_{**}(\bar
x_{[\omega]},\bar x_{\bar d[\bold y]},\bar x_{\bar c[\bold y]})$,
we shall prove that $\bar e,\bar d,\bar c$ are as required in 
$\boxtimes$; and let
$\bar e^+ = \bar e \char 94 \bar d \char 94 \bar c$.  That is, we have three
demands on $\bar e^+$, i.e. on $\bar e,\bar d,\bar c$ 
in $\boxtimes$, noting $\bar e^+ \in {}^u(M_{\bold x})$ recalling $u$
is from $(*)_0(e)$, in other words, $\bar e,\bar d,\bar c$ are
sequences of elements of $M_{\bold x}$ of the right lengths; let $u =
\ell g(\bar e^+) = \ell g(\bar e^+_*)$.

First, \underline{clause (c)} there 
says``realizing $\tp(\bar e^+_*,A_*) = \tp(\bar e_* \char 94 \bar
d_{\bold y} \char 94 \bold c_y,A_*)$" recalling $(*)_9$; 
we shall show that moreover 
$\bar e^+$ realizes $\tp(\bar e^+_*,A_{**})$.
This sufices as $A_* \subseteq A_{**}$.  Why ``$\bar e^+$ realizes
$\tp(\bar e^+_*,A_{**})$" holds?   Just recall $q_{**}(\bar y_{[u]}) 
= \tp(\bar e^+_*,M_{\bold x})$ by the choice of
$\bar e^+_*$ (see $(*)_7(a) + (*)_8 + (*)_9$), by 
the definition of $\Lambda$ and the
choice of $\cD_*$ and $\bar e^+$ above.  We shall deal with 
\underline{clause (b)}
in $\boxplus_5$ below and with \underline{clause (a)} in $\boxplus_7$ below.

Now
\medskip

\noindent
\begin{enumerate}
\item[$\boxplus_2$]  $\bar e^+_* \rest u_n,\bar e^+ \rest u_n$ 
and all $\bar e^+_A \rest u_n$ for $A \in \cS'_{1,n}$ realize the 
same complete $\Delta_n$-type (can
add same complete type) over $B^+_{\bold y,v_n,h_*}$.
\end{enumerate}
\medskip
\noindent
[Why?  First, recaling $h_*$ is from $(*)_6$
there is $p_n \in \bold S^{u(n)}_{\Delta_n}
(B^+_{\bold y,v_n,h_*})$ such
that $\bar e^+_A$ (equivalently $\bar e^+_A \rest u_n$),
 realizes $p_n$ when $A \in \cS'_{1,n}$ by
$(*)_6(d)$ recalling $(*)_2(d)$ and $\cS'_{1,n} \subseteq \cS_{1,n}$,
see after $\boxplus_0$.  For $\bar e^+_*$ by 
its choice in $(*)_7$, the choice of $\cD$ and the previous sentence;
lastly, for $\bar e^+$ it realizes $q_{**}(\bar x_{[u]})$ and recall that
$B^+_{\bold y} \cup A_* \subseteq C_n \subseteq \Dom(q_{**})$, the definition
of $q_{**}$ and the previous sentence.]
\medskip
\noindent
\begin{enumerate}
\item[$\boxplus_3$]   $\bar e^+_* \rest u_n,\bar e^+ \rest u_n$ and $\bar e^+_A
 \rest u_n$ (for $A \in \cS'_{1,n+1}$) all realize the 
same $\Delta_n$-type over $\bar c_{\bold y} + B^+_{\bold y,v_n,h_*}$.
\end{enumerate}
\medskip
\noindent
[Why?  Compared to $\boxplus_2$ we add $\bar c_{\bold y}$.  First,
  all the $\bar e^+_A$ for $A \in \cS'_{1,n}$ realizes the same
  $\Delta_n$-type over $\bar c_{\bold y} + B^+_{y,v_n,h_*}$ as this
  holds for $B^+_{\bold y,v_n,h_*}$ and recalling \ref{c34}(4) the type
$\tp(\bar c_{\bold y},M_{\bold x})$ is finitely 
satisfiable in $B^+_{\bold y,v_n,h_*}$ and all
  those sequences are from $M_{\bold x}$.  Second, the equality for 
$\bar e^+_*$ and $\bar e^+_A$'s as in the proof of
$\boxplus_2$.  Third, for $\bar e^+$ and the $\bar e^+_A$'s as 
$\tp(\bar c_{\bold y,v_n},M_{\bold x})$ is
finitely satisfiable in $B^+_{\bold y,v_n}$ when 
$\iota(\bold x) = 2$ (also it is locally does not split over
$B^+_{\bold y,v_n}$ but have to do more earlier if $\iota_{\bold x} =1$)
and those sequences are from $M_{\bold x}$, using $\boxplus_2$ of course.]
\mn
\begin{enumerate}
\item[$\boxplus_4$]  if $\varphi = \varphi(\bar x_{\bar d[\bold x]},\bar
x_{\bar c[\bold x]},\bar y'_{[u]})$ and $\Lambda_\varphi = \{p \in
\Lambda:\gC \models \varphi[\bar d_{\bold x},\bar c_{\bold x},\bar
e^+_{A(p)}]\}$ belongs to $\cD_*$ \then \, $\gC \models \varphi[\bar
d_{\bold x},\bar c_{\bold x},\bar e^+]$.
\end{enumerate}
\mn
This will take awhile.

Recalling $(*)_{11}$ note that
\mn
\begin{enumerate}
\item[$\oplus_{4.1}$]  $(a) \quad \varphi = \varphi(\bar x_{\bar
  d[\bold x]},\bar x_{\bar c[\bold x]},\bar y'_{[u]}) \in
  \Gamma^1_{\bold x}$
\sn
\item[${{}}$]  $(b) \quad \varphi_1 = \varphi_1(\bar x_{\bar d[\bold x]},
\bar x_{\bar c[\bold x]},\bar y_{[\omega]}) := \psi_{2,\varphi}$
\sn
\item[${{}}$]  $(c) \quad \varphi_2 = \varphi_2(\bar x_{\bar d[\bold x]},
\bar x_{\bar c[\bold x]},\bar y_{[\omega]}) = \psi_{2,\varphi_1}$.
\end{enumerate}
\mn
[Why?  For clause (a) just reflect, for clause (b) and (c) recall $(*)_5$.]

Now:
\mn
\begin{enumerate}
\item[$\oplus_{4.2}$]   $(a) \quad \gC \models ``\varphi_1[\bar
d_{\bold x},\bar c_{\bold x},\bar e_{C_n}]"$ for $n < \omega$
\smallskip

\noindent
\item[${{}}$]   $(b) \quad \varphi_1(\bar x_{\bar d[\bold x]},\bar c_{\bold x},
\bar e_{C_n}) \vdash \varphi(\bar x_{\bar d[\bold x]},\bar c_{\bold x},\bar
e^+_{A(p)})$ for $p \in \Lambda_\varphi,n < \omega$.
\end{enumerate}
\mn
[Why?  Clause (a) holds by the choice of $\varphi_1$ (in
$\oplus_{4.1}$), the choice of $\bar\psi_2$ (in $(*)_5$) and the choice
of the $\bar e_A$'s, in particular, $\bar e_{C_n}$.

Clause (b) holds by the choice of $\varphi_1$
the choice of $\bar\psi_2$ (in $(*)_5$) and the choice of $\bar e_{C_n}$
recalling $\bar e^+_{A(p)} \subseteq C_n$ for $p \in \Lambda$ by the
choice of $C_n$ after $\boxplus_1$.]

Hence letting $\vartheta_2 = \vartheta_2(\bar x_{\bar c[\bold x]},\bar
y''_{[\omega]},\bar y'_{[u]}) := (\forall \bar x_{\bar d[\bold x]})
(\varphi_1(\bar x_{\bar d[\bold x]},
\bar x_{\bar c[\bold x]},\bar y''_{[\omega]}) 
\rightarrow \varphi(\bar x_{\bar d[\bold x]},\bar x_{\bar c[\bold x]},
\bar y'_{[u]}))$, we have
\mn
\begin{enumerate}
\item[$\oplus_{4.3}$]   $(a) \quad \gC \models ``\varphi_1[\bar
d_{\bold x},\bar c_{\bold x},\bar e_{C_n}]"$ for $n < \omega$
\sn
\item[${{}}$]   $(b) \quad \gC \models ``\vartheta_2[\bar c_{\bold x},
\bar e_{C_n},\bar e^+_{A(p)}]$" for $p \in \Lambda_\varphi,n < \omega$.
\end{enumerate}
\mn
[Why?  As $\bar e^+_{A(p)} \subseteq C_n$ for $p \in \Lambda,n < \omega$.]

By \ref{d10}, that is by $(*)_6(d)$, 
recalling $A(p) \in \cS'_{1,|p|} \subseteq \cS_{1,|p|,A}$ from
$\oplus_{4.3}(b)$ as $\tp(\bar c_{\bold x},M_{\bold x})$ is finitely
satisfiable in $B^+_{\bold x}$ it follows that for some $n_1 < \omega$ 
(recalling $\Delta^1_\omega = \cup\{\Delta^1_n:
n < \omega\},\Delta^1_n \subseteq \Delta^1_{n+1}$ by $(*)_1$) we have
\mn
\begin{enumerate}
\item[$\oplus_{4.4}$]   if $\bar c \in {}^{\ell g(\bar c[\bold x])}
(B^+_{\bold x})$ and $p,q \in \Lambda_{\ge n_1}$ and $n < \omega$
\then \, $\gC \models ``\vartheta_2[\bar c',\bar e_{C_n},\bar
e^+_{A(q)}] \equiv \vartheta_2[\bar c',\bar e^+_{C_n},\bar
e^+_{A(p)}]"$.
\end{enumerate}
\mn
Hence by the choice of $\bar e^+$ (after $\boxplus_1$)
\mn
\begin{enumerate}
\item[$\oplus_{4.5}$]   if $\bar c' \in {}^{\ell g(\bar c[\bold x])}
(B^+_{\bold x})$ and $p \in \Lambda_{\ge n_1}$ and $n < \omega$
\then \, $\gC \models ``\vartheta_2[\bar c',\bar e_{C_n},\bar e^+] 
\equiv \vartheta_2[\bar c',\bar e_{C_n},\bar e^+_{A(p)}]"$.
\end{enumerate}
\mn
As $\tp(\bar c_{\bold x},M_{\bold x})$ is finitely satisfiable in
$B^+_{\bold x} \subseteq A_{**}$, clearly
\mn
\begin{enumerate}
\item[$\oplus_{4.6}$]   if $p \in \Lambda_{\ge n_1}$ and $n < \omega$
\then \, $\gC \models ``\vartheta_2[\bar c_{\bold x},\bar e_{C_n},\bar e^+] 
\equiv \vartheta_2[\bar c_{\bold x},\bar e_{C_n},\bar e^+_{A(p)}]"$.
\end{enumerate}
\mn
By $\oplus_{4.6}$ and $\oplus_{4.3}(b)$ we get
\mn
\begin{enumerate}
\item[$\oplus_{4.7}$]   $\gC \models ``\vartheta_2[\bar c_{\bold x},
\bar e_{C_n},\bar e^+]"$
\end{enumerate}
\mn
which means
\mn
\begin{enumerate}
\item[$\oplus_{4.8}$]   $\varphi_1(\bar x_{\bar d[\bold x]},\bar c_{\bold x},
\bar e_{C_n}) \vdash \varphi(\bar x_{\bar d[\bold x]},\bar c_{\bold x},\bar e^+)$.
\end{enumerate}
\mn
By $\oplus_{4.3}(a)$ we have
\mn
\begin{enumerate}
\item[$\oplus_{4.9}$]  $\gC \models \varphi_1[\bar d_{\bold x},\bar
c_{\bold x},\bar e_{C_n}]$ for $n < \omega$.
\end{enumerate}
\mn
By $\oplus_{4.8} + \oplus_{4.9}$ we have
\mn
\begin{enumerate}
\item[$\oplus_{4.10}$]   $\gC \models ``\varphi[\bar d_{\bold x},
\bar c_{\bold x},\bar e^+]"$
\end{enumerate}
\mn
So $\boxplus_4$ has been proved indeed.
\mn
\begin{enumerate}
\item[$\boxplus_5$]  clause (b) of $\boxtimes$ holds (for our choice
of $\bar e^+$).
\end{enumerate}
\mn
[Why?  We have to check clauses $(\alpha),(\beta),(\gamma)$ Definition
  \ref{c5}(1)(f).  For every $A \in I$ the 
pair $(\bar d_A,\bar c_A)$ solve $(\bold m_1,A_*)$ 
hence $\bar d_{\bold x} \char 94 \bar c_{\bold x} \char 94
\bar d_A \char 94 \bar c_A$ realizes $r_1$ so recalling $e^+_A = \bar e_A
\char 94 \bar d_A \char 94 \bar c_A$ by $\boxplus_4$ also
$(\bar d_{\bold x},\bar c_{\bold x},\bar d,\bar c)$ realizes $r_1$ so 
clause $(\beta)$ there holds.
Also as in the proof of $\boxplus_2$ just easier, $\bar d \char 94
\bar c$ realizes $\tp(\bar d_{\bold x} \char 94 \bar c_{\bold x},
A_{**})$ and is from $M_{\bold x}$ by the choices of $\bar e,\bar
d,\bar c$ after $\boxplus_1$, so clause $(\alpha)$ there holds.  
As in the proof of $\boxplus_4$ easily $\bar d \char 94
\bar c$ and $\bar d_{A(p)} \char 94 \bar c_{A(p)}$ realize the same
type over $\bar c_{\bold x} + A_{**}$ for $p \in \Lambda$.

Let $\varphi = \varphi(\bar x_{\bar d[\bold x]},\bar x_{\bar c[\bold x]},
\bar x'_{\bar d[\bold x]},\bar x'_{\bar c[\bold x]},\bar y) \in 
\Gamma^2_{\bar\psi_1}$ and $\bar
b \in {}^{\ell g(\bar y)}(A_*)$ be such that $\models \varphi[\bar
  d_{\bold x},\bar c_{\bold x},\bar d,\bar e,\bar b]$ so $\psi_{1,\varphi} =
\psi_{1,\varphi}(\bar x_{\bar d[\bold x]},\bar x_{\bar c[\bold x]},
\bar x'_{\bar d[\bold x]},\bar x'_{\bar c[\bold x]}) \in r_1$ and
$\psi_{1,\varphi}(\bar x_{\bar d[\bold x]},\bar c_{\bold x},\bar
d_{A(p)},\bar c_{A(p)}) \vdash \varphi(\bar x_{\bar d[\bold x]},
\bar c_{\bold x},\bar d_{A(p)},\bar c_{A(p)},\bar b)$.

As $\psi_{1,\varphi} \in r_1$ necessarily $\models
``\psi_{1,\varphi}[\bar d_{\bold x},\bar c_{\bold x},\bar d_{A(p)},\bar
  c_{A(p)}] \wedge \psi_{1,\varphi}[\bar d_{\bold x},\bar c_{\bold
    x},\bar d,\bar c]"$ and as $\bar d_{A(p)} \char 94 \bar
c_{A(p)},\bar d \char 94 \bar c$ realize the same type over $A_*$ and
as $\bar b \subseteq A_*$ and the end of the previous sentence,
$\psi_{1,\varphi}(\bar x_{\bar d[\bold x]},\bar c_{\bold x},\bar
d,\bar c) \vdash \varphi(\bar x_{\bar d[\bold x]},\bar c_{\bold
  x},\bar d,\bar c,\bar b)$.

The last sentence says that clause $(\gamma)$ from \ref{c5}(1)(f)
holds.  Together indeed we have proved that $\bar d,\bar c$ satisfies clause
(b) of $\boxtimes$, i.e. $(\bar d,\bar c)$ solves 
$(\bar{\bold m}_1,A_{**})$ as promised.]

We are left with clause (a) of $\boxtimes$. For the rest of the proof
let $\bar x_{\bar d} = \bar x_{\bar d[\bold y]},\bar x_{\bar c} = 
\bar x_{\bar c[\bold y]}$.
\mn
\begin{enumerate}
\item[$\boxplus_6$]  for $\varphi = \varphi(\bar x_{\bar d},\bar
x_{\bar c},\bar z)$ let $\varphi_0 = 
\varphi,\varphi_1 = \varphi_1(\bar x_{\bar d},\bar x_{\bar c},\bar
y_{[\omega]}) = \psi_{2,\varphi}$ where $\bar\psi_2$ is from $(*)_5$ 
and $\varphi_2 = \varphi_2(\bar x_{\bar
d},\bar x_{\bar c},\bar y_{[\omega]}) := \psi_{2,\varphi_1}$ so
$\varphi_2 = \psi^*_{2,\varphi}$, see $(*)_{11}$.
\end{enumerate}

Now to finish the proof of $\boxtimes(a)$ hence of the theorem, it
suffices to show:
\medskip

\noindent
\begin{enumerate}
\item[$\boxplus_7$]   $\gC \models \varphi_2[\bar d_{\bold y},\bar
c_{\bold y},\bar e]$ and $\varphi_2(x_{\bar d},\bar c_{\bold
y},\bar e) \vdash \varphi(\bar x_{\bar d},\bar c_{\bold y},\bar b)$
\underline{when} (if we add $\bar e$ in $\varphi$ we have a problem in
$\oplus_{7.10}$ as the $\bar e$ is changed):
\begin{enumerate}
\item[$(a)$]   $\varphi = \varphi(\bar x_{\bar d},\bar x_{\bar c};\bar
z) \in \Gamma_1$
\smallskip

\noindent
\item[$(b)$]  ${\gC} \models \varphi[\bar d_{\bold y},\bar c_{\bold y},
\bar b]$ and $\bar b \in {}^{\ell g(\bar z)}(A_*)$.
\end{enumerate}
\end{enumerate}
\medskip

\noindent
Why?  So assume clauses (a),(b) of $\boxplus_7$ and eventually 
we shall prove the desired conclusions of
$\boxplus_7$.  The first part, $\gC \models \varphi_2[\bar
d_{\bold y},\bar c_{\bold y},\bar e]$ holds by $\boxplus_4$ by
$\oplus_{7.2}(a)$ and the second part by $\oplus_{7.10}$ below. 

Recalling the choice of $\bar\psi_2$ in $(*)_5$ recalling $\varphi_0 =
\varphi$ and the formula 
$\varphi_1(\bar x_{\bar d},\bar x_{\bar c},\bar y_{[\omega]})$ is
equal to $\psi_\varphi$, clearly we have: letting $I_{\bar b}
:= \{A \in I:\bar b \in {}^{\ell g(\bar z)} A\}$ and $J = \{\bar b \in
{}^{\ell g(\bar b)}(M_{\bold x}): \models \varphi_0[\bar d,\bar c,\bar b]\}$
\medskip

\noindent
\begin{enumerate}
\item[$\oplus_{7.1}$]   $(a) \quad \gC \models \varphi_1[\bar d_{\bold y},
\bar c_{\bold y},\bar e_A]$ when $A \in I$
\smallskip

\noindent
\item[${{}}$]  $(b) \quad \varphi_1(\bar x_{\bar d},\bar c_{\bold y},
\bar e_A) \vdash \varphi_0(\bar x_{\bar d},\bar c_{\bold y},\bar b)$
when $A \in I_{\bar b} \cap J$
\sn
\item[${{}}$]  $(c) \quad \varphi_1(\bar x_{\bar d},\bar c_{\bold y},
\bar e_*) \vdash \varphi_0(\bar x_{\bar d},\bar c_{\bold y},\bar b)$
when $A \in J$.
\end{enumerate}

\noindent
Hence by $(*)_7(b)$, recalling $\varphi_2 = \psi_{2,\varphi_1}$ 
\medskip

\noindent
\begin{enumerate}
\item[$\oplus_{7.2}$]   $(a) \quad \gC \models \varphi_2[\bar d_{\bold y},
\bar c_{\bold y},\bar e_A]$ for $A \in I$ hence
$\gC \models \varphi_2[\bar d_{\bold y},\bar c_{\bold y},
\bar e_*]$ so 

\hskip25pt $\{\varphi_2(\bar d_{\bold y},
\bar c_{\bold y},\bar y_{[\omega]})\} \in \Lambda$ pedentically
$\{\varphi'_2(\bar d_{\bold y},\bar c_{\bold y},\bar y_{[u]})\} \in
\Lambda$ where

\hskip25pt  $\varphi'_2(\bar x_{\bar d},x_{\bar c},y_{[u]}) =
\varphi_2(\bar x_{\bar d},x_{\bar c},\bar y_{[u]})$
\smallskip

\noindent
\item[${{}}$]  $(b) \quad \varphi_2(\bar x_{\bar d},\bar c_{\bold y},
\bar e_*) \vdash \varphi_1(\bar x_{\bar d},\bar c_{\bold y},\bar
e_{A_{*,n}})$ for $n < \omega$ 
\end{enumerate}
\mn
moreover
\mn
\begin{enumerate}
\item[$\oplus_{7.3}$]   $(a) \quad$ if $p \in \Lambda$ satisfies
$\varphi_2(\bar d,\bar c,\bar y_{[\omega]}) \in p$ then
$\gC \models \varphi_2[\bar d,\bar c,\bar e_{A(p)}]$
\sn
\item[${{}}$]  $(b) \quad$ if $A \in I$ then 
$\varphi_2(\bar x_{\bar d},\bar c,\bar e_A) \vdash \{\varphi_1
(\bar x_{\bar d},\bar c,\bar e'):\gC \models \varphi_1[\bar d,\bar c,
\bar e']$ and

\hskip25pt  $\bar e' \in {}^\omega A\}$
\sn
\item[${{}}$]  $(c) \quad$ like (a) replacing $\bar e_{A(p)}$ by $\bar e_*$.
\end{enumerate}
\mn
So letting $\Lambda_* = \{p \in \Lambda:\varphi_2(\bar d_{\bold
y},\bar c_{\bold y},\bar y_{[\omega]}) \in p\}$ we have $\Lambda_* \in
\cD_*$ and let

\[
\vartheta_1 = \vartheta_1(\bar x_{\bar c},\bar y'_{[\omega]},
\bar y_{[\omega]}) := (\forall \bar x_{\bar d})(\varphi_2
(\bar x_{\bar d},\bar x_{\bar c},\bar y'_{[\omega]}) 
\rightarrow \varphi_1(\bar x_{\bar d},\bar x_{\bar c},\bar y_{[\omega]})).
\]

\mn
So we have, by $\oplus_{7.1}(b)$ and $\oplus_{7.3}(c)$ as $\bar
e_{A_*,n} \subseteq A_{**} \subseteq A(p)$ by $\boxplus_1$
\medskip

\noindent
\begin{enumerate}
\item[$\oplus_{7.4}$]   $\gC \models \vartheta_1[\bar c_{\bold y},
\bar e_{A(p)},\bar e_{A_*,n}]$ for $p \in \Lambda_*,n < \omega$.
\end{enumerate}

\noindent
But by \ref{d10}, that is, by $(*)_6(d)$ recalling $\bar e_{A(p)} \in
\cS'_{1,|p|} \subseteq \cS_{1,|p|,A_*}$ it follows that
for some $n_1 < \omega$ (as $\cup\{\Delta^1_n:n < \omega\} = \Delta^1_\omega$)
we have:
\mn
\begin{enumerate}
\item[$\oplus_{7.5}$]   if $\bar c' \in {}^{\ell g(\bar c[\bold y])}
(B^+_{\bold x})$ and $p,q \in \Lambda_{\ge n_1}$ \then \, 

$\gC \models ``\vartheta_1[\bar c',\bar e_{A(p)},
\bar e_{A_*,n}] \equiv \vartheta_1[\bar c',\bar e_{A(q)},\bar
e_{A_*,n}]$" for $n < \omega$.
\end{enumerate}

\noindent
Hence by the choice of $\bar e$
\medskip

\noindent
\begin{enumerate}
\item[$\oplus_{7.6}$]   if $\bar c' \in {}^{\ell g(\bar c[\bold
y])}(B^+_{\bold x})$ and $p \in \Lambda_{\ge n_1}$ and $n < \omega$,
\then

$\gC \models \vartheta_1[\bar c',\bar e,\bar e_{A_{*,n}}] \equiv 
\vartheta_1[\bar c',\bar e_{A(p)},\bar e_{A_*,n}]$.
\end{enumerate}

\noindent
As $\tp(\bar c_{\bold y},M_{\bold x})$ is finitely satisfiable in 
$B^+_{\bold x,h_*}$, clearly
\medskip

\noindent
\begin{enumerate}
\item[$\oplus_{7.7}$]   if $p \in \Lambda_{\ge n_1}$ and $n < \omega$ 
then $\gC \models ``\vartheta_1[\bar c_{\bold y},\bar e,
\bar e_{A_*,n}] \equiv \vartheta_1[\bar c_{\bold y},\bar e_{A(p)},\bar
e_{A_*,n}]"$.
\end{enumerate}
\mn
Next
\mn
\begin{enumerate}
\item[$\oplus_{7.8}$]   if $n < \omega$ 
then $\gC \models \vartheta_1[\bar c_{\bold y},\bar e,\bar e_{A_{*,n}}]$.
\end{enumerate}
\mn
[Why?  By $\oplus_{7.4} + \oplus_{7.7}$ because there is $p \in \Lambda_{\ge
n_1} \cap \Lambda_*$ which holds as $\Lambda_{\ge n_1} \in \cD_*$ and
$\Lambda_* \in \cD_*$ and $\cD_*$ is an ultrafilter on $\Lambda$.]

So by the choice of $\vartheta_1$
\medskip

\noindent
\begin{enumerate}
\item[$\oplus_{7.9}$]   $\varphi_2(\bar x_{\bar d},\bar c_{\bold y},
\bar e) \vdash \varphi_1(\bar x_{\bar d},\bar c_{\bold y},\bar e_{A_*,n})$.
\end{enumerate}

\noindent
By $\oplus_{7.9} + \oplus_{7.1}(b)$ applied to $A = A_{*,n}$ we have (recall
$\bar b$ is from $\boxplus_7$ hence $\bar b \subseteq A_* \subseteq A_{**}$)
\mn
\begin{enumerate}
\item[$\oplus_{7.10}$]   $\varphi_2(\bar x_{\bar d},\bar c_{\bold y},
\bar e) \vdash \varphi_0(\bar x_{\bar d},\bar c_{\bold y},\bar b)$.
\end{enumerate}
\mn
This proves the second clause in the
desired conclusion of $\boxplus_7$.

So we are done proving $\boxplus_7$.

As said above (before $\boxplus_7$) proving $\boxplus_7$ finish the proof.

\noindent
2) We repeat the proof above with some changes.  In $(*)_5(a)$ we
replace $\qK^\odot_{\kappa,\bar\mu,\theta}$
by $\uK^\odot_{\kappa,\bar\mu,\theta}$ respectively.  We change $(*)_8
+ (*)_9$ naturally and also the rest should be clear.
\end{PROOF}

Now we get a ``density of $\tK_{\kappa,\mu,\theta}$ in ZFC" for
$\theta = \aleph_0$ and some pairs $\kappa,\mu$.
\begin{conclusion}
\label{d39}  
If $T$ is countable, $\theta = \aleph_0,\mu$ is strong limit and 
$(\mu > \cf(\mu) \ge \aleph_1 \wedge \kappa = \mu^+)$ 
or ($\mu = \cf(\mu) = \kappa$), \then \, for
every $\bold m \in \rK^\oplus_{\kappa,\mu,\theta}$ there is
$\bold n \in \tK^\oplus_{\kappa,\mu,\theta}$ such that
$\bold m \le_1 \bold n$.
\end{conclusion}

\begin{remark}
1) Do we need $\cf(\mu) > 2^\theta$?  No, see \ref{b22}(1) and
   \ref{b30} but ``$\mu$ is strong limit" is assumed above.

\noindent
2) Recall that $\rK^\oplus_{\kappa,\mu,\theta}$ means
   $\rK^\oplus_{\kappa,\bar\mu,\theta}$ with $\bar\mu =
   (\mu_2,\mu_1,\mu_0) = (\kappa,\mu,\mu)$. 

\noindent
3) This is enough for the recounting of types for $\kappa$ strongly
inaccessible.  Also for $\kappa = \mu^+,\mu$ strong limit singular
 of uncountable cofinality, but only if $\mu = \aleph_\mu$ we can deduce the
 correct upper bound on the number of types up to conjugacy in $\bold
 S(M),M \in \EC_{\lambda,\lambda}(T)$, still if $\mu < \aleph_\mu$ the
 upper bound is $\mu$, smaller than the value for independent $T$.
\end{remark}

\begin{PROOF}{\ref{d39}}
 We choose $\bold m_n$ by induction on $n < \omega$ such that
\medskip

\noindent
\begin{enumerate}
\item[$\boxplus_1$]   $(a) \quad \bold m_n \in \rK_{\kappa,\mu,\theta}$
\smallskip

\noindent
\item[${{}}$]   $(b) \quad \bold m_0 = \bold m$
\smallskip

\noindent
\item[${{}}$]   $(c) \quad r[\bold m_{n+1}]$ is complete
\smallskip

\noindent
\item[${{}}$]  $(d) \quad \bold m_m \le^+_1 \bold m_n$ when $n = m+1$.
\end{enumerate}

\noindent
Why can we carry the induction?

For $n=0$, by clause (b) this is trivial.

For $n=m+1$ by \ref{b22}(1) there is $\bold y_m$ such that 
$\bold x_{\bold m_m} \le_2 \bold y_m \in \qK'_{\kappa,\mu,\theta}$
hence to $\qK'_{\mu,\mu,\theta}$ hence by \ref{b22}(3) using our use
of $(\mu,\mu,\theta)$ rather than $(\kappa,\mu,\theta)$ we have $\bold
y_{\bold m} \in \qK_{\mu,\mu,\theta}$;
as $\cf(\mu) > \aleph_0$ by \ref{b24}(1) for some $\bar\psi_m$, we have
$(\bold y_m,\bar\psi_m,\emptyset) \in \qK^\odot_{\mu,\mu,\theta}$.  Hence
$\bold y_m \in \qK_{\kappa,\mu,\theta}$, why?  
If $\kappa = \mu$ trivially
and if $\kappa = \mu^+$ by \ref{b24}(2), which is O.K. by
the assumptions $\theta < \cf(\mu) < \mu$.  As $\cf(\mu) > \theta$ we have
$|B^+_{\bold y}| < \mu$ and as $\mu$ is strong limit we have
$\beth_k(|B_{\bold f}| + \theta) < \mu$ for $k < \omega$, 
so we can apply \ref{d33}.

By \ref{d33} there is $\bold n_m \in \rK^\oplus_{\kappa,\mu,\theta}$ such
that $\bold m_m \le^+_1 \bold n_m$ such that $\bold y_m \le_1 \bold
x_{\bold n_m}$.

Lastly, as $\kappa$ is regular $> 2^\theta$, by 
Observation \ref{d36}(1A) there is a complete $r_n \supseteq
r[\bold n_n]$ such that $\bold m_n := (\bold x_{\bold
n_m},\bar\psi_{\bold n_m},r_n) \in \rK^\oplus_{\kappa,\mu,\theta}$ 
is $\le_1$-above $\bold n_m$.  So $\bold m_n$ is as required.

Now $\bold n = \text{ lim}\langle \bold m_n:n < \omega\rangle$ is as
required by \ref{c70}.  
\end{PROOF}

\begin{remark}
\label{d42}
Also \ref{d39} is enough for the ``generic pair conjecture" for the
relevant cardinals.
\end{remark}

\newpage

\section {Stronger Density}
\bigskip

\subsection{More density of $\tK$} \

The following will help us to prove density of 
$\tK_{\kappa,\mu,\theta}$ replacing $\kappa$ by $\kappa^{+n}$ in
\ref{d39}.  Unfortunately, we are stuck in $\kappa^{+\omega}$, still
this gives more cases for the recounting of types. 

\begin{claim}
\label{p2}
\underline{Crucial Claim}   
There is an indiscernible sequence $\bold I = \langle \bar
 a_\alpha:\alpha < \lambda\rangle$ in $M_{\bold x}$ such that letting
 $\bar a_\lambda$ realizes $\Av(\bold I,M_{\bold x} + \bar c_{\bold x})$,
 the types of $\bar a_\lambda$ and $\bar d_{\bold x}$
 over $M_{\bold x} + \bar c_{\bold x}$ are not weakly orthogonal \when
\,:
\mn
\begin{enumerate}
\item[$(a)$]  $\kappa$ is regular
\sn
\item[$(b)$]   $\bold x \in \pK_{\kappa,\bar\mu,\theta}$ and 
$\iota_{\bold x} = 2$ 
\sn
\item[$(c)$]   $\lambda = \ntr_{\lc}(\bold x)$ is regular, see
  Definition \ref{b43}(2)
\sn
\item[$(d)$]   $(\alpha) \quad u_{\bold x}$ is finite 

\underline{or} just
\sn
\item[${{}}$]  $(\beta) \quad \alpha < \lambda
\Rightarrow |\alpha|^{|T|} < \lambda$
\sn
\item[$(e)$]   $\kappa > \lambda  \ge \mu_1 > |B_{\bold x}|$
\sn
\item[$(f)$]   $\lambda > 2^{|B_{\bold x}|+\theta}$ and 
$\kappa \ge \beth_\omega(|B_{\bold x}| + \theta)$.
\end{enumerate}
\end{claim}

\begin{discussion}
\label{p4}  
1) Recall $\ntr(\bold x)$ is regular or is $\le \theta$, see
 Definition \ref{b43}, Observation \ref{b47}(1) but this is not necessarily so
 for $\ntr_{\lc}(\bold x)$, on it we know only that 
it is regular or its cofinality is $\le \theta$.

\noindent
2) Why above ``$u_{\bold x}$ is finite"?  Otherwise in \ref{p2}
there is a problem.  The reason is a pcf one: maybe
$\lambda \in \pcf({\ga}_{\bold x,< \lambda})$ where we let
${\ga}_{\bold x,<\lambda} = \{\kappa_{\bold x,i}:i \in u_{\bold x}$ and
$\kappa_{\bold x,i} < \lambda\}$, even the case $\lambda \in \{\kappa_{\bold
x,i}:i \in u_{\bold x}\}$ need care.

Even under G.C.H., if $\lambda = \chi^+,\cf(\chi) \le \theta$ we
have a problem.  The problem is in fixing the ``essential" type of $\bar
e_\alpha$ for $\alpha \in [\alpha_{\varepsilon},\alpha_{\varepsilon +1})$
over $\bold x$; which has more information than its type over
$B_{\bold x}$  but less than its type over 
$B^+_{\bold x}$ and is preserved if we
replace $\bold x$ by a very similar $\bold x'$, we can use just $\bold
x_{[h]}$ which is smooth see Definition
\ref{b27}(1),(2) and \ref{b30} and \ref{b35}.

The first idea for saving the day was to get $\langle \bar
e^*_\varepsilon:\varepsilon < \lambda\rangle$ tree indiscernible for
some $\bar g = \langle g_\alpha \in \Pi{\ga}_{\bold x,\lambda}:\alpha
< \lambda\rangle <_{J_{< \lambda}
[{\ga}_{\bold x,\lambda}]}$-increasing 
and cofinal which is ``nice".  Did not seem to work.

Second is a weaker version: demand something on $\langle \bar
e^*_{\varepsilon_0},\dotsc,\bar e^*_{\varepsilon_{n-1}}\rangle$ only
when $g_{\varepsilon_0} < g_{\varepsilon_1} < \ldots$.

The second is not good enough to classify $f^{\aut}_{T,\theta}(-)$.  
Still, when $\kappa = \mu^{+n},\mu$
regular, $\mu = \mu^\theta$ this may help but we prefer not too, when
we can.

The solution is to do it locally, i.e. to deal with local density for
$\qK$ (in $\pK$),
deal with one $\varphi$, then pretend you have no $\varphi$ and deal
with the case $u_{\bold x}$ is finite, i.e. \ref{d10} whose original
aim was to help \ref{d33}.

\noindent
3) The proof serves also for a related more local result, \ref{p7},
   there we just replace stage A; it also serves \S(5B).

\noindent
4) We may use normal $\bold x$ so $\bar x_{\bar c}$ disappears.
\end{discussion}

\begin{PROOF}{\ref{p2}}
\underline{Stage A}:  By Claim \ref{b35}, without loss of generality
\mn
\begin{enumerate}
\item[$\otimes_0$]   $\bold x$ is smooth (see Definition
\ref{b27}) so $\langle \bold I_{\bold x,\kappa}:\kappa \in 
{\ga}_{\bold x}\rangle$ are well defined.
\end{enumerate}
\mn
As we are assuming $\ntr_{\lc}(\bold x) = \lambda$ is regular, so, in
 particular, of cofinality $> \theta$, see Definition \ref{b43}(2) and 
by \ref{b47}(3), there are $\bar\psi,\varphi_*$ such that
(see Definition \ref{b43}(5) on $\lambda$-illuminate, \ref{b5}(8) for
$\Gamma^1_{\bold x}$, \ref{b9}(3A),(3B) on illuminate):
\mn
\begin{enumerate}
\item[$\otimes_1$]  $(a) \quad \bar \psi$ does $\lambda$-illuminate
 $\bold x$ so $\Gamma^1_{\bar\psi} = \Gamma^1_{\bold x}$, 
\sn
\item[${{}}$]   $(b) \quad \varphi_* = \varphi_*(\bar x_{\bar d},\bar
x_{\bar c},\bar y) \in \Gamma^1_{\bold x}$
\sn
\item[${{}}$]   $(c) \quad \psi_{\varphi_*}$ does not
$\lambda^+$-illuminate\footnote{In the present proof, we can demand
 that no $\psi$ does $\lambda^+$-illuminate $(\bold x,\varphi)$} 
$(\bold x,\varphi_*)$
\sn
\item[${{}}$]  $(d) \quad$ \wilog \,
$\psi_{\neg\varphi_*} = \psi_{\varphi_*}$
\sn
\item[${{}}$]   $(e) \quad$ if $A \subseteq M_{\bold x}$ has
cardinality $< \lambda$ then there is $\bar e \in {}^\theta(M_{\bold
x})$ such that

\hskip25pt  $\tp(\bar d_{\bold x},\bar c_{\bold x} \char 94 \bar e)
\vdash \tp(\bar d_{\bold x},\bar c_{\bold x} \char 94 \dotplus A)$
according to $\bar\psi$; follows by (a).
\end{enumerate}
\mn
Hence for some $A$
\mn
\begin{enumerate}
\item[$\otimes_2$]   $(a) \quad A \subseteq M_{\bold x}$ has cardinality
$\lambda$
\sn
\item[${{}}$]   $(b) \quad$ for no $\bar e \in {}^\theta(M_{\bold
x})$ does $\psi_{\varphi_*}(\bar x_{\bar d},\bar x_{\bar c},\bar e)$
solves $(\bold x,A,\varphi_*)$, see \ref{b9}(1A)
\sn
\item[${{}}$]   $(c) \quad$ let $\langle a_\alpha:\alpha <
\lambda\rangle$ list $A$
\sn
\item[$\otimes_3$]  $(a) \quad$ let 
$\varphi_0 = \varphi_*(\bar x_{\bar d},\bar x_{\bar c},\bar y_0)$
\sn
\item[${{}}$]   $(b) \quad \varphi^{\bold t}_0 = \varphi^{\iif(\bold
  t)}_0$ is $\varphi_0$ if
$\bold t = 1$ and $\neg \varphi_0$ if $\bold t = 0$, so
$\vartheta_{\varphi^{\bold t}_0}$ are well defined
\sn
\item[${{}}$]   $(c) \quad$ let $\varphi_1 = 
\varphi_1(\bar x_{\bar d},\bar x_{\bar c},\bar x_{[\theta]}) \in
\Gamma^1_{\bar\psi}$ be $\psi_{\varphi_0}= \psi_{\varphi_*}$ 
\sn
\item[${{}}$]   $(d) \quad$ let $\varphi_2 = 
\varphi_2(\bar x_{\bar d},\bar x_{\bar c},\bar x'_{[\theta]}) \in
\Gamma^1_{\bar\psi}$ be $\psi_{\varphi_1}$
\sn
\item[${{}}$]  $(e) \quad$ let $\varphi_3 = 
\varphi_3(\bar x_{\bar d},\bar x_{\bar c},\bar x''_{[\theta]}) \in
\Gamma^1_{\bar\psi}$  
be $\psi_{\varphi_2}$
\sn
\item[$\otimes_4$]  $(a) \quad$ let $\Delta_0 = 
\{\vartheta'_{\varphi_1}(\bar x'_{[\theta]},
x''_{[\theta]},\bar x_{\bar c})\}$ where
$\vartheta'_{\varphi_1} = \vartheta_{\varphi_1}(\bar x_{\bar c},\bar
x'_{[\theta]},\bar x_{[\theta]})$,

\hskip25pt  see \ref{b9}(1C)
\sn
\item[${{}}$]   $(b) \quad$ let $\Delta_1 = 
\{\vartheta'_{\varphi_2}(\bar x'_{[\theta]},
\bar x''_{[\theta]},\bar x_{\bar c})\}$ where
$\vartheta'_{\varphi_2} = \vartheta_{\varphi_2}
(\bar x_{\bar c},\bar x''_{[\theta]},\bar x'_{[\theta]})$
\sn
\item[${{}}$]   $(c) \quad$ let $\Delta_2 \subseteq 
\bbL(\tau_T)$ be finite large enough such that clause (i) of \ref{d10} 

\hskip25pt holds with $(\Delta_n,\Delta^1_n)$
 there standing for $(\Delta_0,\Delta_2)$ 

\hskip25pt here and for $(\Delta_1,\Delta_2)$ here.
\end{enumerate}
\mn
Note that
\mn
\begin{enumerate}
\item[$\otimes_5$]  $\Pi\{\lambda_i:i \in 
u_{\bold x}$ and $\lambda_i < \lambda\} < \lambda$.
\end{enumerate}
\mn
[Why it holds?  By clause (d) of the assumption; 
important for \ref{p7}.]
\bigskip

\noindent
\underline{Stage B}:  Let $I = ([M_{\bold x}]^{< \kappa},\subseteq)$.
Recall that for every $\bar e \in {}^\theta(M_{\bold x})$ 
for some $h \in \Pi\{\kappa_{\bold x,i}:i \in
u_*\}$ the pair $(B_{\bold x} + \bar e,\bar{\bold I}_{\bold x,h})$ is a
$(\bar\mu,\theta)$-set.  

Now let $\bar e_A$ be as guaranteed by $\otimes_1(e)$ above 
for $A \in I$.  Let ${\cD}_I$ be the club filter on $I$.

We apply\footnote{We could use $\Delta_n$ with union $\bbL(\tau_T)$
if $\theta = \aleph_0 = |T|$, if so we do not have to care in choosing
$\Delta_0$.} Theorem \ref{d10} with $1,M_{\bold x},(B_{\bold x},
\bar{\bold I}_{\bold x}),\langle \bar e_A:A \in I \rangle,\langle \bar
e_A:A \in I\rangle,\langle \Delta_0\rangle,\langle\Delta_2\rangle,
{\cD}_I,\cD_I$ here standing for $\bold k,M,\bold f,\langle \bar
e^1_A:A \in I\rangle,\langle \bar e^2_t:
t \in I\rangle,\langle \Delta_n:n < \bold k\rangle,
\langle \Delta^1_n:n < \bold k\rangle,{\cD}_1,\cD_2$ 
there.  We get $h^*_0,q_0,\cS'_1,\cS_2$ here standing for
$h_*,q,\cS_{1,0},\cS_2$ there.

Next we apply \ref{d10} again with $1,M_{\bold x},
(B_{\bold x},\bar{\bold I}_{\bold x}),\langle \bar e_A:A \in
I\rangle,\langle \bar e_A:A \in I\rangle,\langle \Delta_1\rangle,
\langle \Delta_2\rangle,{\cD}_I \rest
\cS'_1,\cD_I \rest \cS'_1$ here standing for $\bold k,
M,\bold f,\langle \bar e^1_A:A \in I\rangle,\langle \bar e^2_s:
s \in I_\ell\rangle,
\langle \Delta_n:n < \bold k\rangle,\langle \Delta^1_n:n < \bold k\rangle,
{\cD}_1,\cD_2$ there.
We get $h^*_1,q_1,\cS_0,\cS_1$ here standing for 
$h_*,q,\cS_{1,0},\cS_2$ there and \wilog \, $\cS_0,\cS_1 
\subseteq \cS'_1$.  Note that $q_\ell \in \bold S^{\theta +
  \theta}_{\Delta_\ell}(B^+_{\bold x})$ for $\ell=0,1$. 

Let $h = \max\{h^*_0,h^*_1\}$ and $B_* = B_{\bold x,u_*,h}$ and let
$S_\ell = \{3 \alpha + \ell:\alpha < \lambda\}$ for $\ell=0,1,2$.

We shall show that there is a quadruple $(\bar N,\bold I,\bar A,B_*)$  
such that:
\mn
\begin{enumerate}
\item[$\boxplus_1$]  $(a) \quad \bar N = \langle N_\alpha:
\alpha < \lambda\rangle$ and $\bold I = \langle \bar e_\alpha:\alpha <
\lambda\rangle$ and $\bar A = \langle A_\alpha:\alpha < \lambda\rangle$
\sn
\item[${{}}$]  $(b) \quad N_\alpha \prec M_{\bold x}$ is $\prec$-increasing
\sn
\item[${{}}$]  $(c) \quad \|N_\alpha\| < \lambda$
\sn
\item[${{}}$]   $(d) \quad B_{\bold x} \subseteq N_0$ and $a_\alpha
\in N_{\alpha +1}$; hence $A +  B_{\bold x} \subseteq N_\lambda :=
\cup\{N_\alpha:\alpha < \lambda\}$
\sn
\item[${{}}$]  $(e) \quad i \in u_{\bold x} \wedge |\bold
I_{\bold x,i}| < \lambda \Rightarrow \bold I_{\bold x,i} \subseteq N_0$ 
\sn
\item[${{}}$]   $(f) \quad$ if 
$i \in u_{\bold x} \wedge \kappa_{\bold x,i} =
\lambda$ then $\bar a_{\bold x,i,\alpha} \subseteq N_{\alpha +1}$ for
$\alpha < \kappa_{\bold x,i}$, hence $\bold I_{\bold x,i} \subseteq N_\lambda$
\sn
\item[${{}}$]  $(g) \quad \bar e_\alpha \in {}^\theta(N_{\alpha +1})$
\sn
\item[${{}}$]  $(h) \quad \tp(\bar d_{\bold x},\bar c_{\bold x} 
+ \bar e_\alpha) \vdash \tp(\bar d_{\bold x},\bar c_{\bold x}
 \dotplus (N_\alpha + B^+_{\bold x}))$ according 
to $\bar\psi$
\sn
\item[${{}}$]  $(i) \quad \bar c_\alpha \char 94 
\bar d_\alpha \subseteq
{}^{\theta^+>}(N_{\alpha +1})$ realize $\tp(\bar c_{\bold x} \char 94
\bar d_{\bold x},N_\alpha + B^+_{\bold x} + \bar e_\alpha)$

\hskip25pt where $\ell g(\bar c_\alpha) = \ell g(\bar c_{\bold
x}),\ell g(\bar d_\alpha) = \ell g(\bar d_{\bold x})$
\sn
\item[${{}}$]  $(j) \quad N_\alpha \subseteq A_\alpha \in
N_{\alpha +1}$
\sn
\item[${{}}$]  $(k) \quad A_\alpha \in \cS_\ell
\Leftrightarrow \alpha \in S_\ell$ and $\bar e_\alpha = \bar e_{A_\alpha}$
\sn
\item[${{}}$]  $(l) \quad$ if $\beta < \alpha$ and $\beta \in
S_0,\alpha \in S_1$ then $q_0 = \tp_{\Delta_0}(\bar e_\beta
\char 94 \bar e_\alpha,B_*)$
\sn
\item[${{}}$]  $(m) \quad$ if $\beta < \alpha$ and $\beta \in
S_1,\alpha \in S_2$ then $q_1 = \tp_{\Delta_1 }(\bar e_\beta
\char 94 \bar e_\alpha,B_*)$.
\end{enumerate}
How?  We shall choose $N_0,\bar e_\alpha,\bar c_\alpha,\bar
d_\alpha,q_\alpha$ by induction on $\alpha$
satisfying the relevant conditions.

In the induction step, first $N_\alpha$ exists as it should just be
$\prec M_{\bold x}$ and include $< \lambda$ specific elements and has
to be of cardinality $< \lambda$.  Second, $A_\alpha$ exists, if
$\alpha \in S_\ell$ it can be any member of
$\cS_\ell$ satisfying $< \lambda$ requirements, each such
requirement is satisfied by a set of $A$'s from ${\cD}_I +
\cS_\ell$ which is a $\lambda$-complete filter.  

Third, $\bar e_\alpha$ exists by the definition of
$\lambda = \ntr_{\lc}(\bold x)$ and choice of $\bar\psi$, 
more exactly by $\otimes_1(e)$.

Fourth, $\bar c_\alpha \char 94 \bar d_\alpha$ exists 
as $M_{\bold x}$ is $\kappa$-saturated and $\kappa >
\lambda$; so we are done carrying the induction.

Let $u_1 = \{i \in u_{\bold x}:\kappa_{\bold x,i} < \lambda\}$ and
$u_2 = u_{\bold x} \backslash u_1$.  For each $\alpha < \lambda$ by
\ref{b35}(5) we choose a function $h = h_\alpha \in
\Pi\{\kappa_{\bold x,i}:i \in u_{\bold x}\}$
such that $(B_{\bold x} + N_\alpha + \bar e_\alpha + \bar d_{\bold x}
+ \bar c_{\bold x},\bar{\bold I}_{\bold x,u_2,h})$ is a
$(\bar\mu,\theta)$-set recalling $\bar{\bold I}_{\bold x,u_2,h} 
= \langle \bold I_{\bold x,i,h(\kappa_i)}:i \in u_2\rangle$ and $\bar{\bold
I}_{\bold x,i,h(\kappa_i)} = \langle \bar a_{i,\alpha}:\alpha \in
[h(\kappa_i),\kappa_{\bold x,i})\rangle$.

So
\mn
\begin{enumerate}
\item[$\boxplus_2$]  $\langle \bar a_{\bold x,\partial,\beta}:\beta
\in I_{\bold x,\partial,h_\alpha(\partial)}\rangle$ is 
indiscernible over $B_{\bold x} + N_\alpha + \bar c_{\bold x} + 
\bar d_{\bold x} + \bar e_\alpha + \cup\{\bar
a_{\bold x,\sigma,\alpha}:\sigma \in {\ga}_{\bold x}
\backslash \{\partial\}$ and $\alpha \in [h_\alpha(\sigma),\sigma)\}$
for every $\partial \in \ga_{\bold x}$.
\end{enumerate}

Note that
\mn
\begin{enumerate}
\item[$\boxplus_3$]  we can replace $\langle (N_\alpha,\bar
c_\alpha,\bar d_\alpha,\bar e_\alpha):\alpha < \lambda\rangle$ by
$\langle(\cup\{N_{f(\beta)+1}:\beta < 1 + \alpha\} \cup N_\alpha),\bar
c_{f(\alpha)},\bar d_{f(\alpha)},\bar e_{f(\alpha)}):\alpha <
\lambda\rangle$ when $f:\lambda \rightarrow \lambda$ is increasing (so
$\alpha \le f(\alpha)$) and $\ell \in \{1,2,3\} \wedge \alpha \in
S_\ell \Rightarrow f(\alpha) \in S_\ell$.
\end{enumerate}
\mn
Hence recalling noting $\Pi(a_{\bold x} \backslash \lambda^+)$ is
$\lambda^+$-directed and $\cf(\Pi(\ga \cap \lambda))) < \lambda$ by
$\otimes_5$, for some $h_* \in \Pi \ga$ we have $\ell \le 2
\Rightarrow \lambda = \sup\{\alpha \in S_\ell:h_\alpha \rest (\ga
\backslash \{\lambda\}) \le h_*\}$ and $h_*(\lambda) = 0$ hence there is $S'
\subseteq \lambda$ such that for every $\alpha < \lambda$ there is
$\beta \in S'$ such that $\otp(S' \cap \beta) = \alpha \wedge
\bigwedge\limits_{\ell} (\alpha \in S_\ell) \equiv \beta \in S_\ell)$
hene shrinking $\langle a_\alpha:\alpha < \lambda\rangle$ by a
subsequence and as we can replace $h_\alpha$ by any
bigger function in $\Pi {\ga}_{\bold x}$, without loss of generality
\mn
\begin{enumerate}
\item[$\boxplus_4$]   $(a) \quad h_\alpha \rest \ga_{\bold x} \backslash
\{\lambda\} = h_*$, i.e. is constant
\sn
\item[${{}}$]  $(b) \quad \langle h_\alpha(\alpha):\alpha
<\lambda\rangle$ is increasing
\sn
\item[$\boxplus_5$]  let $B_* = B^+_{\bold x,h_*}$ recalling $h_* \in
\Pi(\ga \backslash \{\lambda\})$ and\footnote{the
difference between $B_*$ and $B^*_\alpha$ is concerning $\bold I_{\bold
x,\lambda}$} $B^*_\alpha = B_{\bold x,h_\alpha}$ for $\alpha <
  \lambda$ so $\alpha < \beta < \lambda \Rightarrow B_* \subseteq
  B^*_\beta \subseteq B^*_\alpha$.
\end{enumerate}
\mn
Also without loss of generality $\langle e_{\alpha,0}:\alpha <
\lambda\rangle$ is with no repetitions.
\bigskip

\noindent
\underline{Stage C}:  Let $N^*_\lambda \prec M_{\bold x}$ be of cardinality
$< \kappa$ such that $B^+_{\bold x} \cup N_\lambda \subseteq N^*_\lambda$.
We choose $N^+_\lambda$ expanding $N^*_\lambda$ such that 
$P^{N_\lambda}_0 =
|N_\lambda|,P^{N^+_\lambda}_1 = \{e_{\alpha,0}:\alpha < \lambda\},
P^{N^+_\lambda}_2 = \{(a,e_{\alpha,0}): a \in N_\alpha,\alpha <
\lambda\}$ and $F^{N^+_\lambda}_i(e_{\alpha,0}) 
= e_{\alpha,i}$ for $i < \theta,P^{N^+_\lambda}_{3 + \ell} =
\{e_{\alpha,0}:\alpha \in S_\ell\}$ for $\ell=0,1,2$ so $N^+_\lambda$
and the vocabulary $\tau(N^+_\lambda)$ are well defined.

We shall choose an increasing sequence 
$\langle \alpha_\varepsilon = \alpha(\varepsilon):\varepsilon <
\lambda\rangle$ enumerating in increasing order 
a thin enough club of $\lambda$.

We shall prove in this stage that there are $N^\oplus_\lambda$ and $\langle
\bar e^*_\varepsilon:\varepsilon < \lambda\rangle$ such that:
\mn
\begin{enumerate}
\item[$\boxplus_6$]  $(a) \quad N^\oplus_\lambda$ is 
an elementary extension of $N^+_\lambda$
\sn
\item[${{}}$]  $(b) \quad \bar e^*_\varepsilon \in
{}^\theta\{a:N^\oplus_\lambda \models 
P_2(a,e_{\alpha(\varepsilon +1),0})\}$
\sn
\item[${{}}$]  $(c) \quad e^*_{\varepsilon,i} = 
F^{N^\oplus_\lambda}_i(e^*_{\varepsilon,0})$ for $i <
\theta,\varepsilon < \lambda$
\sn
\item[${{}}$]  $(d) \quad e^*_{\varepsilon,0} \in 
P^{N^\oplus_\lambda}_4$ and $\neg P^{N^\oplus_\lambda}_2
(e^*_{\varepsilon,0},e_{\alpha(\varepsilon),0})$
\sn
\item[${{}}$]  $(e) \quad \langle \bar e^*_\varepsilon:
\varepsilon < \lambda\rangle$ is an
indiscernible sequence in $N^\oplus_\lambda \rest \tau_T$
\sn
\item[${{}}$]  $(f)(\alpha) \quad \bar e^*_\varepsilon(\varepsilon <
\lambda),\bar e_\alpha(\alpha \in S_1)$ realize the same 
$\bbL(\tau_T)$-type over $B_*$
\sn
\item[${{}}$]  $\quad (\beta) \quad$ if $\beta < \lambda$
then all $\bar e^*_\varepsilon$ such that 
$\alpha(\varepsilon) < \beta$ and $\bar e_\alpha$ such that

\hskip40pt    $\alpha \in S_1 \cap \beta$
realize the same $\bbL(\tau_T)$-type over $B^*_\beta$
\sn
\item[${{}}$]  $\quad (\gamma) \quad$ if $\alpha \le
  \alpha_\varepsilon,\varepsilon < \lambda$ and $\alpha \in S_0$ \then
  \, $\bar e_\alpha \char 94 \bar e^*_\varepsilon$ realizes $q_0$
\sn
\item[${{}}$]  $\quad (\delta) \quad$ if $\varepsilon <
  \lambda,\alpha_{\varepsilon +1} \le \alpha$ and $\alpha \in S_2$
  \then \, $\bar e^*_\varepsilon \char 94 \bar e_\alpha$ realizes $q_1$.
\end{enumerate}
\mn
First note:
\mn
\begin{enumerate}
\item[$\boxplus_{6.1}$]  there is an increasing continuous sequence $\langle
\alpha_\varepsilon:\varepsilon < \lambda\rangle$ of limit ordinals 
$< \lambda$ such that: for every $n < \omega$,
finite\footnote{Alternatively use finite $u \subseteq \ell g(\bar
e_0),\Delta$ finite $\subseteq \bbL(\tau_T)$ and get $\langle \bar
e_{\beta_0} \rest u,\dotsc,\bar e_{\beta_1} \rest u)$ is
$\Delta$-indiscernible.} $\Delta \subseteq \cup\{\Gamma_{(\theta)_n}:m <
\omega\} = \{\varphi:\varphi =
\varphi(\bar x_{\bar e_0},\bar x_{\bar e_1},\dotsc,\bar x_{\bar
e_{m-1}})$ and $m < \omega,\varphi \in \bbL(\tau_T)\}$ and for every
$0 = \varepsilon_0 < \varepsilon_1 < \ldots < \varepsilon_n$ 
we can find $\beta_\ell \in
[\alpha_{\varepsilon_\ell},\alpha_{\varepsilon_{\ell +1}}) \cap 
S_1$ for $\ell <n$
such that $\langle \bar e_{\beta_0},\bar e_{\beta_1},\dotsc,\bar
e_{\beta_{n-1}}\rangle$ is a $\Delta$-indiscernible sequence (in
$M_{\bold x}$).
\end{enumerate}
\mn
[Why?  For each such pair $(\Delta,n)$ define a game
$\Game_{\Delta,n}$ with $n$ moves, in the $m$-th move the antagonist
chooses an ordinal $\beta_m < \lambda$ which is $> \sup\{\gamma_k:k <
m\}$ and the protagonist chooses $\gamma_m \in [\beta_m,\lambda) \cap
S_1$.
In the end of a play the protagonist wins the play when $\langle \bar
e_{\gamma_0},\dotsc,\bar e_{\gamma_{n-1}}\rangle$ is a
$\Delta$-indiscernible sequence.  This game is determined so we choose
a winning strategy {\bf st}$_{\Delta,n}$ for the winner.  Let $E =
\{\delta < \kappa:\alpha < 1 + \delta \Rightarrow h_\alpha(\alpha) + 1
< \delta$ and $\delta$ is closed under {\bf st}$_{\Delta,n}$ for
every pair $(\Delta,n)$ as above$\}$.  As the number of pairs
$(\Delta,n)$ as above is $\le \theta < \lambda = 
\text{ cf}(\lambda)$, clearly $E$ is a club of $\kappa$
and let $\bar \alpha = \langle \alpha_\varepsilon:\varepsilon <
\kappa\rangle$ list $E$ in increasing order.

It is enough to prove that the protagonist wins every such game
$\Game_{\Delta,n}$.
Now for each pair $(\Delta,n)$ if this fails the sequence $\langle 
\bar e_{\alpha(k)+1}:k <
\omega\rangle$ has an infinite subsequence $\langle \bar e_{k(i)}:i <
n\rangle$ of length $n$ which is
$\Delta$-indiscernible, hence the protagonist wins in least in one
play of $\Game_{\Delta,n}$ in which the antagonist uses the strategy {\bf
st}$_{\Delta,n}$, i.e. when he chooses $\gamma_m = \alpha_{k(m)}$, it is
legal by the choice of $E$, so {\bf st}$_{\Delta,n}$ is not a winning
strategy for the antagonist hence it is for the protagonist
$\Game_{n,\Delta}$.  So easily $\bar\alpha$ is as required in 
$\boxplus_{5.1}$.].

Now let $N^+_\lambda$ be a $\|N_\lambda\|^+$-saturated elementary
extension of $N^+_\lambda$, by $\boxplus_{5.1}$ we can find in it a
sequence $\langle \bar e^*_\varepsilon:\varepsilon < \lambda\rangle$
as promised in $\boxplus_6$.  Note that clause (f) of $\boxplus_4$
can be gotten by thinning the sequence $\langle
e^*_\varepsilon:\varepsilon < \lambda\rangle$.
\bigskip

\noindent
\underline{Stage D}:  There is $N^\oplus_\lambda$ such that
\mn
\begin{enumerate}
\item[$\boxplus_7$]  $(a) \quad N^\oplus_\lambda$ 
has cardinality $\lambda$ (by the LST theorem)
\sn
\item[${{}}$]  $(b) \quad N^\oplus_\lambda \rest \tau_T \prec M_{\bold x}$ 
(by renaming, possible as $M_{\bold x}$ is $\kappa$-saturated while 

\hskip25pt $\kappa > \lambda = \|N^\oplus_\lambda\|$)
\sn
\item[${{}}$]   $(c) \quad$ if $\bar b \in {}^{\omega >}
(N_\lambda),\varepsilon < \lambda,\vartheta \in \bbL(\tau_T)$ and
  $\alpha \in (\alpha_\varepsilon,\alpha_{\varepsilon +1}) \cap S_1
  \Rightarrow$

\hskip25pt $\gC \models \vartheta[\bar e_\alpha,\bar b,
\bar c_{\bold x}]$ then $\models \vartheta[\bar e^*_\varepsilon,\bar b,
\bar c_{\bold x}]$.
\end{enumerate}
\mn
[Why?  As $\tp(\bar c_{\bold x},M_{\bold x})$ is finitely satisfiable
  in $B^*_0$ and $B^*_0 \subseteq M_{\bold x},|B^*_0| <
  \kappa,M_{\bold x}$ is $\kappa$-saturated.]
\mn
\begin{enumerate}
\item[${{}}$]  $(d) \quad$ let 
$N'_\varepsilon = M_{\bold x} \rest \{a \in N^\oplus_\lambda:
N^\oplus_\lambda \models P_2(a,e^*_{\varepsilon,0})\}$ for any
$\varepsilon < \lambda$.
\end{enumerate}
\mn
Now recall 
\mn
\begin{enumerate}
\item[$\odot_1$]   ${\gC} \models ``\vartheta_{\varphi^{\bold t}_0}
[\bar c_{\bold x},\bar e_{\alpha(\varepsilon)},\bar b]"$ when $\bold t
\in \{0,1\},\bar b \in {}^{\ell g(\bar y_0)}(N_{\alpha(\varepsilon)})$ and
${\gC} \models \varphi^{\bold t}_0[\bar d_{\bold x},
\bar c_{\bold x},\bar b]$
\end{enumerate}
\mn
[Why?  By the choice of $\bar e_{\alpha(\varepsilon)}$ and of
  $\vartheta_{\varphi^{\bold t}_0}$ recalling $\psi_{\varphi^{\bold
      t}_0} = \psi_{\varphi_*}$ by $\otimes_3$.]
\mn
\begin{enumerate}
\item[$\odot_2$]  ${\gC} \models \varphi_1[\bar d_{\bold x},
\bar c_{\bold x},\bar e_{\alpha(\varepsilon)}$.]
\end{enumerate}
\mn
[Why?  As $\varphi_1 = \psi_{\varphi_0}$, see $\otimes_3(c)$ 
and the choice of the $\bar e_{\alpha(\varepsilon)}$.]
\mn
\begin{enumerate}
\item[$\odot_3$]  ${\gC} \models ``\varphi_{\ell +1}[\bar d_{\bold x},
\bar c_{\bold x},\bar e_\alpha] \wedge \vartheta_{\varphi_\ell}[\bar
c_{\bold x},\bar e_\alpha,\bar e_{\alpha(\varepsilon)}]"$ when
$\alpha \in (\alpha_\varepsilon,\alpha_{\varepsilon +1})$ and $\ell
  \in \{1,2\}$.
\end{enumerate}
\mn
[Why?  As $\varphi_{\ell +1} = \psi_{\varphi_\ell}$, see 
$\otimes_3(d)$ and $\odot_2$.]
\mn
\begin{enumerate}
\item[$\odot_4$]   $q_0 = \tp_{\Delta_0}(\bar e_{\alpha(\varepsilon)} 
\char 94 \bar e_\alpha,B^*_{\alpha_{\varepsilon +1}})$ when
$\alpha_\varepsilon \in S_0$ and $\alpha
\in (\alpha_\varepsilon,\alpha_{\varepsilon +1}) \cap S_1$.
\end{enumerate}
\mn
[Why?  See clause $(\ell)$ of $\boxplus_1$ above, recall that $q_0$ is
  from the application of Theorem \ref{d10} in the beginning of Stage B.]
\mn
\begin{enumerate}
\item[$\odot_5$]  $q_0 \rest B^*_\beta = 
\tp_{\Delta_0}(\bar e_{\alpha(\varepsilon)} 
\char 94 \bar e^*_\varepsilon,B^*_\beta)$ when $\alpha_{\varepsilon
  +1} \le \beta < \lambda$.
\end{enumerate}
\mn
[Why?  By $\odot_4$ and the choice of $\bar e^*_\varepsilon$,
i.e. $\boxplus_6(f)(\gamma)$.]
\mn
\begin{enumerate}
\item[$\odot_6$]  ${\gC} \models \vartheta_{\varphi_1}[\bar
c_{\bold x},\bar e^*_\varepsilon,\bar e_{\alpha(\varepsilon)}]$ iff
$\vartheta_{\varphi_1}[\bar c_{\bold x},\bar e_\alpha,\bar
e_{\alpha(\varepsilon)}]$ when $\varepsilon < \lambda \wedge
\alpha \in (\alpha_\varepsilon,\alpha_{\varepsilon +1}) \cap S_1$;
recalling $\alpha_\varepsilon \in S_0$ being a limit ordinal.
\end{enumerate}
\mn
[Why?  By $\boxplus_6(f)(\beta)$ and as $\vartheta_{\varphi_1}(\bar c_{\bold
x},\bar x_{[\theta]},\bar e_{\alpha(\varepsilon)})$ is a
$\Delta_0$-formula, see by $\otimes_4(a)$ recalling $\tp(\bar c_{\bold
x},M_{\bold x})$ is finitely satisfiable in $B^*_\beta$ as $B^*_\beta
  = B^+_{\bold x,h_\beta}$.]
\mn
\begin{enumerate}
\item[$\odot_7$]  ${\gC} \models \varphi_2[\bar d_{\bold x},
\bar c_{\bold x},\bar e^*_\varepsilon]$.
\end{enumerate}
\mn
[Why?  Clearly $\varphi_3 = \psi_{\varphi_2}$, see $\otimes_3(e)$ and 
$\vartheta_{\varphi_2}$ are well defined and by $\odot_3$ we know that 
${\gC} \models
\varphi_2[\bar d_{\bold x},\bar c_{\bold x},\bar e_\alpha]$ for every
$\alpha \in (\alpha_\varepsilon,\alpha_{\varepsilon +1})$.  Let $\beta
\in S_2 \backslash \alpha(\varepsilon +1)$.
Again by $\odot_3$ we have 
${\gC} \models ``\psi_{\varphi_2}[\bar d_{\bold x},\bar c_{\bold
x},\bar e_\beta] \wedge \vartheta_{\varphi_2}[\bar c_{\bold x},
\bar e_\beta,\bar e_\alpha]"$ for $\alpha \in
(\alpha_\varepsilon,\alpha_{\varepsilon +1})$.

But by $\boxplus_1(m)$ 
for every $\alpha \in (\alpha_\varepsilon,\alpha_{\varepsilon +1})
\cap S_1$ we have $q_1 = \tp_{\Delta_1}(\bar e_\alpha \char 94
\bar e_\beta,B_*)$.  So for every $\bar c \in
{}^{\ell g(\bar c_{\bold x})}(B^*_\beta)$ we have ${\gC} \models
``\vartheta_{\varphi_2}[\bar c,\bar e_\beta,\bar e_\alpha]"$ 
iff $\vartheta_{\varphi_2}(\bar c,\bar x''_{[\theta]},\bar
x'_{[\theta]}) \in q_1$.
Hence for every $\bar c \in {}^{\ell g(\bar c_{\bold x})}(B^*_\beta)$ we have
${\gC} \models ``\vartheta_{\varphi_2}[\bar c,\bar e_\beta,
\bar e^*_\varepsilon]"$ iff
$\vartheta_{\varphi_2}(\bar c,\bar x'_{[\theta]},\bar x''_{[\theta]})
\in q_1$.  As this holds for every such $\bar c$ and $\tp(\bar c_{\bold
x},M_{\bold x})$ is finitely satisfiable in $B_*$, clearly ${\gC}
\models ``\vartheta_{\varphi_2}[\bar c_{\bold x},\bar e_\beta,
\bar e_\alpha] \equiv \vartheta_{\varphi_2}[\bar c_{\bold x},
\bar e_\beta,\bar e^*_\varepsilon]"$ for every $\alpha \in
(\alpha_\varepsilon,\alpha_{\varepsilon +1}) \cap S_1$.

By the conclusions of the last two paragraphs
${\gC} \models \vartheta_{\varphi_2}[\bar c_{\bold x},
\bar e_\beta,\bar e^*_\varepsilon]$ and by the conclusion of the first
of them ${\gC} \models \psi_{\varphi_2}[\bar d_{\bold x},\bar
c_{\bold x},\bar e_\beta]$.  Together recalling the definition of
$\vartheta_{\varphi_2}$ we get
${\gC} \models \varphi_2[\bar d_{\bold x},\bar c_{\bold x},
\bar e^*_\varepsilon]$, i.e. we are done proving $\odot_7$.]
\mn
\begin{enumerate}
\item[$\odot_8$]  ${\gC} \models ``\vartheta_{\varphi_1}[\bar c_{\bold
x},\bar e^*_\varepsilon,\bar e_{\alpha(\varepsilon)}]"$.
\end{enumerate}
\mn
[Why?  Note $\gC \models ``\vartheta_{\varphi_1}[\bar c_{\bold x},
\bar e_\alpha,\bar e_{\alpha(\varepsilon)}]"$ for $\alpha \in
(\alpha_\varepsilon,\alpha_{\varepsilon +1}) \cap S_1$.  Let $\bold I
= \{\bar c' \in {}^{\ell g(\bar c[\bold x])}(B_*):\gC \models 
``\vartheta_{\varphi_1}[\bar c',\bar e_\alpha,\bar
e_{\alpha(\varepsilon)}]"$ but clearly $\bold I = \{\bar c' \in
{}^{\ell g(\bar c[\bold x])}(B_*):\vartheta_{\varphi_2}[\bar c_{\bold x},
\bar x'_{[\theta]},\bar x''_{[\theta]}) \in q_0\}$ so does not depend
on $\alpha$ hence, second, $\bold I = \{\bar c' \in {}^{\ell g(\bar
c[\bold x])}(B_{\bold x}):\gC \models \vartheta_{\varphi_1}[\bar c'],
\bar e^*_\varepsilon,\bar e_{\alpha(\varepsilon)}]\}$.  As $\tp(\bar
c_{\bold x},M_{\bold x})$ is finitely satisfiable in $B_*$, we get
$\gC \models \vartheta_{\varphi_1}[\bar c_{\bold x},
\bar e^*_\varepsilon,\bar e_{\alpha(\varepsilon)}]$.]
\mn
\begin{enumerate}
\item[$\odot_9$]  ${\gC} \models ``\psi_{\varphi_1}[\bar d_{\bold
x},\bar c_{\bold x},\bar e^*_\varepsilon]",
{\gC} \models ``\varphi^{\bold t[\bar b]}_0[\bar d_{\bold x},
\bar c_{\bold x},\bar b]"$ and when $\bar b \in{}^{\ell g(\bar y)}
(N_{\alpha(\varepsilon)})$ and $\bold t[\bar b]$ is chosen such that
${\gC} \models \vartheta_{\varphi_0^{\bold t[\bar b]}}[\bar c_{\bold x},
\bar e^*_\varepsilon,\bar b]$.
\end{enumerate}
\mn
[Why?  Recalling $\psi_{\varphi_1} = \varphi_2$, note that 
${\gC} \models \psi_{\varphi_1}[\bar d_{\bold x},
\bar c_{\bold x},\bar e^*_\varepsilon]$ by $\odot_7$, i.e. the first
conclusion of $\odot_9$ holds.  By $\odot_7$ we have ${\gC} \models 
\vartheta_{\varphi_1}[\bar c_{\bold x},\bar e^*_\varepsilon,\bar
e_{\alpha(\varepsilon)}]$ which means that $\psi_{\varphi_1}(\bar x_{\bar
d},\bar c_{\bold x},\bar e^*_\varepsilon) \vdash \varphi_1(\bar
x_{\bar d},\bar c_{\bold x},\bar e_{\alpha(\varepsilon)})$.  But
$\varphi_1 = \psi_{\varphi_0} = \psi_{\neg \varphi_0}$ so by
$\otimes_1(d) + \otimes_3(a) + \otimes_3(b)$ 
we have $\varphi_1(\bar x_{\bar d},\bar c_{\bold x},
\bar e_{\alpha(\varepsilon)}) \vdash \varphi^{\bold t}_0(\bar x_{\bar d},
\bar c_{\bold x},\bar b)$.

As $\vdash$ is transitive we have $\psi_{\varphi_1}(\bar x_{\bar d},
\bar c_{\bold x},\bar e^*_\varepsilon) \vdash 
\varphi^{\bold t}_0(\bar x_{\bar d},
\bar c_{\bold x},\bar b)^{\bold t}$ which by $\odot_7$ means ${\gC} \models
\vartheta_{\varphi_0}[\bar c_{\bold x},\bar e^*_\varepsilon,\bar b]$, 
i.e. the second conclusion of $\odot_9$ so we are done.]
\bigskip

\noindent
\underline{Stage E}:  By the choice of $\varphi_* = \varphi_0$ and
letting $\psi = \psi_{\varphi_*}$  and
of the set $A$ see $\otimes_2$ we can find
\mn
\begin{enumerate}
\item[$(*)_1$]    an ultrafilter $D$ on ${}^{\ell g(\bar z)}
(N_\lambda)$ such that for every $\bar e' \in {}^\theta(M_{\bold x})$
and $\psi'$ satisfying $\psi'(\bar x_{\bar d},\bar c_{\bold x},\bar e') 
\in \tp(\bar d_{\bold x},\bar c_{\bold x} + M_{\bold x})$ 
the\footnote{can use ``$\{b \in {}^{\ell g(\bar z)}(A):\ldots\}$"}
 set $\{\bar b \in {}^{\ell g(\bar z)}(N_\lambda):
{\gC} \models (\exists \bar x_{\bar d})
(\psi'(\bar x_{\bar d},\bar c_{\bold x},\bar e') \wedge \varphi_*
(\bar x_{\bar d},\bar c_{\bold x},\bar b)) \wedge (\exists \bar
x_{\bar d})(\psi'(\bar x_{\bar d},\bar c_{\bold x},\bar e') \wedge \neg
\varphi_*(\bar x_{\bar d},\bar c_{\bold x},\bar b))\}$ belongs to $D$;
\end{enumerate}
\mn
this is as in the proof of \ref{q5} or \ref{b22}(2),
i.e. \cite[2.10=tp25.36]{Sh:900}
\mn
\begin{enumerate}
\item[$(*)_2$]   For $\bar b \in {}^{\ell g(\bar y)}
(N_\lambda)$ let $\varepsilon(\bar b) = \min\{\varepsilon < \lambda:
\bar b \subseteq N_{\alpha(\varepsilon)}\}$ and
\sn
\item[$(*)_3$]   let $\bold t(*)$ be
such that $\{\bar b \in {}^{\ell g(\bar y)}(N_\lambda):\bold t[\bar
b] = \bold t(*)\} \in D$ recalling $\bold t[\bar b] \in \{0,1\}$ is such that
${\gC} \models \varphi^{\bold t[\bar b]}_0[\bar d_{\bold x},
\bar c_{\bold x},\bar b]$.
\end{enumerate}
\mn
Note that
\mn
\begin{enumerate}
\item[$(*)_4$]  for every $\varepsilon < \lambda$ the set $\{\bar b
\in {}^{\ell g(\bar y)}(N_\lambda):\varepsilon(\bar b) \ge \varepsilon\}$
belongs to $D$.
\end{enumerate}
\mn
[Why?  Otherwise $\psi(\bar x_{\bar d},\bar c,\bar e_\varepsilon)$
contradicts the choice of $D$.]

We use ultrapower to get $(\bar b',\bar e')$ in ${\gC}$ realizing
$p(\bar y,\bar x_{[\theta]}) = 
\{\vartheta(\bar d_{\bold x},\bar c_{\bold x},\bar y,
\bar x_{[\theta]},\bar a):\vartheta(\bar x_{\bar d},\bar x_{\bar c},\bar
y,\bar x_{[\theta]},\bar z) \in \bbL(\tau_T)$ and 
$\bar a \in {}^{\ell g(\bar z)}(M_{\bold x})$ and 
$\{\bar b \in {}^{\ell g(\bar z)}(N_\lambda):{\gC} \models 
\vartheta[\bar d_{\bold x},\bar c_{\bold x},\bar b,\bar
e^*_{\varepsilon(\bar b)},\bar a]\} \in D\}\}$.

Now
\mn
\begin{enumerate}
\item[$(*)_5$]   $\langle \bar e^*_\varepsilon:\varepsilon < 
\lambda\rangle \char 94 \langle \bar e'\rangle$ is an
indiscernible sequence.
\end{enumerate}
\mn
[Why?  By $(*)_6$ below.]
\mn
\begin{enumerate}
\item[$(*)_6$]   $\bar e'$ realizes $\Av(\langle 
\bar e^*_\varepsilon:\varepsilon < \lambda\rangle,M_{\bold x} +
\bar c_{\bold x})$.
\end{enumerate}
\mn
[Why?  As $\langle \bar e^*_\varepsilon:\varepsilon < \kappa\rangle$
is an indiscernible sequence (and $T$ is dependent), the average is
well defined.  Now recall $(*)_4$ and the choice of $\bar e'$.]
\mn
\begin{enumerate}
\item[$(*)_7$]  $\tp(\bar d_{\bold x},M_{\bold x} + \bar c_{\bold x}),
\tp(\bar b',M_{\bold x} + \bar c_{\bold x})$ are not weakly
orthogonal.
\end{enumerate}
\mn
[Why?  By the choice of $\bar b'$ and of $\cD$, i.e. as 
witnessed by $\varphi_*$.]

But (an important point for Claim \ref{p7}) we need a 
more effective version of $(*)_7$.

Let $p_1(\bar x_{\bar d}) = \tp(\bar d_{\bold x},
M_{\bold x} + \bar c_{\bold x})$ and $p_2(\bar y) = 
\tp(\bar b',M_{\bold x} + \bar c_{\bold x}),
p_3(\bar x_{[\theta]}) = \tp(\bar e',M_{\bold x} + \bar c_{\bold x})$
\mn
\begin{enumerate}
\item[$(*)^+_7$]  $p_1(\bar x_{\bar d}) \cup p_2(\bar y)
\cup \{\varphi^{\bold t}_*(\bar x_{\bar d},\bar c_{\bold x},\bar y)\}$ is
consistent for $\bold t = 0,1$.
\end{enumerate}
\mn
[Why?  By the choice of the ultrafilter $D$ and of the sequence $\bar b'$.]
\mn
\begin{enumerate}
\item[$(*)_8$]  ${\gC} \models \psi_{\varphi_*}[\bar d_{\bold
x},\bar c_{\bold x},\bar e']$ and ${\gC} \models
\vartheta_{\varphi_*^{\bold t(*)}}[\bar c_{\bold x},\bar e',\bar b']$.
\end{enumerate}
\mn
[Why?  Because by $(*)_4$ and $\odot_8$ 
for every $\bar b \in {}^{\ell g(\bar y)}(N_\lambda)$
we have ${\gC} \models \psi_{\varphi_1}[\bar d_{\bold x},\bar c_{\bold
x},\bar e_{\varepsilon(\bar b)}]$ and ${\gC} \models
\vartheta_{\varphi^{\bold t[\bar b]}_*}
[\bar c_{\bold x},\bar e_{\varepsilon(\bar b)},\bar b]$.]
\mn
\begin{enumerate}
\item[$(*)_9$]   $p_1(\bar x_{\bar d}) \cup p_3(\bar x_{[\theta]})
\cup \{\pm \psi_{\varphi_*}(\bar x_{\bar d},\bar c_{\bold x},\bar
x_{[\theta]})\}$ are consistent.
\end{enumerate}
\mn
[Why?  First, clearly ${\gC} \models \psi_{\varphi_*}[\bar d_{\bold x},\bar
c_{\bold x},\bar e']$ by $(*)_8$ hence $p_1(\bar x_{\bar d}) \cup
p_3(\bar x_{[\theta]}) \cup \{\psi_{\varphi_*}(\bar x_{\bar d},\bar
c_{\bold x},\bar x_{[\theta]})\}$  being realized by $\bar d_{\bold x}
\char 94 \bar e'$ is consistent.

By $(*)^+_7$
for some $\bar d'$ realizing $p_1(\bar x_{\bar d})$ and $\bar b''$
realizing $p_2(\bar y)$ recalling $\bold t(*)$ is from $(*)_3$
 we have ${\gC} \models 
\neg \varphi^{\bold t(*)}_*[\bar d',\bar c_{\bold x},\bar b'']$; 
as $p_1(\bar x_{\bar d}) = \tp(\bar d_{\bold x},
M_{\bold x} + \bar c_{\bold x})$ \wilog \,
$\bar d' = \bar d_{\bold x}$.  As $\tp(\bar b'',\bar c_{\bold x} +
M_{\bold x}) = p_2(\bar y) = \tp(\bar b',M_{\bold x} 
+ \bar c_{\bold x})$ for some $\bar e''$ we
have $\tp(\bar b'' \char 94 \bar e'',M_{\bold x} + \bar c_{\bold x}) =
\tp(\bar b' \char 94 \bar e',M_{\bold x} + \bar c_{\bold x})$.

Note that ${\gC} \models \vartheta_{\varphi_*^{\bold t(*)}}
[\bar c_{\bold x},\bar e',\bar b']$ hence by $(*)_8$ we have
${\gC} \models \vartheta_{\varphi_*^{\bold t(*)}}
[\bar c_{\bold x},\bar e'',\bar b'']$.

Now if ${\gC} \models \psi_{\varphi_*}[\bar d_{\bold x},\bar
c_{\bold x},\bar e'']$ then by the definition of
$\vartheta_{\varphi_*,\bold t}$, see $\odot_9$ and the last
sentence, ${\gC} \models \varphi_*^{\bold t(*)}[\bar d_{\bold x},
\bar c_{\bold x},\bar b'']$ but $\bar d_{\bold x} \char
94 \bar e''$ realizes $p_1(\bar x_{\bar d}) \cup p_3(\bar
x_{[\theta]})$ we have a contradiction to
the choice of $\bar b''$ hence ${\gC} \models \neg \psi_{\varphi_*}[\bar
d_{\bold x},\bar c_{\bold x},\bar e'']$ thus finishing the proof of
$(*)_9$.]

As $p_3(\bar x_{[\theta]})$ was defined as $\tp(\bar e',M_{\bold x} +
\bar c_{\bold x})$ by $(*)_6 + (*)_9$ we are done. 
\end{PROOF}
\bigskip

We shall not use \ref{p2} as stated but a variant which the proof
gives (as mentioned in the proof).
\begin{claim}
\label{p7}  Assume $\bold x \in \pK_{\kappa,\mu,\theta},{\cP} 
\subseteq \{u \subseteq
v_{\bold x}:u \cap u_{\bold x}$ is finite$\}$ is $\subseteq$-directed
with union $v_{\bold x}$ and $\kappa > \lambda = \ntr_{\lc}(\bold x)$
and $i \in v_{\bold x} \backslash u_{\bold x} \Rightarrow \lambda
\ge \beth_\omega(|B_i| + \theta)$ and $\lambda$ is regular, this is similar to
\ref{p2} but omitting the assumption ``$u_{\bold x}$ is finite".  We still
can find $(\bar\psi,A,u,\varphi_*)$ such that
\mn
\begin{enumerate}
\item[$(a)$]  $\bar\psi$ illuminates $(\bold x,\lambda,\Gamma^1_{\bold x})$
\sn
\item[$(b)$]   $A \subseteq M_{\bold x}$ is of cardinality $\lambda,u
\in {\cP}$
and $\varphi_* = \varphi_*(\bar x_{\bar d},\bar x_{\bar c,u},\bar y)$ such
that no $\psi'(\bar x_{\bar d},\bar c_{\bold x,u},\bar e') \in 
\tp(\bar d_{\bold x},\bar c_{\bold x} \dotplus M_{\bold x})$
solve $(\bold x,\varphi,A)$ [i.e. for no finite $u \subseteq \ell
  g(\bar c_{\bold x})$ and $\psi' = \psi'(\bar x_{\bar d},
\bar x_{\bar c,u},\bar z)$ and $\bar e' \in {}^{\ell g(\bar z)}
(M_{\bold x})$ do we have $\gC \models ``\psi'[\bar d_{\bold x},
\bar c_{\bold x,u},\bar e']"$ and
$\psi'(\bar x_{\bar d},\bar c_{\bold x,u},\bar e') \vdash \{\varphi_*(\bar
x_{\bar d},\bar c_{\bold x},\bar b):\bar b \in {}^{\ell g(\bar y)}A$
and $\varphi_*(\bar x_{\bar d},\bar c_{\bold x},\bar b) \in 
\tp(\bar d_{\bold x},A + \bar c_{\bold x})\}]$; let
$\psi_{\varphi_*} = \psi_{\varphi_*}(\bar x_{\bar d},\bar x_{\bar
c},\bar x_{[\theta]}) \equiv \psi_{\varphi_*}(\bar x_{\bar d},\bar
x_c,\bar x_{[v]}),v \subseteq \theta$ finite
\sn
\item[$(c)$]   we get the result of \ref{p2} with 
$\bar c_{\bold x}$ replaced by $\bar c_{\bold x,u} = \bar c_{\bold
x,u} \rest u = \langle \bar c_{\bold x,i}:i \in u\rangle$, i.e.:
\begin{enumerate}
\item[$(*)$]  there is an indiscernible sequence $\bold I = \langle \bar
b_\alpha:\alpha < \lambda\rangle$ in $M_{\bold x},\ell g(\bar
b_\alpha) = \ell g(\bar y)$ and $\bar b^0,\bar b^1 \in {}^v{\gC}$ realizing 
$\Av(\bold I,M_{\bold x} + \bar c_{\bold x,u})$ 
and ${\gC} \models \varphi[\bar d_{\bold x},\bar c_{\bold x,u},
\bar b^{\bold t}]^{\bold t}$ for $\bold t = 0,1$.
\end{enumerate}
\end{enumerate}
\end{claim}

\begin{PROOF}{\ref{p7}}
We can find $\bar\psi,\varphi_*,A$ as in $\otimes_1$ in the proof of
\ref{p2}, so $\bar\psi$ is as in clause (a) there.  Then we find
$\varphi_\ell,\vartheta'_\ell$ as in $\oplus_3 + \oplus_4$ in the proof
of \ref{p2}, Stage A.  Next let $u_* \in {\cP}$ be such
that $i \in u_*$ \Iff \, $x_{\bar c_{\bold x,i}}$ is not dummy in
$\varphi_0$ or in $\varphi_1$ or in $\varphi_2$ or in $\varphi_3$.  Now use the
proof of \ref{p2} from Stage B on, however not on $\bold x$ but $\bold
x_{[u_*]}$, see Definition \ref{b5}(10).  
\end{PROOF}

\begin{conclusion}
\label{p9}  
Assume $T$ is countable, $\theta =
\aleph_0,\mu$ strong limit of uncountable cofinality and $\mu \le
\kappa = \cf(\kappa) < \mu^{+ \omega}$.  \Then \, for every $\bold m \in 
\rK^\oplus_{\kappa,\mu,\theta}$ with $u_{\bold m}$ finite 
there is $\bold n \in \tK^\oplus_{\kappa,\mu,\theta}$ 
such that $\bold m \le_1 \bold n$.
\end{conclusion}

\begin{remark}
If we assume $\cf(\mu) = \aleph_0$ and $\kappa = \mu^+$, then we can get a
weaker version of density of $\tK_{\kappa,\mu,\theta}$.
\end{remark}
\bigskip

\begin{PROOF}{\ref{p9}}
  Without loss of generality $\ell g(\bar d) = \omega$.

Let $\langle \varphi_n(\bar x_{[\omega + \omega]},\bar y_{[\omega]},
\bar z_n):n <\omega\rangle$ list all formulas of such form, 
each appearing infinitely many times.  \Wilog \,
$\varphi_n = \varphi_n(\bar x_{[w_n]},\bar y_{[n]},\bar z_n),w_n :=
[0,n) \cup [\omega,\omega +n)$.  We choose $\bold m_n$ 
by induction on $n < \omega$ such that:
\mn
\begin{enumerate}
\item[$\boxplus$]  $(a) \quad \bold m_n \in 
\rK^\oplus_{\kappa,\mu,\theta}$ and $u(\bold m_n)$ is finite,
moreover $\in [n,\omega)$; we may 

\hskip25pt add $v(\bold m_n)$ finite (as we can assume $v_{\bold m}$ finite)
\sn
\item[${{}}$]  $(b) \quad \bold m_n = \bold m$
\sn
\item[${{}}$]  $(c) \quad$ if $n = m+1$ then $\bold m_m \le_1 \bold
  m_{n+1}$ and $r_{\bold m_n}$ is complete
\sn
\item[${{}}$]  $(d) \quad$ if $n=m+1$ and there is $\bold m' \in 
\rK^\oplus_{\kappa,\mu,\theta}$ satisfying

\hskip25pt  $(\bold m_m \le_1 \bold n \wedge (\bold m' \text{ is } 
(\varphi_m(\bar x_{[w_m]},\bar y_{[m]},\bar z_m),i)$ active 
for some

\hskip25pt  $i \in v(\bold m') \backslash v(\bold m_m))$ \then \,
$\bold m_n$ satisfies this
\sn
\item[${{}}$]  $(e) \quad$ if $n=m+1$ and the assumption in clause (d)
fails, but there is

\hskip25pt  $\bold m' \in \rK^\oplus_{\kappa,\mu,\theta}$ 
satisfying $\bold m_n \le_1 \bold m'$ and $\varphi_m \in 
\Gamma^2_{\bar\psi[\bold m']}$

\hskip25pt  \then \, $\bold m_{n+1}$ satisfies this.
\end{enumerate}
\mn
We can carry the induction for clauses (d) + (e) because if there is
such $\bold m'$ we can find $\bold m''$ such that $\bold m_m \le_1
\bold m'' \le_1 \bold m'$ such that $u[\bold m'']$ is finite, and the
demand ``$r_{\bold m_{m+1}}$ is complete" is not a problem by \ref{d36}(1A).

Having carried the induction let $\bold n = \lim\langle \bold m_n:n <
\omega\rangle$, we have to show that $\bold n \in
\rK^\oplus_{\kappa,\mu,\theta}$ and more.  
If $\Gamma_{\bar\psi[\bold n]} = \Gamma^2_{\bold x[\bold n]}$
we shall be done by \ref{c70}, 
so toward contradiction assume $\varphi = \varphi(\bar
x_{\bar d[\bold n]},\bar x_{\bar c[\bold n]},\bar y) \in
\Gamma^2_{\bold x[\bold n]} \backslash \Gamma_{\bar\psi[\bold n]}$,
let $k$ be such that $\varphi = \varphi_k$ hence $u =
\{n:\varphi_n = \varphi\}$ is infinite.  By clause (d) and
\ref{b20}(2), i.e. \cite{Sh:900} the set $\{m:\varphi_m =
\varphi_k$ and the assumption in $\boxplus(d)$ holds$\}$ is finite.
So choose $m$ such that $\varphi_m = \varphi_k$ but the assumption of
$\boxplus(d)$ fails.  By \ref{p2} more exactly \ref{p7}, the assumption
of $\boxplus(e)$ holds; why?  the point is that 
$\ntr_{\lc}(\bold x_{\bold n}) = \{\mu^+,\mu^{+2},\ldots,\mu^{+n}\}$.
So the conclusion of $\boxplus(e)$ holds, contradiction.

Lastly, $\bold n = \cup\{\bold m_n:n < \omega\}$ is well defined and
by Claim \ref{c70}, using (c)$'$ there, $\bold m = \bold m_0
\le_1 \bold n \in \tK^\oplus_{\kappa,\mu,\theta}$.
\end{PROOF}
\bigskip

\noindent
\centerline{$* \qquad * \qquad *$}
\bigskip

\noindent
\begin{discussion}
\label{p10}
We may like to cover every $\kappa = \kappa^{< \kappa} \ge
\beth_\omega$; at least and/or when for countable $T$ as in \ref{d39}, 
G.C.H. holds).  For this, we are still left with the case $\cf(\mu) 
= \aleph_0$, for this we
have to redo some previous definitions and claims, so this is
presently delayed.
\end{discussion}

\begin{conclusion}
\label{p13}  Assume G.C.H. and $T$ is countable and $\theta =
\aleph_0$ and $\mu$ is strong limit of cofinality $> \aleph_0$ and
$\kappa = \cf(\kappa) \in (\mu,\mu^{+ \omega})$.

\noindent
1) For every $\kappa$-saturated $M$ of cardinality $\kappa$ and $\bar d
\in {}^{\theta^+ >}{\gC}$ there is $\bold x \in 
\tK_{\kappa,\kappa,\theta}$ with $\bar d \trianglelefteq \bar d_{\bold x}$.

\noindent
2) Hence $M \in \EC_{\kappa,\kappa}(T) \Rightarrow |\bold S^\theta(M)/
   \equiv_{\aut}| \le \kappa$.

\noindent
3) If $M$ is a saturated model of $T$ of cardinality $\kappa$ \then \, 
$\gS^\theta_{\aut}(M)$ has
cardiality $\le \mu$.
\end{conclusion}
\bigskip

\begin{remark}
\label{p15}
1) For $\vK$ this is easier.

\noindent
2) When $\cf(\mu) = \aleph_0$, maybe see more in \cite{Sh:F973}.
\end{remark}
\bigskip

\subsection {Density of $\vK$; Exact recounting of types and $\vK$} \

We prove the density of $\vK_{\kappa,\bar\mu,\theta}$.  We use
\ref{q3}(1) but not \ref{q3}(2)-(6).

Recall that we have difficulties when $\ntr_{\lc}(\bold x)$ was
singular.  This motive defining relatives in \ref{n4},\ref{n8} and
investigating them.  This succeeds but not applicable to $\rK$ only to $\vK$.

\begin{convention}
\label{n1}
We here tend to use $\varphi \in \Gamma^1_{\bold x}$ as $\varphi(\bar
x_{d,\rho},\bar x_{\bar c,\varrho},\bar y)$ where $\rho \in {}^n(\ell
g(\bar d_{\bold x}))),\varrho \in {}^m(\ell g(\bar c_{\bold x}))$ for
some $n,m$.
\end{convention}

\begin{definition}
\label{n4}
Assume for $\bold x \in \pK_{\kappa,\bar\mu,\theta}$ and $\varphi =
\varphi(\bar x_{\bar d,\rho},\bar x_{\bar c,\varrho},\bar y)$, so
$\varphi$ determines $\rho,\varrho$ and 
$\psi = \psi(\bar x_{d,\rho},\bar c_{\bold x,\varrho_0}) \in
\tp(\bar d_{\bold x,\rho},\bar c_{\bold x,\varrho_0})$.  Below we
   may omit $\psi$ when $\rho = \langle \rangle = \varrho_0$, the
   role of $\psi$ in (1) is minor.

\noindent
1) Let $\bold k(\varphi,\psi,\bold x)$ be the maximal $n$ such that there is an
 increasing sequence $\varrho_1 \in {}^n(v_{\bold x})$ 
which witness it, which means (note that $\psi$ has a role only via
 $\varrho_0$):
\mn
\begin{enumerate}
\item[$\bullet$]  $\ell < \ell g(\varrho) \wedge k < \ell g(\varrho_1)
  \Rightarrow \varrho(\ell) <_{v_{\bold x}} \varrho_1(k)$
\sn
\item[$\bullet$]  $\ell < \ell g(\varrho_0) \wedge k < \ell g(\varrho_1)
\Rightarrow \varrho_0(\ell) <_{v_{\bold x}} \varrho_1(\ell)$
\sn
\item[$\bullet$]  $\bar c^*_{\ell,0},\bar c^*_{\ell,1}$ are
  subsequences of $\bar c_{\bold x,\varrho_1(\ell)}$ realizing the same
  type over $\bar c_{\bold x,< \varrho_1(\ell)} + M_{\bold x}$
\sn
\item[$\bullet$]  $\gC \models ``\varphi[\bar d_{\bold x,\rho},
\bar c_{\bold x,\varrho},\bar c^*_{\ell,1}] \wedge \neg \varphi[\bar
  d_{\bold x,\rho},\bar c_{\bold x,\varrho},\bar c^*_{\ell,0}]$.
\end{enumerate}
\mn
2) We define $\bold k_{\du}(\varphi,\psi,\bold x)$ as the maximal $n$ such
that some $\eta$ witness it which means; $\du$ stands for duplicate:
\mn
\begin{enumerate}
\item[$\bullet$]  $\eta$ is an increasing sequence in $w_{\bold x}$ of
  length $\ell g(\rho)$
\sn
\item[$\bullet$]  $\ell g(\bar d_{\bold x,\rho(\ell)})= \ell g(\bar
  d_{\bold x,\eta(\ell)})$ for $\ell < \ell g(\rho)$
\sn
\item[$\bullet$]  $\psi(\bar x_{\bar d,\eta},\bar c_{\bold
  x,\varrho_0}) \in \tp(\bar d_{\bold x,\eta},\bar c_{\bold x,\varrho_0})$
\sn
\item[$\bullet$]  $\tp_\varphi(\bar d_{\bold x,\eta},\bar
  c_{\bold x,\varrho} \dotplus M_{\bold x}) = \tp_\varphi(\bar d_{\bold
  x,\rho},\bar c_{\bold x,\varrho} \dotplus M_{\bold x})$
\sn
\item[$\bullet$]  $n = \bold k(\psi(\bar x_{\bar d,\eta},\bar c_{\bold
  x,\varrho_0}),\varphi(\bar x_{\bar d,\eta},\bar x_{\bar c,
\varrho},\bar y),\bold x)$.
\end{enumerate}
\end{definition}

\begin{claim}
\label{n6}
1) In Definition \ref{n4}, $\bold k(\varphi,\psi,\bold x)$ is well defined
and $< \ind(\varphi)$.

\noindent
2) Also $\bold k_{\du}(\varphi,\psi,\bold x)$ is well defined and 
$< \ind(\varphi)$.

\noindent
3) If $(\varphi,\psi,\bold x)$ is as in \ref{n4} and $\bold x \le_1
   \bold y \in \pK_{\kappa,\bar\mu,\theta}$ \then \,:
\mn
\begin{enumerate}
\item[$\bullet$]  $(\varphi,\psi,\bold y)$ is as in \ref{n4}
\sn
\item[$\bullet$]  $\bold k(\varphi,\psi,\bold x) \le 
\bold k(\varphi,\psi,\bold y)$
\sn
\item[$\bullet$]  $\bold k_{\du}(\varphi,\psi,\bold x) \le 
\bold k_{\du}(\varphi,\psi,\bold y)$.
\end{enumerate}
\mn
4) If $\bold x \in \pK_{\kappa,\bar\mu,\theta}$ \then \, there is $\bold
y$ such that (see \ref{n19}):
\mn
\begin{enumerate}
\item[$\bullet_1$]  $\bold x \le_1 \bold y \in \pK_{\kappa,\mu,\theta}$
\sn
\item[$\bullet_2$]  if $(\varphi,\psi,\bold x)$ is as in \ref{n4} and 
$\bold y \le \bold z \in \pK_{\kappa,\bar \mu,\varphi}$ then 
$\bold k(\varphi,\psi,\bold y) = \bold k(\varphi,\psi,\bold z)$ and 
$\bold k_{\du}(\varphi,\psi,\bold y) = \bold k_{\du}(\varphi,\psi,\bold z)$.
\end{enumerate}
\mn
5) Like (4) but in $\bullet_2$ it applies to $(\varphi,\psi)$ such
that $(\varphi,\psi,\bold y)$ is as in \ref{n4}.
\end{claim}

\begin{PROOF}{\ref{n6}}
1),2) By the definition of $\ind(\varphi)$ it is always finite as $T$
  is dependent, see \cite{Sh:900} or \ref{q3}(1).

\noindent
3) Read the definition.

\noindent
4) By parts (1),(2),(3) as
\mn
\begin{enumerate}
\item[$(*)_1$]  $(\bold K_{\kappa,\bar\mu,\theta},\le_1)$ is a partial order
\sn
\item[$(*)_2$]  in this partial order an increasing sequence $\langle
\bold x_i:i < \delta\rangle$ of length $< \theta^+$ has a
  $\le_1$-upper bound $\bold x_\delta$; moreover is the union so if
  $(\varphi,\psi,\bold x_\delta)$ is as in \ref{n4} then for some $i <
\delta,(\varphi,\psi,\bold x_i)$ is as in \ref{n4}
\sn
\item[$(*)_3$]  for any $\bold x \in \pK_{\kappa,\bar\mu,\theta}$,
  there $\le \theta$ relevant pairs $(\varphi,\psi)$.
\end{enumerate}
\mn
5) Similarly using (4).
\end{PROOF}

\begin{claim}
\label{n5}
If (A) then (B) where
\mn
\begin{enumerate}
\item[$(A)$]  $(a) \quad \bold x,\varphi = \varphi(\bar x_{\bar d,\rho},
\bar x_{\bar c,\rho},\bar y),\psi(\bar x_{\bar d,\rho},
\bar x_{\bar c,\varrho_0})$ are as in Definition \ref{n4}
\sn
\item[${{}}$]  $(b) \quad n = \bold k_{\du}(\varphi,\psi,\bold x)$ and $\rho_1$ witness it
\sn
\item[${{}}$]  $(c) \quad$ let $\varphi' = \varphi(\bar x_{\bar d,\rho_1},\bar
  x_{\bar c,\varrho},\bar y)$
\sn
\item[${{}}$]  $(d) \quad$ let 
$\psi' = \psi(\bar x_{\bar d,\rho_1},\bar c_{\bold
  x,\varrho_0})$ hence $\in \tp(\bar d_{\bold x,\rho_1},\bar c_{\bold
  x,\varrho_0})$
\sn
\item[${{}}$]  $(e)\quad$ let 
$\varrho_1$ witness $\bold k(\varphi',\psi',\bold x)$ 
with $\langle (c^*_{\ell,0},\bar c^*_{\ell,1}):\ell < \ell
  g(\varrho_1)\rangle$ as in 

\hskip25pt  Definition \ref{n4}(1)
\sn
\item[${{}}$]  $(f) \quad 
\psi''(\bar x_{\bar d,\rho_1},\bar c_{\bold x,\varrho_0
  \char 94 \varrho_1}) = \psi''(\bar x_{\bar d,\rho_1},\bar c_{\bold
  x,\varrho_0 \char 94 \varrho_1}) \wedge \bigwedge\limits_{\ell <
  \ell g(\varrho_1)} (\varphi(\bar x_{\bar d,\rho},\bar c_{\bold
  x,\varrho},\bar c^*_{\ell,1}) \wedge$

\hskip25pt $\neg \varphi(\bar x_{\bar d,\rho},
\bar c_{\bold x,\varrho},\bar c^*_{\ell,0}))$
\sn
\item[$(B)$]  $(a) \quad (\varphi',\psi'',\bold x)$ are as in \ref{n4}
\sn
\item[${{}}$]  $(b) \quad \bold k_{\du}(\varphi',\psi'',\bold x) = 0$.
\end{enumerate}
\end{claim}

\begin{PROOF}{\ref{n5}}
Straightforward.
\end{PROOF}

\begin{definition}
\label{n8}
1) For $\bold x \in \pK_{\kappa,\bar\mu,\theta}$ and $\varphi_* =
   \varphi_*(\bar x_{\bar d},\bar x_{\bar c},\bar y) \in
   \Gamma^1_{\bold x}$ and $\psi_* = \psi_*(\bar x_{\bar d},\bar
   c_{\bold x},\bar b_*) \in \tp(\bar d_{\bold x},\bar c_{\bold x}
   \dotplus M_{\bold x})$ let $\ntr_w(\varphi_*,\psi_*,\bold x)$ be
the maximal $\lambda$ such that: if $A \subseteq M_{\bold x}$ and
   $|A| < \lambda$ \then \, 
for some finite $p \subseteq \tp_{\pm \varphi_*}(\bar
x_{\bar d},\bar c_{\bold x} \dotplus M_{\bold x})$ we have $p
   \cup \{\psi_*(\bar x_{\bar d},\bar c_{\bold x},b_*)\} \vdash
\tp_{\pm \varphi_*}(\bar d_{\bold x},\bar c_{\bold x} \dotplus A)$.

\noindent
2) For $\bold x \in \pK_{\kappa,\bar\mu,\theta}$ let $\ntr_w(\bold x)
   = \min\{\ntr_w(\varphi_*,\psi_*,\bold x):\varphi_*,\psi_*$ as above$\}$.
\end{definition}

\begin{claim}
\label{n10}
1) For $\bold x,\varphi_*,\psi_*$ as in Definition \ref{n8}(1) the
   cardinal $\ntr_w(\varphi_*,\psi_*,\bold x)$ is a regular (infinite)
   cardinal.

\noindent
2) For $\bold x \in \pK_{\kappa,\bar\mu,\theta}$ the cardinal
   $\ntr_w(\bold x)$ is a regular (infinite) cardinal.

\noindent
3) If $\lambda := \ntr_w(\varphi_*,\psi_*,\bold x)$ is $> \aleph_0$
\then \, for some $m$ we can replace ``$p \subseteq ...$ finite" by: for
   some fix $n$ and $\eta \in {}^n 2$ we have $``p \subseteq 
\tp_{\pm \varphi}(\bar d_{\bold x},\bar c_{\bold x} 
+ M_{\bold x})$ has the form
   $\{\varphi(\bar x,\bar a_\ell)^{\iif(\eta(\ell))}:\ell < n\}$.
\end{claim}

\begin{PROOF}{\ref{n10}}
1) Toward contradiction assume $\lambda =
   \ntr_w(\varphi_*,\psi_*,\bold x)$ is singular and $A \subseteq
   M_{\bold x}$ has cardinality $\lambda$.  We shall prove that for some finite
   $p \subseteq \tp_{\varphi_*}(\bar d_{\bold x},\bar c_{\bold x}
   \dotplus M_{\bold x})$ we have $p \cup \{\psi_*(\bar x_{\bar
   d},\bar c_{\bold x})\} \vdash \tp_{\varphi_*}(\bar d_{\bold x},\bar
   c_{\bold x} \dotplus M_{\bold x})$, this suffices.

Let $\langle A_\varepsilon:\varepsilon < \cf(\lambda)\rangle$ be a
$\subseteq$-increasing sequence of subsets of $A$ with union $A$ with
each $A_\varepsilon$ having cardinality $< \lambda$.  For each
$\varepsilon < \cf(\lambda)$ there is a finite $p_\varepsilon \in
\tp_{\varphi_*}(\bar d_{\bold x},\bar c_{\bold x} \dotplus (M_{\bold
  x})$ such that $p_\varepsilon  \cup \{\psi_*(\bar x_{\bar d},\bar
c_{\bold x})\} \vdash \tp_{\varphi_*}(\bar d_{\bold x},\bar c_{\bold
  x} \dotplus A_\varepsilon)$.  As $p_\varepsilon$ is finite also
$B_\varepsilon := \Dom(p_\varepsilon)$ is finite hence the cardinality
of $B := \cup\{B_\varepsilon:\varepsilon < \cf(\lambda)\}$ is $\le
\cf(\lambda)$.  As $B \subseteq M_{\bold x}$ there is a finite $q
\subseteq \tp_{\varphi_*}(\bar d_{\bold x},\bar c_{\bold x} \dotplus
M_{\bold x})$ such that $q \cup\{\psi_*(\bar x_{\bar d},\bar c_{\bold
  x})\} \vdash \tp_{\varphi_*}(\bar d_{\bold x},
\bar c_{\bold x} + M_{\bold x})$.

Now $q$ is as required.

\noindent
2) Follows from part (1).
\end{PROOF}

\noindent
The following is a replacement of \ref{p2} of \S(5A).
\begin{cc}
\label{n13}
If (A) then (B) \underline{where}:
\mn
\begin{enumerate}
\item[$(A)$]  $(a) \quad \lambda < \kappa$ is regular $\ge \mu$
\sn
\item[${{}}$]  $(b) \quad \bold x \in \pK_{\kappa,\bar\mu,\theta}$ and
$\iota_{\bold x} = 2$
\sn
\item[${{}}$]  $(c) \quad \varphi_* = \varphi_*(\bar x_{\bar d},\bar
  x_{\bar c,\varrho},\bar y)$ 
\sn
\item[${{}}$]  $(d) \quad \psi = \psi(\bar x_{\bar d},\bar x_{\bar c},\bar z)$
\sn
\item[${{}}$]  $(e) \quad \psi_*(\bar x_{\bar d}) = \psi(\bar x_{\bar
  d},\bar c,\bar b) \in \tp(\bar d_{\bold x},
\bar c_{\bold x} \dotplus M_{\bold x})$
\sn
\item[${{}}$]  $(f) \quad u_{\bold x}$ is finite, but see \ref{n15}
\sn
\item[${{}}$]  $(g) \quad \lambda = \cf(\lambda) = \min\{|A|:A \subseteq
  M_{\bold x}$ and there is no finite 

\hskip25pt $p \subseteq \tp_{\pm \varphi_*}(\bar d_{\bold x},
\bar c_{\bold x} \dotplus M_{\bold x})$
  such that $p \cup \{\psi_*(\bar x_{\bar d})\} \vdash \tp_{\pm
  \varphi_*}(\bar d_{\bold x},\bar c_{\bold x} \dotplus A)$
\sn
\item[${{}}$]  $(h) \quad$ if $A \subseteq M_{\bold x},|A| < \lambda$
  \then \, there is $\bar b_A \subseteq M$ such that $\varphi_*(\bar
  x_{\bar d},\bar c_{\bold x},\bar b_A) \in$

\hskip25pt $\tp(\bar d_{\bold x},\bar c_{\bold x} 
\dotplus M_{\bold x})$ and $\{\varphi_*(\bar x_{\bar d},
\bar c_{\bold x},\bar b_A)\} \cup \{\psi_*(\bar x_{\bar d})\} \vdash
\tp_{\pm \varphi_*}(\bar d_{\bold x},\bar c_{\bold x} \dotplus A)$
\sn
\item[$(B)$]  there are $\bold I,\bar d'$ such that
\sn
\begin{enumerate}
\item[$(a)$]  $\bold I = \langle \bar a_{\alpha,0} \char 94 \bar
  a_{\alpha,1}:\alpha < \lambda\rangle$ is an indiscernible sequence
  in $M_{\bold x}$
\sn
\item[$(b)$]  $\ell g(\bar a_{\alpha,\ell}) = \ell g(\bar y)$
\sn
\item[$(c)$]  $\bar a_{\lambda,0} \char 94 
\bar a_{\lambda,1}$ realizes $\Av(\bold I,\bar c_{\bold x} + M_{\bold x})$
\sn
\item[$(d)$]  $\tp(\bar a_{\lambda,0},\bar c_{\bold x} \dotplus
  M_{\bold x}) = \tp(\bar a_{\lambda,1},\bar c_{\bold x} \dotplus
  M_{\bold x})$
\sn
\item[$(e)$]  $\bar d'$ realizes $\{\psi_*(\bar x_{\bar d})\} \cup \tp_{\pm
  \varphi_*}(\bar d_{\bold x},\bar c_{\bold x} \dotplus M_{\bold x})
  \cup \{\varphi_*(\bar x_{\bar d},\bar a_{\lambda,1}),\neg
  \varphi_*(\bar x,\bar c_{\bold x},\bar a_{\lambda,0})\}$.
\end{enumerate}
\end{enumerate}
\end{cc}

\begin{remark}
\label{n14}
1) Note $\psi_*(\bar x_{\bar d})$ correspond to $q_*$ in \S(5C), so we
can restrict its form if necessary, see \S(5C).

\noindent
2) How will we justify clause (A)(h)?
\mn
\begin{enumerate}
\item[$(a)$]  we can manipulate $\varphi_*$ such that
  $\{\varphi_*(M,\bar a):\bar a\} = \{\neg \varphi_*(M,\bar a):\bar a\}$
  and $\emptyset,{}^{\ell g(\bar x)}(M)$ belongs to it, as in the
  proof of \ref{k4}
\sn
\item[$(b)$]  replacing $\varphi(\bar x,\bar y)$ by
  $\bigwedge\limits_{\ell < m} \varphi(\bar x,\bar y_\ell)$ change little.
\end{enumerate}
\mn
3) We could have weakened clause (B)(c) to $\Delta$-types for $\Delta$
derived from \ref{q3}, in fact $\Delta = \Lambda$ with $n = 
\bold k(\varphi_*,\bold x)$. 

\noindent
4) So \ref{q3}(2)-(6) is what is really required but we do not need it.
\end{remark}

\begin{observation}
\label{n15}
We can omit (A)(f) of \ref{n17}.
\end{observation}

\begin{PROOF}{\ref{n15}}
Let $\varphi_* = \varphi(\bar x_{\bar d,\rho},\bar x_{\bar
  c,\varrho},\bar y)$ with $\rho,\varrho$ as in convention \ref{n1}
  and work with $\bold y$ which is like $\bold x$ but $\bar d_{\bold
  y} = \bar d_{\bold x} \rest \Rang(\rho),\bar c_{\bold y} =  \bar
  c_{\bold x} \rest \Rang(\varrho)$.

Now reflect.
\end{PROOF}

\begin{PROOF}{\ref{n13}}

\noindent
\underline{Proof of \ref{n13}}

We repeat the proof of Claim \ref{p2}, making minor changes in Stages
(A)-(D) and replacing stage (E) as follows:
\medskip

\noindent
\underline{Stage (A)-(D)}:

We omit $\otimes_1(b)-(e)$, using clauses of (A) of the claim when
quoted.

In $\otimes_2(b)$ the set $A$ exemplify (A)(g) of the claim

In $\otimes_3$ let $\varphi_\ell = \varphi_*$ for $\ell=0,1,2$ 
(or just omit and replace $\varphi_\ell$ by $\varphi_*$ when used,
justified by clause (A)(h) hence $\bar e_\alpha,\bar e^*_\varepsilon
\in {}^{\ell g(\bar y)})(M_{\bold x})$

In $\otimes_4(a)$ add ``$+ \otimes_1$".

We replace $\bar x_{[\theta]}$ by $\bar x_{\ell g(\bar y)}$ or $\bar
y$ recalling $\varphi_* = \varphi(x_{\bar d},\bar x_{\bar c},\bar y)$.
\bigskip

\noindent
\underline{Stage E}:  Let $M^+ \prec \gC$ be such that $\bar
d_{\bold x} + \bar c_{\bold x} + M_{\bold x} \subseteq M^+$ and let $I
= \{(\bar b,\varepsilon,\zeta,\bar d):\bar b \in {}^{\ell g(\bar
y)}(N_\lambda),\varepsilon < \zeta < \lambda$ and $\bar d \in M^+\}$.

For every $\xi < \lambda$ and $p \in \cP = \{p:p$ is finite and $p
\subseteq \tp_{\pm \varphi_*}(\bar d_{\bold x},\bar c_{\bold x} \dotplus
M_{\bold x})\}$ we let $I_{p,\xi}$ be the set of $(\bar
b,\varepsilon,\zeta,\bar d)) \in I$ such that
\mn
\begin{enumerate}
\item[$\bullet$]  $\bar b \in {}^{\ell g(\bar y)}(N_\lambda)$
\sn
\item[$\bullet$]  $\bar d$ realizes $p$
\sn
\item[$\bullet$]  $\varepsilon < \zeta$ are from $[\xi,\lambda)$
\sn
\item[$\bullet$]  $\gC \models \psi_*[\bar d]$
\sn
\item[$\bullet$]  $\gC \models \varphi_*[\bar d,\bar c_{\bold
x},\bar e^*_\varepsilon]$
\sn
\item[$\bullet$]  $\gC \models \neg \varphi_*[\bar d,\bar c_{\bold
x},\bar e^*_\zeta]$.
\end{enumerate}
\mn
Now note that
\mn
\begin{enumerate}
\item[$(*)_1$]  if $p_1 \subseteq p_2 \subseteq \tp_{\pm \varphi_*}(\bar
d_{\bold x},\bar c_{\bold x} \dotplus M_{\bold x})$ are finite and $\xi_1 <
\xi_2 < \lambda$ then $I_{p_2,\xi_2} \subseteq I_{p_1,\xi_1}$.
\end{enumerate}
\mn
[Why?  Read the definition.]
\mn
\begin{enumerate}
\item[$(*)_2$]  if $(p,\xi) \in \cP \times \lambda$ then $I_{p,\xi}
\ne \emptyset$.
\end{enumerate}
\mn
[Why?  As $\models \varphi[\bar d_{\bold x},\bar c_{\bold x},\bar
e^*_\xi]$ clearly $p_1 := p \cup \{\varphi_*(\bar x_{\bar d},\bar
c_{\bold x},\bar e^*_\varepsilon)\}$ belongs to $\cP$, hence by clause (A)(g)
there is $\bar b \in {}^{\ell g(\bar y)}(N_\lambda)$ such that
$p_1(x_{\bar d}) \cup \{\psi_*(\bar x_{\bar d})\} \cup \{\varphi_*(\bar
x_{\bar d},\bar c_{\bold x},\bar b)^{\iif(\bold t)}\}$ is consistent for
$\bold t=0,1$.  Hence recalling $\gC \models \varphi_*[\bar d_{\bold
x},\bar c_{\bold x},\bar b]^{\iif(\bold t(\bar b_j))}$ there is $\bar
d$ in $\gC$ realizing $p_1(\bar x_{\bar d}) \cup \{\psi_*(\bar x_{\bar
d})\}$ and $\gC \models \varphi_*[\bar d_{\bold x},\bar c_{\bold
x},\bar b]^{\iif(1-\bold t(\bar b)}$.  Next choose $\zeta < \lambda$
such that $\zeta > \varepsilon$ and $\bar b \in N_\zeta$.  Now
$\varphi_*(\bar x_{\bar d},\bar c_{\bold x},\bar e^*_\zeta) \vdash
\varphi_*[\bar x_{\bar d},\bar c_{\bold x},\bar b]^{\iif(\bold t(\bar
b))}$ so necessarily $\gC \models ``\neg \varphi_*[\bar d,\bar
c_{\bold x},\bar e^*_\zeta]^{\iif(\bold t[\bar b])}"$.

Clearly $(\bar b,\varepsilon,\zeta,\bar d) \in I_{p,\xi}$ so $I_{p,\xi}
\ne \emptyset$ as promised.]
\mn
\begin{enumerate}
\item[$(*)_3$]  choose an ultrafilter $\cD$ on $I$ such that
$(p,\xi) \in \cP \times \lambda \Rightarrow I_{p,\xi} \in \cD$.
\end{enumerate}
\mn
[Why?  As $I_{p,\xi} \subseteq I$ using $(*)_1$ and $(*)_2$
above.]
\mn
\begin{enumerate}
\item[$(*)_4$]  $p(\bar y,\bar y',\bar y'',\bar x_{\bar d})$ is the
  following complete type over $\bar c_{\bold x} + M_{\bold x}$; where $\bar
y',\bar y''$ has length $\ell g(\bar y)$:
\newline
$p(\bar y,\bar y',\bar y'',\bar
x_{\bar d}) = \{\vartheta(\bar c_{\bold x},\bar y,\bar y',\bar
y'',x_{\bar d},\bar e):\vartheta = \vartheta(\bar x_{\bar c[\bold x]},\bar
y,\bar y',\bar y'',\bar x_{\bar d},\bar z) \in \bbL(\tau_T),\bar e \in
{}^{\ell g(\bar z)}(M_{\bold x})$ and the set $\{(\bar
b,\varepsilon,\zeta,\bar d) \in I:\gC \models \vartheta[\bar c_{\bold
x},\bar b,\bar e^*_\varepsilon,\bar e^*_\zeta,\bar d,\bar e]\}$
belongs to $\cD\}$.
\end{enumerate}
\mn
[Why?   As $\cD$ is an ultrafilter on $I$.]
\mn
\begin{enumerate}
\item[$(*)_5$]  choose  $(\bar b',\bar a_0,\bar a_1,\bar d')$ in $\gC$
realizing $p(\bar y,\bar y',\bar y'',\bar x_{\bar d})$.
\end{enumerate}
\mn
Now note
\mn
\begin{enumerate}
\item[$(*)_6$]   $(a) \quad \bar a_1$ realizes $\Av(\langle \bar
  e^*_\varepsilon:\varepsilon < \lambda\rangle,\bar c_{\bold x} +
  M_{\bold x})$
\sn
\item[${{}}$]   $(b) \quad \bar a_0$ realizes $\Av(\langle 
e^*_\varepsilon:\varepsilon < \lambda\rangle,\bar a_0 + \bar c_{\bold x} +
  M_{\bold x})$
\sn
\item[${{}}$]   $(c) \quad \bar a_0 \char 94 \bar a_1$ 
realizes $\Av(\langle \bar e^*_{2\varepsilon} \char 94 \bar e^*_{2
  \varepsilon +1}:\varepsilon < \lambda,\bar c_{\bold x} +   M_{\bold
  x}\rangle)$.
\end{enumerate}
\mn
[Why?  Think.]
\mn
\begin{enumerate}
\item[$(*)_7$]  $(a) \quad \bar d'$ realizes $\tp_{\varphi_*}(\bar
  d_{\bold x},\bar c_{\bold x} \dotplus M_{\bold x})$
\sn
\item[${{}}$]   $(b) \quad \bar d'$ realizes $\varphi_{\bar d},
\bar c_{\bold x},\bar a_0)$ and $\neg \varphi_*(\bar x_{\bar d},\bar
c_{\bold x},\bar a_1)$.
\end{enumerate}
\mn
So clearly we are done.
\end{PROOF}

\begin{claim}
\label{n17}
Assume $\bold x \in \pK_{\kappa,\bar \mu,\theta},\mu \ge
\beth_\omega,\varphi_* = \varphi(\bar x_{\bar d,\rho},\bar x_{\bar
  c,\varrho},\bar y)$ and\footnote{May add parameters from $M_{\bold x}$,
  but can use trivial members of $\bar c_{\bold x}$, i.e. $\bar
  c_{\bold x,i} \subseteq M_{\bold x}$.} $\psi_* = \psi_*(\bar
x_{d,\rho},\bar x_{\bar c,\varrho_0}) \in \tp(\bar d_{\bold x,\rho},\bar
c_{\bold x,\varrho_0})$ and $\bold k(\varphi_*,\psi_*,\bold x)=0$.  
If $\lambda :=
\ntr_w(\varphi_*,\psi_*,\bold x)$ is $< \kappa$ \then \, there is
$\bold y$ such that $\bold x \le_1 \bold y \in
\pK_{\kappa,\bar\mu,\theta}$ and $0 < \bold k(\varphi_*,\psi_*,\bold y)$.
\end{claim}

\begin{PROOF}{\ref{n17}}
Let $A_* \subseteq M_{\bold x}$ witnessing
$\ntr_w(\varphi_*,\psi_*,\bold x) = \lambda$.

Let $\bar A_* = \langle A_\varepsilon:\varepsilon < \lambda\rangle$,
etc.
\bigskip

\noindent
\underline{Case 1}:  $\lambda < \mu$.

As in \ref{b20} that is \cite[2.8=tp25.33]{Sh:900} and see Definition
\cite[2.6=tp25.32]{Sh:900} but here we elaborate.
\mn
\begin{enumerate}
\item[$\boxplus_1$]  Let $J$ be the set of pairs $(q,\Gamma)$ such that:
\sn
\begin{enumerate}
\item[$(a)$]  $q = q(\bar x_{\bar d,\rho}) \subseteq 
\tp_{\pm \varphi_*}(\bar d_{\bold x,\rho},\bar c_{\bold x,\varrho}
\dotplus M)$ is finite
\sn
\item[$(b)$]  $\Gamma = \Gamma(\bar y)$ is a finite subset of
$\Lambda = \{\vartheta(\bar y,\bar c):\vartheta(\bar y,\bar z) \in
  \bbL(\tau_T)$ and $\bar c \in {}^{\ell g(\bar z)}(M_{\bold x})\}$
\end{enumerate}
\item[$\boxplus_2$]  for a pair $(q,\Gamma) \in J$ we say 
$(\bar c_0,\bar c_1)$ does $A_*$-exemplifies $(q,\Gamma)$ \when \, :
\sn
\begin{enumerate}
\item[$(a)$]  $\bar c_0,\bar c_1 \in {}^{\ell g(\bar y)}(A_*)$
\sn
\item[$(b)$]  $\gC \models ``\vartheta[\bar c_0] \equiv \vartheta[\bar
c_1]"$ when $\vartheta(\bar y) \in \Gamma(\bar y)$
\sn
\item[$(c)$]  $\{\psi(\bar x_{\bar d})\} \cup q(\bar x_{\bar d,\eta}) \cup
\{\varphi(\bar x_{\bar d,\rho},\bar c_1),\neg \varphi(\bar x_{\bar
d,\rho},\bar c_0)\}$ is consistent
\end{enumerate}
\item[$\boxplus_3$]  the family $\{\{(\bar c_0,\bar c_1):(\bar
c_0,\bar c_1)$ does $A_*$-exemplifies $(q,\Gamma)\}:(q,\Gamma) \in J\}$ has the
finite intersection property.
\end{enumerate}
\mn
[Why does $\boxplus_3$ hold?  Otherwise we can find 
$(q_\ell,\Gamma_\ell) \in J$ for $\ell < n$
such that no $(\bar c_0,\bar c_1)$ does $A_*$ exemplify $(p,\Gamma_\ell)$ for
every $\ell < n$.  Define the two-place relation $E$ on ${}^{\ell
g(\bar y)}(A_*)$:
\mn
\begin{enumerate}
\item[$\odot_{3.1}$]  $\bar c_0 E \bar c_1$ iff $\bar c_0,\bar c_1 \in
{}^{\ell g(\bar y)}(A_*)$ and $\gC \models ``\vartheta[\bar c_0] \equiv
\vartheta[\bar c_1]"$ for every $\vartheta(\bar y) \in
\bigcup\limits_{\ell < n} \Gamma_\ell$
\end{enumerate}
\mn
clearly $E$ is an equivalence
relation with finitely many equivalence classes.  Let $\langle \bar
c^*_\ell:\ell < \ell(*)\rangle$ be a set of representatives and let
$q_* = (\bigcup\limits_{\ell < n} q_\ell) \cup \{\varphi(\bar x_{\bar
d,\rho},\bar c^*_\ell)^{\bold t}:\bold t \in \{0,1\}$ and $\gC
\models \varphi[\bar d_{\bold x,\rho},\bar c^*_\ell]^{\bold t}\}$.

So
\mn
\begin{enumerate}
\item[$\odot_{3.2}$]  $q_*$ is a finite subset of 
$\tp_{\pm \varphi}(\bar d_{\bold x,\rho},\bar c_{\bold x,\varrho}
  \dotplus M_{\bold x})$ and
\sn
\item[$\odot_{3.3}$]  $q_* \vdash \tp_{\pm \varphi}
(\bar d_{\bold x,\rho},\bar c_{\bold x,\varrho} \dotplus A_*)$.
\end{enumerate}
\mn
But $\odot_{3.2}$ 
contradicts the choice of $A_*$ so $\boxplus_3$ holds indeed.]

So
\mn
\begin{enumerate}
\item[$\boxplus_4$]  there is an ultrafilter on ${}^{2 \ell g(\bar
y)}(A_*)$ extending the family from $\boxplus_3$.
\end{enumerate}
\mn
Choose such an ultrafilter $D$.

Let $(\bar c'_0,\bar c'_1)$ realizes $\Av(D,\bar d_{\bold x} +
\bar c_{\bold x} + M_{\bold x})$, so clearly
\mn
\begin{enumerate}
\item[$\boxplus_5$]  the following set of formulas is finitely
  satisfiable in $\gC$:
\[
\{\psi_*(\bar x_{\bar d})\} \cup \tp_{\pm \varphi_*}(\bar d_{\bold
  x,\rho},\bar x_{\bar c,\varrho} \dotplus M_{\bold x}) \cup 
\{\varphi_*(\bar x_{\bar d,\rho},\bar c_{\bold x,\varrho},\bar c'_1),
\neg \varphi_*(\bar x_{\bar d,\rho},\bar c_{\bold x,\varrho},\bar c'_\rho)\}.
\]
\end{enumerate}
\mn
So let $\bar d'$ realize the type from $\boxplus_5$ and define $\bold
y \in \pK_{\kappa,\mu,\theta}$ by
\mn
\begin{enumerate}
\item[$\boxplus_6$]  $(a) \quad M_{\bold y} = M_{\bold x}$
\sn
\item[${{}}$]  $(b) \quad w_{\bold y} = w_{\bold x} + \{t_*\}$
\sn
\item[${{}}$]  $(c) \quad \bar d_{\bold y} \rest w_{\bold x} = \bar
  d_{\bold x}$ and $\bar d_{\bold y,t_*} = \bar d'$
\sn
\item[${{}}$]  $(d) \quad v_{\bold y} = v_{\bold x} + \{s_*\}$ and
  $u_{\bold y} = u_{\bold x}$
\sn
\item[${{}}$]  $(e) \quad \bar c_{\bold y} \rest v_{\bold x} = \bar
  c_{\bold x}$ and $\bar c_{\bold y,s_*} = \bar c'_1 \char 94 \bar c'_1$
\sn
\item[${{}}$]  $(f) \quad \bold I_{\bold y} = \bold I_{\bold x}$
\sn
\item[${{}}$]  $(g) \quad \bar B_{\bold y,s}$ is equal to $B_{\bold
  x,s}$ if $s \in v_{\bold x} \backslash u_{\bold x}$ and is equal
  to $A_*$ if $s=s_*$.
\end{enumerate}
\mn
Clearly $\bold y$ is as required.
\bigskip

\noindent
\underline{Case 2}:  $\lambda \ge \mu$ is singular.

Impossible by \ref{n10}.
\bigskip

\noindent
\underline{Case 3}:  $\lambda \ge \mu$ regular
\mn
\begin{enumerate}
\item[$\oplus_1$]  \wilog \,
\sn
\begin{enumerate}
\item[$(a)$]  for some $\bar c_1,\bar c_0 \in {}^{\ell g(\bar
  y)}(M_{\bold x})$ we have
\newline
$\gC \models ``(\forall \bar x_{\bar d})
[\varphi_*(\bar x_{\bar d},\bar c_{\bold x},\bar c_1) \wedge \neg
  \varphi_*(\bar x_{\bar d},\bar c_{\bold x},\bar c_0)]"$
\sn
\item[$(b)$]  for every $\bar c_1 \in {}^{\ell g(\bar y)}(M_{\bold
  x})$ for some $\bar c_0 \in {}^{\ell g(\bar y)}(M_{\bold x})$ we
  have
\newline
$\gC \models ``(\forall \bar x_{\bar d})[\varphi_*(\bar x_{\bar d},
\bar c_{\bold x},\bar c_1) \equiv \neg \varphi_*(\bar x_{\bar d},
\bar c_{\bold x},\bar c_2)]"$.
\end{enumerate}
\end{enumerate}
\mn
[Why?  We can use $\varphi'(\bar x_{\bar d},\bar x_{\bar c},\bar y
  \char 94 (\langle y_0,y_1,y_2 \rangle) = [y_0 = y_2 \rightarrow
    \varphi(\bar x_{\bar d},\bar x_{\bar c},y)] \wedge [y_0 \ne y_2
    \wedge y_1 = y_2 \rightarrow \neg \varphi(\bar x_{\bar d},
\bar x_{\bar c},y)]$ so if $B \subseteq M_{\bold x},|B| \ge 2$ 
we can use for $y_0,y_1,y_2$ members of $B$; see more in the proof of 
\ref{k4}.]
\mn
\begin{enumerate}
\item[$\oplus_2$]  $(a) \quad$ let $n_*$ be as $m$ in \ref{n10}(3)
\sn
\item[${{}}$]  $(b) \quad$ let $\varphi_{**} = \varphi_{**}(\bar
  x_{\bar d,\rho}, \bar c_{\bold x,\varrho},\bar y_{**})$ where $\bar y = \bar
  y_0 \char 94 \ldots \char 94 \bar y_{n(*)-1},\ell g(\bar y_\ell) =
  \ell g(\bar y)$ 

\hskip25pt  and $\varphi_{**} =
  \bigwedge\limits_{\ell < n(*)} \varphi_*(\bar x_{\bar d,\rho},\bar
  c_{\bold x,\varrho},\bar y_\ell)$.
\end{enumerate}
\mn
Now
\begin{enumerate}
\item[$\oplus_3$]  $\lambda = \ntr_w(\varphi_{**},\psi,\bold x)$.
\end{enumerate}
\mn
[Why?  Think.]

Now we shall use \ref{n13} + \ref{n15} for $\varphi_{**},\psi$ getting
$d',\bar a_0,\bar a_1,\bold I = \langle (\bar a_{\alpha,0},\bar
a_{\alpha,1}):\alpha < \lambda\rangle$.

Why this suffice?  We choose $\bold y$ by
\mn
\begin{enumerate}
\item[$\oplus_4$]  $(a) \quad M_{\bold y} = M_{\bold x}$
\sn
\item[${{}}$]  $(b) \quad \bar d_{\bold y} = \bar d_{\bold x} \char 94
  \langle \bar d' \rangle$ i.e. $w_{\bold y} = w_{\bold x} \cup
  \{s\},w_{\bold x}, p <_{\bold y} s,\bar d_{\bold y,s} = \bar d'$.
\sn
\item[${{}}$]  $(c) \quad \bar c_{\bold y} = \bar c_{\bold x} \char 94
  (\bar a_0 \char 94 \bar a_1)$, i.e. $\bar a_0 \char 94 \bar a_1 =
  \bar c_{\bold y,t},v_{\bold y} = v_{\bold x} + \{t\}$
\sn
\item[${{}}$]  $(d) \quad \bold I_{\bold y,t} = \bold I$.
\end{enumerate}
\mn
This is this possible?  We just have to check that the relevant
condition in \ref{n13}, i.e. the clauses in (A) holds which is straight.
\end{PROOF}

\begin{conclusion}
\label{n19}
1) For every $\bold x \in \pK_{\kappa,\mu,\theta}$ \underline{there}
   is $\bold y$ such that $\bold x \le_1 \bold y \in
   \pK_{\kappa,\mu,\theta}$ and: if $\varphi_*(\bar x_{\bar
   d,\rho},\bar x_{\bar c,\varrho},\bar y) \in \Gamma^1_{\bold
   x},\psi_*(\bar x_{\bar d,\rho},\bar x_{\bar c,\varrho}) \in
   \tp(\bar d_{\bold x,\rho} \char 94 \bar c_{\bold
   x,\varrho},\emptyset)$ and $\bold y \le_1 \bold z \in
\pK_{\kappa,\mu,\theta}$ then $\bold k_{\du}(\varphi_*,\psi_*,\bold z) 
= \bold k_{\du}(\varphi_*,\psi_*,\bold y)$.

\noindent
2) Above $\bold y \in \uK_{\kappa,\mu,\theta}$, see Definition \ref{c23}(3C).
\end{conclusion}

\begin{PROOF}{\ref{n19}}
By \ref{n13}, recalling \ref{n5}, as in \ref{n10}(4).
\end{PROOF}

\begin{conclusion}
\label{n23}
If $\kappa > \mu \ge \beth_\omega + \theta^+,\theta \ge |T|$ \then \, 
$\vK^\otimes_{\kappa,\mu,\theta}$ is $\le_1$-dense in 
$\sK^\oplus_{\kappa,\mu,\theta}$.  Moreover, if $\bold m \in
\sK^\oplus_{\kappa,\mu,\theta}$ then for some $\bold n$ we
have $\bold m \le_1 \bold n \in \uK^\otimes_{\kappa,\mu,\theta}$.
\end{conclusion}

\begin{PROOF}{\ref{n23}}
Assume $\bold m \in \sK^\oplus_{\kappa,\mu,\theta}$ we apply
\ref{n6}(5) to $\bold x$ getting $\bold y$ as there.  By  \ref{n17} we
have
\mn
\begin{enumerate}
\item[$(*)_1$]  if $\varphi = \varphi(\bar x_{\bar d,\rho},\bar
  x_{\bar c,\varrho},\bar y)$ and $\psi = \psi(\bar x_{\bar
  d,\rho},\bar x_{\bar c,\varrho_0}) \in \tp(\bar d_{\bold
  x,\rho},\bar c_{\bold x,\rho_0})$ and 
$\bold k_{\du}(\varphi,\psi,\bold y) = 0$ 
\then \, $\ntr_w(\varphi,\psi,\bold y) \ge \kappa$.
\end{enumerate}
\mn
Clearly we can find $\bold n$ such that
\mn
\begin{enumerate}
\item[$(*)_2$]  $(a) \quad \bold m \le_1 \bold n \in
  \sK^\oplus_{\kappa,\mu,\theta}$
\sn
\item[${{}}$]  $(b) \quad \bold x_{\bold n} = \bold y$
\sn
\item[${{}}$]  $(c) \quad$ if $\varphi,\psi$ are as in $(*)_1$ then
  $\varphi \in \Gamma^3_{\bar\psi_{\bold n}}$.
\end{enumerate}
\mn
By \ref{n5} + \ref{n13} clearly $\bold n \in
\vK^\oplus_{\kappa,\mu,\theta}$.

The ``moreover" is proved similarly.
\end{PROOF}

\begin{theorem}
\label{n25}

\noindent
\underline{The recounting theorem}  Assume $\kappa = \kappa^{< \kappa}
= \aleph_\alpha = \mu + \alpha \ge \mu \ge 
\beth_\omega + \theta^+,\theta \ge |T|$.

\Then \, for any $M \in \EC_{\kappa,\kappa}(T)$ 
the cardinality of $\bold S^\theta(M)/\equiv_{\aut}$ 
is $\le 2^{< \mu} + |\alpha|^{\theta + |T|}$.
\end{theorem}

\begin{PROOF}{\ref{n25}}
By \ref{n23} and \ref{c39}(2).
\end{PROOF}
\bigskip

\subsection{Exact recounting of types and $\vK$} \
\bigskip

The following analysis look more carefully at decomposition and
$\varphi \in \Gamma^1_{\bold x}$: eventually it was not used in
proving the density of $\vK$.

Here we use $\ind(\varphi)$.
\begin{definition}
\label{q3}  
1) For $\varphi = \varphi(\bar x,\bar y',\bar y)$ let $\ind(\varphi) 
= \ind_T(\varphi) = \min\{n:$ the set
$\{\varphi(\bar x_\eta,\bar y',\bar y_k)^{[\eta(k)]}:\eta \in {}^n 2$
and $k<n\}$ is inconsistent with $T\}$; compare with \ref{b5}(5), \ref{b22}(1).

\noindent
2) Above if $\bar y'$ is the empty sequence we may omit it; we may
 ignore the case $\ind(\varphi)=1$; it is always $\ge 1$.

\noindent
3) For $\bold x \in \pK_{\kappa,\mu,\theta}$ and $\varphi =
\varphi(\bar x_{\bar d},\bar x_{\bar c},\bar y) \in 
\Gamma^1_{\bold x}$ and $k < \ind_T(\varphi)$ let

\noindent
$\Lambda^1_{\bold x,\varphi,k} = \{\psi:\psi = \psi(\bar y_k,
\bar y^+_0 \char 94 \ldots \char 94 \bar y^+_{k-1},y^+_{k+1} \char 94
\ldots \char 94 \bar y^+_{\text{ind}(\varphi)-1};\bar x_{\bar c})$ and for 
each $\ell \in \text{ ind}(\varphi) \backslash \{k\},\bar y_{\ell,0}$ 
or $\bar y_{\ell,1}$ is a dummy in $\psi\}$ 
where we fix $y^+_m = \bar y_{m,0} \char 94 \bar y_{m,1},\ell g
(\bar y_{m,0}) = \ell g(\bar y_m) = \ell g(\bar y_{m,1})$
and in $\bar y_k \char 94 \bar y^+_0 \char 94 \ldots \char 94 
\bar y^+_{\ind(\varphi)-1}$ there is no repetitions.

\noindent
4) $\Lambda^0_{\bold x,\varphi,k} = \{\psi_{\bold
   x,\varphi,\eta,\nu,k}:\psi_{\bold x,\varphi,\eta,\nu,k} =
\psi_{\bold x,\varphi,\eta,\nu,k}(\bar y_k,\bar y^+_0 \char 94
   \ldots \char 94 \bar y^+_{k-1},\bar y^+_{k+1} \char 94 \ldots \char
94 \bar y^+_{\ind(\varphi)-1};\bar x_{\bar c}) =
   (\exists \bar x_{\bar d})(\bigwedge\limits_{m < \ind(\varphi),m
   \ne k} (\varphi(\bar x_{\bar d},\bar x_{\bar c},\bar
   y_{m,\eta(m)})^{[\nu(m)]} \wedge \varphi(\bar x_{\bar d},\bar
   x_{\bar c},\bar y_k)^{[\nu(k)]}))$ where $\eta,\nu 
\in {}^{\ind(\varphi)}2\}$.

\noindent
5) Let $\Omega^0_{\bold x,\varphi,k,\bar c_*} =
\{\psi_{\bold x,\varphi,\eta,\nu,k}(\bar y,\bar c_*,\bar e,
\bar c_{\bold x}):\bar e = \bar e_{k+1} \char 94 \ldots \char 94 
\bar e_{\ind(\varphi)-1}$ and $\bar e_m \in 
{}^{2 \ell g(\bar y)}(M_{\bold x})$ for each $m\}$ where:
$\bar c_* = \bar c^*_{0,0} \char 94 \bar c^*_{0,1} \char 94 \ldots \char 94
\bar c^*_{k-1,0} \char 94 c^*_{k-1,1},\ell g(\bar c^*_{\ell,0}) =
\ell g(\bar c^*_{\ell,1}) = \ell g(\bar y)$.

\noindent
6) For $\Lambda \subseteq \Lambda^1_{\bold x,\varphi,k}$ and we let 

\noindent
$\Omega^1_{\bold x,\varphi,\Lambda,k,\bar c_*} =
\{\psi(\bar y,\bar c_*,\bar e,\bar c_{\bold x}):\psi = \psi(\bar y,
\bar y^+_0 \char 94 \ldots \char 94 \bar y^+_{k-1},\bar y^+_{k+1}
\char 94 \ldots \char 94 \bar y^+_{\ind(\varphi)-1},
\bar x_{\bar c}) \in \Lambda$ and 
$\bar e = \bar e_{k+1} \char 94 \ldots \char 94 
\bar e_{\ind(\varphi)}$ and $\bar e_m \in {}^{2 \ell g(\bar y)}
(M_{\bold x})$ for each $m\}$.
\end{definition}

A relative of \ref{b20}(1) = \cite[2.8=tp25.33]{Sh:900} imitating
$\vK$ is
\begin{definition}
\label{q6}  
Let $\bold x \in \pK_{\kappa,\bar\mu,\theta}$ be
 normal\footnote{This indicates we may forget $\bar c_{\bold x}$ and
 instead have a set of sequences some $\bar d_i$'s which function as
 $\bar c_i$'s, so we have $B_\eta$ or $\bold I_\eta$ but even if $\ell
 g(\eta_\ell)  = \eta_\ell +1,\eta_1(n_1) = \eta_2(n_2)$ we still may
 have $B_{\eta_1} \ne B_{\eta_2}$, etc.} and 
$\varphi = \varphi(\bar x_{\bar d},\bar x_{\bar c},\bar y) \in 
\Gamma^1_{\bold x}$, really $\varphi = 
\varphi(\bar x_{\bar d,\rho},\bar x_{\bar c,\varrho},\bar y)$ for some
 $\rho \in {}^{\omega >}(\ell g(\bar d_{\bold x})$ and $\varrho \in
 {}^{\omega >}(\ell g(\bar c_{\bold x}))$ and $n = \ind(\varphi)$, 
see \ref{q3}(1).

We call $\bold w$ an $(\bold x,\varphi)$-witness when $\bold w =
\langle(\bar c_{k,0},\bar c_{k,1}):\ell < n\rangle = \langle
\bar c_{\bold w,k,0},\bar c_{\bold w,k,1}:k < \bold n_{\bold
w}\rangle$, their concatanation is denoted by $\bar c_{\bold w}$  
and there is $\rho_1$ exemplifying it such that
\mn
\begin{enumerate}
\item[$(a)$]  let $\bar c_k = \bar c_{k,0} \char 94 \bar c_{k,1}$
and $\bar c_{< k} = (\bar c_0,\dotsc,\bar c_{k-1})$
\sn
\item[$(b)$]  $\bar c_{k,0},\bar c_{k,1}$ are 
finite subsequences of some $\bar c_{\bold x,i(k)}$
with $\langle i(k):k < n\rangle$ increasing and $\varrho(\ell) < i(k)$
for $\ell < \ell g(\varrho),k < n$
\sn
\item[$(c)$]  $\bar c_{k,0},\bar c_{k,1}$ satisfies the same
formulas from $\Omega^0_{\bold x,\varphi,k,\bar c_{<k}}$
\sn
\item[$(d)$]   $\rho_1 \in{}^{\ell g(\rho)}\ell g(\bar d_{\bold x})$
and
\sn
\item[$(e)$]  $\bar d_{\bold x,\rho}$ 
and $\bar d_{\bold x,\rho_1}$ satisfies the same formulas 
from $\{\varphi(\bar x_{\bar d,\rho},\bar c_{\bold x,\varrho},\bar b):
\bar b \in {}^{\ell g(\bar y)}(M_{\bold x})\}$
\sn
\item[$(f)$]  $\bar d_{\bold x,\rho_1}$ realizes $q_{\bold w} := 
\{\varphi(\bar x_{\bar d,\rho_1},\bar c_{\bold x,\varrho}),
\bar c_{k,1}) \wedge \neg \varphi(\bar x_{\bar d,\rho_1},\bar c_{\bold
x,\varrho},\bar c_{k,0}): k < n\}$
\sn
\item[$(g)$]  (nec?) $\tp(\bar c_{i(k)},\bar c_{< i(k)} + M_{\bold x})$
is\footnote{First, we can use just $\tp_{\Delta_k}$ for $\Delta_k$
large enough.  Second, does clause (g) follows from the earlier ones?}
finitely satisfiable in $M_{\bold x}$.
\end{enumerate}
\end{definition}

\begin{observation}
\label{q8}
Above in Definition \ref{q6}, $\ell g(\bold w) < \ind_T(\varphi)$.
\end{observation}

\begin{definition}
\label{q10}  In Definition \ref{q6}

\noindent
0) We say $\bold w$ is a maximal $(\bold x,\varphi)$-witness \when \,
it is an $(\bold x,\varphi)$-witness and there is no 
$(\bold x,\varphi)$-witness $\bold w_1$ such that $\bold w \triangleleft
\bold w_1$.

\noindent
1) We say $\bold w$ is a successful $(\bold x,\varphi)$-witness
\when \, it is an $(\bold x,\varphi)$-witness and for every
$(\bold x_1,\bold w_1)$ satisfying $\bold x \le_1 \bold x_1 \in 
\pK_{\kappa,\bar\mu,\theta}$ and $\bold w_1$ an 
$(\bold x_1,\varphi)$-witness $\bold w \trianglelefteq
\bold w_1$ we have $\bold w = \bold w_1$.

\noindent
2) We say $\bold x \in \pK_{\kappa,\bar\mu,\theta}$ is full for
   $(\kappa,\bar\mu,\theta)$ \when \, $\bold x$ is normal and 
for every $\varphi \in \Gamma^1_{\bold x}$ there is a
successful $(\bold x,\varphi)$-witness.
\end{definition}

\begin{remark}
In Definition \ref{q10} we may consider ``every maximal 
$(\bold x,\varphi)$-witness is successful".
\end{remark}

\begin{definition}
\label{q13} 
1) Let $\ntr_{\ve}(\varphi(\bar x_{\bar d,\eta},\bar x_{\bar c,\varrho},
\bar y),\bold w,\bold x)$ where $\bold x \in 
\pK_{\kappa,\bar\mu,\theta}$ and $\bold w$ is a $\varphi(\bar x_{\bar
d,\eta},\bar x_{\bar c,\varrho},\bar y)$-witness, be the 
minimal cardinal $\lambda$ such
that, recalling $q_{\bold w}$ is from Definition \ref{q6}(f);
\mn
\begin{enumerate}
\item[$(*)$]  for every $A \subseteq M_{\bold x}$ of cardinality $<
\lambda$ there is a finite $q_A = q_A(\bar x_{\bar d,\eta}) 
\subseteq \tp_{\pm \varphi}(\bar d_{\bold x,\eta},\bar c_{\bold x}
\dotplus M_{\bold x})$ such that 
$q_A(\bar x_{\bar d,\eta}) \cup q_{\bold w} \vdash \tp_\varphi(\bar
d_{\bold x,\eta},\bar c_{\bold x} \dotplus A)$.
\end{enumerate}
\mn
2) Let $\ntr_{\ve}(\bold x) = \min\{\ntr_{\ve}
(\varphi(\bar x_{\bar d,\eta},\bar x_{\bar c,\varrho},\bar y),
\bold w,\bold x):\varphi = \varphi(x_{\bar d,\eta},\bar x_{\bar
c,\varrho},\bar y)$ and $\bold w$ is a maximal $(\bold
   x,\varphi)$-witness$\}$.  It is regular (see case 2 in the proof of
\ref{q15} below and we can replace finite by ``of cardinality $< n_*$" if
   $\lambda > \aleph_0$, see case 3 there.
\end{definition}

\begin{discussion}
\label{q14}
1) The point is that looking for $q \subseteq 
\tp_\varphi(\bar d_{\bold x,\eta},\bar c_{\bold x} 
\dotplus M_{\bold x})$ enables us to deal with singular $\ntr_{\ve}(\varphi)$.

\noindent
2) Do we really have to change $\bar d_{\bold x,\eta}$ to $\bar
   d_{\bold x,\rho}$ in the definition of $\ntr_{\varphi_v}(\bold
   x,\varphi,\bold w)$? when we succeed, i.e. is it $\kappa$?

Part is a finite subset $q$ of $\tp^\pm_\varphi(\bar d_{\bold
x,\eta},M_{\bold x})$ so $\eta,\rho $ are not distinguished.  But we
have $q_{\bold w}$ is a $\pm \varphi$-type on $\bar c^{\bold
w}_{k,\ell}(k < n_{\bold w},\ell=0,1\}$ and 
$\tp(\bar c_{\bold w},M_{\bold x})$ is definable.
\end{discussion}
\newpage

\section {Indiscernibles} 

\begin{hypothesis}
\label{e0}  
$T$ dependent.
\end{hypothesis}

\begin{theorem}
\label{e1}
Assume $\kappa = \kappa^{< \kappa} > \mu = \beth_\omega +
\theta,\theta \ge |T|$ and $M \in \EC_{\kappa,\kappa}(T)$.

If $\varepsilon < \theta^+$ and $p \in \bold S^\varepsilon(M)$ \then
\, there is an indiscernible sequence $\bold I = \langle \bar
a_\alpha:\alpha < \kappa\rangle$ of $\varepsilon$-tuples from $M$,
i.e. $\bar a_\alpha \in {}^\varepsilon M$ for $\alpha < \kappa$ such
that $p = \Av(\bold I,M)$.
\end{theorem}

\begin{PROOF}{\ref{e1}}
Let $\bar d$ realize $p$ hence for some $\bold x \in
\rK_{\kappa,\kappa,\theta}$ we have $\bar d_{\bold x} = \bar d,\bar
c_{\bold x} = \langle \rangle,v_{\bold x} = \emptyset$.  Let $\bold m =
(\bold x,\langle \rangle,\langle \rangle)$, it $\in
\rK^\oplus_{\kappa,\kappa,\theta}$.   By the density of
$\uK^\oplus_{\kappa,\mu,\theta}$ there is $\bold n \in
\uK^\otimes_{\kappa,\mu,\theta}$ such that $\bold m \le_1 \bold n$,
hence $\bar d \trianglelefteq \bar d_{\bold x[\bold m]}$.  
By \ref{e6} below we are done.
\end{PROOF}
\bigskip

\subsection {Indiscernibility and materializing $\bold m$} \
\bigskip

\begin{claim}
\label{e6}  
1) Assume $\bold m = (\bold x,\bar\psi,r) \in 
\tK^\oplus_{\kappa,\bar\mu,\theta}$ and $M_{\bold x}$ has cardinality
 $\kappa$, \then \, for some $\langle \bar c_\alpha \char 94 
\bar d_\alpha:\alpha \le \kappa + \omega\rangle$ we  have:
\mn
\begin{enumerate}
\item[$(a)$]   $\bar c_\alpha,\bar d_\alpha$ are from $M_{\bold x}$
 and $\bar c_{\bold x} \char 94 \bar d_{\bold x} \char 94 \bar
 c_\alpha \char 94 \bar d_\alpha$ realizes $r$ for $\alpha < \kappa$
\sn
\item[$(b)$]   $\bold I = \langle \bar c_\alpha \char 94 
\bar d_\alpha:\alpha \le \kappa + \omega\rangle$ is an indiscernible
sequence and $(\bar c_\kappa,\bar d_\kappa) = (\bar c_{\bold x},\bar
d_{\bold x})$ 
\sn
\item[$(c)$]  $\tp(M,\cup\{\bar c_\alpha \char 94 
\bar d_\alpha:\alpha < \kappa\} + \bar c_\kappa + \ldots 
+ \bar c_{\kappa + n_*}) \vdash \tp(M,\cup\{\bar d_\alpha:
\alpha \le \kappa + n_*\})$
\sn
\item[$(d)$]   if $A \subseteq M_{\bold x}$ is finite\footnote{we
can say of cardinality $< \kappa$, but for \ref{e8} sake we use this form}
 and $\alpha < \kappa$ is large enough \then \, 
$\tp(\bar c_\alpha \char 94 d_\alpha,A) = 
\tp(\bar c_{\bold x} \char 94 \bar d_{\bold x},A)$ and 
$\tp(\bar d_{\bold x},\bar c_{\bold x} + \bar d_\alpha) \vdash 
\tp(\bar d_{\bold x},\bar c_{\bold x} + A + \bar c_\alpha
 \char 94 \bar d_\alpha)$  according to $\bar\psi$.
\end{enumerate}
\mn
2) For $(\bold m,\bold w) \in \vK^\otimes_{\kappa,\bar\mu,\theta}$ similarly 
\underline{but} replace clause (d) by
\mn
\begin{enumerate}
\item[$(d)'$]  if $A \subseteq M_{\bold x}$ is finite and $\alpha <
\kappa$ is large enough then 
$(\bar c_\alpha,\bar d_\alpha)$ solve $(\bold m,A + 
\cup\{\bar c_\beta \char 94 \bar d_\beta:\beta < \alpha\}$ in the 
$\rK^\oplus_{\kappa,\bar\mu,\theta}$-sense-see Definition \ref{c5}(f).
\end{enumerate}
\end{claim}

\begin{definition}
\label{e8}  
1) We say an indiscernible sequence 
$\bold I = \langle \bar c_s \char 94 \bar d_s:s \in I\rangle$
materialize $\bold m \in \tK^\oplus_{\kappa,\mu,\theta}$ when
in the linear order $I$ there is no last element and for some 
$\bar c_n,\bar d_n$ for $n <\omega$ the sequence $\langle
\bar c_s \char 94 \bar d_s:s \in I + \omega\rangle$ satisfies (a)-(d) of
Claim \ref{e6}, and $(\bar c_n,\bar d_n)$ here standing for $\bar
c_{\kappa +n},\bar c_{\kappa +n}$ there.

\noindent
1A) Similarly for $\vK^\otimes_{\kappa,\bar\mu,\theta}$.

\noindent
2) $\bold I$ is also said to materialize $\bold x$ \when \, this 
holds for some $\bold m$ with $\bold x = \bold x_{\bold m}$.

\noindent
3) We say that $D$ is the ultrafilter of $\bold m \in 
\tK^\oplus_{\kappa,\mu,\theta}$ or just $\bold m \in 
\tK^\otimes_{\kappa,\mu,\theta}$ and we may (see \ref{e13}(3)),
$D = D_{\bold m}$, \when \, $D \in \uf(\bar c_{\bold x[\bold m]} 
\char 94 \bar d_{\bold x[\bold m]},M_{\bold x[\bold m]})$ 
satisfies: for every $A \subseteq M_{\bold m[\bold x]}$ of 
cardinality $< \kappa$ and sequence 
$\bar c' \char 94 \bar d'$ from $M_{\bold m[\bold x]}$, the sequence
   realizes $\tp(D,A)$ \Iff \, $\bar c' \char 94 \bar d'$ realizes
$\tp(\bar c_{\bold x[\bold m]} \char 94 \bar d_{\bold x[\bold
   m]},A)$ and $\bar c_{\bold x[\bold m]} \char 94 \bar d_{\bold
   x[\bold m]} \char 94 \bar c' \char 94 \bar d'$ realizes $r_{\bold
   m}$, recalling Definition \ref{a31}(7).
\end{definition}

\begin{PROOF}{\ref{e6}}
\underline{Proof of \ref{e6}}  
1) Let $\langle a_\alpha:\alpha < \kappa\rangle$ list
 $M_{\bold x}$ and choose $(\bar c_\alpha,\bar d_\alpha)$ in $M$ by
 induction on $\alpha < \kappa$ which solves
$(\bold x,\bar\psi,r)$ over $A_\alpha := \{a_\beta + \bar c_\beta +
\bar d_\beta:\beta < \alpha\}
\cup B^+_{\bold x}$, see clause (f) Definition \ref{c5}.

Next, let $(\bar c_\kappa,\bar d_\kappa) = (\bar c_{\bold x},\bar
d_{\bold x})$.  By \ref{c42} for each $\alpha < \kappa$ 
the sequence $\langle \bar c_\beta \char 94 \bar d_\beta:\beta \in 
[\alpha,\kappa]\rangle$ is indiscernible over $A_\alpha$
and choose $(\bar c_{\kappa +n},\bar d_{\kappa +n})$ for
$n \in [1,\omega)$ such that $\langle \bar c_\beta \char 94 \bar d_\beta:\beta
\in [\alpha,\kappa + \omega)\rangle$ is an indiscernible sequence over
$A_\alpha$ for every $\alpha < \kappa$, possible by compactness, so
clauses (a),(b),(d) of \ref{e6}(1) hold.

We are left with clause (c).  By clause (d) we have $\tp(\bar d_{\bold
x},c_{\bold x},\bar d_\alpha + \bar c_\alpha) \vdash \tp(\bar
d_{\bold x},\bar c_{\bold x} + \bar d_\alpha + \bar c_\alpha +
A_\alpha)$.  Now for stationarily many $\alpha < \kappa$ we have
$\tp(\bar d_{\bold x},\bar c_{\bold x} + \bar d_\alpha + \bar c_\alpha) 
\vdash \tp(\bar d_{\bold x},\bar c_{\bold x} + \bar d_\alpha +
\bar c_\alpha + A_\alpha + \sum\limits_{n < \omega} \bar c_{\kappa +
n})$, otherwise by Fodor Lemma we get contradiction to \ref{b20}(2).
So by indiscernibility we get, for $n < \omega,\beta < \kappa + n$
that $\tp(\bar d_{\kappa + n +1},c_{\kappa +n} + \bar d_\beta + \bar
c_\beta) \vdash \tp(\bar d_{\kappa +n},c_{\kappa +n} +
\sum\limits_{i < \beta} \bar c_i \char 94 \bar d_i + A_{\beta,\kappa}
+ \sum\limits_{m \ge n} \bar c_{\kappa +n})$.

Hence for $n < \omega$ we have $\tp(M_\kappa,
\sum\limits_{\alpha < \kappa +n} \bar c_\alpha \char 94
\bar d_\alpha + \sum\limits_{n < \omega} c_{\kappa +m}) \vdash 
\tp(M_\kappa,\sum\limits_{\alpha \le \kappa +n} \bar c_\alpha \char 94
d_\alpha + \sum\limits_{m < \omega} \bar c_{\kappa +n})$ hence we get
the desired conclusion.

\noindent
2) Similarly.
\end{PROOF}

A variant of \ref{e6}
\begin{claim}
\label{e10}  
If $\bold m = (\bold x,\bar\psi,r) \in \tK^\oplus_{\kappa,\mu,\theta}$ 
with $M_{\bold x}$ of cardinality
$\kappa$ and $I_2 = I_1 \times \bbZ$ ordered lexicographically of
course, $I_1$ is a saturated model 
of $\Th(\bbQ,<)$ of cardinality $\kappa$, \then \, we can 
find $\langle \bar c_s \char 94 \bar d_s:
s \in I_2 + \{\kappa\}\rangle$ such that:
\mn
\begin{enumerate}
\item[$(a)$]  $(\bar c_\kappa,\bar d_\kappa) = 
(\bar c_{\bold x},\bar d_{\bold x})$ 
\sn
\item[$(b)$]   $\langle \bar c_s \char 94 \bar d_s:s \in I_2 +
\{\kappa\}\rangle$ is an indiscernible sequence
\sn
\item[$(c)$]   $M_{\bold x}$ is $|T|^+$-atomic over $\cup\{\bar c_s
\char 94 \bar d_s:s \in I\}$
\sn
\item[$(c)'$]   if $J_2 = J_1 \times \bbZ$, where $J_1$ 
(is a linear order which) extends $I_1$ and 
$\bar c_s,\bar d_s$ for $s \in J_2 \backslash I_2$ are 
such that $\langle \bar c_s \char 94 \bar d_s:s \in J\rangle$ is
indiscernible, \then \, $\tp(M,\cup\{\bar c_s \char 94 \bar d_s\}:
s \in I\}) \vdash \tp(M,\cup\{\bar c_s \char 94 \bar d_s:s \in J\})$
\sn
\item[$(d)$]   if $s \in I_2$ then $\bar d_{\bold x} \char 94 \bar
c_{\bold x} \char 94 \bar d_s \char 94 \bar c_s$ realizes $r$ and for
every $A \in [M_{\bold x}]^{< \kappa}$ for every large enough $t
\in I_2$ we have $\tp(\bar d_{\bold x},\bar c_{\bold x} + \bar d_t
+ c_t) + \tp(d_{\bold x},\bar c_{\bold x} + \bar d_t + c_t +
\sum \limits_{s<t} \bar c_s \char 94 \bar d_s +A)$.
\end{enumerate}
\end{claim}

\begin{remark}
1) In \ref{e10} we cannot use $I_2$ a 
saturated model of $\Th(\bbQ,<)$ as then some 
$b \in M_{\bold x}$ may induce a cut with both cofinalities $> |T|$.

\noindent
2) In \ref{e10} we can replace $\bbZ$ by any linear order
with at least two elements but $< \lambda$.

\noindent
3) Note that if $\bold m \in \tK^\oplus_{\kappa,\bar\mu,\theta}
   \subseteq \vK^\oplus_{\kappa,\bar\mu,\theta}$, then also
$(\bold m,\bold w) \in \vK^\otimes_{\kappa,\bar\mu,\theta}$
for $\bold w$ the ``identity" on $\Gamma^2_{\bold x_{\bold m}}$,
 see Definition \ref{c23}(4C).
\end{remark}

\begin{PROOF}{\ref{e10}}
Let $\langle \bar a_\alpha:\alpha < \kappa\rangle$ list the
finite sequences of $M_{\bold x}$ each appearing stationarily many times.

Let $\langle t_\alpha:\alpha < \kappa\rangle$ list the elements of $I_1$
without repetitions and for technical reasons $\langle t_n:n
<\omega\rangle$ is increasing.

Now we choose $J_{1,\alpha},J_{2,\alpha},
\langle \bar c_s \char 94 \bar d_s:s \in J_{2,\alpha}\rangle$ 
by induction on $\alpha < \lambda$ such that
\mn
\begin{enumerate}
\item[$(a)$]  $J_{1,\alpha}$ is a subset of $I_1$ of cardinality $<
\lambda,\subseteq$-increasing continuous
\sn 
\item[$(b)$]  $J_{2,\alpha} \subseteq J_{1,\alpha} \times \bbZ$ ordered
lexicographically and contains $J_{1,\alpha} \times \{0\}$
\sn
\item[$(c)$]  $\{t_\beta:\beta < \alpha\} \subseteq J_{1,\alpha}$ and
$\{t_\beta:\beta < \alpha\} \times \bbZ \subseteq J_{2,\alpha}$
\sn
\item[$(d)$]  $\bold I_\alpha = \langle \bar c_t \char 94 \bar
d_t:t \in J_{2,\alpha}\rangle$ is indiscernible for $\alpha \ge \omega$
\sn
\item[$(e)$]  $(\bar c_t,\bar d_t)$ solves $(\bold
x,\bar\psi,\cup\{\bar c_s \char 94 \bar d_s:s <_{J_{2,\alpha}} t\}$
for each $t \in J_{2,\alpha}$
\sn
\item[$(f)$]  if $\alpha < \omega$ or $\alpha = 1$ mod $3$
let $\beta(\alpha) < \lambda$ be minimal such 
that $t_{\beta(\alpha)}$ is $<_I$-above $J_{1,\alpha}$ 
\then \, $\bar c_{(t_{\beta(\alpha)},0)} \char 94 \bar
d_{t_{(\beta(\alpha)},0)}$ solves $(\bold x,\bar\psi,r,A_\alpha)$ 
where $A_\alpha := \cup\{\bar c_t \char 94 \bar d_t:
t \in J_{2,\alpha} \cup \{(t_{\beta(\alpha)},m):m < n\}\}$ 
and $J_{1,\alpha +1} = J_{1,\alpha}
\cup\{t_{\beta(\alpha)}\},J_{2,\alpha +1} = J_{2,\alpha} \cup
\{(t_{\beta(\alpha)},n):n \in \bbN\}$
\sn
\item[$(g)$]  if $\alpha = 2 \mod 3,\alpha \ge \omega$ and
let $J_{1,\alpha +1} = J_{1,\alpha} \cup\{t_\gamma:\beta \le
\alpha,t_\beta < t_{\beta(\alpha)}\},
J_{2,\alpha +1} = J_{1,\alpha +1} \times \bbZ$ and choose
$\bar a_t$ for $t \in J_{2,\alpha +1} \backslash J_{2,\alpha}$ (such
that (d) + (e) holds)
\sn
\item[$(h)$]  if $\alpha = 3 \beta \ge \omega$, then we choose
$J_{1,\alpha +1},J_{2,\alpha +1},
\langle(\bar c_s,\bar d_s):s \in J_{2,\alpha +1}
\backslash J_{2,\alpha}\rangle$ such that\footnote{$J_{2,\alpha +1}$
has an infinite end segment included in $J_{2,\alpha}$}, 
if possible, for some finite $I \subseteq I_1 \backslash J_{1,\alpha}$
we have $J_{1,\alpha +1} = J_{1,\alpha}
\cup I,J_{2,\alpha} = J_{1,\alpha} \times \bbZ$ and defining $I^\alpha
\in K_p$ as $(I \times \bbZ,P_s)_{s \in I},P_s = \{s\} \times \bbZ$, 
the sequence $\langle \bar c_s \char 94 \bar d_s:s \in I^\alpha\rangle$ is not
indiscernible over $\bar a_\beta$.
\end{enumerate}
\mn
It is easy to carry the induction.

The main point is to verify clause (c) hence (c)$'$.  By
\cite[3.4]{Sh:715} or see \S(1C), if $\bar a \in {}^n(M_{\bold x})$
and $\varphi = \varphi(\bar x_{\bar d[\bold x]},\bar x_{\bar c[\bold
x]},\bar z_{[n]})$ \then \,
there is an expansion of $I_2$ to $I^+_2 =
(I^+_2,P_0,\dotsc,P_n)$ each $P_\ell$ a (non-empty) convex subset of
$I_2$ such that $\langle \bar c_s \char 94 d_s:s \in I_2)$ is
$\{\varphi\}$-indiscernible over $\bar a$.

\Wilog \, if $t \in I_1,\ell \le n$ and $(\{t\} \times \bbZ) \cap
P_\ell \ne \emptyset \wedge (\{t\} \times \bbZ) \backslash P_\ell \ne
\emptyset$ \then \, $P_\ell \subseteq \{t\} \times \bbZ$ and let
$I^1_{\bar a}$ be the set of such $t$'s.
Let $\alpha < \kappa$ be such that $\ell \le n \Rightarrow P_\ell \cap
J_{2,\alpha} \ne \emptyset$, \wilog \, $J_{2,\alpha} = J_{1,\alpha}
\times \bbZ,\alpha = \omega \alpha$.

By clause (h) of the construction we get that $\tp_\varphi(\bar
a,\{\bar c_s \char 94 \bar d_s:s \in I^1_{\bar a} \times \bbZ\})
\vdash \tp_\varphi(\bar a,\{\bar c_s \char 94 \bar d_s:s \in
I_2\})$, treating $\bar c_s \char 94 \bar d_s$ are singletons, of course.

As this holds for any such $\varphi$ we are done.
\end{PROOF}

\begin{observation}
\label{e12}  
1) If $\bold I = \langle \bar c_s \char 94 \bar d_s:
s \in I\rangle$ materializes $\bold m \in 
\tK^\oplus_{\kappa,\mu,\theta}$ \then \, we can replace $\bold I$ by
$\bold I \rest J$ for any $J \subseteq I$ cofinal in $I$ and $\cf(I) \ge
\kappa$.

\noindent
2) If $\|M_{\bold x}\| = \kappa$ then $\cf(I) = |I|$ is necessarily $\kappa$.
\end{observation}

\begin{remark}
Recall that if $T$ is stable (or just $\bold I$ is an indiscernible
set), necesssarily we get that $\bar d_s$ is algebraic over $\bar c_s$.
\end{remark}

\begin{PROOF}{\ref{e12}}
Straightforward.   
\end{PROOF}

\begin{claim}
\label{e13}  
1) If $\bold m = (\bold x,\bar\psi,r) \in 
\tK^\oplus_{\kappa,\bar\mu,\theta}$ or $\bold m = (\bold
x,\bar\psi,r,\bold u) \in \vK^\otimes_{\kappa,\mu,\theta}$, 
\then \, any two materializations 
$\bold I_1,\bold I_2$ of $\bold m$ are equivalent, see
Definition \ref{a67}(5).

\noindent
2) If $\bold x \in \uK^\otimes_{\kappa,\bar\mu,\theta}$ or
 $\bold x \in \vK^\otimes_{\kappa,\bar\mu,\theta}$ and
$M_{\bold x}$ has cardinality  
$\kappa$, the number of materializations of $\bold x$ up to 
equivalence is $\le 2^\theta$.

\noindent
3) If $\bold m \in \tK^\oplus_{\kappa,\mu,\theta}$ \then \,
there is one and only one $D = D_{\bold m}$, the ultrafilter of
$\bold m$, see \ref{e8}(3).
\end{claim}

\begin{PROOF}{\ref{e13}}
  1) Suppose $\bold I_\ell = \langle \bar c_{\ell,s} \char 94  
\bar d_{\ell,s}:s \in I_\ell\rangle$ is a materialization of $\bold m$
and $\bar c^\ell_n \char 94 \bar d^\ell_n$ be as in Definition \ref{e8},
or see \ref{e6}, for $\ell=1,2$.  We can replace $I_\ell$ by any
cofinal sequence hence \wilog \, $\otp(I_\ell) = \kappa_\ell = 
\cf(\kappa_\ell)$, so by \ref{e12} $\kappa_\ell \ge \kappa > |T|$.
Without loss of generality $\kappa_1 \le \kappa_2$, now we let $I_\ell
= \{t_{\ell,\varepsilon}:\varepsilon < \kappa_\ell\}$ with
$t_{\alpha,\varepsilon}$ being $<_{I_\ell}$-increasing with $\varepsilon$.  

First assume $\bold m \in \tK^\oplus_{\kappa,\mu,\theta}$, so for every
$\alpha < \kappa_1$ for some $h_1(\alpha) < \kappa_2$ we have:
\mn
\begin{enumerate}
\item[$(*)_\alpha$]   $\tp(\bar d_{\bold x},\bar c_{\bold x} + 
\bar d_{2,t_2,\beta}) \vdash \tp(\bar d_{\bold x},
\bar c_{\bold x} + \bar d_{2,t_2,\beta} + \{\bar c_{1,t_1,\varepsilon}
\char 94 \bar d_{1,t_1,\varepsilon}:\varepsilon < \alpha\})$ 
if $\beta \in [h_1(\alpha),\kappa_2)$.
\end{enumerate}
\medskip

\noindent
\underline{Case 1}:  $\kappa_1 < \kappa_2$

Then $\beta(*) = \sup\{h_1(\alpha):\alpha <
\kappa_1\}$ is $< \kappa_2$, so applying $(*)_\alpha$ for every
$\alpha < \kappa_1$, for $\beta = \beta(*)$ we get that $\bar
d_{2,t_2,\beta(*)+1} \char 94 \bar c_{2,t_2,\beta(*)+1}$ realizes
$\tp(\bar d_{\bold x} \char 94 \bar c_{\bold x} \cup \{\bar
c_{1,t_1,\beta} \char 94 \bar d_{1,t_1,\beta}:\beta < \alpha\}$ which is
realized in $M_{\bold x}$ so
we get contradiction.
\medskip

\noindent
\underline{Case 2}:   $\kappa_1 = \kappa_2$

So $h_2:\kappa_2
\rightarrow \kappa_1$ can be defined similarly and let $E = \{\delta <
\kappa_1:\delta$ a limit ordinal such that $\alpha < \delta
\Rightarrow h_1(\alpha) < \delta \wedge h_2(\alpha) < \delta\}$, it is
a club of $\kappa_1$.

Now for any $h:E \rightarrow \{1,2\}$ the sequence $\bold I_h =
\langle \bar a_{t_{h(\alpha),\alpha}}:\alpha \in E\rangle$ is an
indiscernible sequence by Claim \ref{c42} or \ref{c47}.

So all the $\bold I_h$'s and $\bold I_1,\bold I_2$ are equivalent.
Second, assume $(\bold m,\bold u) \in 
\vK^\otimes_{\kappa,\bar\mu,\theta}$, easy too.  The case of $\vK$ is similar.

\noindent
2) For $\bold x \in \tK^\oplus_{\kappa,\bar\mu,\theta}$ 
the number of pairs $(\bar\psi,r)$ such that $(\bold x,\bar\psi,r)
\in \tK^\oplus_{\kappa,\bar\mu,\theta}$ is $\le 2^\theta$, and
   now apply part (1).  Similarly, if $\bold x \in 
\vK_{\kappa,\bar\mu,\theta}$ then the number of triples
$(\bar\psi,r,\bold u)$ such that $(\bold x,\bar\psi,r,\bold u) \in
 \vK^\oplus_{\kappa,\bar\mu,\kappa}$.

\noindent
3) E.g. force by Levy$(\kappa,\|M_{\bold x}\|)$ and use absoluteness.
\end{PROOF}

\begin{definition}
\label{e17}  
Assume $M \in \EC_{\kappa,\kappa}(T)$ and $p \in \bold S^\sigma(M)$.  

Let
\mn
\begin{enumerate}
\item[$(a)$]   $\bbI_p = \{\bold I:\bold I$ is an (endless) indiscernible 
sequence in $M$ with $\Av(\bold I,M) = p\}$
\sn
\item[$(b)$]  $\bbI^\chi_p = \{\bold I \in 
\bbI_p:\bold I\text{ has length } \chi\}$
\sn
\item[$(c)$]  $\bbI^*_p = \bbI^\kappa_p$.
\end{enumerate}
\end{definition}

\begin{definition}
\label{e19}  Assume $M \prec {\gC},\bold m =  
(\bold x,\bar\psi,r) \in \tK^\oplus_{\kappa,\bar\mu,\theta}$ or
$\bold m = (\bold x,\bar\psi,r,\bold u) \in 
\vK^\oplus_{\kappa,\bar\mu,\theta}$ and 
$\gamma < \ell g(\bar d_{\bold x})$ and 
$p=p(\bar x) \in \bold S^\gamma(M_{\bold x})$.
We say that $\bold I$ materializes the quadruple
$(p,\bold x,\bar\psi,r)$ or $(p,\bold x,\bar\psi,r,\bold w)$ 
or in $(p,\bold m)$ \when \,:
\mn
\begin{enumerate}
\item[$(a)$]  $\bold I = \langle \bar b_s \char 94 \bar c_s 
\char 94 \bar d_s:s \in I\rangle$ is an 
indiscernible sequence in $M_{\bold x}$ 
\sn
\item[$(b)$]  $\langle \bar c_s \char 94 \bar d_s:s \in I\rangle$ 
materialize $\bold m = (\bold x,\bar\psi,r)$
\sn
\item[$(c)$]   $\langle \bar b_s:s \in I\rangle \in \bbI_p$
\sn
\item[$(d)$]  \underline{Case 1}:  $\bold m \in 
\tK^\oplus_{\kappa,\bar\mu,\theta}$ for every finite $A \subseteq M_{\bold x}$ 
for every large enough $s \in I$ we have $\tp(\bar d_{\bold x},
\bar c_{\bold x} + \bar d_s) \vdash 
\tp(\bar d_{\bold x},\bar c_{\bold x} \dotplus 
(A + \bar b_s + \bar c_s + \bar d_s))$ according to $\bar\psi$
\bigskip

\noindent
\underline{Case 2}:  $\bold m \in \vK^\oplus_{\kappa,\bar\mu,\theta}$: 
for every finite $A \subseteq M_{\bold x}$ for every large enough 
$s \in I$, the pair $(c_s,d_s)$ solves $(\bold m,A + \bar b_s)$.
\end{enumerate}
\end{definition}

\begin{claim}
\label{e21}  
Assume $M,\bold m,\gamma,p$ are as 
in Definition \ref{e19} and $\|M_{\bold x}\| = \kappa$.

If $\bold I_* = \langle \bar b^*_\alpha:\alpha < \kappa\rangle \in
\bbI^*_p$, \then \, there is $\bold I = \langle \bar b_\alpha
\char 94 \bar c_\alpha \char 94 \bar d_\alpha:\alpha < \kappa\rangle$ 
which materialize $(p,\bold x,\bar\psi,r)$ such that
the sequences $\bold I_*$ and $\langle\bar b_\alpha:
\alpha < \kappa\rangle$ are equivalent (even are equal on
   a stationary set of indices).
\end{claim}

\begin{PROOF}{\ref{e21}}
Let $\langle a_\alpha:\alpha < \kappa\rangle$ list
the members of $M_{\bold x}$.  Now repeat the proof of \ref{e6}
 before choosing $(\bar c_\alpha,\bar d_\alpha)$ in stage $\alpha$, choose
minimal $\gamma(\alpha) < \lambda$ such that 
$\bar b_{\gamma(\alpha)}$ realizes $\Av(\bold I_*,A'_\alpha)$ where
$A'_\alpha := \cup\{\langle a_\beta\rangle \char 94 
\bar b_{\gamma(\beta)} \char 94 \bar c_\beta \char 94 
\bar d_\beta:\beta < \alpha\}$ and choose 
$(\bar c_\alpha,\bar d_\alpha)$ as a solution of $(\bold x,\bar\psi,r)$ over
$A'_\alpha + \bar b_{\gamma(\alpha)}$.

As $\langle \bar b_\alpha:\alpha <\kappa \rangle,\langle \bar c_\alpha
\char 94 \bar d_\alpha:\alpha < \kappa\rangle$ are indiscernible sets,
for some type $r_\alpha$ we have $(\forall^\kappa \beta <
\kappa)(\forall^\kappa \gamma < \beta)[\tp(\bar b_\beta \char 94
\bar c_\gamma \char 94 \bar d_\gamma,\bar c_{\bold x} + \bar d_{\bold
x} + A'_\alpha) = r_\alpha]$, and clearly
$\bar b_{\gamma(\alpha)} \char 94 \bar c_\alpha \char 94 \bar
d_\alpha$ realizes the type and $r_\alpha$ increases with $\alpha$.  

So again by \ref{c42} or \ref{c47}, the sequence $\langle (\bar
b_{\gamma(\alpha)},\bar c_\alpha,\bar d_\alpha):\alpha <
\lambda\rangle$ is indiscernible and also the rest should be clear.
\end{PROOF}

\begin{discussion}
\label{e24}
Recall \ref{a37}(4).
If we replace the type by its $\omega$-th iteration, see \cite{Sh:93},
 i.e. if $\langle \bar d_n:n < \omega\rangle$ is an indiscernible
 sequence witnessing $D \in \uf(\tp(\bar a,M))$ \then \,
$\tp(a,\bar d_0 \char 94 \bar d_1 \char 94 \ldots \char 94 D)$ determine $D$.
\end{discussion}

\begin{definition}
\label{e25} 
1) For $M \prec {\gC}_T$, 
ultrafilter $D$ on ${}^\zeta M$ and $\bold I_D = \langle 
\bar a_n:n < \omega\rangle$ based on $D$ (see Definition \ref{a31}(6)), 
let ${\bold T}_D$ be the set of sequences
$\langle (A_s,\bar a_s,\Delta_s):s \in u\rangle$ such that:
\mn
\begin{enumerate}
\item[$(a)$]   $u$ is an inverted tree with root being maximal
\sn
\item[$(b)$]   $A_s \subseteq M$ finite, increases with $s \in u$
\sn
\item[$(c)$]   $\bar a_s \in {}^\zeta M$
\sn
\item[$(d)$]   finite $\Delta_s \subseteq \Gamma_\zeta$ which is
$\subseteq$-increasing
with $s \in u$ recalling $\Gamma_\zeta = \{\varphi(\bar x_0,
\dotsc,\bar x_n;\bar y):\bar x_\ell = \langle
x_{\ell,\varepsilon}:\varepsilon < \zeta\rangle$ and $y = \langle
y_\ell:\ell < n\rangle$ for some $n\}$ 
\sn
\item[$(e)$]  $\langle \bar a_t\rangle \char 94 \bold I_D$ is
$\Delta_t$-indiscernible over $\cup\{A_s \cup (\bar a_s \rest w):s
<_u t$ and $w \subseteq \gamma$ is the finite set of places not dummy
in $\Delta_s\}$.
\end{enumerate}
\mn
2) For ${\gn} \in {\bold T}_D$ let ${\gn} = \langle A_{\gn,s},
\bar a_{{\gn},s},\Delta_{{\gn},n}):s \in u_{\gn}
\rangle,u[{\gn}] = u_{\gn}$ and let $\max(\gn)$ be the
$\le_u$-maximal member (= root) of $u_{\gn}$ and 
$(A_{\gn},\bar a_{\gn},\Delta_{\gn}) = 
(A_{{\gn},\max({\gn})},\bar a_{{\gn},\max({\gn})},
\Delta_{{\gn},\max({\gn})})$.  Lastly, 
$\bar a_{\gn}$ is $\bar a_{{\gn},\max({\gn})}$.

\noindent
4) If ${\gn} \in \bold T_D$ and $s \in u_{\gn}$ let $u_{\gn}
 \rest (\le s)$ be $u \rest \{s_1 \in u_{\gn}:s_1 
\le_{u_{\gn}} s_1\}$ as a partial order and let ${\gn} 
\rest (\le s)$ be $\langle (A_{{\gn},s_1},\bar a_{{\gn},s_1},
\Delta_{{\gn},s_1}):s_1 \in u \rest (\le s)\rangle$.

\noindent
5) We say ${\bold w} = \langle {\gn}_t:t \in I\rangle$ is a witness for $D$
\when \,:
\mn
\begin{enumerate}
\item[$(a)$]   $I$ is a directed partial order
\sn
\item[$(b)$]  $\bold n_t \in \bold T_D$ for every $t \in I$
\sn
\item[$(c)$]  if $t_1 <_I t_2$ then for some $s \in 
u[{\gn}_{t_2}]$ we have ${\gn}_{t_1} = {\gn}_{t_2} \rest (\le s)$.
\end{enumerate}
\mn
6) In part (5) let $(A_t,\bar a_t,\Delta_t) = (A_{\gx,t},
\bar a_{{\gx},t},\Delta_{{\gx},t})$ denote $(A_{{\gn}_t,
\max({\gn}_t)},\bar a_{{\gn}_t,\max({\gn}_t)},
\Delta_{{\gn}_t,\max({\gn}_t)})$. 
\end{definition}

\begin{claim}
\label{e27}  
1) If $D$ is an ultrafilter on ${}^\zeta M$ \then \, there 
is a witness ${\gx} = \langle {\gn}_t:t \in I\rangle$ for 
$D$, let ${\gn}[t] = {\gn}_t$.

\noindent
2) If $\Delta \subseteq \Gamma_\zeta(\tau_T)$ is finite and $A
 \subseteq M$ is finite \then \, for some $t_0 \in I$, if $k <
\omega$ and $t_0 <_I \ldots < t_k$ then $\langle 
\bar a_{{\gn}[t_\ell]}:\ell \le k\rangle$ is 
$\Delta$-indiscernible over $A$.
\end{claim}

\begin{PROOF}{\ref{e27}}
See \cite[\S1]{Sh:715} or an exercise.
\end{PROOF}
\bigskip

\subsection {Indiscernible existence from bounded directionality} \
\bigskip

We affirm here the conjecture from \S(1C) for the case $k=1$,
for dependent theory $T$ of bounded directionality.  
We state the more informative version
(see Defintion \ref{a80}(1)).

\begin{claim}
\label{e28}
\underline{The Strong Indiscernible Existence Theorem}  
1) Let $T$ be of finite directionality, see Definition \ref{a41}.
Assume $\kappa = \cf(\kappa) > \theta = |\gamma| + |T|,\bar
d_\alpha \in {}^\gamma{\gC}$ for $\alpha \in \kappa$ and $\langle
\bar d_\alpha:\alpha < \kappa\rangle$ is a type-increasing sequence,
see \ref{e31}(0) below, \then \, for some $I \in K_{q,\theta}$
expanding $(\kappa,<)$ the sequence $\langle \bar d_s:s \in I\rangle$
is mod club locally-indiscernible, see Definition \ref{a80}.  

\noindent
2)   Let $T$ be of bounded directionality.  \Then \, we get a similar
 result for $I \in K_{\reg,\theta}$, see Definition \ref{e30d} below.
\end{claim}

\begin{PROOF}{\ref{e29}}
  1) By \ref{e33} + \ref{e36}(2) below.

\noindent
2) By \ref{e33} + \ref{e41} below.
\end{PROOF}

\begin{definition}
\label{e30d}  
1) $K_{\reg,\zeta}$ is the
class of structures $I = (I,<_I,P^I_i,F^I_j)_{i < \theta,j < \theta}$
such that $(I,<)$ is a linear order, $P^I_i$ a unary relation and
$F^I_j$ is a unary function such that $F^I_j(t) \le_I t$; reg stands for
regressive.
 
\noindent
2) Assume $\kappa$ is regular uncountable and $I \in
   K_{\reg,\theta}$ expand $(\kappa,<)$.  We say the sequence
   $\langle \bar d_\alpha:\alpha \in I\rangle$ is mod club locally
   indiscernible \when \, $(\bar d_s \in {}^\zeta{\gC},\zeta < \theta$
   and) for some club $E$ of $\kappa$, 
for every $\eta \in {}^\theta 2,\nu \in {}^\theta(\kappa +1)$ and
   finite $\Delta \subseteq \bbL(\tau_T)$, 
we have: if $S_{\eta,\nu} = \{\alpha \in E:I \models
P_i(\alpha)^{[\eta(i)]}$ for every $i \in u$ and 
$F^I_j(\alpha) = \nu(j) \vee (F^I_j(\alpha) = \alpha \wedge 
\nu(j) = \kappa)$ for every 
$j\in v\}$ is unbounded in $\kappa$ then $\langle \bar
   d_\alpha:\alpha \in S_{\eta,\nu}\rangle$ is $\Delta$-indiscernible.
\end{definition}

\begin{definition}
\label{e31} 
0) We say that $\langle \bar d_\alpha:\alpha < \beta\rangle$ 
is type-increasing over $B$ \when \, $\tp(\bar d_\beta,
\cup\{\bar d_\alpha:\alpha < \beta\} \cup B)$ 
is $\subseteq$-increasing with $\alpha$; if $B = \emptyset$ we may omit it.

\noindent
1) Let $\aK_{\lambda,\kappa,\theta}$ be the class of $\bold x$ consisting of
\mn
\begin{enumerate}
\item[$(a)$]  $\bar M = \langle M_\alpha:\alpha \le \kappa\rangle$,
which is $\prec$-increasing, $\alpha < \kappa \Rightarrow
\|M_\alpha\| < \lambda$
\sn 
\item[$(b)$]  $\bold I = \langle \bar d_\alpha:\alpha \le
\kappa\rangle$ and $\bar d = \bar d_\kappa$ is of length $< \theta^+$
\sn
\item[$(c)$]  $\bar d_\alpha \in {}^{\ell g(\bar d)}(M_\kappa)$
realizes $\tp(\bar d,M_\alpha)$
\sn
\item[$(d)$]  $M_\kappa = \cup\{M_\alpha:\alpha < \kappa\}$
\sn
\item[$(e)$]  $\bar d_{\bold x} = \bar d_\kappa$.

\end{enumerate}
\mn
2) Let $\eK_{\lambda,\kappa,\theta}$, be the class of $\bold x$ consisting of
\mn
\begin{enumerate}
\item[$(a)$]  $\bar M$ as above
\sn 
\item[$(b)$]  $\bar{\bold I} = \langle \bold I_\alpha:\alpha < \kappa
\rangle$ and $\bar d = \bar d_{\bold x}$
\sn
\item[$(c)$]  each $\bar d' \in \bold I_\alpha$ belongs to
   $M_\kappa$ and realizes tp$(\bar d,M_\alpha)$
\sn
\item[$(d)$]  $M_\kappa = \cup\{M_\alpha:\alpha < \kappa\}$.
\end{enumerate}
\mn
3) Let $\aK_{\kappa,\theta} = \aK_{\kappa,\kappa,\theta}$ and
$\eK_{\kappa,\theta} = \eK_{\kappa,\kappa,\theta}$.
\end{definition}

\begin{observation}
\label{e32}  
If $\bold x \in \aK_{\lambda,\kappa,\theta}$ \then \, 
for a unique $\bold y \in \eK_{\lambda,\kappa,\theta}$ we have 
$M_{\bold y,\alpha} =M_{\bold x,\alpha}$ for $\alpha 
\le \kappa$ and $\bold I_{\bold y,\alpha} =
\{\bar a_{\bold x,\alpha}\}$.
\end{observation}

\begin{claim}
\label{e33}  If $\langle \bar d_\alpha:\alpha \le
\kappa\rangle$ is type increasing 
and $\kappa = \cf(\kappa) > |T|$, \then \, there is $\bold
x \in \aK_{\kappa,\theta}$ such that for a club of $\alpha <
\kappa$ we have $\bar d_{\bold x,\alpha} = \bar d_\alpha$.
\end{claim}

\begin{PROOF}{\ref{e33}}
Let $C_\alpha = \cup\{\bar d_\beta:\beta < \alpha\}$ for
$\alpha \le \kappa +1$.
We can find a sequence $\langle a_\alpha:\alpha < \kappa\rangle$ such
that
\mn
\begin{enumerate}
\item[$(a)$]   $\tp(a_\alpha,A_\alpha)$, where $A_\alpha := C_\kappa
\cup\{a_\beta:\beta < \alpha\})$, does not split over some $B_\alpha
\subseteq A_\alpha$ of cardinality $\le |T|$;
\sn
\item[$(b)$]   every finite type
over $A_\alpha,\alpha < \kappa$ is realized by some $a_\beta,\beta <
\kappa$.
\end{enumerate}
\mn
This is possible by \cite[III,7.5,pg.140]{Sh:c} or see 
\cite[4.24=np4.10]{Sh:715}.  So
$A_\kappa$ is well defined and is the universe of a model $M \prec
{\gC}$.

As $\kappa$ is regular, for every $\alpha$ for some $\beta_\alpha <
\kappa$ we have: the type $\tp(\bar d_\beta,A'_\alpha)$ where $A'_\beta
= \cup\{\bar d_\gamma \char 94 \langle a_\gamma\rangle:\gamma <
\beta\})$ is the same for all $\beta \in [\beta_\alpha,\kappa)$, just
consider the definition of non-splitting.

Hence \wilog \, this holds for $\beta \in [\beta_\alpha,\kappa +1)$,
too.  Clearly ${\cU} := \{\delta < \kappa:A'_\delta$ is universe of an
elementary submodel of $M\}$ is a club of $\kappa$.

Define $\bold x$ by letting $M_{\bold x,\kappa} = M$ and for $\alpha <
\kappa$ letting $M_{\bold x,\alpha} = M
\rest A'_{\min}(\cU \backslash \alpha)$ and $\bar d_\alpha = \bar
d_{\min(\cU \backslash \alpha)}$.  Clearly $\bold x$ is as
  required.  
\end{PROOF}

\begin{claim}
\label{e36}  
Assume $T$ is of finite directionality.
Assume $\bold x \in \aK_{< \lambda,\kappa,\theta},\kappa = 
\cf(\kappa) > \theta \ge |T|$ and $\zeta = \ell g(\bar d_{\bold x})$.

\noindent
1) If $\kappa > 2^\theta$ \then \, for some club ${\cU}$
of $\kappa$ and partition $\langle S_i:i < 2^\theta\rangle$ of
${\cU}$, letting $I = ({\cU},<,S_i)_{i < 2^\theta}$ the
sequence $\langle \bar d_{\bold x,\alpha}:\alpha \in I\rangle$ is an 
indiscernible sequence.

\noindent
2) If $\Delta \subseteq \Gamma_\zeta := 
\{\varphi(\bar x_0,\dotsc,\bar x_{n-1};\bar y):\bar x_\ell = \langle
x_{\ell,\varepsilon}:\varepsilon < \zeta\rangle$ and $\bar y = \langle
   y_\ell:\ell < n\rangle$ for some $n\}$ is finite \then \, for
 some club ${\cU}_\Delta$ of $\kappa$ and finite partition
 $\langle P_{\Delta,\ell}:\ell < \ell_\Delta\rangle$ of 
${\cU}_\Delta$ we have:
 $\langle \bar d_{\bold x,\alpha}:\alpha \in ({\cU}_\Delta,
<,P_{\Delta,\ell})_{\ell < \ell_\Delta}\rangle$ is 
$\Delta$-indiscernible in the sense of \ref{e39} below.
\end{claim}

Before we prove, similarly:
\begin{claim}
\label{e38}  
Assume $T$ is of finite directionality.
As in \ref{e36} for $\eK_{\kappa,\theta}$.

In full: assume $\bold x \in \eK_{< \lambda,\kappa,\theta},
\kappa = \cf(\kappa) > \theta
\subseteq |T|,\zeta$ and finite $\Delta \subseteq \Gamma$ where
$\Gamma$ is as in \ref{e36}(2), there are functions $F_n:\cup\{\bold
I_\alpha:\alpha < \kappa\} \rightarrow \kappa$ for
$n,n_\Delta(*),n_\Delta$ large enough (i.e. $\varphi(\bar
x_0,\dotsc,\bar x_{n-1};\bar y) \in \Delta \Rightarrow n < n(*)$ such
that if $\varphi(\bar x_0,\dotsc,\bar x_{n-1};\bar y) \in \Delta,m <
n$ and $\kappa > \alpha^\iota_0 > \ldots > \alpha^\iota_{n-1} \ge
\gamma$ and for $\iota = 1,2$ and $\bar d^\iota_\ell \in
I_{\alpha^\iota_\ell}$ for $\ell < m,\iota = 1,2$ and $k<m \Rightarrow
\alpha^\iota_k >
F_{m-k}(\alpha^\iota_{k(1)},\dotsc,\alpha^\iota_{m-1})$ and $\bar b
\in {}^\zeta(M_{\bold x,\gamma}),\bar d^*$ and $\bar d^*_m,\dotsc,\bar
d^*_{n-1} \in {}^\zeta(M_{\bold x,\gamma})$ \underline{then}

\[
{\gC} \models \varphi[\bar d_{\alpha^1_0},\dotsc,\bar
d_{\alpha^1_{m-1}},\bar d^*_m,\dotsc,\bar d^*_{n-1},\bar b]
\]
\mn
\Iff \,

\[
{\gC} \models \varphi[\bar d_{\alpha^2_0},\dotsc,
\bar d_{\alpha^2_{m-1}},\bar d^*_m,\dotsc,\bar d^*_{n-1},\bar b].
\]
\end{claim}

\begin{definition}
\label{e39}  
For $\Gamma$ as in \ref{e36}(2) and $\Delta
\subseteq \Gamma$ we say $\langle \bar d_\alpha:\alpha <
\alpha(*)\rangle$ is $\Delta$-indiscernible over $A$ \when \, if $m
\le n,\alpha(*) > \alpha_0 > \ldots > \alpha_{m-1}$ and $\alpha(*) >
\beta_0 > \ldots > \beta_{m-1}$ and $\bar b \in {}^{\ell g(\bar y)}A$
and $\bar d^*_\ell \in {}^\zeta A$ for $\ell = m,\dotsc,n-1$ then

\[
{\gC} \models \varphi[\bar d'_{\alpha_0},\dotsc,
\bar d'_{\alpha_{m-1}},\bar d^*_m,\dotsc,\bar d^*_{n-1};\bar b]
\]
\mn
\Iff \,

\[
{\gC} \models \varphi[\bar d'_{\beta_0},\dotsc,
\bar d'_{\beta_{m-1}},\bar d^*_m,\dotsc,\bar d^*_{n-1}].
\]
\end{definition}

\begin{discussion}
\label{e39u}  
Even for singletons we cannot replace
 ``finite" by one in \ref{e36}, because even for 
$T = \Th(\bbQ,<)$, a cut has two cofinalities in general.
\end{discussion}

\begin{claim}
\label{e41}  Let $T$ be of bounded directionality.

\noindent
1) Assume $\bold x \in \aK_{\kappa,\theta}$ and $\ell g(\bar
   d_\alpha) = \zeta$ and $\Delta \subseteq \Gamma_\zeta$ is finite.
\Then \, we can find a club $E$ of $\kappa$ and a regressive
   function $f$ on $E$ such that for every $\gamma <  \kappa$ the
   sequence $\langle \bar d_\alpha:\alpha \in E$ and $f(\alpha) =
   \gamma\rangle$ is $\Delta$-indiscernible or is empty.

\noindent
2) Parallel for \ref{e38}.
\end{claim}

\begin{PROOF}{\ref{e36}}
\underline{Proof of \ref{e36}}  
1) It follows from part (2) as if ${\cU},\langle
P_{\Delta,i}:i < \ell_\Delta,\Delta \subseteq \Gamma$
finite$\rangle$ is as gotten there, we let $E =
\{(\alpha,\beta):\alpha,\beta \in {\cU}$ and $\alpha \in
P_{\Delta,i} \Leftrightarrow \beta \in P_{\Delta,i}$ for every
 finite $\Delta \subseteq \bbL(\tau_T)$ and
$\langle P_i:i < i(*) \le 2^\theta\rangle$ list the $E$-equivalence
classes then $({\cU},<,P_i)_{i<i(*)}$ is as required.

\noindent
2) We prove here also \ref{e21}(2).
We call $\Delta \subseteq \Gamma_\zeta$ cyclically 
closed \when \,: if
   $\varphi(\bar x_0,\dotsc,\bar x_{n-1};\bar y) \in \Delta$ then some
   $\varphi'(\bar x_0,\dotsc,\bar x_{n-1};\bar y) \in \Delta$ is
   equivalent to $\varphi(\bar x_1,\dotsc,\bar x_{n-1},\bar x_0;\bar
   y)$.  Clearly it suffices to deal only with cyclically closed $\Delta$'s.

Let $\bbD_\Delta = \{D \cap \Deef^\zeta_\Delta(M_{\bold x}):D
\in \uf(\tp(\bar d,M_{\bold x,\kappa})\}$ says 
$\langle D_{\Delta,i}:i < \lambda_\Delta\rangle$ list it.
So by \ref{a41} and \ref{a42} we have:
\bigskip

\noindent
\underline{Case 1}:  $T$ of finite directionality.  Then $\bbD_\Delta$ is
finite so $\lambda_D$ is finite.
\bigskip

\noindent
\underline{Case 2}:  Not Case 1 but $T$ of bounded directionality.  
Then $\lambda_D \le \kappa$.

Let $\bar d_\alpha = \bar d_{\bold x,\alpha},\bar d = d_{\bold x,\kappa}$.
Now fix $\zeta = \ell g(\bar d_{\bold x}),\Delta \subseteq
\Gamma_\zeta,\Delta$ finite (cyclically closed).
For each $i < \lambda_\Delta$ choose $D_i \in \uf(\tp
(\bar d,M_{\bold x,\kappa}))$ such that $D_i \cap 
\Deef^\zeta_\Delta(M_{\bold x}) = D_{\Delta,i}$.  Let $\gamma_*$
be $\kappa + \omega$ if $T$ is of finite dimensionality and $\kappa +
\kappa$ otherwise.

Let $J = ([\kappa,\gamma_*),<,P^J_i)_{i < \lambda_\Delta} \in
K_{p,\lambda_\Delta}$ be such that each $P^J_i$ is unbounded.  For $\gamma
\in [\kappa,\gamma_*)$ let $\bold i(\gamma)$ be such that $\gamma$
belongs to $P^J_{\bold i(\gamma)}$ and let $J_\gamma = J
\rest \{\gamma\}$.   Let $J_{\alpha,i} =
(\{\alpha\},<,P^{J_{\alpha,i}}_j)_{j < \lambda_\Delta} \in
K_{p,\lambda_\Delta}$ be such that $P^{J_{\alpha,i}}_i = \{\alpha\}$.

We can choose $\bar d_\gamma$ for $\gamma \in [\kappa,\gamma_*)$ so
  redefining $\bar d_{\bold x}$ such that
$\tp(\bar d_\gamma,\cup\{\bar d_\gamma:\gamma \in [\gamma,\gamma_*)\} 
\cup M_{\bold x,\kappa})$ is equal to $\Av(D_{\bold i(\beta)},
\cup\{\bar d_\beta:\beta \in (\gamma,\gamma_*)\} \cup M_{\bold
  x,\kappa})$.  How?  For any finite $u \subseteq [\kappa,\gamma_*)$
  we can use downward induction and now use general compactness.

For $i < \lambda_\Delta$, let 
$S_i = \{\alpha < \kappa$: the sequence $\langle \bar
d_\varepsilon:\varepsilon \in J_{\alpha,i} + J\rangle$ is
$\Delta$-indiscernible over $M_{\bold x,\alpha}$ and for simplicity
$\alpha \notin \cup\{S_j:j < i\}\}$.

For $\alpha < \kappa$ let $\bold i(\alpha)=
i \Leftrightarrow \alpha \in S_i$ and let $J_\alpha =
J_{\alpha,\bold i(\alpha)}$, so $\bold i(\alpha)$ may be undefined.
 
Now
\mn
\begin{enumerate}
\item[$\boxplus_1$]   if $i < \lambda_\Delta$ and $\gamma < \kappa$ then the
sequence $\langle d_\alpha:\alpha \in S_i \backslash \gamma\rangle
\char 94 \langle \bar d_\gamma:\gamma \in [\kappa,\gamma_*)$ and 
$\bold i(\gamma) =
i\rangle$ is $\Delta$-indiscernible over $M_{\bold x,\gamma}$.
\end{enumerate}

Moreover
\mn
\begin{enumerate}
\item[$\boxplus_2$]   Assume $\gamma < \kappa$, let
${\cW}_{\Delta,\gamma} = \bigcup\limits_{i < \lambda_\Delta} 
S_i \cup[\kappa,\gamma_*)
\backslash \gamma,J_{\Delta,\gamma} = \Sigma\{J_\alpha:\alpha \in
{\cW}_{\Delta,\gamma}\} + J \in K_{p,\lambda_\Delta}$ let $J_\Delta =
J_{\Delta,0},I_{\Delta,\gamma} = \Sigma\{J_\alpha:\alpha \in 
{\cW}_{\Delta,\gamma}\}$ and $I_\Delta = I_{\Delta,0}$.
\Then \, the sequence $\langle \bar d_\alpha:\alpha \in
J_{\Delta,\gamma}\rangle$ is $\Delta$-indiscernible 
over $M_{\bold x,\gamma}$.
\end{enumerate}
\mn
[Why?  Without loss of generality consider only $\gamma \in 
{\cW}_{\Delta,0}$.
Let $\varphi = \varphi(\bar x_0,\dotsc,\bar x_{n-1};\bar y) \in
\Delta$.  We now prove by
induction on $m$, the statement for $m$ simultaneously for all $\gamma
< \kappa$.  That is
\mn
\begin{enumerate}
\item[$(*)$]    if $n \ge m$ and $\varphi(\bar x_0,\dotsc,\bar
  x_{n-1},\bar y) \in \Delta$ and 
$J_\Delta \models \alpha > \alpha^\iota_0 > \ldots
> \alpha^\iota_{m-1} \ge \gamma$ for $\iota = 1,2$ and $k< m \wedge i
  < \ell_\Delta \Rightarrow 
[\alpha^1_k \in P^{J_\Delta}_i \leftrightarrow \alpha^2_k \in
P^{J_\Delta}_i]$ and $\bar b \in {}^{\ell g(\bar y)}(M_{\bold x,\gamma})$ and
$\bar d^*_m,\dotsc,\bar d^*_{n-1} \in {}^\zeta(M_{\bold x,\gamma})$
then ${\gC} \models \varphi[\bar d_{\alpha_0},\dotsc,
\bar d_{\alpha^1_{m-1}},\bar d^*_m,\dotsc,\bar d^*_{n-1},
\bar b^*]$ iff ${\gC} \models \varphi[\bar d_{\alpha^2_0},\dotsc,
\bar d_{\alpha^2_{m-1}},\bar d^*_m,\dotsc,\bar d^*_{n-1},\bar b]$.
\end{enumerate}
\mn
We prove this by induction on $|\{\alpha^\iota_k:k<m$ and $\iota =
1,2\} \cap \kappa|$.  If it is zero this should be clear by the
use of ultrafilters.  If not, let $(\iota,k)$ be such that
$\alpha^\iota_k \notin \kappa$ and $\iota + 2k$ maximal.

Let $\beta \in [\kappa,\gamma_*)$ be such that $\bold
i(\alpha^\iota_k) = \bold i(\beta)$.  Easily $\langle \bar
d_\beta:\beta \in J\rangle$ is indiscernible over $M_{\bold x}$ so
\wilog \, $\{\alpha^{\iota(1)}_{k(1)}:\iota(1) \in \{1,2\}$ and $k(1)
< m\} \cap J$ is\footnote{would be easier if we choose $J$ with no
first member} disjoint to $[\kappa,\beta +1)$.  But now
 note that replacing $\alpha^\iota_k$ by $\beta$ 
does not change the truth value.  So $\boxplus_2$ hence $\boxplus_1$
indeed holds.]

Clearly $\boxplus_1 + \boxplus_2$ are nice but will say nothing if,
e.g. each $S_i$ is empty.
\mn
\begin{enumerate}
\item[$\boxplus_3$]  the set $S := \kappa 
\backslash \cup\{S_i \backslash i:i < \lambda_\Delta\}$ is non-stationary.
\end{enumerate}
\mn
[Why?  Toward contradiction assume $S$ is a stationary subset of
$\kappa$. For each $\delta \in S$ and $i < \lambda_\Delta$ we know $\delta
\notin S_i$, hence there are $n_{\delta,i} = n(\delta,i)$ and 
formula $\varphi_{\delta,i}(\bar x_0,\dotsc,\bar x_{n_{\delta,i}-1};
\bar y_{\delta,i}) \in \Delta$ and 
$\bar d^*_{\delta,i,1},\dotsc,\bar d^*_{\delta,i,n-1} \in 
\cup\{\bar d_\gamma:\gamma \in [\kappa,\gamma_*)\}$ and
$\bar b_{\delta,i} \in 
{}^{\ell g(\bar y_{\delta,i})}(M_{\bold x,\delta})$ such that 
${\gC} \models \neg \varphi_{\delta,i}[\bar d_\delta,\bar d^*_{\delta,i,1},
\dotsc,\bar d^*_{\delta,i,n-1},\bar b_{\delta,i}]$ but
$\varphi^*_{\delta,i}(\bar x_0,\bar d^*_{\delta,i,1},\dotsc,
\bar d^*_{\delta,i,n-1};\bar b_{\delta,i}) \in \Av(D_i,
\cup\{\bar d_\gamma:\gamma \in [\kappa,\delta_*)\} + M_{\bold
x,\delta})$.  Only finitely many of the
members of $\bar d^*_{\delta,i,\ell}$ matter 
say $\bar d^*_{\delta,i,\ell} \rest v_{\delta,i,\ell},v_{\delta,i,\ell}$ 
finite.
\bigskip

\noindent
\underline{Case 1}:  $\lambda_\Delta$ is finite 

Let $C_\delta = \cup\{\Rang(\bar d^*_{\delta,i,\ell} \rest
v_{\delta,i,\ell}:\ell < n_{\delta,i}$ and $i < \lambda_\Delta\} \cup
\{\Rang(\bar b_{\delta,i,\ell})):i < \ell_\Delta\}$ so it is finite.

Also $C_\delta \subseteq \cup\{\bar d_\gamma:\gamma \in
[\kappa,\gamma_*)\} \cup M_{\bold x,\delta}$ and $\cup\{\bar
  d_\gamma:\gamma \in [\kappa,\gamma_*)\}$ has cardinality $\le |T|$.

Hence by Fodor lemma for some $C_* \subseteq M_{\bold x}$ 
the set $S' = \{\delta \in S:C_\delta = C_*\}$ is a 
stationary subset of $\kappa$.  The number of
possibilities for $\langle(n_{\delta,i},\varphi_{\delta,i})\rangle \char
94 \langle \bar d_{\delta,i,\ell} \rest v_{\delta,i,\ell}:
i < \lambda_\Delta,\ell < n\rangle$ is $\le |T|$ and the number of
possible $\langle \bar b_{\delta,i}:i < \lambda_\Delta\rangle$ is
finite so for some stationary $S'' \subseteq S'$ for every $\delta
\in S''$ we have $n_{\delta,i} = n_{*,i},\varphi_{\delta,i} =
\varphi_{*,i,v_{\delta,i,\ell}} = v_{*,i,\ell},
\bar d_{\delta,i,\ell} \rest v_{\delta,i,\ell} = \bar d_{*,i,\ell},
\bar b_{\delta,i} = \bar b_{*,i}$.

Let $\cD$ be an ultrafilter on $\kappa$ to which $\kappa$ and every
club of $\kappa$ belong as well as $S''$.

Let $D' = \{\bold I \subseteq {}^\zeta(M_{\bold x,\kappa})$: the set
$\{\alpha < \kappa:\bar d_\alpha \in \bold I\}$ belongs to $\cD\}$,
clearly $D'$ is an ultrafilter on ${}^\zeta(M_{\bold x,\kappa})$.  As
tp$(\bar d,M_{\bold x,\kappa}) = \cup\{\tp(\bar d_\alpha,
M_{\bold x,\alpha}):\alpha < \kappa\}$ clearly $D' \in 
\uf(\tp(\bar d_{\bold x},M_{\bold x,\kappa}))$ or pedantically $D' \cap
\Deef({}^\zeta(M_{\bold x})) \in \uf(\tp(\bar d_{\bold x},M_{\bold
x,\kappa}))$ hence we can find $i < \lambda_\Delta$ 
such that $\cD_i \cap \Deef_\Delta(M_{\bold x,\kappa}) = 
\cD' \cap \Deef(M_{\bold x,\kappa})$.  But this is a
contradiction to $\{\bar d_\delta:\delta \in S''\} \in \cD$ and the
choice of $\varphi_{\delta,i}(\bar x_0,\bar d^*_{*,i,1},\dotsc,\bar
d^*_{*,i,n_{*,i}-1},\bar b_{*,i})$.
\bigskip

\noindent
\underline{Case 2}:  $\lambda_\Delta$ is infinite.

For $\alpha < \kappa +1$ let $M^+_{\bold x,\alpha}$ be 
$(M_{\bold x,\alpha})_{[\bar d]}$, also $E := \{\delta < \kappa:M^+_{\bold
x,\delta} \prec M^+_{\bold x,\kappa}\}$ is a club of $\kappa$.  For
each $\delta \in E$ choose $D_\delta \in \uf(\tp(\bar
d_{\bold x},M_{\bold x,\delta}))$ and choose $\bar d_{\delta,n} \in
{}^\zeta{\gC}$ for $n < \omega$ such that $\bar d_{\delta,n}$
realizes $\Av(D_\delta,\cup\{\bar d_{\delta,m}:m \in (n,\omega) + M_{\bold
x,\delta}\}$.  As $T$ has bounded directionality 
for each $\varphi = \varphi(\bar
x_0,\dotsc,\bar x_{n(\varphi)-1};\bar y_\varphi) \in \bbL(\tau_T)$
and $\delta \in E$ there are formulas $\psi_\delta(\bar
y_\varphi,\bar z_{\varphi,\delta}) \in \bbL(\tau_{M^+_{\bold
x,\kappa}})$ and $\bar c_{\varphi,\delta} \in 
{}^{\ell g(\bar z_{\varphi,\delta})}(M_{\bold x,\delta})$ such that
\mn
\begin{enumerate}
\item[$\odot$]   for $\bar b \in {}^m(M_{\bold x,\delta})$ we have:
${\gC} \models \varphi[\bar d_{\delta,0},\dotsc,\bar
d_{\delta,n-1},\bar b]$ \Iff \, $M_{\bold x,\delta} \models
\psi_{\varphi,\delta}[\bar b,\bar c_{\varphi,\delta}]$.
\end{enumerate}

Note that if: $\varphi(\bar x_0,\dotsc,\bar x_{n_\varphi-1};\bar y)
\in \Delta$ and $m < n_\varphi$ \then \, 
\newline
$\varphi(\bar x_0,\dotsc,\bar
x_{m-1};\bar x_m \char 94 \ldots \char 94 \bar x_{n(\varphi)-1} \char
94 \bar y) \in \Delta$.

For transparency \wilog \, $\tau_T$ is countable, $\zeta < \omega$ let
$\langle \Delta_n:n < \omega\rangle$ be $\subseteq$-increasing with
union $\Gamma_\zeta$ and $\Delta_0 = \Delta$ and each $\Delta_n$
finite.  For induction on $n$,
for some stationary $S'_n \subseteq S \cap E$ we have $\delta \in S_n
\wedge \varphi \in \Delta_n \Rightarrow 
\psi_{\varphi_\delta}(\bar y_\varphi,\bar
z_{\varphi,\delta}) = \psi_{\varphi,*}(\bar y_\varphi,\bar
z_{\varphi,*})$ and $m < n \Rightarrow S_n \subseteq S_m$.  Let $\cD$ be
a uniform ultrafilter on $\kappa$ such that $n < \omega \Rightarrow
S_n \in \cD_*$, let $\langle \bar d_{*,n}:n < \omega\rangle$ realize
$\Av(\cD_*,\left<\langle \bar d_{\delta,n}:n < \omega\rangle:\delta \in
S\right>,M_{\bold x})$.  Easily $\langle \bar d_{*,n}:n < \omega\rangle$ is
indiscernible over $M_{\bold x}$, each $\bar d_{*,n}$ realizes
$\tp(\bar d_{\bold x},M_{\bold x})$ and $\tp(\bar d_{*,n},\cup\{\bar
d_{*,m},m \in (n,\omega)\} + M_{\bold x})$ is finitely satisfiable in
$M_{\bold x}$, so it is based on some $D \in \uf(\tp(\bar
d_{\bold x},M_{\bold x}))$, so for some $i < \lambda_\Delta$ we have
$D_i \cap \Deef^\zeta_\Delta(M_{\bold x}) = D \cap 
\Deef^\zeta_\Delta(M_{\bold x})$.  We easily get a contradiction to the
choice of $S$ as disjoint to $(S_i \backslash (i+1))$ and $S_1$.
\end{PROOF}
\newpage

\section {Applications} 
\bigskip

\subsection {The generic pair conjecture/On uniqueness 
of $(\kappa,\sigma)$-limit models} \
\bigskip

We now return to the $(\kappa,\sigma)$-limit model conjecture and
generic pair conjecture for $\kappa$.

We shall not deal with the first, only represent it.  The second, the
generic pair conjecture was solved in \cite{Sh:900} for $\kappa > |T|$
measurable.  Here we solve it for $\kappa = \kappa^{<
\kappa} > |T| + \beth_\omega$, it is the case $\xi = 1$ in Definition 
\ref{e71}.

Note that even under GCH the picture is somewhat cumbersome when: 
$\kappa = \chi^+ = 2^\chi >
|T| + \beth_\omega$ and $\chi$ strong limit singular.  It is natural to
restrict ourselves to $S^{\kappa^+}_{\ged}$ (see \cite{Sh:108}).
We may still like to deal with $|T| < \kappa < \beth_\omega$.

Presently, the proof is complete only for $\xi =1$, i.e. the
generic pair conjecture.

Now we rephrase the conjecture; the use of $2^\lambda = \lambda^+$ (in
addition to $\lambda = \lambda^{< \lambda}$) is for transparency only
as an equivalent version without it is absolute under forcing with
Levy$(\lambda^+,2^\lambda)$, see \S1.

\begin{definition}
\label{e71} 
1) We say that $T$ satisfies the uniqueness of limit models above $\mu$ 
\when \, for any $\mu$-complete forcing
notion $\bbQ$ in $\bold V^{\bbQ}$ and $\xi < \lambda$
we have $(A) \Rightarrow (B)_\xi$, see below.  Omitting $\mu$ means
$\mu = |T| + \beth_\omega$.

\noindent
2) For regular $\lambda > |T|$ and ordinal $\xi < \lambda$
we say that $T$ satisfies the
$(\lambda,\xi)$-limit uniqueness \when \, for every $\lambda$-complete
 forcing notion $\bbQ$ such that $\bold V^{\bbQ} \models
 ``\lambda = \lambda^{< \lambda} \wedge 2^\lambda = \lambda^+"$ \
 clause $(B)_\xi$ holds.

\noindent
3) We can add above ``for the trivial $\bbQ$" or other
restrictions.  Instead ``for the trivial $\bbQ$" we may say
``presently" 

\noindent
\underline{where}
\mn
\begin{enumerate}
\item[$(A)$]  $(a) \quad \lambda = \lambda^{<\lambda}$ 
and $2^\lambda = \lambda^+ \ge \mu$
\sn
\item[${{}}$]  $(b) \quad$ density for 
$\vK^\otimes_{\lambda,\lambda,\theta}$ holds for every $\theta <
\lambda$, see \S5 
\sn
\item[${{}}$]  $(c) \quad \langle M_\alpha:\alpha < \lambda^+\rangle$ is a
$\prec$-increasing continuous chain of

\hskip25pt  models of cardinality $\lambda$ with union $M$, 
a saturated model of

\hskip25pt  cardinality $\lambda^+$
\sn
\item[$(B)_\xi$]  for some club ${\cU}$ of $\lambda^+$, if 
$\langle \alpha_{\ell,\varepsilon}:\varepsilon \le \xi\rangle$ is an
increasing continuous

\hskip25pt sequence of ordinals from ${\cU}$ for $\ell=1,2$ such that 

\hskip25pt [$\varepsilon < \xi$ non-limit $\Rightarrow
\alpha_{\ell,\varepsilon}$ of cofinality $\lambda$] \then \, there is

\hskip25pt an isomorphism $\pi$ from $M_{\alpha_{1,\xi}}$ onto
$M_{\alpha_2,\xi}$ mapping

\hskip25pt  $M_{\alpha_{1,\varepsilon}}$ onto
$M_{\alpha_{2,\varepsilon}}$ for every $\varepsilon \le \xi$.
\end{enumerate}
\end{definition}

\noindent
We now translate the relevant questions represented in \S0 to this definition.
\begin{observation}
\label{e73}  
Assume $T$ is dependent.

\noindent
0) If $|T| < \lambda = \lambda^{< \lambda}$, \then \, $T$ has
$(\lambda,0)$-uniqueness (even for the trivial forcing).

\noindent
1) Assume $|T| < \lambda = \lambda^{< \lambda}$ and 
$2^\lambda = \lambda^+$.  \Then \, $T$ has
$(\lambda,1)$-uniqueness, for trivial forcing \Iff \, $T$ satisfies
the generic pair conjecture \Iff \, in $(B)_1$ of \ref{e71} above, if
$\alpha_1 < \beta_1,\alpha_2 < \beta_2$ are all from ${\cU}$ and
has cofinality $\lambda$ then $(M_{\beta_1},M_{\alpha_1}) \approx
(M_{\beta_2},M_{\alpha_2})$.

\noindent
2) Assume $|T| < \lambda = \lambda^{< \lambda}$ and
$2^\lambda = \lambda^+$ and $\sigma = \cf(\sigma) \in
[\aleph_0,\lambda]$.  Then $T$ has uniqueness of
$(\lambda,\sigma)$-model \Iff \, $T$ has
$(\lambda,\sigma)$-limit-uniqueness for the trivial forcing.
\end{observation}

\begin{theorem}
\label{e74}  
$T$ satisfies the generic pair conjecture for $\lambda$ \when \, 
$\lambda = \lambda^{< \lambda} > |T|^+ + \beth^+_\omega$.
\end{theorem}

\begin{remark}
\label{e77}
This is closed to the proof from \cite{Sh:900} as we could restrict
ourselves to $\bold x$ with $u_{\bold x} = \emptyset$.
\end{remark}

\begin{PROOF}{\ref{e74}}
By older works, we can assume $T$ is dependent.
  Without loss of generality $2^\lambda = \lambda^+$ by
absoluteness, see \cite{Sh:877}.

So let $\langle M_\alpha:\alpha < \lambda^+\rangle$ be
given, $M = \cup\{M_\alpha:\alpha < \lambda^+\}$.  Let $E$ be the
set of limit $\delta < \lambda^+$ such that:
\mn
\begin{enumerate}
\item[$\circledast_\delta$]  $(a) \quad$ for every $\alpha <
   \delta$ for some $\beta \in (\alpha,\delta)$ the model $M_\beta$ is
   saturated
\sn
\item[${{}}$]  $(b) \quad$ if $\alpha < \beta < \delta,\zeta <
  \lambda,\{\bar b_1,\bar b_2\} \subseteq {}^\zeta(M_\beta)$ and there
  is an automorphism $g$ of $M$ such that $g(M_\alpha) =
  M_\alpha,g(\bar b_1) = (\bar b_2)$ \then \, there is such $g$
  mapping $M_\delta$ onto itself.
\end{enumerate}
\mn
So
\mn
\begin{enumerate}
\item[$(*)_0$]  $(a) \quad E$ is a club of $\lambda^+$
\sn
\item[${{}}$]  $(b) \quad$ if $\alpha \in E$ has cofinality $\lambda$
\then \, $M_\alpha$ is saturated.
\end{enumerate}
\mn
[Why?  As $M$ is saturated and $\lambda = \lambda^{< \lambda}$.] 
\mn
\begin{enumerate}
\item[$(*)_1$]   if $\alpha < \lambda^+$ and
$M_\alpha$ is saturated and $\bold m_1,\bold m_2 \in 
\vK^\otimes_{\lambda,\lambda,<\lambda}$ satisfies
$M_{\bold m_1} = M_\alpha = M_{\bold m_2}$ and 
$\bold m_1 \le_1 \bold m_2$ and $\bar c_{\bold m_1} \char 94 
\bar d_{\bold m_1}$ is from $M$
\then \, there is an automorphism $g$ of $M$ over $B^+_{\bold m_2}$
mapping $M_\alpha$ onto itself such that
$g(\bar c_{\bold m_1} \char 94 \bar d_{\bold m_1}) = 
\bar c_{\bold m_1} \char 94 \bar d_{\bold m_1}$.
\end{enumerate}
\mn
[Why?  See uniqueness of $M_{[\bold x]}$ in \ref{c32}, see Definition
  \ref{b5}(6).]

Fix $\alpha_1 < \beta_2,\alpha_2 < \beta_2$, all from $E$ and of
cofinality $\lambda$ and we have to prove just that
$(M_{\beta_1},M_{\alpha_1}) \cong (M_{\beta_2},M_{\alpha_2})$.
Let AP be the set of triples $(\bold m_1,\bold m_2)$ satisfying:
\mn
\begin{enumerate}
\item[$(*)_2$]  $(a) \quad \bold m_\ell \in 
\vK^\otimes_{\lambda,\lambda,<\lambda}$, and $r_{\bold m_\ell}$ is
complete
\sn
\item[${{}}$]  $(b) \quad M_{\bold x[\bold m_\ell]} =
M_{\alpha_\ell}$
\sn
\item[${{}}$]  $(c) \quad \bar c_{\bold x[\bold m_\ell]} \char 94
\bar d_{\bold x[\bold m_\ell]} \subseteq M_{\beta_\ell}$
\sn
\item[${{}}$]  $(d) \quad g$ is an elementary mapping with domain 
$B^+_{\bold x[\bold m_\ell]} + \bar c_{\bold x[\bold m_2]} 
+ \bar d_{\bold x[\bold m_1]}$
\sn
\item[${{}}$]  $(e) \quad g$ maps $\bold m_1$ onto $\bold m_2$.
\end{enumerate}
\mn
Let the two place relation $\le_{\AP}$ on $\AP$ be
\mn
\begin{enumerate}
\item[$(*)_3$]   $(\bold m_1,\bold m_2,g) \le_{\AP} (\bold
n_1,\bold n_2,h)$ \Iff \, both triples are from $\AP$, and $g \subseteq
h$ and $\bold m_1 \le_1 \bold n_1,\bold m_2 \le_1 \bold n_2$.
\end{enumerate}
\mn
Now
\mn
\begin{enumerate}
\item[$(*)_4$]  $\AP \ne \emptyset$.
\end{enumerate}
\mn
[Why?  Use $\bold m_\ell$ which is empty except $M_{\bold m_\ell} =
M_{\alpha_\ell}$, see \ref{c27}(3).]
\mn
\begin{enumerate}
\item[$(*)_5$]  if the sequence $\langle(\bold m_{1,\varepsilon},
\bold m_{2,\varepsilon},g_\varepsilon):\varepsilon <
\zeta\rangle$ is $\le_{\AP}$-increasing and $\zeta$ is a limit
ordinal $< \lambda$ \then \, this sequence has a
$\le_{\AP}$-lub, its union $(\bold m_{1,\zeta},
\bold m_{2,\zeta},g_\zeta)$, i.e.
\begin{enumerate}
\item[$(a)$]  $\bold x_{\bold m_\ell,\zeta} = \cup\{\bold
x_{\bold m_{\ell,\varepsilon}}:\varepsilon < \zeta\}$ for $\ell=1,2$
\sn
\item[$(b)$]   similarly for $\bar\psi_{\bold m_{\ell,\varepsilon}}$
\sn
\item[$(c)$]  similarly for r$_{\bold m_{\ell,\varepsilon}}$
\sn
\item[$(d)$]  $g_\zeta = \cup\{g_\varepsilon:\varepsilon < \zeta\}$.
\end{enumerate}
\end{enumerate}
\mn
[Why?  See \ref{c76}.]

The main point is
\mn
\begin{enumerate}
\item[$(*)_6$]  if $(\bold m_1,\bold m_2,g) \in \AP$ and
$\ell \in \{1,2\}$ and $A \subseteq M_{\beta_\ell}$ has cardinality $<
\lambda$ \then \, for some $(\bold n_1,\bold n_2,h) \in \AP$
which is $\le_{\AP}$-above $(\bold m_1,\bold m_2,g)$ we have $A
\subseteq B_{\bold n_\ell} + \bar c_{\bold n_\ell} + \bar d_{\bold
n_\ell}$.
\end{enumerate}
\mn
[Why?  By symmetry we can assume $\ell=1$.  Now trivially we can find
$\bold x \in \pK_{\lambda,\lambda,<\lambda}$ such that $\bold
x_{\bold m_1} \le_1 \bold x$ and $A \subseteq 
\Rang(\bar d_{\bold x[\bold m_1]})$.  By \ref{n23} there is $\bold n'_1 \in
\rK^\oplus_{\lambda,\lambda,<\lambda}$ such that $\bold m_1 \le_1
\bold n'_1$ and $\bold x \le_1 \bold x_{\bold n'_1}$.  

Let $C_\ell = \bar d_{\bold x[\bar{\bold m}_\ell]}
 + \bar c_{\bold x[\bold m_\ell]} + B^+_{\bold x[\bold m_\ell]}$.  
Now recall that by \ref{c32}, the model
$(M_{\alpha_1})_{[C_1]}$ is $(\lambda,\bold D_\ell)$-sequence
homogeneous and moreover $g$ induces a
mapping from $\bold D_1$ onto $\bold D_2$, because 
$g$ maps $\bold m_1$ to $\bold m_2$.  
So there is an isomorphism $f$ from $M_{\alpha_1}$ onto 
$M_{\alpha_2}$ such that $f \cup g$ is an elementary mapping 
(of ${\gC}$), hence it can be extended to an automorphism 
$f^+$ of $M$.  Now $(\bold n_1,f^+(\bold n_1),
f^+ \rest (B_{\bold n_1} + \bar c_{\bold n_1} +
\bar d_{\bold n_1}))$ is almost as required but $f(\bar c_{\bold n_1}
\char 94 \bar d_{\bold n_1})$ is $\subseteq M$ rather than $\subseteq
M_{\beta_2}$.  But $\beta_2 \in E$ hence the definition of $E$ we can
finish.]

Now by $(*)_4 + (*)_5 + (*)_6$ we can find a
$\le_{\AP}$-increasing sequence $\langle (\bold
m_{1,\varepsilon},\bold m_{2,\varepsilon},g_\varepsilon):\varepsilon <
\lambda\rangle$ such that: for any $A_1 \subseteq M_{\beta_1},A_2
\subseteq M_{\beta_2}$ of cardinality $< \lambda$ 
for some $\varepsilon < \lambda$ we have
$A_\ell \subseteq B_{\bold m_{\ell,\varepsilon}} + \bar c_{\bold
m_{\ell,\varepsilon}} + \bar d_{\bold m_{\ell,\varepsilon}}$ for
$\ell=1,2$.

So $g_\lambda = \cup\{g_\varepsilon:\varepsilon < \lambda\}$ is an
isomorphism as required.  
\end{PROOF}

\begin{discussion}
\label{e75}  
1) So we know that $T_2 = \Th(M_{\alpha_0},M_{\alpha_1})$ for
  every $\alpha_0 < \alpha_1$ of cofinality $\lambda$ from ${\cU}$,
is a complete theory and does not depend in 
$(\alpha_0,\alpha_1)$ and even on $\lambda$.
  But we may like to understand it better, see Kaplan-Shelah \cite{KpSh:946}.

\noindent
2) Still $M_{[\alpha_0,\alpha_1]} = (M_{\alpha_0},M_{\alpha_1})$ is
close to being sequence-homogeneous.
So this leads us to deal with dependent finite diagrams $\bold D$.  
Because if we like to deal with
$(\lambda,\zeta)$-uniqueness we have to look at
$(M_{\alpha_0},M_{\alpha_1})$ for any $\bar a \in {}^{\lambda >}
(M_{\lambda^+})$.
\end{discussion}
\bigskip

\noindent
\centerline {$* \qquad * \qquad *$}
\bigskip

\subsection {Criterion for saturativity} \
\bigskip

\begin{claim}
\label{e80} 
Assume $\sigma > \mu = (2^{|T|})^+ + \beth^+_\omega$.

Then $M$ is $\sigma$-saturated \Iff \,
\mn
\begin{enumerate}
\item[$(a)$]  $M$ is $\mu$-saturated
\sn
\item[$(b)$]  if $\kappa \in [\mu,\sigma)$ and $\langle
a_\alpha:\alpha < \kappa\rangle$ is an indiscernible sequence in $M$
\then \, for some $a \in M$ the sequence $\langle a_\alpha:\alpha <
\kappa\rangle \char 94 \langle a \rangle$ is indiscernible
\sn
\item[$(c)$]  if $\kappa \in [\mu,\sigma)$ is regular, $\langle
a_s:s \in I_1 + I_2\rangle$ is an indiscernible sequence in $M$ where
$I_1 \cong (\kappa,<),I_2 \cong (\alpha,>)$ for some $\alpha \le
\kappa +1$ \then \, for some $a \in M$ the sequence $\langle a_s:s
\in I_1 \rangle \char 94 \langle a \rangle \char 94 \langle a_t:t \in
I_2\rangle$ is an indiscernible sequence.
\end{enumerate}
\end{claim}

\begin{PROOF}{\ref{e80}}
The ``only if" implication is obvious.  For the ``if" direction assume
(a),(b),(c) and we prove that $M$ is $\kappa^+$-saturated by induction
on $\kappa \in [\mu,\sigma)$; clearly this suffices.  By clause (a) of
the assumption the model $M$ is $\mu$-saturated.
So by the induction
hypothesis $M$ is $\kappa$-saturated.  Let $A_* \subseteq M$ be of
cardinality $\kappa,p_* \in \bold S(A_*)$ and we should prove that
$p_*$ is realized in $M$.  Let $\bold x \in 
\pK_{\kappa,\mu,\theta},\theta = |T|,\bar d_{\bold x} = \langle
d\rangle$ where $d$ realizes $p_*,M_{\bold x} = M,v_{\bold x} = 0$.  

Now
\mn
\begin{enumerate}
\item[$\boxplus_1$]  if $m=1,\bold I = \langle \bar a_\alpha:\alpha
< \kappa \rangle$ is an indiscernible sequence of $m$-tuples from $M$
and $A \subseteq M$ have cardinality $\le \kappa$ \then \, the type
Av$(\bold I,A)$ is realized in $M$.
\end{enumerate}
\mn
[Why?  Choose $b_\alpha \in M$ for $\alpha < \kappa$ such that $A
\subseteq \{b_\alpha:\alpha < \kappa\}$ and let $A_\alpha = \cup\{\bar
a_\beta \char 94 \langle \bar b_\beta \rangle:\beta < \alpha\}$ for
$\alpha \le \kappa$.  Let $\{\varphi_\varepsilon(x,\bar c_\varepsilon):
\varepsilon < \kappa\}$ list the type $q = 
\Av(\bold I,A_\kappa)$ and for $\bar a \in {}^m{\gC}$ define
$\varepsilon(\bar a)$ as min$\{\varepsilon \le \kappa$: if
$\varepsilon < \kappa$ then ${\gC} \models \neg
\varphi_\varepsilon[\bar a,\bar c_\varepsilon]\}$.  We try to choose $\bar
a'_\alpha,\varepsilon_\alpha$ by induction on $\alpha < \kappa$ such that
\mn
\begin{enumerate}
\item[$(*)$]  $(a) \quad \bar a'_\alpha$ realizes $p_\alpha := 
\Av(\bold I,\cup\{\bar a'_\beta:\beta < \alpha\})$
\sn
\item[${{}}$]  $(b) \quad$ if $\alpha$ is even then $\varepsilon(\bar a')$ is
minimal, i.e. $\varepsilon(\bar a') \le \varepsilon(\bar a'')$
 whenever 

\hskip25pt $\bar a''_\alpha$ realizes $\Av(\bold I,\cup\{\bar a'_\beta:
\beta < \alpha\})$
\sn
\item[${{}}$]  $(c) \quad$ if $\alpha$ is odd then ${\gC} \models
\varphi[\bar a'_\alpha,\bar c_{\varepsilon(\bar a'_{\alpha -1})}]$.
\end{enumerate}
\mn
We can choose $\bar a'_\alpha$ satisfying clause (a) as $p_\alpha$ is
an $m$-type in $M$ of cardinality $< \kappa$ and 
$M$ and is $\kappa$-saturated. 

If $\alpha$ is even it follows trivially that we can satisfy clause
(b), too.  If $\alpha$ is odd, and $\varepsilon(\bar a'_{\alpha-1}) =
\kappa$ then $\bar a'_{\alpha -1}$ is as required, i.e. so we are done
proving $\boxplus_1$, so assume $\varepsilon(\bar a'_{\alpha -1}) <
\kappa$ hence also $p_\alpha
\cup\{\varphi_{\varepsilon(\bar a'_\alpha)}(\bar x,\bar
c_{\varepsilon(\bar a'_{\alpha -1})})\}$ is well defined and being a
subset of $q$ it is an $m$-type in $M$ hence is realized in $M$, 
and any $\bar a'_\alpha$ realizing it is O.K.

Having carried the definition, clearly $\langle \bar a_\alpha:\alpha <
\kappa\rangle$ is an indiscernible sequence; also
by clause (b) of the theorem there is
$\bar a'_\kappa \in {}^m M$ such that $\langle \bar a'_\alpha:\alpha
\le \kappa\rangle$ is an indiscernible sequence.  If $\bar a'_\kappa$
realizes $q$ we are done, if not choose $\varepsilon < \kappa$ such
that ${\gC} \models \neg\varphi_\varepsilon[\bar a'_\kappa,\bar
c_\varepsilon]$.  So for every even $\alpha < \kappa,\bar a'_\kappa$
satisfies clause (a) hence $\varepsilon(\bar a'_\alpha) \le
\varepsilon$.  So for some $\zeta \le \varepsilon$ the set $\{\alpha <
\kappa:\alpha$ even and $\varepsilon(\bar a'_\alpha) = \zeta\}$ is
infinite.  But by $(*)(b)+(c)$ this is a contradiction to
``$\langle \bar a'_\alpha:\alpha < \kappa\rangle$ is an indiscernible
sequence" from the beginning of the paragraph and 
``$T$ is dependent".  So $\boxplus_1$ holds.]
\mn
\begin{enumerate}
\item[$\boxplus_2$]  if $B \subseteq A \subseteq M,|B| < \kappa,|A|
\le \kappa,m=1$ and $p$ is an $m$-type over $A$ which is finitely
satisfiable in $B$, \then \, $p$ is realized in $M$.
\end{enumerate}
\mn
[Why?  Let $D$ be an ultrafilter on ${}^m B$ such that $p \subseteq
\Av(D,A)$.  Let $\langle A_\alpha:\alpha < \kappa\rangle$ be
$\subseteq$-increasing with union $A$ such that $|A_\alpha| \le
|\alpha|$.  Choose $\bar a_\alpha \in {}^m M$ by induction on $\alpha
< \kappa$ such that $\bar a_\alpha$ realizes $\Av(D,\{\bar a_\beta:
\beta < \alpha\} \cup A_\alpha \cup B)$.  So $\bold I = \langle
\bar a_\alpha:\alpha <\kappa\rangle$ is an indiscernible sequence
(\cite[I,\S1]{Sh:c} or see \cite[\S1]{Sh:715}) and
$p \subseteq \Av(\bold I,A)$, hence by $\boxplus_1$ we are done
proving $\boxplus_2$.]
\mn
\begin{enumerate}
\item[$\oplus_1$]    if $\kappa$ is singular, $p_*$ is realized
in $M$.
\end{enumerate}
\mn
[Why?  Let $\langle A^*_\varepsilon:\varepsilon < 
\cf(\kappa)\rangle$  be $\subseteq$-increasing with union $A_*$ such
that $|A^*_\varepsilon| < \kappa$ for $\varepsilon < \kappa$.  As $M$
is $\kappa$-saturated for each $\varepsilon < \kappa$ there is
$a_\varepsilon \in M$ realizing $p_* \rest A_\varepsilon$.  Let $B =
\{a_\varepsilon:\varepsilon < \cf(\kappa)\}$ and let $D$ be an
ultrafilter on $B$ such that $\varepsilon < \kappa \Rightarrow
\{a_\zeta:\zeta \in (\varepsilon,\cf(\kappa)\} \in D$.  Clearly
$p_* \subseteq \Av(D,A_*)$ hence by $\boxplus_2$ we are done.]
\mn
\begin{enumerate}
\item[$\oplus_2$]    if $\kappa$ is regular, \then \, $p_*$ is realized in $M$.
\end{enumerate}
\mn
[Why?  Let $\langle A^*_\alpha:\alpha < \kappa\rangle$ be 
$\subseteq$-increasing with union $A_*$ such
that $|A^*_\alpha| < \kappa$ for $\alpha < \kappa$.  Let $\bold m = (\bold
y,\bar \psi,r,\bold u) \in \vK^\otimes_{\kappa,\mu,\theta}$ be such
that $\bold x \le_1 \bold y$, exists by \ref{n23}.  We can choose
$(\bar d_\alpha,\bar c_\alpha)$ by induction on
$\alpha < \kappa$ such that it solves $(\bold m,A_\alpha \cup \{\bar d_\beta
\char 94 \bar c_\beta:\beta <  \alpha\} \cup B^+_{\bold y})$.  By
\ref{c42} the sequence $\langle \bar c_\alpha \char 94
\bar d_\alpha:\alpha < \kappa\rangle$ is an indiscernible sequence over
 $B^+_{\bold y}$.

Let $d'_\alpha = d_{\alpha,0}$ so $\bold I = \langle d'_\alpha:\alpha
< \kappa\rangle$ is an indiscernible sequence and $d'_\alpha$ realizes
$p_* \rest A_\alpha$.  Hence $\Av(\bold I,A_*)$ is equal to $p_*$.  So
by $\boxplus_1$ we are done.]  

By $\oplus_1 + \oplus_2$ we are done.
\end{PROOF}

\noindent
Another result of interest is (compare with \ref{e1})
\begin{conclusion}
\label{e83}
Assume $\vK_{\kappa,\bar\mu,\theta}$ is dense and $\varepsilon <
\theta^+$.

If $M$ is a $\kappa$-saturated model \then \, for any $p \in \bold
S^\varepsilon(M)$ there is a $\kappa$-complete filter on
${}^\varepsilon|M|$ which is an ultrafilter when restricted to
$\Deef_\varepsilon(M)$, see Definition \ref{a31}. 
\end{conclusion}
\newpage

\section {Concluding Remark}

Another relative of \ref{d21} is
\begin{claim}
\label{d23}
1) Assume $\varphi_n(\bar x_{[\zeta]},\bar y_n)$ is a formula for $n
   <n_*$.  If $\cD$ is a filter on $I$ and $\bar a_t \in {}^\zeta \gC$
   for $t \in I$, \then \, there is $\bar{\cS}$ such that
\mn
\begin{enumerate}
\item[$(a)$]  for some $k_*,\bar{\cS} = \langle \cS_k:k < k_*\rangle$ is a
partition of $I$
\sn
\item[$(b)$]  $\cS_k \in \cD^+$
\sn
\item[$(c)$]  if $\ell< n_*$ and $\bar b \in {}^{\ell g(\bar
y_\ell)}\gC$ \then \, for some truth value $\bold t$ and $k < k_*$ we have
$\{t \in \cS_k:\gC \models \varphi_\ell[\bar a_t,\bar b]^{\iif(\bold t)}\} =
\cS_k \mod \cD$.
\end{enumerate}
\mn
1A) Above we find $\bar S$ such that
\mn
\begin{enumerate}
\item[$(a),(b)$]  are as there and
\sn
\item[$(c)$]  if $\bar b_n \in {}^{\ell g(\bar y_\ell)}\gC$ for
  $n<n_*$ \then \, for some $k$ we have
\sn
\begin{enumerate}
\item[$\bullet$]  for each $n<n_*$ for some truth value $\bold t$, the
  set $\{t \in \cS_k:\gC \models \varphi_n[\bar a_t,\bar
  b_n]^{\iif(\bold t)}\}$ is $= \cS_k \mod \cD$.
\end{enumerate}
\end{enumerate}
\mn
2) Assume $\cD_\ell$ is a filter on $I_\ell$ for $\ell=0,1$ and
$C \subseteq \gC_T,\Delta \subseteq \bbL(\tau_T)$ are finite 
and $\bar a_{\ell,t} \in {}^{m(i)}\gC$ for $t
\in I_\ell,\ell < 2$.  \Then \, we can find $\cS_\ell  \in D^+_\ell$ for
$\ell < 2$ such that for some $q$ we have $(\forall^{\cD_0} s_0 \in
\cS_0)(\forall^{\cD_1} s_1 \in \cS_1)
[q = \tp_\Delta(\bar a_{0,s_0} \char 94 \bar
a_{1,s_1},C)]$.

\noindent
3) Like part (2) for $\langle (I_\ell,\cD_\ell):\ell < n_*\rangle$. 
\end{claim}

\begin{PROOF}{\ref{d23}}
1) We try to choose $n_\ell,\bar b_\ell,\bar{\cS}_\ell$ by induction
   on $\ell \in \bbN$ such that
\mn
\begin{enumerate}
\item[$\boxplus$]  $(a) \quad n_\ell < n_*$
\sn
\item[${{}}$]  $(b) \quad \bar b_\ell$ has length 
$\ell g(\bar y_{n_\ell})$
\sn
\item[${{}}$]  $(c) \quad \bar{\cS}_\ell = \langle \cS_\eta:\eta \in
{}^{\ell +1}2\rangle$ is a partition of $I$
\sn
\item[${{}}$]  $(d) \quad \cS_\eta \in \cD^+$ for $\eta \in {}^\ell 2$
\sn
\item[${{}}$]  $(e) \quad \cS_\eta = \{t \in I$: if $k < 
\ell g(\eta)$ then $\gC \models ``\varphi_{n_k}
[\bar a_t,\bar b_k]^{\iif(\eta(k))}"\}$.
\end{enumerate}
\mn
We stipulate $\bar{\cS}_{-1} = \langle \cS_{<>}\rangle,\cS_{<>} = I$.

If we succeed, we get a contradiction to ``$T$ is dependent".  Arriving to
$\ell$, clearly $\bar{\cS}_{\ell -1}$ has been defined, and if we cannot
choose $n_\ell,\bar b_\ell$ are required, the conclusion of part (1)
holds. 

\noindent
1A)  Similarly; e.g. \wilog \, $\zeta$ is finite from a failure we get that for
every $k$ we can find $A,|A| \le (\Sigma\{\ell g(\bar y_n):n < n_*\}
\times k,|\bold S^\zeta_{\{\varphi_n:n < n_*\}}(A)| \ge 2^k$,
contradiction to ``$T$ dependent" (see \ref{h9}(b)).

\noindent
2) Let $\Delta = \{\varphi^1_n(\bar x_{[m(0)]},\bar y_{[m(1)]},\bar
   z_n):n < n_*\}$ and $\Phi = \{\varphi^1_n(\bar x_{[m(0)]},\bar
   y_{[m(1)]},\bar c_k):n <n_*,\bar c \in {}^{\ell g(\bar z_n)}C)$,
   it is finite and clearly it suffices to deal with one pair
   $(\varphi^1_n,\bar c)$, as we can replace $\cD_\ell$ by $\cD_\ell +
\cS$ when $\cS \in \cD^+_\ell$ and let
$\varphi^2_n = \varphi^1_n(\bar y_{[m(1)]},\bar x_{[m(0)]},\bar z_n) =
\varphi^1_n(\bar y_{[m(0)]},\bar y_{[m(1)]},\bar z_n)$.  We apply part
   (1) with $m(n),1,\varphi^2_n,\langle \bar a_t \char 94 \bar c:
t \in I_\ell\rangle,\cD_1$ here standing for
   $\zeta,n_*,\varphi_n,\langle \bar a_t:t \in I\rangle,\cD$ there and
   get $\bar{\cS}_1 = \langle \cS_{1,k}:k < k_*\rangle$ as there.  We
   define a funtion $h:I_0 \rightarrow \{0,\dotsc,k_*-1\}$, by $h(s) =
\min\{k:(\forall^{\cD_1} t \in I_1)\varphi^2_n(\bar a_{1,t},\bar
   a_{0,s},\bar c)$ or $(\forall^{\cD_1} t \in I_1)(\neg
   \varphi^2_n(\bar a_{1,t},\bar a_{0,s},\bar c))\}$.  By the choice
of $\bar{\cS}$, this is a well defined function.  Clearly for some
   $k$ and $\bold t$, the set 
$\cS_0 := \{s \in I_1:h(s)=k$ and $(\forall^{\cD_1} t \in
   I_1)[\varphi^2_n(\bar a_{1,1},\bar a_{0,s},\bar c)^{\iif(\bold
   t)}]\}$ belongs to $\cD^+_0$
   and let $\cS_1 = \cS_{1,k}$, clearly we are done.
\end{PROOF}

Here we look again at decomposition as in \cite{Sh:900}, i.e. with
$u_{\bold x} = \emptyset$.
\begin{claim}
\label{k2}
Assume $\Delta = \{\varphi_*(\bar x,\bar y),\neg \varphi_*(\bar x,\bar
y)\},m = \ell g(\bar x)$ and $m < \omega$ and $n_* = \ind(\varphi_*)$.  
For any $A(\subseteq \gC)$ and $p \in \bold
S^m_\Delta(A)$ is consistent with $r_*(\bar x)$ and $\mu > \aleph_0$
we can find the following objects:
\mn
\begin{enumerate}
\item[$(A)$]  $(a) \quad \bar d \in {}^m \gC$ realizing $p$
\sn
\item[${{}}$]  $(b) \quad n_1 < \ind(\varphi(\bar x,\bar y))$ 
\sn
\item[${{}}$]  $(c) \quad A_n \subseteq A$ has cardinality $< \mu$ for $n<n_1$
\sn
\item[${{}}$]  $(d) \quad \bar b_{n,0},\bar b_{n,1} \in {}^{\ell
  g(\bar y)} \gC$ for $n < n_1$
\sn
\item[${{}}$]  $(e) \quad D_n$ is an ultrafilter on ${}^{\ell g(\bar
  y)}(A_n)$
\sn
\item[${{}}$]  $(f)(\alpha) \quad \bar b_{n,0} \char 94 \bar b_{n,1}$ realizes
  $\Av(D_n,\{\bar b_{k,0},\bar b_{k,1}:k<n\} + A$
\sn
\item[${{}}$]  $\qquad (\beta) \quad$ if $\ell \le n$ then 
$\bar b_{\ell,0},\bar b_\ell$ hence realizes the same type over

\hskip35pt  $\{b_{k,\iota}:k \le n,k \ne \ell,\iota < 2\} +A$ 
\sn
\item[${{}}$]  $\qquad (\gamma) \quad$ if $\eta,\nu \in
  {}^{n+1}2,\langle \bar b_{0,\eta(0)},\bar b_{1,\eta(1)},\dotsc,\bar
  b_{n,\eta(n)}\rangle,\langle \bar b_{0,\nu(0)},b_{1,\nu(1)},
\ldots\rangle$ 

\hskip35pt realizes the same type over $A$
\sn
\item[${{}}$]  $(g) \quad p \cup r^*_{n_1}$ is consistence where
  $r^*_k(\bar x) = \{\varphi(\bar x,\bar b_{\ell,1}) \equiv \neg
\varphi(\bar x,\bar b_{\ell,0}):\ell < k\}$
\sn
\item[${{}}$]  $(h) \quad \bar d$ realizes the type from (g)
\sn
\item[$(B)$]  $(a) \quad$ if $q \subseteq p$ has cardinality $< \mu$
  then for some finite $r \subseteq p$ we have 

\hskip25pt $r \cup r_* \vdash q$
\sn
\item[${{}}$]  $(b) \quad$ for some $n_2$ depending on $p$ only, we can demand

\hskip25pt  $|r| = n_2$ so $\{r_* \cup r:r \subseteq p,|r| = n_2\}$ is a
$\mu$-directed partial ordered  

\hskip25pt by $r_1 \le r_2 \Leftrightarrow (r_2 \vdash r_1)$.
\end{enumerate}
\end{claim}

\begin{remark}
\label{k3}
1) This is a relative of a claim from \cite{Sh:900}.  We lose not
   fixing $\bar d$ a priori but can use e.g. finite $\Delta$.

\noindent
2) We can chose $D_n$ such that if $B \supseteq A_n$ and $\bar b'_0
   \char 94 \bar b_1$ realizes $\Av(D_n,B)$ then $\tp(\bar b_0,B) =
   \tp(\bar b_1,B)$, moreover the two natural projections of $D_n$ to
   an ultrafilter on ${}^{\ell g(\bar y)}(A_n)$ are equivalent.

\noindent
3) If we are analyzing $\tp(\bar d,A)$ and already have $\bar c$ as in
   decompositions, w can work in $\gC^*_{\bar c} =
   (C,c_\alpha)_{\alpha < \ell g(\bar c)}$ and use $\varphi' =
   \varphi(\bar x_{\bar d},\bar c,\bar y),A'=A$ and apply the claim.

\noindent
4) This may be used in \S5.
\end{remark}

\begin{PROOF}{\ref{k2}}
We try to choose $(A_n,D_n,\bar b_{n,0},\bar b_{n,1})$ by induction on $n$ such
that
\mn
\begin{enumerate}
\item[$\boxplus$]  clauses $(c),(d),(e),(f)(\alpha),(\beta)$ of (A) 
of the assumption holds as well as (g) of (A), i.e. 
$p \cup r^*_{n+1}$ is consistent where $r^*_{n+1}$ is from clause (A)(g).
\end{enumerate}
\mn
Note that $p \cup r^*_0$ is consistent by an assumption.
\bigskip

\noindent
\underline{Case 1}:  We can carry the induction for $n <
\ind(\varphi)$.  

We get a contradiction to the definition of $\ind(-)$
as in \cite{Sh:900}.
\bigskip

\noindent
\underline{Case 2}:  We are stuck in $n_1$ (i.e. cannot choose for
$n_1$)
\mn
\begin{enumerate}
\item[$\oplus_1$]  clause (B)(a) holds.
\end{enumerate}
\mn
Why?  Toward contradiction, let 
$q(\bar x) \subseteq p(\bar x)$ be of cardinality $< \mu$ be a
counterexample so let $q(\bar x) = \{\varphi_\alpha(\bar x,\bar
b_\alpha):\alpha < \mu_*\}$ where $\mu_* = |q(\bar x)|$.

For any finite $r \subseteq p$ let

\[
\cU_r = \{\alpha < \mu_*:r(\bar x) \cup r^*_{n_1} \nVdash
``\varphi_\alpha(\bar x,\bar b_\alpha)" \text{ and } r(\bar x) \cup
r^*_{n_1} \nVdash ``\neg \varphi_\alpha(\bar x,\bar b_\alpha)"\}
\]

\[
\cU^1_r = \{(\alpha,\beta) \in \mu_*: r(\bar x) \cup r^*_{n_1} \cup
\{\varphi_\beta(\bar x,\bar b_\beta),\neg \varphi_\alpha(\bar x,\bar
b_\alpha)\} \text{ is consistent}\}.
\]

\mn
Clearly
\mn
\begin{enumerate}
\item[$(*)_1$]  $\cU_r \ne \emptyset$.
\end{enumerate}
\mn
[Why?  Otherwise recalling 
$r \subseteq p,r$ is as promised in (B) of the claim.] 

So $\varphi_\alpha(\bar x,\bar b_\alpha) = \varphi_*(\bar x,\bar
b_\alpha)^{\iif(\bold t(\alpha)}$ for some truth value $\bold t(\alpha)$.

Let $\bar y_\ell = \langle y_{\ell,k}:k < \ell g(\bar y)\rangle$.

Let
\mn
\begin{enumerate}
\item[$(*)_2$]  $\Gamma_* = \{\psi(\bar y_1,\bar c) \equiv \psi(\bar
  y_2,\bar c):\psi = \psi(\bar y,\bar z) \in \bbL(\tau_T)$ and $\bar c
  \in {}^{\ell g(\bar z)}(\Sigma\{b_{n,\iota}:n < n_1,\iota < 2\}+A)\}$.
\end{enumerate}
\mn
Now
\mn
\begin{enumerate}
\item[$(*)_3$]  for any finite $\Gamma \subseteq \Gamma_*$ let
\newline
$\cU^2_\Gamma = \{(\alpha,\beta) \in \mu_* \times \mu_*:(\bar
b_\alpha,\bar b_\beta)$ realizes $\Gamma$ and $\bold t(\alpha) = \bold
t(\beta)\}$.
\end{enumerate}
\mn
Now
\mn
\begin{enumerate}
\item[$(*)_4$]  if $\Gamma \subseteq \Gamma_*$ is finite then 
$\cU^2_\Gamma$ is an equivalence relation on $\mu_*$ with $\le
  2^{|\Gamma|+1}$ equivalence class.
\end{enumerate}
\mn
[Why?  By inspection.]
\mn
\begin{enumerate}
\item[$(*)_5$]  if $r(\bar x) \subseteq p(\bar x)$ and 
$\Gamma \subseteq \Gamma_*$ are finite then 
$\cU^1_r \cap \cU^2_\Gamma \ne \emptyset$.
\end{enumerate}
\mn
[Why?  As $\cU^2_\Gamma$ has $\le 2^{|\Gamma|+1}$ equivalence classes, we can
find a sequence $\langle \alpha(j):j < 2^{|\Gamma|+1}\rangle$ of
ordinals $< \mu_*$ represent all the $\cU^2_\Gamma$-equivalence
classes.  Let $r_1(\bar x) = r(\bar x) \cup
\{\varphi_{\alpha_j}(\bar x,\bar b_{\alpha(j)}:j < 2^{|\Gamma|+1}\}$ as
$q(\bar x) \subseteq p(\bar x)$, necessarily $r_1(\bar x)$ is a subset
of $p(\bar x)$ and of course it is finite.  So $\cU_{r_1} \ne
\emptyset$ and choose $\beta \in \cU_{r_1}$ and let $j <
2^{|\Gamma|}$ be such that $\alpha_j,\beta$ are
  $\cU^r_\Gamma$-equivalent.  Recalling $\varphi_*(\bar x,\bar
  b_{\alpha(j)})^{\iif(\bold t(\alpha,j))} \in p$ so in particular
  $r_1(\bar x) \cup \{\varphi_*(\bar x),\bar
  b_{\alpha(j)})^{\iif(1-\bold t(\alpha(j)))}\}$ is consistent.

Let $\bar d'$ realize it then the pairs 
$(\alpha(j),\beta),(\beta,\alpha(j))$ belongs
to $\cU^2_p$ and at least one of them belongs to $\cU^1_r$.  So $(*)_5$
holds indeed.]
\mn
\begin{enumerate}
\item[$(*)_6$]  If $r_1,r_2 \subseteq p(\bar x)$ and
$\Gamma_1,\Gamma_2 \subseteq \Gamma_*$ are finite then $\cU^1_{r_1
\cup r_2} \cap \cU^2_{\Gamma_1 \cup \Gamma_2} \subseteq (\cU^1_{r_1}
\cap \cU^2_{\Gamma_1}) \cap (\cU^1_{r_2} \cap \cU^2_{\Gamma_2})$.
\end{enumerate}
\mn
[Why?  By inspection.]

By $(*)_5 + (*)_6$ clearly
\mn
\begin{enumerate}
\item[$(*)_7$]  there is an ultrafilter $\cD_{n_1}$ on $\mu_* \times
\mu_*$ such that:

if $r(\bar x) \subseteq p(\bar x)$ is finite and $\Gamma \subseteq
\Gamma_*$ is finite then $\cU^1_{r(\bar x)} \cap \cU^2_\Gamma \in
\cD_{n_1}$
\sn
\item[$(*)_8$]  let $\bold t$ be such that $\{(\alpha,\beta) \in \mu_*
  \times \mu_*:\bold t(\alpha) = \bold t(\beta)$ is equal to $\bold
  t\}$ belongs to $D$.
\end{enumerate}
\mn
[Why well defined?  As $\cU^2_\emptyset \in \cD_{n_1}$ by $(*)_7$ and
  see $(*)_3$.]

Let $\bar b_{n_1,0},\bar b_{n_1,1} \in {}^{\ell g(\bar y)}\gC$ be
such that $\bar b_{n_1,0} \char 94 \bar b_{n_1,1}$ realize
$\Av(\cD_{n_1},\langle \bar b_\alpha \char 94 \bar b_\beta:(\alpha,\beta) \in
\mu_* \times \mu_*\rangle,\Sigma\{\bar b_{k,\iota}:k < n_1$ and $\iota
< 2\} +A)$.  Clearly $\cD_{n_1}$ satisfies clause (A)(e).

Let $A_{n_1} := \cup\{\Rang(\bar b_\alpha):\alpha < \mu_*\}$ so clause
(A)(i) holds.  Now $(\bar b_{n_1,0},\bar b_{n_1+1})$ satisfies clauses
(A)(d) and (A)(f)$(\alpha),(\beta)$ and $r^*_{n_*+1}$ is well defined.

Lastly, concerning clause (e), the set $p(\bar x) \cup r^*_{n_1+1}$ is well
defined and consistent because for any finite $r(\bar x)
\subseteq p(\bar x)$, for the $\cD_{n_1}$-majority of $(\alpha,\beta)
\in \mu_* \times \mu_*,p(\bar x) \cup r^*_{n_1} \cup \{\varphi_*(\bar
x,\bar a_\beta)^{\iif(\bold t)},\neg \varphi_*(\bar x,\bar
a_\alpha)^{\iff(\bold t)}\}$ is inconsistent,
contradiction to $D$ assumptions.  So indeed $(A_{n_1},D'_n,
\bar b_{n_1,0},\bar b_{n_1,1})$ are as required.

Contradiction to the case assumption so really to ``$\oplus_1$ fail".
So indeed $\oplus_1$, i.e. clause $(B)(a)$ holds.
\mn
\begin{enumerate}
\item[$\oplus_2$]  choose $\bar d \in {}^n \gC$ realizing $p(\bar x)
\cup r^*_{n_1+1}$ so clauses (A)(a),(b) hold.
\end{enumerate}
\mn
[Why possible?  As $p(\bar x) \cup r^*_{n_1+1}$ is consistent by the
  induction assumption, i.e. clause (A)(g), see above.]
\mn
\begin{enumerate}
\item[$\oplus_3$]  clause (A)(f)$(\gamma)$ holds.
\sn
\item[$\oplus_4$]  clause (B)(b) holds.
\end{enumerate}
\mn
[Why?  Otherwise for every $n$ there is $q_n(\bar x) \subseteq p(\bar x)$ of
cardinality $< \mu$ for which in clause (B)(a) there is no $r(\bar x)
\subseteq p(\bar x)$ with $n$ elements such that $r(\bar x) \cup
r^*_{n_1}(\bar x) \vdash q_n(\bar x)$.  Still there is a finite
$r_n(\bar x) \subseteq p(x)$ such that $r_n(\bar x) \cup
r^*_{n_1}(\bar x) \vdash q_n(\bar x)$.
  Let $q(\bar x) =
\cup\{q_n(\bar x):n \in \bbN\}$, by (B)(a) there is a finite $r(\bar
x) \subseteq p(\bar x)$ such that $r(\bar x) \cup r^*_{n_1} \vdash
q(\bar x)$; let $n = |r(\bar x)|$ and we get a contradiction to the
choice of $q_n(\bar x)$.

Together by $\oplus_1 - \oplus_4$ and the induction hypothesis
$\boxplus$ we are done.
\end{PROOF}

\begin{claim}
\label{k4}
Assume $\Delta \subseteq \{\varphi:\varphi = \varphi(\bar x_{[m]},\bar
y) \in \bbL(\tau_T)\}$ is finite and closed under negation (well we
stipulate $\neg \neg \varphi = \varphi$).  Then \ref{k2} holds.
\end{claim}

\begin{PROOF}{\ref{k4}}
We may repeat the proof.  Alternatively we can in
\cite[Ch.II]{Sh:c} manipulate $\Delta$ to one formula 
$\varphi_*$, i.e. let $\Delta =
\{\varphi_\ell(\bar x,\bar y_\ell):\ell < n_*\}$ and we can consider only
$A$ with at least two members.  Let $\ell g(\bar y_\ell) =
k_\ell,(\forall \ell < n_*)(k_\ell \le k_0)$ let

\begin{equation*}
\begin{array}{clcr}
\varphi_*(\bar x,\bar y_0 \char 94 \langle z_0,z_1,z_2,&z_{2n_*+1}) =
\bigwedge \limits_{\ell < n_*} (z_{2n_* +1} = z_\ell \wedge
\bigwedge\limits_{k< \ell} z_{2n_* +1} \ne z_\ell \rightarrow
\varphi_\ell(\bar x,\bar y_0 \rest k_\ell)) \\
  &\wedge \bigwedge\limits_{\ell < n_*} (z_{2n_*+1} = z_{n_*+\ell}
\wedge \bigwedge\limits_{k < n_* + \ell} z_{2n_*+1} \ne z_k \rightarrow \neg
\varphi_\ell(\bar x,\bar y_\ell \rest k_\ell)) \\
  &\wedge (\bigvee\limits_{\ell < 2n_*+1} z_{2n_*+1} = z_\ell).
\end{array}
\end{equation*}

\mn
So
\mn
\begin{enumerate}
\item[$(*)_1$]  for any $\bar c \in {}^{(k_0 + 2n_*+2)}A$ one of the
following cases occurs:
\sn
\begin{enumerate}
\item[$(a)$]  for some $\ell < n_*$ and $\bar b \triangleleft \bar c$
and truth value $\bold t$ we have $(\forall \bar x)[\varphi_*(\bar x,\bar c)
\equiv \varphi_\ell(\bar x,\bar b)^{\iif(\bold t)}]$
\sn
\item[$(b)$]  $(\forall \bar x)\varphi(\bar x,\bar c)$
\sn
\item[$(c)$]  $(\forall \bar x)(\neg \varphi(\bar x,\bar c)$
\end{enumerate}
\item[$(*)_2$]  if $a^*_0 \ne a^*_1,\ell < n_*$ and $b \in {}^{\ell
g(\bar y)}\gC$ and $\bold t$ a truth value \then \, for some $\bar c
\subseteq (\Rang(\bar b) \cup \{a_0,a_1\})$ we have $(\forall \bar x)
[\varphi_*(\bar x,\bar c) \equiv \varphi_\ell(\bar x,\bar b)^{\iif(\bold t)}]$
\sn
\item[$(*)_3$]  if $a_0 \ne a_1$ then for some $\bar c_0,\bar c_1 \in 
{}^{2n_*+2)}\{a_0,a_1\}$ we have $\gC \models (\forall \bar
x)(\varphi(\bar x,\bar c_1)) \wedge (\forall \bar x)(\neg \varphi(\bar
x,\bar c_0))$.
\end{enumerate}
\end{PROOF}
\newpage


\end{document}